\pdfoutput=1
\RequirePackage{ifpdf}
\ifpdf 
\documentclass[pdftex]{sigma}
\else
\documentclass{sigma}
\fi

\usepackage{amscd,amsfonts,psfrag}
\usepackage[active]{srcltx}
\usepackage[all]{xy}

\numberwithin{equation}{section}

\newtheorem{teo}{Theorem}[section]
\newtheorem{lem}[teo]{Lemma}
\newtheorem{cor}[teo]{Corollary}
\newtheorem{prop}[teo]{Proposition}

\theoremstyle{definition}
\newtheorem{exa}[teo]{Example}

\newtheorem{defi}[teo]{Definition}

\newtheorem{Remark}[teo]{Remark}

\newcommand{\Oq}{\Oo_{q}}

\newcommand{\mUq}{\mathbb{U}_q}

\newcommand{\e}{\epsilon}
\newcommand{\ep}{{\epsilon'}}

\newcommand{\mr}{\mathbb{R}}
\newcommand{\mc}{\mathbb{C}}
\newcommand{\mz}{\mathbb{Z}}

\newcommand{\mn}{\mathbb{N}}

\newcommand{\Aa}{{\mathcal A}}

\newcommand{\Ff}{{\mathcal F}}
\newcommand{\Gg}{{\mathcal G}}

\newcommand{\Ll}{{\mathcal L}}
\newcommand{\Mm}{{\mathcal M}}

\newcommand{\Oo}{{\mathcal O}}

\newcommand{\Rr}{{\mathcal R}}

\newcommand{\ra}{\rightarrow}

\newcommand{\smallbullet}{%
\vcenter{\hbox{\,\scalebox{0.5}{$\bullet$}\,}}
}

\begin{document}
\allowdisplaybreaks

\newcommand{\arXivNumber}{2106.04136}

\renewcommand{\PaperNumber}{047}

\FirstPageHeading

\ShortArticleName{Unrestricted Quantum Moduli Algebras, II}

\ArticleName{Unrestricted Quantum Moduli Algebras, II:\\ Noetherianity and Simple Fraction Rings at Roots of~1}

\Author{St\'ephane BASEILHAC and Philippe ROCHE}

\AuthorNameForHeading{S.~Baseilhac and P.~Roche}

\Address{IMAG, Univ Montpellier, CNRS, Montpellier, France}
\Email{\href{mailto:stephane.baseilhac@umontpellier.fr}{stephane.baseilhac@umontpellier.fr}, \href{mailto:philippe.roche@umontpellier.fr}{philippe.roche@umontpellier.fr}}

\ArticleDates{Received May 11, 2023, in final form May 07, 2024; Published online June 06, 2024}

\Abstract{We prove that the quantum graph algebra and the quantum moduli algebra associated to a punctured sphere and complex semisimple Lie algebra $\mathfrak{g}$ are Noetherian rings and finitely generated rings over $\mc(q)$. Moreover, we show that these two properties still hold on $\mc\big[q,q^{-1}\big]$ for the integral version of the quantum graph algebra. We also study the specializations $\Ll_{0,n}^\e$ of the quantum graph algebra at a root of unity $\e$ of odd order, and show that $\Ll_{0,n}^\e$ and its invariant algebra under the quantum group $U_\e(\mathfrak{g})$ have classical fraction algebras which are central simple algebras of PI degrees that we compute.}


\Classification{16R30; 17B37; 20G42; 57M27; 57R56; 81R50}

\section{Introduction}\label{INTRO}

This paper is the second part of our work, initiated in \cite{BR}, on the quantum graph algebra~$\Ll_{g,n}(\mathfrak{g})$ and the quantum moduli algebra $\Mm_{g,n}(\mathfrak{g})$, which are associated to a surface $\Sigma_{g,n+1}$ of genus $g$ with $n+1$ punctures and a complex semisimple Lie algebra $\mathfrak{g}$. As in \cite{BR}, we focus in this paper on punctured spheres ($g=0$, $n\geq 1$). From now on we fix $\mathfrak{g}$, and when no confusion may arise we omit it from the notations of the various algebras.

The algebras $\Ll_{g,n}$ and $\Mm_{g,n}$ are defined over the field $\mc(q)$. They were introduced in the mid 90's by Alekseev--Grosse--Schomerus \cite{AGS1,AGS2} and Buffenoir--Roche \cite{BR1,BR2} by a method called {\it combinatorial quantization}. By this method, the pair formed by $\Ll_{g,n}$ and $\Mm_{g,n}$ appear naturally as a $q$-deformation of the Fock--Rosly \cite{FR} lattice model of the algebra of functions on the ``classical'' moduli space $\Mm_{g,n}^{cl}$ of flat $\mathfrak{g}$-connections on the surface $\Sigma_{g,n+1}$.

In \cite{BR}, we showed that both $\Ll_{0,n}$ and $\Mm_{0,n}$ have integral forms $\Ll_{0,n}^A$ and $\Mm_{0,n}^A$ defined over the ring \smash{$A= \mc\big[q,q^{-1}\big]$} (in fact we could have taken \smash{$\mathbb{Q}\big[q,q^{-1}\big]$} or \smash{${\mathbb Z}\big[q,q^{-1}\big]$} as ground ring, see Section \ref{statement}). One can thus consider the specializations of these algebras at \smash{$q=\epsilon \in \mc^\times$}, which we denote by \smash{$\Ll_{0,n}^{\e}$} and \smash{$\Mm_{0,n}^{A, \e}$} respectively. The algebra \smash{$\Ll_{0,n}^A$} is endowed with an action of the Lusztig integral form \smash{$U_A^{\rm res} =U_A^{\rm res}(\mathfrak{g})$} of the quantum group \smash{$U_q=U_q(\mathfrak{g})$}, and $\Mm_{0,n}^A$ is the subalgebra of invariant elements under this action. Therefore,
\[\Mm_{0,n}^A := \big(\Ll_{0,n}^A\big)^{U_A^{\rm res}},\qquad \Mm_{0,n} := \Ll_{0,n}^{U_q} = \Mm_{0,n}^A \bigotimes_A \mc(q).\]
The definition of $\Ll_{0,n}^A$ is based on the original combinatorial quantization method, together with twists of module-algebras and Lusztig's theory of canonical bases of quantum groups \cite{Lusztig}. This allows us to address the structure and representation theory of $\Ll_{0,n}^A$ and $\Mm_{0,n}^A$ by means of quantum groups, following ideas of classical invariant theory. In particular, we obtained that $\Ll_{0,n}$ and~$\Ll_{0,n}^\e$ have no nontrivial zero divisors (and therefore do as well the subalgebras~\smash{${\Mm_{0,n}}$},~\smash{$\Ll_{0,n}^A$},~\smash{$\Mm_{0,n}^A$}, and \smash{$\bigl(\Ll_{0,n}^\e\bigr)^{U_\e^{\rm res}}$}, where $U_\e^{\rm res}$ is the specialization of $U_A^{\rm res}$ at $q=\e$). Also, by extending the quantum coadjoint action of De Concini--Kac--Procesi \cite{DC-K,DC-K-P1,DC-K-P2}, we described in the ${\mathfrak{sl}_2}$ case an action by derivations of the center $\mathcal{Z}\bigl(\Ll_{0,n}^\e\bigr)$ of $\Ll_{0,n}^\e$ on $\Ll_{0,n}^\e$, and we defined a~subalgebra
\smash{$\mathcal{Z}\bigl(\Ll_{0,n}^\e\bigr)^\Gg \subset\mathcal{Z}\bigl(\Ll_{0,n}^\e\bigr)$},
 which is a finite extension of the ring of regular functions on the character variety of the sphere with $(n+1)$ punctures (see \cite[Corollary~7.20 and Theorem~8.8]{BR}). Moreover, from these results we derived an action by derivation of \smash{$\mathcal{Z}\bigl(\Ll_{0,n}^\e\bigr)^\Gg$} on~\smash{$\Mm_{0,n}^{A,\e}({\mathfrak{sl}_2})$}.

Representations of a quotient (the semisimplification) of $\Mm_{g,n}^{A, \e}$ were already constructed and classified in \cite{AG3}; they involve only the irreducible representations of the finite-dimensional ``small'' quantum group $\mathfrak{u}_{\e}(\mathfrak{g})$. Moreover, \cite{AG3} deduced from these representations of $\Mm_{g,n}^{A, \e}$ a family of representations of the mapping class groups of surfaces, that is equivalent to the one associated to the Witten--Reshetikhin--Turaev TQFT \cite{RT,Wi}. Recently, representations of another, larger quotient of $\Mm_{g,n}^{A, \e}$, and the corresponding representations of the mapping class groups of surfaces, were constructed in \cite{Faitg,Faitg2}. These representations are equivalent to those previously obtained by Lyubashenko--Majid \cite{LM}, and are associated to the TQFT defined in \cite{dRGGP-MR,dRGP-M}. In the ${\mathfrak{sl}_2}$ case, they involve the irreducible and also the principal indecomposable representations of the small quantum group $\mathfrak{u}_{\e}({\mathfrak{sl}_2})$. The related link and $3$-manifold invariants coincide with those of~\cite{BBG,Mu}.

In general, the representation theory of \smash{$\Mm_{g,n}^{A, \e}$} is by now far from being understood. Because $\Mm_{g,n}^{A, \e}$ deforms the classical moduli space $\Mm_{g,n}^{cl}$, it is natural to expect that its representation theory provides $(2+1)$-dimensional TQFTs for $3$-manifolds endowed with general flat $\mathfrak{g}$-connections, extending the known TQFTs based on quantum groups (where purely topological ones correspond to the trivial connection). A family of such invariants, called {\it quantum hyperbolic invariants}, has already been defined for $\mathfrak{g}={\mathfrak{sl}_2}$ by means of certain $6j$-symbols, {\it Deus ex machina}; they are closely connected to classical Chern--Simons theory, provide generalized volume conjectures, and contain quantum Teichm\"uller theory (see \cite{BB-3,BB-2,BB-1,BB0,BB1,BB2,BB3}). It is part of our present program, initiated in \cite{B}, to shed light on these invariants and to generalize them to arbitrary $\mathfrak{g}$ by developing the representation theory of $\Mm_{g,n}^{A, \e}$.

 The quantum moduli algebras have also been recognized as central objects from the viewpoints of factorization homology \cite{B-BZ-J}, multiplicative quiver varieties \cite{GJS19} and (stated) skein theory~\mbox{\cite{BFR,BFK,CL,Faitg4}}. In another direction, one may expect that the equivalence proved in~\cite{Meus} between combinatorial quantisation for the Drinfeld double $D(H)$ of a finite-dimensional semisimple Hopf algebra $H$, and Kitaev's lattice model in topological quantum computation, can be extended to the setup of quantum moduli algebras.

In the present paper, we study $\Ll_{0,n}$, its integral form $\Ll_{0,n}^A$, and the specialization~$\Ll_{0,n}^\e$ of~$\Ll_{0,n}^A$ at~$q=\e$ a primitive root of unity of odd order. We study also the subalgebras of invariant elements $\Mm_{0,n} = \Ll_{0,n}^{U_q}$ and $\bigl(\Ll_{0,n}^\e\bigr)^{U_\e}$. Here, $U_\e$ is the specialization of $U_A$ at $q=\epsilon$, where $U_A$ is the De Concini--Kac integral form of $U_q$ (see Section \ref{statement}). Our results hold for every complex semisimple Lie algebra $\mathfrak{g}$. The main ones are proofs that $\Ll_{0,n}$, $\Ll_{0,n}^A$ and $\Mm_{0,n}$ are Noetherian and finitely generated rings (see Theorem \ref{HN}), and that the classical fraction algebras of $\Ll_{0,n}^\e$ and \smash{$\bigl(\Ll_{0,n}^\e\bigr)^{U_\e}$} are central simple algebras of PI degrees \smash{$l^{nN}$} and \smash{$l^{N(n-1)-m}$} respectively (see Theorem~\ref{degree}). Here, $m$ and $N$ are the rank and the number of positive roots of $\mathfrak{g}$.

In the sequel \cite{BFR} to this paper, in collaboration with M.~Faitg, we extend Theorem~\ref{HN} to the algebras $\Ll_{g,n}$ and $\Mm_{g,n}$, associated to arbitrary finite type surfaces (arbitrary genus and number of punctures). Also, we show that $\Mm_{g,n}$ is isomorphic to the $\mathfrak{g}$-skein algebra of $\Sigma_{g,n+1}$, and $\Ll_{g,n}$ to the stated skein algebra of the compact surface $\overline{\Sigma}_{g,n+1}$ with one boundary component and one marked point on the boundary component. This was proved for $\mathfrak{g}={\mathfrak{sl}_2}$ in \cite{Faitg4}. In this specific case $\mathfrak{g}={\mathfrak{sl}_2}$, the fact that the stated skein algebra of any finite type surface is Noetherian and finitely generated was proved in \cite{LY}. Still in the ${\mathfrak{sl}_2}$ case, for related results, e.g., on non-zero divisors and computation of PI degrees, see \cite{BW,BW1,FKL,KK,K0,K1,KQ,Le}. For recent results on~${\mathfrak{g} ={\mathfrak{sl}_n}}$, see \cite{LeSik, Wang}.

By using the analysis developed in the present paper for $\Ll_{0,n}^A$, one can define the integral form~\smash{$\Ll_{g,n}^A$} as well, and show that it is a Noetherian and finitely generated ring. We do not have a proof yet of these properties for the algebra \smash{$\Mm_{0,n}^{A}$}, which seems to be much more difficult to handle. We note that there is a strict inclusion \smash{$\Mm_{0,n}^{A,\e} \subset \bigl(\Ll_{0,n}^\e\bigr)^{U_\e}$}. This is discussed after Theorem~\ref{Llibre}. In \cite{BFR2}, we study further properties of $\bigl(\Ll_{g,n}^\e\bigr)^{U_\e}$, and we consider also the subalgebra~\smash{$\Mm_{g,n}^{A,\e}$}.

\subsection{Statement of results}\label{statement} Let us recall a few notations and facts from \cite{BR}. Let $U_q$ be the simply-connected quantum group of $\mathfrak{g}$, defined over the field $\mc(q)$. From $U_q$ one can define a $U_q$-module algebra $\Ll_{0,n}$, called (quantum, daisy) {\it graph algebra}, where $U_q$ acts by means of a right coadjoint action. The set of invariant elements of $\Ll_{0,n}$ for this action is an algebra; we denote it \smash{$\Mm_{0,n}:=\Ll_{0,n}^{U_q}$} and call it {\it quantum moduli algebra}. As a $\mc(q)$-module $\Ll_{0,n}$ is just $\mathcal{O}_q^{\otimes n}$, where $\mathcal{O}_q = \mathcal{O}_q(G)$ is the standard quantum function algebra of the connected and simply-connected Lie group $G$ with Lie algebra~$\mathfrak{g}$. The product of $\Ll_{0,n}$ is obtained by twisting both the product of each factor~$\Oo_q$ and the product between them. It is equivariant with respect to a (right) coadjoint action of~$U_q$, which defines the structure of $U_q$-module of $\Ll_{0,n}$.

The module algebra $\Ll_{0,n}$ has an {integral form} $\Ll_{0,n}^A$, which is defined over $A=\mc\big[q,q^{-1}\big]$, and endowed with an (coadjoint) action of the Lusztig \cite{Lusztig2} integral form $U_A^{\rm res}$ of $U_q$. It is obtained by replacing $\mathcal{O}_q$ in the construction of $\Ll_{0,n}$ with the restricted dual $\mathcal{O}_A$ of the integral form~$U_A^{\rm res}$, or equivalently with the restricted dual of the integral form $\Gamma$ of $U_q$ defined by De Concini--Lyubashenko \cite{DC-L}. Since $U_A^{\rm res}$ contains the De Concini--Kac \cite{DC-K} integral form $U_A$, and $U_A$ has the same set of invariant elements in $\Ll_{0,n}^A$, we systematically denote the latter
\[\Mm_{0,n}^A := \big(\Ll_{0,n}^A\big)^{U_A} \qquad \bigl(=\big(\Ll_{0,n}^A\big)^{U_A^{\rm res}}\bigr).\]
We call $\Mm_{0,n}^A$ the {\it integral} quantum moduli algebra.

A cornerstone of the theory of $\Mm_{0,n}$ is a map $\Phi_n$ originally due to Alekseev \cite{A}, building on works of Drinfeld \cite{Dr} and Reshetikhin and Semenov-Tian-Shansky \cite{RSTS}. In \cite{BR}, we showed that~$\Phi_n$ eventually provides isomorphisms of module algebras and algebras respectively,
\[\Phi_n\colon\ \Ll_{0,n} \ra \big(U_q^{\otimes n}\big)^{\rm lf},\qquad \Phi_n\colon\ \Mm_{0,n} \ra \big(U_q^{\otimes n}\big)^{U_q},
\]
where $U_q^{\otimes n}$ is endowed with a right adjoint action of $U_q$, and $\big(U_q^{\otimes n}\big)^{\rm lf}$ is the subalgebra of locally finite elements with respect to this action. When $n=1$ the algebra \smash{$U_q^{\rm lf}$} has been studied in great detail by Joseph--Letzter \cite{Jos,JL,JL2}; we will use simplified proofs of their results, obtained in~\cite{VY}.

All the material we need about the results discussed above is described in \cite{BR}, and recalled in Sections \ref{ALGdef} and~\ref{INTEG}.

Our first result, proved in Section \ref{HNsec}, is the following.

\begin{teo}\label{HN} $\Ll_{0,n}$, $\Mm_{0,n}$ and the integral form $\Ll_{0,n}^A$ are Noetherian rings, and finitely generated rings.
\end{teo}
It follows immediately from the theorem that the specializations $\Ll_{0,n}^\e$, $\e \in \mc^\times$, are Noetherian and finitely generated rings as well. In \cite{BR} we proved that all these algebras (and therefore $\Mm_{0,n}^A$ and $\Mm_{0,n}^{A,\e}$) have no nontrivial zero divisors.

Because the construction of the integral form \smash{$\Ll_{0,n}^A$} is based on the Kashiwara--Lusztig theory of canonical bases, we could have defined \smash{$\Ll_{0,n}^A$} over the ground ring $\mathbb{Z}\big[q,q^{-1}\big]$, and Theorem~\ref{HN} for \smash{$\Ll_{0,n}^A$} holds true as well in this generality. Since we are mainly interested in the representation theory of the specializations \smash{$\Ll_{0,n}^\e$} and \smash{$\Mm_{0,n}^{A,\e}$}, which will be addressed in \cite{BFR2}, the choice of~${A=\mc\big[q,q^{-1}\big]}$ is natural. Note however that the proof of Proposition \ref{L0NAfgfree} uses that~${\mc\big[q,q^{-1}\big]}$ is a~PID.

We describe the background material on canonical bases in Section \ref{canbasemodqg}; we have tried to make the exposition pedestrian and self-contained, so as to be more accessible to non experts.

After we finished this work, we discovered that \cite{DomLen} already proved that $\Ll_{0,1}(\mathfrak{gl}(n))$ and $\Ll_{0,n}(\mathfrak{gl}(2))$ are Noetherian and finitely generated rings. Our approach here is completely different. For $\Ll_{0,n}$, we adapt the proof given by Voigt--Yuncken \cite{VY} of a result of Joseph \cite{Jos}, which asserts that \smash{$U_q^{\rm lf}$} is a Noetherian ring (see Theorem~\ref{LonNoeth}). For~$\Mm_{0,n}$, we deduce the result from the one for $\Ll_{0,n}$, by following a line of proof of the Hilbert--Nagata theorem in classical invariant theory (see Theorem \ref{MonNoeth}).

At present, we do not have a proof that $\Mm_{0,n}^A$ is a Noetherian and finitely generated ring for arbitrary $\mathfrak{g}$ and $n\geq 1$, though it is natural to expect it is the case. Indeed, when ${\mathfrak{g}= {\mathfrak{sl}_2}}$, $\Mm_{0,n}^A({\mathfrak{sl}_2})$~is isomorphic to the skein algebra of a sphere with $n+1$ punctures (see \cite[Theorem~8.6]{BR}), which is finitely generated and Noetherian by results of \cite{Bul} and \cite{PS}. In our general situation, key arguments in the proof of Theorem~\ref{HN} for~$\Mm_{0,n}$ depend on the existence of a~Reynolds operator on the $U_q$-module~$\Ll_{0,n}$, and one can easily show there is no Reynolds operator on~$\Ll_{0,n}^A$. This follows from the corresponding fact for the integral quantum coordinate ring~$\Oo_A$ (see Remark~\ref{noHaarRoverA}).

From Section \ref{ORDER}, we consider the specializations $\Ll_{0,n}^\e$ of $\Ll_{0,n}^A$ at $q=\e$, a primitive root of unity of odd order $l$ (and coprime to $3$ if $\mathfrak{g}$ has $G_2$ components). In~\cite{DC-L}, De Concini--Lyubashenko introduced a central subalgebra $\mathcal{Z}_0(\Oo_\e)$ of $\Oo_\e$ isomorphic to the coordinate ring~$\Oo(G)$, and proved that the $\mathcal{Z}_0(\Oo_\e)$-module $\Oo_\e$ is projective of rank $l^{\dim \mathfrak{g}}$. As observed by Brown--Gordon--Stafford~\cite{BGS}, Bass' cancellation theorem in $K$-theory and the fact that $K_0(\Oo(G))\cong \mz$, proved by Marlin~\cite{Marlin}, imply that this module is free. Alternatively, this follows also from the fact that~$\Oo_\e$ is a cleft extension of $\Oo(G)$ by the dual of the Frobenius--Lusztig kernel $\mathfrak{u}_\e(\mathfrak{g})$, as proved by Andruskiewitsch--Garcia (see \cite[Remark~2.18\,(b)]{AnGa}, and also \cite[Section~3.2]{BC}; this argument was explained to us by K.A.~Brown).

The Section \ref{ORDER} proves the analogous property for $\Ll_{0,n}^\e$. Namely:
\begin{teo}\label{Llibre} $\mathcal{Z}_0(\Oo_\e)^{\otimes n}$ is a central subalgebra of $\Ll_{0,n}^\e$, and $\Ll_{0,n}^\e$ is a free $\mathcal{Z}_0(\Oo_\e)^{\otimes n}$-module of rank \smash{$l^{n.\dim \mathfrak{g}}$}, isomorphic to the \smash{$\mathcal{Z}_0(\Oo_\e)^{\otimes n}$}-module \smash{$\Oo_\e^{\otimes n}$}.
\end{teo}
In the sequel we systematically denote $\mathcal{Z}_0\bigl(\Ll_{0,n}^\e\bigr):= \mathcal{Z}_0(\Oo_\e)^{\otimes n}$.
We prove the first and third claims of Theorem \ref{Llibre} in Proposition \ref{Z0L0n}. The arguments use results of De Concini--Kac \cite{DC-K}, De Concini--Procesi \cite{DC-K-P1,DC-K-P2}, and De Concini--Lyubashenko \cite{DC-L}, that we recall in Sections \ref{intdual}--\ref{UOUlf}. Let us stress that the algebra structures of $\Ll_{0,n}^\e$ and $\Oo_\e^{\otimes n}$ are completely different.

Since $\mathcal{Z}_0(\Oo_\e)\cong \Oo(G)$, we may deduce the second claim of Theorem \ref{Llibre} from the first and third claims together with the results of \cite{DC-L,Marlin}, or \cite{AnGa}, recalled above. Nevertheless, we give a self-contained proof that $\Ll_{0,1}^\e$ is finite projective of rank $l^{\dim \mathfrak{g}}$ over $\mathcal{Z}_0\bigl(\Ll_{0,1}^\e\bigr)$, by adapting the original arguments of De Concini--Lyubashenko \cite[Theorem~7.2]{DC-L}. In particular, we study the coregular action of the braid group of $\mathfrak{g}$ on $\Ll_{0,1}^\e$; by the way, in the appendix, we provide different proofs of some technical facts shown in \cite{DC-L}. Of course, it remains an exciting problem to describe the centralizing extension $\mathcal{O}(G)^{\otimes n} \subset \Ll_{0,n}^\e$ (and similarly $\mathcal{O}(G)^{\otimes n} \subset \bigl(\Ll_{0,n}^\e\bigr)^{U_\e}$ below), aiming at generalizing the results of~\cite{AnGa} and finding a direct, more structural proof of freeness in Theorem~\ref{Llibre}.
Also, we note that bases of $\Ll_{0,n}^\e$ over $\mathcal{Z}_0\bigl(\Ll_{0,n}^\e\bigr)$ are complicated. The only case we know is for~$\Oo_\e({\mathfrak{sl}_2})$, described in \cite{DRZ}, and it is far from being obvious (see \eqref{basisSL2}).

In Section \ref{DEGREEsec}, we turn to fraction rings. As mentioned above $\Ll_{0,n}^\e$ has no nontrivial zero divisors. Therefore, its center $\mathcal{Z}\bigl(\Ll_{0,n}^\e\bigr)$ is an integral domain. Denote by $Q(\mathcal{Z}\bigl(\Ll_{0,n}^\e\bigr))$ its fraction field. Denote by \smash{$\bigl(\Ll_{0,n}^\e\bigr)^{U_\e}$} the subring of $\Ll_{0,n}^\e$ formed by the invariant elements of $\Ll_{0,n}^\e$ with respect to the right coadjoint action of $U_\e$. The center $\mathcal{Z}(\Ll_{0,n}^\epsilon)$ of $\Ll_{0,n}^\e$ is contained in~\smash{$\bigl(\Ll_{0,n}^\e\bigr)^{U_\e}$} (this follows~from \cite[Proposition~6.19]{BR}). Note also that we trivially have an in\-clusion \smash{$\Mm_{0,n}^{A,\e}\subset \bigl(\Ll_{0,n}^\e\bigr)^{U_\e}$}, and these two algebras are distinct in general. For instance, when $n=1$, we~have $\bigl(\Ll_{0,1}^\e\bigr)^{U_\e}=\mathcal{Z}(\Ll_{0,1}^\epsilon)$, which is a finite extension of $\mathcal{Z}_0(\Oo_\e) \cong \Oo(G)$ (see Lemma~\ref{ZfiniteZ0}). On another hand, \smash{$\Mm_{0,1}^{A,\e}$} is the specialization at $q=\e$ of \smash{$\mathcal{Z}\bigl(\Ll_{0,1}^A\bigr)$}, a polynomial algebra in ${\rm rk}(\mathfrak{g})$ variables, which may be identified via $\Phi_1$ with the center $\mathcal{Z}(U_A)$ of the integral form~$U_A$.

Consider the rings
\[Q\bigl(\Ll_{0,n}^\e\bigr) = Q\bigl(\mathcal{Z}\bigl(\Ll_{0,n}^\e\bigr)\bigr) \bigotimes_{\mathcal{Z}(\Ll_{0,n}^\e)} \Ll_{0,n}^\e,\qquad Q\big(\bigl(\Ll_{0,n}^\e\bigr)^{U_\e}\big) = Q\bigl(\mathcal{Z}\bigl(\Ll_{0,n}^\e\bigr)\bigr) \bigotimes_{\mathcal{Z}(\Ll_{0,n}^\e)} \bigl(\Ll_{0,n}^\e\bigr)^{U_\e}.\]
In general, given a ring $A$ with center $\mathcal{Z}(A)$ an integral domain we reserve the notation $Q(A)$ to~the central localization of $A$, i.e., \smash{$Q(A):= Q(\mathcal{Z}(A)) \bigotimes_{\mathcal{Z}(A)} A$}. Though the center \smash{$\mathcal{Z}\big(\bigl(\Ll_{0,n}^\e\bigr)^{U_\e}\big)$} of \smash{$\bigl(\Ll_{0,n}^\e\bigr)^{U_\e}$} is larger than $\mathcal{Z}\bigl(\Ll_{0,n}^\e\bigr)$, the notation \smash{$Q\big(\bigl(\Ll_{0,n}^\e\bigr)^{U_\e}\big)$} is valid, for \smash{$\mathcal{Z}\big(\bigl(\Ll_{0,n}^\e\bigr)^{U_\e}\big)$} is an integral domain finite over $\mathcal{Z}\bigl(\Ll_{0,n}^\e\bigr)$, and hence the central localization of \smash{$\bigl(\Ll_{0,n}^\e\bigr)^{U_\e}$} coincides with \smash{$Q\big(\bigl(\Ll_{0,n}^\e\bigr)^{U_\e}\big)$} as defined above. Throughout the paper, unless we mention it explicitly, we follow the conventions of McConnell--Robson \cite{MC-R} as regards the terminology of ring theory; in particular, for the notions of central simple algebras and PI degrees, see in \cite[Sections~5.3 and~13.3.6--13.6.7]{MC-R}.

Denote by $m$ the rank of $\mathfrak{g}$, and by $N$ the number of its positive roots. In Section \ref{DEGREEsec}, we prove the following.

\begin{teo}\label{degree} \quad
\begin{itemize}\itemsep=0pt
\item[$(1)$] $Q\bigl(\Ll_{0,n}^\e\bigr)$ is a division algebra and a central simple algebra of PI degree \smash{$l^{nN}$}.
\item[$(2)$] \smash{$Q\big(\bigl(\Ll_{0,n}^\e\bigr)^{U_\e}\!\big)$}, $n\!\geq\! 2$, is a division algebra and a central simple algebra of PI degree~\smash{$l^{N(n{-}1){-}m}.$\!}
 \end{itemize}
\end{teo}

The second claim of (1) means that $Q\bigl(\Ll_{0,n}^\e\bigr)$ is a complex subalgebra of a full matrix algebra ${\rm Mat}_d(\mathbb{F})$, where $d=l^{nN}$ and $\mathbb{F}$ is a finite extension of $Q(\mathcal{Z}\bigl(\Ll_{0,n}^\e\bigr))$ such that
\[\mathbb{F}\bigotimes_{Q(\mathcal{Z} (\Ll_{0,n}^\e ))} Q\bigl(\Ll_{0,n}^\e\bigr) = {\rm Mat}_d(\mathbb{F}).\]
That $Q\bigl(\Ll_{0,n}^\e\bigr)$ is a division algebra and a central simple algebra follows from Theorem \ref{Llibre} and the fact that $\Ll_{0,n}^\e$ has no nontrivial zero divisors (proved in \cite{BR}). The computation of \smash{$d=l^{nN}$} uses a lower bound coming from the representation theory of $U_\e$, and a lower bound for the degree of $Q\bigl(\mathcal{Z}\bigl(\Ll_{0,n}^\e\bigr)\bigr)$ as a field extension of $Q\bigl(\mathcal{Z}_0\bigl(\Ll_{0,n}^\e\bigr)\bigr)$, obtained by using specializations to~${q=\e}$ of certain central elements in $\Ll_{0,n}$ (for $q$ generic). In this computation a main role is played by results of De Concini--Kac~\cite{DC-K}.

We deduce (2) from (1), the double centralizer theorem for central simple algebras, a few results of \cite{BR,DC-L}, and Theorem~\ref{Llibre} again.

\subsection{Basic notations} Given a ring $R$, we denote by $\mathcal{Z}(R)$ its center. When $R$ is commutative and has no nontrivial zero divisors, $Q(R)$ denotes its fraction field.

Given a Hopf algebra $H$ with product $m$ and coproduct $\Delta$, we denote by $H^{\rm cop}$ (resp.~$H_{\rm op}$) the Hopf algebra with the same algebra (resp.\ coalgebra) structure as $H$ but the opposite coproduct $\Delta^{\rm cop}:=\sigma \circ {\Delta}$ (resp.\ opposite product $m\circ \sigma$), where $\sigma(x\otimes y) = y\otimes x$, and the antipode ${S}^{-1}$. We use Sweedler's coproduct notation, $ \Delta(x) = \sum_{(x)}x_{(1)} \otimes x_{(2)}$, $x\in H$, and we set \smash{$\Delta^{(1)} := {\rm id}$}, \smash{$\Delta^{(2)} := \Delta$}, and \smash{$\Delta^{(n)} := (\Delta \otimes {\rm id}) \Delta^{(n-1)}$} for $n\geq 3$ (this is not the convention used in \cite{BR}).

The results of this paper hold true for any finite-dimensional complex semisimple Lie algebra~${\mathfrak g}$, but unless we state it differently, we will assume ${\mathfrak g}$ is simple. We will denote its rank by~$m$, and its Cartan matrix by~$(a_{ij})$. We fix a Cartan subalgebra $\mathfrak{h}\subset {\mathfrak g}$ and a basis of simple roots~${\alpha_i \in \mathfrak{h}_{\mr}^{*}}$, and denote by $\mathfrak{b}_\pm$ the Borel subalgebras associated to it. We denote by $N$ the number of positive roots of ${\mathfrak g}$, and by $\rho$ half the sum of the positive roots.

We denote by $d_1,\dots , d_m$ the unique coprime positive integers such that the matrix $(d_ia_{ij})$ is symmetric, and $(\ ,\ )$ the unique inner product on $\mathfrak{h}_{\mr}^{*}$ such that $d_ia_{ij} = (\alpha_i,\alpha_j)$. For any root~$\alpha$, the coroot is $\alpha\check{}=2\alpha/(\alpha,\alpha)$; in particular $\alpha\check{}_i=d_i^{-1}\alpha_i$. The root lattice $Q$ is the $\mz$-lattice in~$\mathfrak{h}_{\mr}^{*}$ defined by $ Q = \sum_{i=1}^m \mathbb{Z} \alpha_i$. The weight lattice $P$ is the $\mz$-lattice formed by all $\lambda\in \mathfrak{h}_{\mr}^{*}$ such that~${(\lambda ,\alpha\check{}_i) \in \mz}$ for every $i=1,\dots,m$. So $ P= \sum_{i=1}^m \mathbb{Z} \varpi_i$, where $\varpi_i$ is the fundamental weight dual to the simple coroot $\alpha\check{}_i$, which satisfies $( \varpi_i,\alpha\check{}_j) = \delta_{i,j}$. Note that $(\lambda,\alpha)\in \mz$ for every $\lambda\in P$, $\alpha\in Q$. We denote by $D$ the cardinality of the quotient lattice $P/Q$. Then $D$ is the smallest positive integer such that $D(\lambda,\mu)\in {\mathbb Z}$ for every $\lambda,\mu\in P$, that is, such that $DP\subset Q$.

We denote by
\[ P_+:= \sum_{i=1}^m \mathbb{Z}_{\geq 0} \varpi_i
\]
 the cone of dominant integral weights, and we put
 \[ Q_+:= \sum_{i=1}^m \mathbb{Z}_{\geq 0} \alpha_i.\]
 Though $Q\subset P$, it is not true that $Q_+\subset P_+$, but we have $DP_+ \subset Q_+$. This last property is not trivial, and follows from the classical fact that the inverse of the Cartan matrix~$(a_{ij})$ has coefficients in $D^{-1}\mn$.

We will use the standard partial order relation $\leq$ on $P$, defined by: $\lambda, \mu\in P$ satisfy $\lambda \leq \mu$ if~${\mu-\lambda \in Q_+}$. In Section \ref{HNsec}, we will also use the partial order relation $\preceq$ on $P$ defined by: $\lambda \preceq \mu$ if $\mu-\lambda \in D^{-1}Q_+$.

We denote by $\mathcal{B}(\mathfrak{g})$ the braid group of $\mathfrak{g}$; we recall its standard defining relations in Appendix~\ref{QW}.

We denote by $G$ the connected and simply-connected algebraic group with Lie algebra $\mathfrak{g}$, and by $T_G$ the maximal torus of $G$ with Lie algebra $\mathfrak{h}$; $N(T_G)$ is the normalizer of $T_G$, $W = N(T_G)/T_G$ is the Weyl group, $B_\pm$ are the Borel subgroups of $G$ with Lie algebra $\mathfrak{b}_\pm$, and $U_\pm \subset B_\pm$ are their unipotent subgroups.

We denote by $\Oo(G)$ the coordinate ring of $G$. It is a commutative Hopf algebra, which can be identified with the restricted dual of the universal enveloping algebra $U(\mathfrak{g})$ (see \cite{Ko,Lusztig3}). Similarly we denote by $\Oo(B_\pm)$ the coordinate ring of~$B_\pm$.

Let $q$ be an indeterminate, let $q^{1/D}$ be such that $\bigl(q^{1/D}\bigr)^D=q$, set $A={\mathbb C}\big[q,q^{-1}\big]$, $q_i=q^{d_i}$, $q_\beta = q^{(\beta,\beta)/2}$ for $\beta\in Q$, and given integers $p$, $k$ with $0\leq k\leq p$, we put
\begin{alignat*}{5}
 & [p]_q = \frac{q^p-q^{-p}}{q- q^{-1}} ,\qquad &&[0]_q! =1,\qquad &&[p]_q! =[1]_q[2]_q\cdots [p]_q ,\qquad &&\left[\begin{matrix} p \\ k \end{matrix}\right]_{q} =\frac{[p]_q! }{[p-k]_q![k]_q!},&\\
& (p)_q = \frac{q^p-1}{q- 1} ,\qquad&& (0)_q! =1 ,\qquad&& (p)_q! =(1)_q(2)_q\cdots (p)_q,\qquad &&\left(\begin{matrix} p \\ k \end{matrix}\right)_{q} =\frac{(p)_q! }{(p-k)_q!(k)_q!}.&
\end{alignat*}
We denote by $\mathcal{A}_0\subset \mc(q)$ the ring of functions regular at $q=0$; this ring is used only in Section~\ref{canbasemodqg}.

We denote by $\e$ a primitive $l$-th root of unity such that $\e^{2d_i}\ne 1$ is also a primitive $l$-th root of unity for all $i\in \{1,\dots,m\}$. This means that $l$ is odd, and coprime to $3$ if $\mathfrak{g}$ is $G_2$. We put $\e_i:=\e^{d_i}$.

In this paper, we use the definition of the unrestricted integral form $U_A(\mathfrak{g})$ given in \cite{DC-L,DC-K-P2}; in \cite{BR} we used the one of \cite{DC-K,DC-K-P1}. The two are (trivially) isomorphic, and have the same specialization at $q=\e$. Also, we denote here by $L_i$ the generators of $U_q(\mathfrak{g})$ we denoted by $\ell_i$ in~\cite{BR}.

In order to facilitate the comparison with the results of~\cite{DC-L}, we note that their generators denoted $K_i$, $E_i$ and $F_i$, that we will denote by $\tilde{K}_i$, $\tilde{E}_i$ and $\tilde{F}_i$, can be written as $K_i$, $K_i^{-1}E_i$ and~$F_iK_i$ in our notations. They satisfy the same algebra relations.

\section{Background results}\label{PREL}

\subsection[On U\_q, O\_q, L\_{0,n}, M\_{0,n}, and Phi\_n]{On $\boldsymbol{U_q}$, $\boldsymbol{\Oo_q}$, $\boldsymbol{\Ll_{0,n}}$, $\boldsymbol{\Mm_{0,n}}$, and $\boldsymbol{\Phi_n}$}\label{ALGdef}
Except when stated differently, we refer to \cite[Sections 2--4 and 6]{BR}, and the references therein for details about the material of this section. We stress that the simply-connected quantum group, that we denote $U_q$ below, was denoted $\tilde U_q$ in \cite{BR}. Also, the adjoint quantum group $U_q^{\rm ad}$ was denoted $U_q$.

The {\it simply-connected} quantum group $U_q = U_q(\mathfrak{g})$ is the Hopf algebra over $\mathbb {C}(q)$ with generators $E_i$, $F_i$, $L_i$, $L_i^{-1}$, $1\leq i \leq m$, and defining relations
\begin{gather*}
L_iL_j=L_jL_i , \qquad L_iL_i^{-1}=L_i^{-1}L_i=1,\qquad L_iE_jL_i^{-1}=q_i^{\delta_{i,j}}E_j ,\qquad L_iF_jL_i^{-1}=q_i^{-\delta_{i,j}}F_j, \\
E_iF_j-F_jE_i=\delta_{i,j}\frac{K_i-K_i^{-1}}{q_i-q_i^{-1}},\\
\sum_{r=0}^{1-a_{ij}} (-1)^r \left[\begin{matrix} 1-a_{ij} \\ r \end{matrix}\right]_{q_i} E_i^{1-a_{ij}-r}E_jE_i^{r} = 0 \qquad {\rm if}\quad i\ne j,\\
\sum_{r=0}^{1-a_{ij}} (-1)^r \left[\begin{matrix} 1-a_{ij} \\ r \end{matrix}\right]_{q_i} F_i^{1-a_{ij}-r}F_jF_i^{r} = 0 \qquad {\rm if}\quad i\ne j,
\end{gather*}
where for $ \lambda=\sum_{i=1}^m m_i \varpi_i \in P$ we set $ K_\lambda=\prod_{i=1}^m L_i^{m_i}$, and $ K_i=K_{\alpha_i}=\prod_{j=1}^m L_j^{a_{ji}}$.
The coproduct $\Delta$, antipode $S$, and counit $\varepsilon$ of $U_q$ are given by
\begin{gather*}
\Delta(L_i)=L_i\otimes L_i ,\qquad \Delta(E_i)=E_i\otimes K_i+1\otimes E_i ,\qquad \Delta(F_i)=F_i\otimes 1 + K_i^{-1}\otimes F_i, \\ S(E_i) = -E_iK_i^{-1},\qquad S(F_i) = -K_iF_i ,\qquad S(L_i) = L_i^{- 1},\\ \varepsilon(E_i) = \varepsilon(F_i)=0,\qquad \varepsilon(L_i)=1.
\end{gather*}
We fix a reduced expression $s_{i_1}\cdots s_{i_N}$ of the longest element $w_0$ of the Weyl group of $\mathfrak{g}$. It induces a total ordering of the positive roots,
\[\beta_1 = \alpha_{i_1},\qquad \beta_2 = s_{i_1}(\alpha_{i_2}), \qquad \dots, \qquad\beta_N = s_{i_1}\cdots s_{i_{N-1}}(\alpha_{i_N}).\]
The root vectors of $U_q$ with respect to such an ordering are defined by
\begin{equation}\label{rootvectdef}
E_{\beta_k} = T_{i_1}\cdots T_{i_{k-1}}(E_{i_k}) ,\qquad F_{\beta_k} = T_{i_1}\cdots T_{i_{k-1}}(F_{i_k}),
\end{equation}
where $T_i$ is the Lusztig algebra automorphism of $U_q$ associated to the simple root $\alpha_i$ \mbox{\cite{Lusztig2,Lusztig}} (see also \cite[Chapter~8]{CP}). The braid group $\mathcal{B}(\mathfrak{g})$ acts on $U_q$ by means of the Lusztig automorphisms. In the appendix, we recall the relation between $T_i$ and the generator $\hat{w}_i$ of the quantum Weyl group, which we will mostly use. Let us just recall here that the monomials $F_{\beta_1}^{r_1}\cdots F_{\beta_N}^{r_N}K_\lambda E_{\beta_N}^{t_N}\cdots E_{\beta_1}^{t_1}$ ($r_i,t_i\in \mn$, $\lambda\in P$) form a basis of $U_q$, the {\it PBW basis}.

$U_q$ is a {\it pivotal} Hopf algebra, with pivotal element \smash{$ \ell :=K_{2\rho}=\prod_{j=1}^m L_j^2$}.
So $\ell$ is group-like, and $S^2(x) = \ell x\ell^{-1}$ for every $x\in U_q$.

The {\it adjoint} quantum group $U_q^{\rm ad} = U_q^{\rm ad}(\mathfrak{g})$ is the Hopf subalgebra of $U_q$ generated by the elements $E_i$, $F_i$ ($i=1,\dots,m$) and $K_\alpha$ with $\alpha\in Q$; so $\ell\in U_q^{\rm ad}$. When $\mathfrak{g}={\mathfrak{sl}_2}$, we simply write the above generators $E=E_1$, $F=F_1$, $L=L_1$, $K=K_1$.

We denote by $U_q(\mathfrak{n}_+)$, $U_q(\mathfrak{n}_-)$ and $U_q(\mathfrak{h})$ the subalgebras of $U_q$ generated respectively by the~$E_i$, the $F_i$, and the $K_\lambda$ ($\lambda\in P$), and by $U_q(\mathfrak{b}_+)$ and $U_q(\mathfrak{b}_-)$ the subalgebras generated by the~$E_i$ and the $K_\lambda$, and by the $F_i$ and the $K_\lambda$, respectively. We do similarly with $U_q^{\rm ad}$, where now $U_q^{\rm ad}(\mathfrak{h})$ is generated by the $K_\lambda$ with~$\lambda\in Q$.

The Hopf algebra $U_q^{\rm ad}$ is not braided in a strict sense, but it has braided categorical completions. Let us recall briefly what this means and implies. For details, we refer to \cite[Sections~2 and~3]{BR} (see also \cite[Section~3.10]{VY}, where $\mUq$ below is formulated in terms of multiplier Hopf algebras).

A $U_q^{\rm ad}$-module $V$ is said {\it of type $1$} if it has finite dimension and the generators $K_i$ are diagonalizable on $V$ with eigenvalues in \smash{$q_i^{\mz}$}. We denote by $\mathcal{C}$ the category of $U_q^{\rm ad}$-modules of type $1$, by~${\rm Vect}$ the category of finite-dimensional $\mc(q)$-vector spaces, and by $F_{\mathcal C}\colon{\mathcal C}\to {\rm Vect}$ the forgetful functor. The category $\mathcal{C}$ is semisimple. The simple objects are highest weight $U_q^{\rm ad}$-modules; we denote by $V_\mu$ the simple module with highest weight $\mu\in P_+$. In the case $\mathfrak{g}={\mathfrak{sl}_2}$, we identify~$P_+$ with $\mathbb{N}$, and denote by $V_n$ the simple module of dimension $n+1$. Note that $V_\mu$ is canonically endowed with a structure of $U_q$-module of type $1$, the generators $L_i$ being diagonalizable with eigenvalues in \smash{$q_i^{\mz/D}$}. The {\it categorical completion} $\mUq^{\rm ad}$ of $U_q^{\rm ad}$ is the set of natural transformations~${F_{\mathcal C}\ra F_{\mathcal C}}$. An element of \smash{$\mUq^{\rm ad}$} is a collection $(a_V)_{V\in {\rm Ob}(\mathcal{C})}$, where $a_V\in {\rm End}_{\mc(q)}(V)$ satisfies~${F_{\mathcal C}(f)\circ a_V=a_W\circ F_{\mathcal C}(f)}$ for any objects $V$, $W$ of $\mathcal{C}$ and any arrow \smash{$f\in {\rm Hom}_{U_q^{\rm ad}}(V,W)$}. It is not hard to see that $\mUq^{\rm ad}$ inherits from $\mathcal{C}$ a natural structure of (completed) Hopf algebra such that the map
\begin{gather}\label{defiota}
\iota\colon\ U_q^{\rm ad} \longrightarrow \mUq^{\rm ad}, \qquad x\longmapsto(\pi_V(x))_{V\in {\rm Ob}(\mathcal{C})}
\end{gather}
is a morphism of Hopf algebras, where $\pi_V\colon U_q^{\rm ad} \ra {\rm End}(V)$ is the representation associated to a module $V$ in $\mathcal{C}$. It is a theorem that this map is injective. From now on, let us extend the coefficient ring of the modules and morphisms in $\mathcal{C}$ to $\mc\big(q^{1/D}\big)$. Put \smash{$\mUq = \mUq^{\rm ad} \bigotimes_{\mc(q)} \mc\big(q^{1/D}\big)$}.
The map $\iota$ above extends to an embedding of $U_q$ in $\mUq$. The category $\mathcal{C}$, with coefficients extended to $\mc\big(q^{1/D}\big)$, is braided and ribbon; we postpone a discussion of that fact to Section \ref{intdual}, where it will be developed. As a consequence, we can regard $\mUq$ as a quasitriangular and ribbon Hopf algebra in a generalized sense (see \cite{BR}). The $R$-matrix of $\mUq$ is the family of morphisms
\[R = (R_{V,W})_{V,W\in {\rm Ob}(\mathcal{C})},\]
where $R_{V,W}\in {\rm End}(V\otimes W)$ is the endomorphism defined by the action of Drinfeld's universal $R$-matrix on $V\otimes W$. The ribbon element of $\mUq$ is defined similarly by Drinfeld's universal ribbon element. One defines the {\it categorical tensor product} \smash{$\mUq^{\hat{\otimes} 2}$} similarly as $\mUq$; in particular it contains all the infinite series of elements of $\mUq^{\otimes 2}$ having only a finite number of non-zero terms when evaluated on a given module $V\otimes W$ of $\mathcal{C}$. There is an expansion of $R$ as an infinite series in~\smash{$\mUq^{\hat{\otimes} 2}$}. Adapting Sweedler's coproduct notation $ \Delta(x)=\sum_{(x)}x_{(1)}\otimes x_{(2)}$, we find convenient to write this series as
\begin{equation}\label{Rcat}
R=\sum_{(R)}R_{(1)}\otimes R_{(2)}.
\end{equation}
We put $R^+:=R$, $R^- := (\sigma\circ R)^{-1}$ where $\sigma$ is the flip map $x\otimes y \mapsto y\otimes x$. We will not use any explicit formula of $R$, but the following factorization formula
\begin{equation}\label{Rmatfact}
R=\Theta \hat{R},
\end{equation}
where
\[\Theta = q^{\sum_{i,j=1}^m (B^{-1})_{ij} H_i\otimes H_j} \in \mUq^{\hat{\otimes} 2},\]
with $B\in M_m(\mathbb{Q})$ the matrix with entries $B_{ij}:=d_j^{-1}a_{ij}$, and
\[\hat{R}=\sum_{(\hat{R})} \hat{R}_{(1)} \otimes \hat{R}_{(2)}\in \mUq(\mathfrak{n}_+)\hat{\otimes} \mUq(\mathfrak{n}_-)\]
(see \cite[Section 3.2]{BR}, and for details, e.g., \cite[Theorem 8.3.9]{CP}, or \cite[Theorem 3.108]{VY}). If $x$,~$y$ are weight vectors of weights $\mu$, $\nu$ respectively, then $\Theta(x\otimes y) = q^{(\mu,\nu)} x\otimes y$. Moreover, $\hat{R}$ has weight $0$ for the adjoint action of $U_q(\mathfrak{h})$; that is, complementary components $\hat{R}_{(1)}$ and $\hat{R}_{(2)}$ have opposite weights.

Recall that we denote by $G$ the connected and simply-connected algebraic group with Lie algebra $\mathfrak{g}$. The {\it quantum function Hopf algebra} $\Oo_q=\Oo_q(G)$ is defined as the restricted dual of~$U_q^{\rm ad}$ with respect to the category $\mathcal{C}$, that is, the set of ${\mathbb C}(q)$-linear maps $f\colon U_q^{\rm ad}\ra \mc(q)$ such that ${\rm Ker}(f)$ contains a cofinite two sided ideal $I$ (i.e., such that $I\oplus M =U_q$ for some finite-dimensional vector space $M$), and $ \prod_{s=-r}^r (K_i - q_i^s) \in I$ for some $r\in \mathbb{N}$ and every $i$ (see, e.g., \cite[Chapter I.7]{BG}).

The space $\Oo_q$ is a Hopf algebra, with structure maps defined dually to those of $U_q^{\rm ad}$. We denote by $\star$ its product. The algebras $\Oo_q(T_G)$, $\Oo_q(U_\pm)$, $\Oo_q(B_\pm)$ are defined similarly, by replacing~$U_q^{\rm ad}$ with $U_q^{\rm ad}(\mathfrak{h})$, $U_q^{\rm ad}(\mathfrak{n}_\pm)$, $U_q^{\rm ad}(\mathfrak{b}_\pm)$, respectively. As a vector space, $\Oq$ is generated by the functionals $x\mapsto w(\pi_V(x)v)$, $x\in U_q^{\rm ad}$, for every object $V\in {\rm Ob}(\mathcal{C})$ and vectors $v\in V$, $w\in V^*$. Such functionals are called {\it matrix coefficients}. Because the morphism $\iota\colon U_q^{\rm ad}\to \mathbb{U}_q$ is injective (see~\eqref{defiota}), the Hopf duality pairing $\langle\cdot,\cdot\rangle \colon \Oq \times U_q^{\rm ad} \ra \mc(q)$ is non degenerate. By extending the coefficient ring from $\mc(q)$ to $\mc\big(q^{1/D}\big)$, we can uniquely extend it to a bilinear pairing
\[\langle\cdot,\cdot\rangle\colon\ \big(\Oq \bigotimes_{{\mathbb C}(q)} {\mathbb C}\big(q^{1/D}\big)\big)\times \mUq \ra \mc\big(q^{1/D}\big)\]
such that the following diagram is commutative:

\centerline{\xymatrix{
 \Oq \otimes U_q^{\rm ad} \ar[r]^{\langle\cdot,\cdot \rangle}\ar[d]_{{\rm id} \otimes \iota} & {\mathbb C}(q) \ar[d]\\
 \bigl(\Oq \bigotimes_{{\mathbb C}(q)} {\mathbb C}\big(q^{1/D}\big) \bigr)\otimes \mUq \ar[r]^{\;\;\;\;\;\;\;\ \;\ \;\;\;\ \langle\cdot,\cdot \rangle}& {\mathbb C}\big(q^{1/D}\big).
 }}

This pairing is defined by $\langle {}_Y\phi{}^w_v , (a_{X})\rangle=w(a_{Y}v)$
for every $(a_{X})\in \mUq$ and ${}_Y\phi{}^w_v\in \Oq$. It is non degenerate.

The maps
\begin{gather}\label{phipm}
\Phi^\pm\colon\ \Oo_q \longrightarrow U_q^{\rm cop},\qquad \alpha\longmapsto (\alpha \otimes {\rm id})\big(R^\pm\big)= \sum_{(R^\pm)} \big\langle \alpha , R_{(1)}^\pm\big\rangle R_{(2)}^\pm
\end{gather}
are well-defined morphisms of Hopf algebras. Here we stress that it is the simply-connected quantum group $U_q^{\rm cop}$ that is the range of $\Phi^\pm$. This will be explained with more details in Section~\ref{intdual}.

Let us make two simple observations, for future reference. Firstly, because $\Oo_q$ is spanned by the matrix coefficients of the objects of $\mathcal{C}$, and $\mathcal{C}$ is semisimple with simple objects the $U_q^{\rm ad}$-modules $V_\mu$, $\mu\in P_+$, there is a decomposition of $U_q$-bimodule
\begin{equation}\label{decompdirectOq}
 \Oo_q = \bigoplus_{\mu\in P_+} C(\mu),
 \end{equation}
 where $C(\mu) = V_\mu^* \otimes V_\mu$, the space of matrix coefficients of $V_\mu$, is endowed with the left action on the factor $V_\mu$ and the right action on $V_\mu^*$, and $\Oo_q$ has the left and right coregular actions $\lhd$ and~$\rhd$, defined by
\[x\rhd \alpha := \sum_{(\alpha)}\alpha_{(1)} \langle \alpha_{(2)},x\rangle,\qquad \alpha \lhd x := \sum_{(\alpha)} \langle \alpha_{(1)},x\rangle \alpha_{(2)}\]
for all $x\in U_q$ and $\alpha\in \Oq$. Here we recall that each $U_q^{\rm ad}$-module $V_\mu$ can be regarded as a~$U_q$-module, so the above expressions make sense. The decomposition \eqref{decompdirectOq} is the {\it Peter--Weyl} decomposition of $\Oo_q$. It will be refined in Section \ref{canbasemodqg}.

Moreover, the algebra $\Oo_q$ is generated by the matrix coefficients of the simple $U_q^{\rm ad}$-modules~$V_{\!\varpi_{\!k}}$ with highest weights the fundamental weights $\varpi_k$, $k=1,\dots,m$; see, e.g., \cite[Proposition~I.7.8]{BG} for a proof. This relies on the standard fact that, for any $\mu,\nu \in P_+$ we have a direct sum decomposition of modules (where $m(\lambda)\in \mn$)
\begin{equation}\label{decomprep}
V_\mu\otimes V_\nu = V_{\mu+\nu} \oplus \bigoplus_{\lambda < \mu+\nu} V_\lambda^{\oplus m(\lambda)}.
\end{equation}
In particular, this decomposition implies that, up to scalar multiples, there is a unique non-zero morphism $V_{\mu+\nu} \ra V_\mu\otimes V_\nu$, which is injective and splits. Dually, this means that, applying the product in $\Oo_q$ followed by the projection onto the subspace $C(\mu+\nu)$ we get a canonical projection map
\begin{equation}\label{projOq}
C(\mu)\otimes C(\nu) \ra C(\mu+\nu).
\end{equation}
The {\it loop algebra} $\Ll_{0,1} = \Ll_{0,1}(\mathfrak{g})$ is defined by twisting the product $\star$ of $\Oo_q$, keeping the same underlying linear space. The new product is equivariant with respect to the right coadjoint action ${\rm coad}^r$ of $U_q$, defined by
\[{\rm coad}^r(x)(\alpha) = \sum_{(x)}{S}(x_{(2)}) \rhd \alpha \lhd x_{(1)}\]
for all $x\in U_q$ and $\alpha\in \Oq$. By equivariant we mean that $\Ll_{0,1}$ is a $U_q$-module algebra. Let us spell out its product and equivariance property. Using the fact that $U_q$ can be regarded as a~subspace of $\mUq$, the actions $\lhd$ and $\rhd$ extend naturally to actions of $\mUq$, and the product of~$\Ll_{0,1}$ is expressed in terms of $\star$ by the formula (see \cite[Proposition 4.1]{BR}):
\begin{equation}\label{mmtilde}
\alpha\beta = \sum_{(R),(R)}(R_{(2')}{S}(R_{(2)}) \rhd \alpha) \star (R_{(1')}\rhd \beta \lhd R_{(1)}),
\end{equation}
where $ \sum_{(R)}R_{(1) }\otimes R_{(2)}$ and $\sum_{(R)}R_{(1') }\otimes R_{(2')}$ are expansions of two copies of $R\in \mUq^{\hat \otimes 2}$. Note that the sum in \eqref{mmtilde} has only a finite number of non-zero terms. By using that
\[R\Delta = \Delta^{\rm cop}R,\]
 this product can equivalently be expressed as
\begin{equation}\label{mmtildebis}
\alpha\beta = \sum_{(R),(R)}(\beta \lhd R_{(1)}R_{(1')}) \star (S(R_{(2)})\rhd \alpha \lhd R_{(2')}).
\end{equation}
This product gives $\Ll_{0,1}$ (like $\Oo_q$) a structure of $U_q$-module algebra for the actions $\rhd$, $\lhd$, but also for ${\rm coad}^r$ (which is not the case of $\Oo_q$). Spelling this out for ${\rm coad}^r$, this means
\[{\rm coad}^r(x)(\alpha\beta)=\sum_{(x)}{\rm coad}^r(x_{(1)})(\alpha){\rm coad}^r(x_{(2)})(\beta).\]
The relations between $\Oo_q$, $\Ll_{0,1}$ and $U_q$ are encoded by the map
\begin{equation}\label{phi1def}
\Phi_1\colon\ \Oo_q\longrightarrow \mUq, \qquad \alpha\longmapsto(\alpha \otimes {\rm id})(RR'),
\end{equation}
where $R' = \sigma\circ R$, and as usual $\sigma\colon x\otimes y \mapsto y\otimes x$. Note that
\begin{equation}\label{RSDmapdef}
\Phi_1 = m\circ \big(\Phi^+\otimes \big(S^{-1}\circ \Phi^-{}\big)\big)\circ \Delta.
\end{equation}
We call $\Phi_1$ the {\it RSD} map, for Drinfeld, Reshetikhin and Semenov-Tian-Shansky introduced it first (see \cite{Dr,Maj,RSTS}). It is a fundamental result of the theory (see \cite{Bau1,Ca,Jos}) that $\Phi_1$ affords an isomorphism of $U_q$-modules
$\Phi_1\colon \Oo_q \ra U_q^{\rm lf}$.
For full details on that result we refer to \cite[Section 3.12]{VY}. Here, $U_q^{\rm lf}$ is the set of \emph{locally finite} elements of $U_q$, endowed with the right adjoint action ${\rm ad}^r$ of $U_q$. It is defined by
\[U_q^{\rm lf} :=\{x\in U_q\mid {\rm rk}_{\mc(q)}({\rm ad}^r(U_q)(x)) < \infty\}\]
and
\[{\rm ad}^r(y)(x) = \sum_{(y)}{S}(y_{(1)}) x y_{(2)}\]
for every $x,y \in U_q$. The action ${\rm ad}^r$ gives in fact $U_q^{\rm lf}$ a structure of right $U_q$-module algebra. It is also a right coideal, that is \smash{$\Delta\big(U_q^{\rm lf}\big) \subset U_q^{\rm lf}\otimes U_q$}. Moreover, $\Phi_1$ affords an isomorphism of $U_q$-module algebras \smash{$
\Phi_1\colon \Ll_{0,1} \ra U_q^{\rm lf}$}.
It is a fact that $\Phi_1$ affords an isomorphism between the centers $\mathcal{Z}(\Ll_{0,1})$ of $\Ll_{0,1}$ and $\mathcal{Z}(U_q)$ of $U_q$ \cite[Proposition~6.24]{BR}. Since $\Phi_1$ is an isomorphism of $U_q$-modules and \smash{$\mathcal{Z}(U_q)=U_q^{U_q}$}, it follows that \smash{$\mathcal{Z}(\Ll_{0,1}) = \Ll_{0,1}^{U_q}$}.

Let us recall a few fundamental results about $U_q^{\rm lf}$ that we will meet again later. Denote by~${T\subset U_q}$ the multiplicative Abelian group formed by the elements $K_{\lambda}$, $\lambda \in P$, and by $T_2\subset T$ the subgroup formed by the elements $K_{\lambda}$, $\lambda \in 2P$. Consider the subset $T_{2-} \subset T_2$ formed by the elements $K_{-\lambda}$, $\lambda \in 2P_+$. Clearly, $T_2 = T_{2-}^{-1}T_{2-}$ and ${\rm Card}(T/T_2) = 2^m$.

\begin{teo} \label{JLteo1} \quad
\begin{itemize}\itemsep=0pt
\item[$(1)$] $U_q^{\rm lf} = \bigoplus_{\lambda \in 2P_+} {\rm ad}^r(U_q)(K_{-\lambda})$.
\item[$(2)$] $U_q = T_{2-}^{-1}U_q^{\rm lf}[T/T_{2}]$, so $U_q$ is a free $T_{2-}^{-1}U_q^{\rm lf}$-module of rank $2^m$.
\item[$(3)$] The ring $U_q^{\rm lf}$ is $($left and right$)$ Noetherian.
\end{itemize}
\end{teo}

These results were proved by Joseph--Letzter in \cite[Theorem~4.10]{JL2}, \cite[Theorem 6.4]{JL}, and \cite[Theorem~7.4.8]{Jos}, respectively (see also \cite[Sections~7.1.6, 7.1.13 and~7.1.25]{Jos}). For~(1) and~(3), we refer also to \cite[Theorems~3.113 and~3.137]{VY}, which provides simpler proofs. For instance, in the ${\mathfrak{sl}_2}$ case, we have
\[U_q({\mathfrak{sl}_2}) = U_q({\mathfrak{sl}_2})^{\rm lf}[K] \oplus U_q({\mathfrak{sl}_2})^{\rm lf}[K].L.\]
The actual values of $\Phi_1$ are complicated in general, however, there is a simple important one, that we describe now. Let $V_{-\lambda}$ be the type $1$ simple $U_q^{\rm ad}$-module of lowest weight $-\lambda\in -P_+$ (i.e., the highest weight \smash{$U_q^{\rm ad}$}-module $V_{-w_0(\lambda)}$ of highest weight $-w_0(\lambda)$, where $w_0$ is the longest element of the Weyl group; note that $-w_0$ permutes the simple roots). Let $v\in V_{-\lambda}$ be a lowest weight vector, and $v^*\in V_{-\lambda}^*$ be such that $v^*(v)=1$ and $v^*$ vanishes on a $U_q^{\rm ad}(\mathfrak{h})$-invariant complement of $v$. Define $\psi_{-\lambda}\in \Oo_q$ by $\langle \psi_{-\lambda},x\rangle = v^*(xv)$, $x\in U_q$. From the definition~\eqref{phi1def}, it is quite easy to see that
\begin{equation}\label{Phi1psi-}
\Phi_1(\psi_{-\lambda}) = K_{-2\lambda}.
\end{equation}
In particular, $\Phi_1(\psi_{-\rho}) = \ell^{-1}$, where as usual $\ell$ is the pivotal element of $U_q$.

\begin{Remark}\label{remsocle} Since $\Ll_{0,1}=\Oo_q$ as a vector space, we still denote by $C(\mu)$, $\mu\in P^+$, the linear subspace generated by the matrix coefficients of $V_\mu$, the $U_q^{\rm ad}$-module of type $1$ and highest weight~$\mu$. It can be proved (see \cite[Section~7.1.22]{Jos}, or \cite[p.~156]{VY}, where different conventions are used) that $\Phi_1$ yields an isomorphism of $U_q$-modules
\begin{equation}\label{restphi1}
\Phi_{1}\colon\ C(-w_0(\mu)) \ra {\rm ad}^r(U_q)(K_{-2\mu}).
\end{equation}
Therefore, the summands in (1) are finite-dimensional \smash{$U_q$}-modules, and the action ${\rm ad}^r$ is completely reducible on \smash{$U_q^{\rm lf}$}. In fact, \smash{$U_q^{\rm lf}$} is the socle of ${\rm ad}^r$ on $U_q$.
\end{Remark}

\begin{Remark}\label{factell} Because $\ell =\prod_{j=1}^m L_j^2$ and $\Phi_1(\psi_{-\rho}) = \ell^{-1}$, a natural question is the factorization of~$\psi_{-\rho}$ in $\Ll_{0,1}$ (see Corollary \ref{Phi1int2}). This question is considered in \cite{JW}, where ${\mathcal L}_{0,1}({\mathfrak g})$ for~${{\mathfrak g}=\mathfrak{gl}(r+1)}$ is analysed and quantum minors are extensively studied. Let us review here some of their results in relation with $\psi_{-\rho}$.

First note that for ${\mathfrak g}=\mathfrak{sl}(r+1)$ the irreducible representation $V_{-\rho}$ of lowest weight~${-\rho}$ is isomorphic to the representation of highest weight $V_{\rho}$ because $-w_0(\rho)=\rho$. By the Weyl formula, the dimension of this representation is
\[ \prod_{\alpha>0}\frac{(2\rho,\alpha)}{(\rho,\alpha)}=2^N.
\]
 In \cite{KS}, a presentation of~${U_q(\mathfrak{gl}(r\!+1))}$ is given, which differs from our presentation of ${U_q(\mathfrak{sl}(r\!+1))}$ only by its subalgebra~${U_q({\mathfrak h})}$, generated by $r+1$ elements ${\mathbb K}_1,\dots , {\mathbb K}_{r+1}$. The inclusion
 \[U_q(\mathfrak{sl}(r+1))\subset U_q(\mathfrak{gl}(r+1))
 \]
 is such that~${K_i={\mathbb K}_i^2{\mathbb K}_{i+1}^{-2}}$, $i=1,\dots,r$. The quantum minors, properly defined in \cite{JW}, of the matrix of matrix elements of the natural representation of $U_q(\mathfrak{gl}(r+1))$ are denoted $\det_q(A_{\geq k})$ for~${k=1,\dots ,r+1}$. We have $\det_q(A_{\geq 1})=1$ in the case of $\mathfrak{sl}(r+1)$. Then \cite{JW} proves that $\det_q(A_{\geq k})=({\mathbb K}_k\cdots{\mathbb K}_{r+1})^2$, and there exists an element $\mathbb{K}\in U_q(\mathfrak{gl}(r+1))$ such that
\[{\mathbb K}^{-2\rho}={\det}_q(A_{\geq 1})^{-r}{\det}_q(A_{\geq 2})\cdots {\det}_q(A_{\geq r+1}).\]
This has to be interpreted as $K_{-2\rho}=\Phi_1(\det_q(A_{\geq 2})\cdots \det_q(A_{\geq r+1}))$ in the case of $\mathfrak{sl}(r+1)$. As a result, this gives the equality
\[
\psi_{-\rho}={\det}_q(A_{\geq 2})\cdots {\det}_q(A_{\geq r+1}).
\]
\end{Remark}

The {\it $($quantum$)$ graph algebra} $\Ll_{0,n} = \Ll_{0,n}(\mathfrak{g})$ is the braided tensor product of $n$ copies of~$\Ll_{0,1}$ (considered as a $U_q$-module algebra). As a linear space and $U_q$-bimodule with actions $\lhd$ and~$\rhd$, it coincides with $\Ll_{0,1}^{\otimes n}$, and thus with $\Oo_q^{\otimes n}$. It is also a right $U_q$-module algebra, with the following action of $U_q$ (extending ${\rm coad}^r$ on $\Ll_{0,1}$):
 \begin{gather}
{\rm coad}_n^r(y)\big(\alpha^{(1)}\otimes \dots \otimes \alpha^{(n)}\big) =
\sum_{(y)} {\rm coad}^r(y_{(1)})\big(\alpha^{(1)}\big) \otimes \dots \otimes {\rm coad}^r(y_{(n)})\big(\alpha^{(n)}\big)\label{actionprod}
\end{gather}
for all $y \in U_q$ and $\alpha^{(1)}\otimes \dots \otimes \alpha^{(n)} \in \Ll_{0,n}$. The product of $\Ll_{0,n}$ can be expressed as follows. For every $1\leq a\leq n$, define $\mathfrak{i}_a\colon \Ll_{0,1}\ra \Ll_{0,n}$ by $\mathfrak{i}_a(x) = 1^{\otimes (a-1)}\otimes x \otimes 1^{\otimes (n-a)}$; $\mathfrak{i}_a$ is an embedding of~$U_q$-module algebras. We will use the notations
\begin{equation}\label{plgtL01a}
\Ll_{0,n}^{(a)}:= {\rm Im}(\mathfrak{i}_a),\qquad (\alpha)^{(a)} := \mathfrak{i}_a(\alpha).
\end{equation}
Take $(\alpha)^{(a)},(\alpha')^{(a)} \in \Ll_{0,n}^{(a)}$ and $(\beta)^{(b)},(\beta')^{(b)} \in \Ll_{0,n}^{(b)}$ with $a<b$. Then the product of $\Ll_{0,n}$ is given by the following formula (see \cite[Section 6]{BR}):
\begin{align}
\big((\alpha)^{(a)}\otimes (\beta)^{(b)}\big) \big((\alpha')^{(a)} \otimes (\beta')^{(b)}\big) ={}& \sum_{(R^1),\dots, (R^4)} \big(\alpha \big(S\big(R^3_{(1)}R^4_{(1)}\big)\rhd \alpha' \lhd R^1_{(1)}R^2_{(1)}\big)\big)^{(a)}\nonumber\\
 &
\otimes \big( \big( S\big(R^1_{(2)}R^3_{(2)}\big) \rhd \beta \lhd R^2_{(2)}R^4_{(2)}\big)\beta'\big)^{(b)},\label{prodL0nform}
\end{align}
where $R^i = \sum_{(R^i)}R_{(1)}^i\otimes R_{(2)}^i$, $i\in \{1,2,3,4\}$, are expansions of four copies of $R\in \mUq^{\hat \otimes 2}$, and on the right-hand side the product is componentwise that of $\Ll_{0,1}$. Later we will use the fact that the product of $\Ll_{0,n}$ is obtained from the standard (componentwise) product of $\Ll_{0,1}^{\otimes n}$ by a process that may be inverted. Indeed, \eqref{prodL0nform} can be rewritten as
\begin{equation}\label{prodL0nform2}
\big((\alpha)^{(a)}\otimes (\beta)^{(b)}\big) \big((\alpha')^{(a)} \otimes (\beta')^{(b)}\big) = \sum_{(F)} (\alpha)^{(a)}\big((\alpha')^{(a)}\cdot F_{(2)}\big)\otimes \big((\beta)^{(b)}\cdot F_{(1)}\big)(\beta')^{(b)},
\end{equation}
where $F=\sum_{(F)} F_{(1)}\otimes F_{(2)} := (\Delta \otimes \Delta)(R')$, and the symbol ``$\cdot$'' stands for the right action of \smash{$\mUq^{\hat \otimes 2}$} on $\Ll_{0,1}$ that may be read from \eqref{prodL0nform}. The tensor $F$ is known as a twist. Then, by replacing $F$ with its inverse $\bar{F} = (\Delta \otimes \Delta)\big(R'{}^{-1}\big)$, one can express the product of
$\Ll_{0,1}^{\otimes n}$ in terms of the product of $\Ll_{0,n}$ by
\begin{equation}\label{prodL0nform3}
(\alpha)^{(a)}(\alpha')^{(a)}\otimes (\beta)^{(b)}(\beta')^{(b)} = \sum_{(\bar{F})}\big((\alpha)^{(a)}\otimes
 \big((\beta)^{(b)}\cdot \bar{F}_{(1)}\big)\big) \big( \big((\alpha')^{(a)}\cdot \bar{F}_{(2)}\big)\otimes (\beta')^{(b)}\big).
\end{equation}
 We call {\it quantum moduli algebra} and denote by $\Mm_{0,n} = {\mathcal L}_{0,n}^{U_q}$ the subalgebra of $\Ll_{0,n}$ formed by the $U_q$-invariant elements.

The map $\Phi_1$ can be extended to $\Ll_{0,n}$ as follows. Consider the following action of $U_q$ on the tensor product algebra $U_q^{{\otimes} n}$, which extends ${\rm ad}^r$ on $U_q$:
\[
{\rm ad}_n^r(y)(x) = \sum_{(y)}\Delta^{(n)}({S}(y_{(1)})) x \Delta^{(n)}(y_{(2)})
\]
for all $y \in U_q$, $x \in U_q^{{\otimes} n}$. This action gives $U_q^{{\otimes} n}$ a structure of right $U_q$-module algebra. In~\cite{A}, Alekseev introduced a morphism of $U_q$-module algebras $\Phi_n\colon \Ll_{0,n} \ra U_q^{\otimes n}$ which extends $\Phi_1$. In~\cite[Proposition~6.7]{BR}, we showed that $\Phi_n$ affords isomorphisms
\begin{equation}\label{Alekseevmap}
\Phi_n\colon\ \Ll_{0,n} \ra \big(U_q^{\otimes n}\big)^{\rm lf},\qquad \Phi_n\colon\ \Mm_{0,n}\rightarrow \big(U_q^{{\otimes} n}\big)^{U_q},
\end{equation}
where $\big(U_q^{\otimes n}\big)^{\rm lf}$ is the set of ${\rm ad}_n^r$-locally finite elements of $U_q^{{\otimes} n}$. We call $\Phi_n$ the {\it Alekseev map}; we do not recall here the definition of $\Phi_n$, for we will not use it. It is a key argument of the proof of~\eqref{Alekseevmap} that the set of locally finite elements of $U_q^{{\otimes} n}$ for $({\rm ad}^r)^{\otimes n}\circ \Delta^{(n)}$ coincides with $\big(U_q^{\rm lf}\big)^{\otimes n}$; this follows from the main result of \cite{KLNY}. Using that the map
\begin{equation}\label{Phipsi}
\psi_n = \Phi_n \circ \big(\Phi_1^{-1}\big)^{\otimes n}\colon\ \big(U_q^{\rm lf}\big)^{\otimes n}\ra \big(U_q^{\otimes n}\big)^{\rm lf}
\end{equation}
intertwines the actions $({\rm ad}^r)^{\otimes n}\circ \Delta^{(n-1)}$ and ${\rm ad}_n^r$, we deduced that ${\rm Im}(\Phi_n) = \big(U_q^{\otimes n}\big)^{\rm lf}$.

\begin{Remark}We have \smash{$\big(U_q^{\rm lf}\big)^{\otimes n} \!\ne\! \big(U_q^{\otimes n}\big)^{\rm lf}$} and in fact there is not even an inclusion.~Indeed, let $\Omega=\big(q-q^{-1}\big)^2FE+qK+q^{-1}K^{-1}$ be the Casimir element of $U_q({\mathfrak{sl}_2})$. We trivially have~\smash{$\Delta(\Omega)\in \big(U_q^{\otimes 2}\big)^{\rm lf}$} but
\[\Delta(\Omega)=\big(q-q^{-1}\big)^2\big(K^{-1}E\otimes FK+F\otimes E\big)+
\Omega\otimes K+ K^{-1}\otimes \Omega-\big(q+q^{-1}\big)K^{-1}\otimes K\]
and therefore \smash{$\Delta(\Omega)\notin \big(U_q^{\rm lf}\big)^{\otimes 2}$}, since \smash{$K\notin U_q^{\rm lf}$} (see, e.g., Theorem \ref{JLteo1}\,(2)). This reflects the fact that \smash{$U_q^{\rm lf}$} is only a right coideal of $U_q$ (and not a subcoalgebra).
\end{Remark}
As in Remark \ref{remsocle}, denote by $C(\mu)$, $\mu\in P^+$, the linear subspace of $\Ll_{0,1}$ generated by the matrix coefficients of $V_\mu$. For every tuple $[\mu] = (\mu_1,\dots,\mu_n)\in P_+^n$ put
\begin{equation}\label{defCmu}
C([\mu]) =C(\mu_1)\otimes \dots \otimes C(\mu_n).
\end{equation}
Then \smash{$ {\mathcal L}_{0,n}=\bigoplus_{[\mu]\in P_+^n} C({[\mu]})$}. Each space $C({[\mu]})$ is a finite-dimensional $U_q$-module under the action ${\rm coad}^r_n$, whence it is completely reducible. Therefore, $
\Ll_{0,n} = \Mm_{0,n} \oplus I$ as $U_q$-modules, where $I$ is the sum of nontrivial isotypical components of $\Ll_{0,n}$. The $\mc(q)$-linear projection map
\begin{equation}\label{Reynoldsdef}
\Rr\colon\ {\mathcal L}_{0,n}\ra \Mm_{0,n}, \qquad {\rm Ker}(\Rr)=I
\end{equation}
is called the {\it Reynolds operator}. For all $\alpha\in \Mm_{0,n}$, $\beta \in \Ll_{0,n}$ it satisfies $
\Rr(\alpha \beta) = \alpha \Rr(\beta)$.
This property will be crucial in the sequel, so let us recall a (classical) proof of it. We can write~${\beta = \Rr(\beta) + \gamma}$ with $\gamma \in I$, and then we have to show $\alpha \gamma \in I$. We can reduce to the case where $\gamma$ is contained in a simple summand $V$ of $I$. Multiplication by the invariant element $\alpha$ yields a surjective map $V \ra \alpha V$, which is a morphism of $U_q$-modules. Since $V$ is simple, it is either the $0$ map, or an isomorphism. In either cases it follows $\alpha V \subset I$ (in fact the first case cannot happen, for $\Ll_{0,n}$ has no nontrivial zero divisors, as explained after \eqref{Zinv}).

We can formulate the Reynolds operator in the following way. Recall that $\Oo_q$ has a unique left (or right, or $2$-sided) Haar integral, that is a linear map $h\colon \Oo_q\ra \mc(q)$ such that
\[
h(1)=1 \qquad \text{and}\qquad ({\rm id}\otimes h)\Delta(\alpha) = h(\alpha)1,\qquad \forall \alpha\in \Oo_q.
\]
(See, e.g., \cite[Proposition 13.3.6]{CP}.) It vanishes on all matrix coefficients except the one of the trivial representation, to which it gives the value $1$. Denote by $\Delta_\Ll\colon {\mathcal L}_{0,n} \ra {\mathcal L}_{0,n}\otimes \Oo_q$ the right coaction dual to the action ${\rm coad}_n^r$ of $U_q$ on ${\mathcal L}_{0,n}$. Then, in analogy with the formula of the averaging operator $\mathcal{C}^\infty(G) \ra \mathcal{C}^\infty(G)^G$, $f\ra [f] = \int_G f\big(g^{-1} \cdot g\big) {\rm d}\mu(g)$, for a locally compact group $G$ with Haar measure ${\rm d}\mu(g)$, it is straightforward that
\begin{equation}\label{formuleRh}
\Rr = ({\rm id}\otimes h)\Delta_\Ll .
\end{equation}
Note that the complete reducibility of $\Ll_{0,n}$ discussed after \eqref{defCmu} follows also from Theorem~\ref{JLteo1}\,(1), since by \eqref{Phipsi} we have an isomorphism of $U_q$-modules
\[
{\mathcal L}_{0,n} \stackrel{\Phi_n}{\longrightarrow} \big(U_q({\mathfrak g})^{\otimes n}\big)^{\rm lf} \stackrel{\psi_n^{-1}}{\longrightarrow}U_q^{\rm lf}({\mathfrak g})^{\otimes n},
\]
where ${\rm lf}$ means respectively locally finite for the action ${\rm ad}_n^r$ of $U_q({\mathfrak g})$ on $U_q({\mathfrak g})^{\otimes n}$, and locally finite for the action ${\rm ad}^r$ of $U_q({\mathfrak g})$ on $U_q({\mathfrak g})$. An explicit basis of $\Mm_{0,n}$ is described in \cite[Pro\-po\-sition~6.22]{BR}.

Finally, let us point out here two important consequences of~\eqref{Alekseevmap}. First, $\Phi_n$ yields isomorphisms between centers, $\mathcal{Z}(\Ll_{0,n}) \cong \mathcal{Z}(U_q)^{\otimes n}$ and \smash{$\mathcal{Z}\big(\Ll_{0,n}^{U_q}\big) \cong \mathcal{Z}\big(\big(U_{q}^{\otimes n}\big)^{U_q}\big)$}, where one can show that \cite[Lemma 6.29]{BR}
\begin{equation}\label{Zinv}
\mathcal{Z}\big(\big(U_{q}^{\otimes n}\big)^{U_q}\big)\cong \Delta^{(n)}(\mathcal{Z}(U_q))\bigotimes_{{\mathbb C}(q)}\mathcal{Z}(U_q)^{\otimes n}.
\end{equation}
Second, $\Ll_{0,n}$ (and therefore $\Mm_{0,n}$) has no nontrivial zero divisors because of the isomorphisms \smash{$\Phi_n\colon \Ll_{0,n} \ra \big(U_q^{\otimes n}\big)^{\rm lf}\subset U_q^{\otimes n}$} and $U_q^{\otimes n}\cong U_q\big(\mathfrak{g}^{\oplus n}\big)$, and the fact that $U_q\big(\mathfrak{g}^{\oplus n}\big)$ has no nontrivial zero divisors (proved, e.g., in \cite{DC-K}).

\subsection{Integral forms and specializations}\label{INTEG} Let $A=\mc\big[q,q^{-1}\big]$. We call integral form of a (Hopf) $\mc(q)$-algebra $H$ a (Hopf) $A$-subalgebra~${}_AH$ such that the canonical map ${}_AH\bigotimes_A \mc(q) \ra H$ is an isomorphism. Note that the standard notion of integral form of $\mc(q)$-algebra uses $\mz\big[q,q^{-1}\big]$ instead of $\mc\big[q,q^{-1}\big]$; our choice is made for simplicity ($\mc\big[q,q^{-1}\big]$ is a principal ideal domain, whereas $\mz\big[q,q^{-1}\big]$ is not).

\subsubsection{Definitions}\label{defintform} The {\it unrestricted} integral form of $U_q$ is the $A$-subalgebra $U_A = U_A(\mathfrak{g})$ introduced by De Concini--Kac--Procesi in \cite[Section 12]{DC-K-P2} (and in a differently normalized form in \cite{DC-K,DC-K-P1}). It is the smallest $A$-subalgebra of $U_q$ which contains the elements ($i=1,\dots,m$)
\begin{equation}\label{normgen0}
\bar{E}_i = \big(q_i-q_i^{-1}\big)E_i,\qquad \bar{F}_i = \big(q_i-q_i^{-1}\big)F_i,\qquad L_i,\qquad L_i^{-1}
\end{equation}
and is stable under the action of $\mathcal{B}(\mathfrak{g})$ given by the Lusztig automorphisms (see \eqref{rootvectdef}). Recall the root vectors $E_{\beta_k}$, $F_{\beta_k}$ defined in \eqref{rootvectdef}. Let us put $q_\beta := q^{(\beta,\beta)/2}$. The algebra $U_A$ is a free $A$-module with basis the monomials $\bar{E}_{\beta_1}^{p_1}\cdots \bar{E}_{\beta_N}^{p_N}K_\lambda \bar{F}_{\beta_N}^{n_N}\cdots \bar{F}_{\beta_1}^{n_1}$, where $\lambda \in P$ and we set
\[\bar{E}_{\beta_k}= \big(q_{\beta_k}-q_{\beta_k}^{-1}\big)E_{\beta_k},\qquad \bar{F}_{\beta_k}= \big(q_{\beta_k}-q_{\beta_k}^{-1}\big)F_{\beta_k}.\]
We denote \smash{$U_A^{\rm lf} := U_A \cap U_q^{\rm lf}$}.
The unrestricted integral form of \smash{$U_q^{\rm ad}$} is defined similarly, as the smallest $A$-subalgebra $U_A^{\rm ad} \subset U_A$ which contains the elements $\bar{E}_i$, $\bar{F}_i$ and $K_i^{\pm 1}$, for $i=1,\dots,m$, and is stable under the Lusztig action of $\mathcal{B}(\mathfrak{g})$.

For $\beta$ a positive root, we define the divided powers
\[E_{\beta}^{(k)} = \frac{E_\beta^k}{[k]_{q_\beta}!} , \qquad F_{\beta}^{(k)} = \frac{F_\beta^k}{[k]_{q_\beta}!}, \qquad k\in \mn.\]
The Lusztig {\it restricted} integral form of $U_q^{\rm ad}$ \cite{Lusztig2,Lusztig} (see also \cite[Chapter 9.3]{CP}) is the $A$-sub\-algebra~$U_A^{\rm res}$ generated by the elements ($i=1,\dots,m$, $k\in \mn^*$)
\[E_i^{(k)}= \frac{E_i^k}{[k]_{q_i}!},\qquad F_i^{(k)}= \frac{F_i^k}{[k]_{q_i}!} ,\qquad K_i ,\qquad K_i^{-1}.\]
The algebra $U_A^{\rm res}$ is a free $A$-module with Poincar\'e--Birkhoff--Witt (PBW) basis
\[
E_{\beta_1}^{(p_1)}\cdots E_{\beta_N}^{(p_N)}\prod_{i=1}^m K_i^{\sigma_i}[ K_i; t_i]_{q_i} F_{\beta_N}^{(n_N)}\cdots F_{\beta_1}^{(n_1)},
\]
where $\sigma_i\in \{0,1\}$, $n_i,p_i, t_i\in \mn$, and we set $[ K_i; 0]_{q_i}:=1$ and
\[[ K_i; t ]_{q_i}= \prod_{s=1}^t \frac{K_iq_i^{-s+1}-K_i^{-1}q_i^{s-1}}{q_i^s-q_i^{-s}}.\]
The integral forms $U_A(\mathfrak{h})$, $U_A(\mathfrak{b}_\pm)$ and $U_A^{\rm res}(\mathfrak{h})$, $U_A^{\rm res}(\mathfrak{b}_\pm)$ associated to the subalgebras $\mathfrak{h}$,~$\mathfrak{b}_\pm \subset \mathfrak{g}$ are the subalgebras of $U_A$ and $U_A^{\rm res}$, respectively, defined in the obvious way. For instance, the~``Cartan'' subalgebra $U_A^{\rm res}(\mathfrak{h}) = U_q(\mathfrak{h})\cap U_A^{\rm res}$ is generated as a $A$-module by the elements $\prod_{i=1}^m K_i^{\sigma_i}[ K_i; t_i]_{q_i}$.

Denote by $\mathcal{C}_A$ the category of $U_A^{\rm res}$-modules of {\it type $1$}, i.e., free $A$-modules of finite rank which have a basis where the elements $K_i$ act diagonally with eigenvalues of the form $q_i^k$, $k\in \mz$ (in general, finiteness of the rank imposes eigenvalues of the form $\pm q_i^k$, $k\in \mz$). The category~$\mathcal{C}_A$ is a~rigid and tensor category. It is not semisimple, and this makes the study of $\mathcal{C}_A$ a~complicated task; for this, see \cite{BR}, and Section \ref{canbasemodqg} below. Every type $1$ finite-dimensional simple $U_q$-module $V_\mu$,~$\mu \in P_+$, has a $U_A^{\rm res}$-invariant full $A$-sublattice, that we denote by ${}_AV_\mu$. These~$U_A^{\rm res}$-modules form the simple objects of $\mathcal{C}_A$. Moreover, $\mathcal{C}_A \otimes \mc\big[q^{1/D},q^{-1/D}\big]$ is a ribbon category (see Section~\ref{intdual}).

The {\it integral} quantum function Hopf algebra $\Oo_A=\Oo_A(G)$ is the (type $1$) restricted dual of~$U_A^{\rm res}$, that is, the $A$-span of the matrix coefficients $x\mapsto v^i(\pi_V(x)v_i)$, $x\in U_A^{\rm res}$, for every module~$V$ in $\mathcal{C}_A$, where $(v_i)$ is an $A$-basis of $V$ and $\big(v^i\big)$ the dual $A$-basis of the dual module~$V^*$ (compare with the definition of $\Oo_q$). We can also regard $\Oo_A$ as the set of $A$-linear maps~${f\colon U_A^{\rm res}\ra A}$ such that ${\rm Ker}(f)$ contains a cofinite two sided ideal $I$, and $ \prod_{s=-r}^r (K_i - q_i^s) \in I$ for some $r\in \mathbb{N}$ and every $i$. Because of the inclusions of $U_A^{\rm res}(\mathfrak{h})$, $U_A^{\rm res}(\mathfrak{n}_\pm)$, $U_A^{\rm res}(\mathfrak{b}_\pm)$ in $U_A^{\rm res}$, there are Hopf epimorphisms from $\Oo_A$ to the $A$-duals of these subalgebras, that we denote by $\Oo_A(T_G)$, $\Oo_A(U_\pm)$ and $\Oo_A(B_\pm)$, respectively.

The algebra $\Oo_A$ has been introduced by Lusztig in \cite{Lusztig2,Lusztig}. It is an integral form of $\Oo_q$, so $\Oo_q = \Oo_A \bigotimes_A \mc(q)$.

$\Oo_A$ is also the restricted dual of the integral form $\Gamma=\Gamma(\mathfrak{g})$ of \smash{$U_q^{\rm ad}$} introduced by De Concini--Lyubashenko in \cite[Sections 2 and~3]{DC-L}; $\Gamma$ is the $A$-subalgebra of \smash{$U_q^{\rm ad}$} generated by the elements ($i=1,\dots,m$)
\[E_i^{(k)}= \frac{E_i^k}{[k]_{q_i}!} ,\qquad F_i^{(k)}= \frac{F_i^k}{[k]_{q_i}!} ,\qquad ( K_i; t )_{q_i}= \prod_{s=1}^t \frac{K_iq_i^{-s+1}-1}{q_i^s-1},\qquad K_i^{-1},\]
where $k\in \mn$, $t\in \mn$ (setting $( K_i; 0)_{q_i}=1$ by convention). Note that the definition of $\Gamma$ is less symmetric than that of $U_A^{\rm res}$. However, $\Gamma$ contains the elements $K_i$, and the commutation relations between the generators \smash{$E_i^{(k)}$}, \smash{$F_i^{(k)}$} imply that the symmetrized elements $[ K_i; t ]_{q_i}$ belong to $\Gamma$. In fact, let us denote $\Gamma(\mathfrak{h}) = U_q(\mathfrak{h}) \cap \Gamma$ and $\Gamma(\mathfrak{b}_\pm) = U_q(\mathfrak{b}_\pm) \cap \Gamma$. It is proved in \cite[Theorem 3.1]{DC-L} that $\Gamma(\mathfrak{h})$ contains $U_A^{\rm res}(\mathfrak{h})$ and that the elements \smash{$ \prod_{i=1}^m K_i^{-\sigma(t_i)}(K_i; t_i)_{q_i}$, $t_i\in \mn$}, where $\sigma(t)$ is the integer part of $t/2$, is an $A$-basis of $\Gamma(\mathfrak{h})$. A PBW basis of $\Gamma$ is formed by the monomials
\[
E_{\beta_1}^{(p_1)}\cdots E_{\beta_N}^{(p_N)}\prod_{i=1}^m K_i^{-\sigma(t_i)}(K_i; t_i)_{q_i} F_{\beta_N}^{(n_N)}\cdots F_{\beta_1}^{(n_1)}.
\]

The inclusion $U_A^{\rm res}\subset \Gamma$ is strict, for the elements $(K_i; t)_{q_i}$, $t\ne 0$, do not belong to $U_A^{\rm res}$. However, the restriction functor $\mathcal{C}_\Gamma \ra\mathcal{C}_A$ is obviously an equivalence, where $\mathcal{C}_\Gamma$ is the category of $\Gamma$-modules of {\it type $1$}, i.e., free $A$-modules of finite rank which have a basis where the elements~$K_i$ act diagonally with eigenvalues of the form $q_i^k$, $k\in \mz$. Therefore, we can identify the two categories, and $\Oo_A$ with the (type $1$) restricted dual of $\Gamma$. We will thus consider the $U_A^{\rm res}$-modules~${}_AV_\mu$, $\mu\in P_+$, equally as $\Gamma$-modules. We will sometimes use $\Gamma$ instead of $U_A^{\rm res}$ in order to make direct the connection with results of De Concini--Lyubashenko about the integral pairings $\pi_A^\pm$ considered in Section~\ref{intdual}.

The {\it integral} form $\Ll_{0,1}^A$ of $\Ll_{0,1}$ is defined as the $U_A^{\rm res}$-module $\Oo_A$ endowed with the product of~$\Ll_{0,1}$. The {\it integral} form \smash{$\Ll_{0,n}^A$} of $\Ll_{0,n}$ is the braided tensor product of $n$ copies of \smash{$\Ll_{0,1}^A$}; in particular, \smash{$\Ll_{0,n}^A = \Oo_A^{\otimes n}$} as $U_A^{\rm res}$-modules. That the products of $\Ll_{0,1}$ and $\Ll_{0,n}$ are well defined over $A$ was shown in \cite[Proposition~6.9]{BR}.

The {\it integral} quantum moduli algebra is
\[\Mm_{0,n}^A := \big(\Ll_{0,n}^A\big)^{U_A^{\rm res}} = \big(\Ll_{0,n}^A\big)^{U_A}.\]

Finally, given $q=\e\in \mc^{\times}$ we define the {\it specializations} $U_\e$, $\Gamma_\e$, $\Oo_\e$, $\Ll_{0,n}^{\e}$ and $\Mm_{0,n}^{A,\e}$ as the~$\mc$-algebras obtained by tensoring $U_A$, $\Gamma$, $\Oo_A$, \smash{$\Ll_{0,n}^A$} and \smash{$\Mm_{0,n}^A$} respectively with $\mc_\e$, the $A$-module~$\mc$ where $q$ acts by multiplication by $\e$. Each one can also be defined as the quotient by the ideal generated by $q-\e$. We find convenient to use the notations
\begin{equation}\label{notambig0}
\big(U_A^{\otimes n}\big)^{U_A}_\e:= \big(U_A^{\otimes n}\big)^{U_A}\bigotimes_A \mc_\e,\qquad \big(U^{\otimes n}\big)^{\rm lf}_\e:= \big(U_A^{\otimes n}\big)^{\rm lf}\bigotimes_A \mc_\e .
\end{equation}
Let us stress here that when $\e$ is a root of unity, taking the locally finite part and taking the specialization at $\e$ are non commuting operations. Indeed, as shown by Theorem \ref{DCKteo1} below, $U_\e$~is finite over $\mathcal{Z}_0(U_\e)$ and therefore all its elements are locally finite for ${\rm ad}^r$; on another hand~\smash{${U_\e^{\rm lf} = U_A^{\rm lf}\bigotimes_A \mc_\e}$} does not contain the elements~$L_i$.

Similarly, taking invariants and taking the specialization at $\e$ are non commuting operations when $\e$ is a root of unity: indeed, it is easily checked that in this case \smash{$\big(U_A^{\otimes n}\big)^{U_A}_\e$} and \smash{$\big(U_\e^{\otimes n}\big)^{U_\e}$}, or \smash{$\Mm_{0,n}^{A,\e} = \Mm_{0,n}^A \bigotimes_A \mc_\e$} and \smash{$\big(\Ll_{0,n}^\e\big)^{U_\e}$}, are distinct spaces. When $\e$ is a root of unity, we will not consider the algebras $\Mm_{0,n}^{A,\e}$ in this paper.

Arguments similar to those mentioned at the end of Section \ref{ALGdef} imply that the algebras $\Ll_{0,n}^A$,~$\Mm_{0,n}^A$ and $\Ll_{0,n}^\ep$, \smash{$\Mm_{0,n}^{A,\ep}$}, $\ep\in \mc^\times$, have no nontrivial zero divisors (see \cite[Propositions~6.11 and~6.30]{BR}).

\subsubsection{Canonical bases and modified quantum groups} \label{canbasemodqg} Because the category $\mathcal{C}_A$ is not semisimple, it is not clear from the above definition of $\Oo_A$ whether or not it is a finitely generated algebra, if $\Mm_{0,n}^A$ is a direct summand of the $A$-module $\Ll_{0,n}^A$, or if the projection map \eqref{projOq} may be refined to a morphism between underlying $A$-modules.

Such properties, using the appropriate formalism developed by Kashiwara--Lusztig, indeed hold true, and will play a key role later. We state them precisely in Proposition \ref{genOA}, Theorem~\ref{teoLuzstigOAinv} and Proposition \ref{teoLuzstigOAscinde}. These results are consequences of the existence of an $A$-basis of $\Oo_A$ with favourable properties, which implies in particular that $\Oo_A$ is a free $A$-module. In order to introduce this $A$-basis it is necessary to consider a variant of $U_q^{\rm ad}$ introduced by Lusztig~\cite{Lusztig}, called {\it modified quantum group}, and use the Kashiwara--Lusztig theory of canonical bases \cite{Kashiwara1,Kashiwara2,Kashiwara3,Lusztig}. We are going to recall the background material step by step.

The Lusztig {\it modified quantum group} is the $\mc(q)$-algebra $\dot{\mathbf{U}}$ obtained by replacing $U_q^{\rm ad}(\mathfrak{h})$ with the direct sum of infinitely many one-dimensional algebras, generated by orthogonal idempotents~$1_\lambda$ indexed by the elements $\lambda$ of the weight lattice $P$ \cite[Chapter~23]{Lusztig}. Namely, as a vector space \smash{$ \dot{\mathbf{U}} = \bigoplus_{\lambda',\lambda''\in P} {}_{\lambda'}\dot{\mathbf{U}}{}_{\lambda''}$}, where
\[{}_{\lambda'}\dot{\mathbf{U}}{}_{\lambda''} = U_q^{\rm ad}\bigg/\bigg(\sum_{\alpha\in Q} \big(K_\alpha -q^{(\alpha,\lambda')}\big)U_q^{\rm ad}+ \sum_{\alpha\in Q} U_q^{\rm ad}\big(K_\alpha -q^{(\alpha,\lambda'')}\big)\bigg).\]
Denote by $\pi_{\lambda',\lambda''}\colon U_q^{\rm ad} \ra {}_{\lambda'}\dot{\mathbf{U}}{}_{\lambda''}$ the canonical projection. The product of $\dot{\mathbf{U}}$ is given by \smash{$\pi_{\lambda_1',\lambda_1''}(s)\pi_{\lambda_2',\lambda_2''}(t) = \pi_{\lambda_1',\lambda_2''}(st)$} if $\lambda_1''=\lambda_2'$ and zero otherwise. Set $1_\lambda := \pi_{\lambda,\lambda}(1)$. The algebra \smash{$\dot{\mathbf{U}}$} has not unit, but the family $(1_\lambda)_{\lambda\in P}$ can be regarded as a substitute of it. Denote by $\Delta$ the collection of maps
\[
\Delta_{\lambda'_1,\lambda'_2,\lambda''_1,\lambda_2''} \colon {}_{\lambda'_1+\lambda'_2} \dot{\mathbf{U}} {}_{\lambda''_1+\lambda''_2} \ra {}_{\lambda'_1} \dot{\mathbf{U}} {}_{\lambda''_1}\otimes {}_{\lambda'_2} \dot{\mathbf{U}} {}_{\lambda''_2}
\]
 such that
\begin{equation}\label{defcoprod}
\Delta_{\lambda'_1,\lambda'_2,\lambda''_1,\lambda''_2} \pi_{\lambda'_1+\lambda'_2,\lambda''_1+\lambda''_2} = (\pi_{\lambda'_1,\lambda''_1}\otimes \pi_{\lambda'_2,\lambda''_2})\Delta_{U_q^{\rm ad}},
\end{equation}
where \smash{$\Delta_{U_q^{\rm ad}}$} is the coproduct of \smash{$U_q^{\rm ad}$}. We can regard $\Delta$ as a (categorically completed) coproduct \smash{$\Delta \colon \dot{\mathbf{U}} \ra \dot{\mathbf{U}} {}^{\hat{\otimes}2}$}. There is a natural structure of $U_q^{\rm ad}$-bimodule on $\dot{\mathbf{U}}$, defined by
\begin{equation}\label{bimodstr}
t'\pi_{\lambda',\lambda''}(s)t'' = \pi_{\lambda'+\nu',\lambda''-\nu''}(t'st'')
\end{equation}
for all \smash{$s\in U_q^{\rm ad}$} and all elements \smash{$t',t''\in U_q^{\rm ad}$} of respective weights $\nu'$, $\nu''$. This structure affords a triangular decomposition of $\dot{\mathbf{U}}$: given bases $\{b^\pm\}$ of \smash{$U_q^{\rm ad}(\mathfrak{n}_\pm)$}, the set of elements $b^+1_\lambda b^-$ (or~$b^-1_\lambda b^+$, or $b^+b^-1_\lambda$), where $\lambda\in P$, is a basis of $\dot{\mathbf{U}}$.

Given any \smash{$U_q^{\rm ad}$}-module $X$ of type $1$, and any weight subspace $X^\lambda\subset X$ of weight $\lambda\in P$, one can define the action of an element $u 1_\lambda \in \dot{\mathbf{U}}$, \smash{$u\in U_q^{\rm ad}$}, on $X$ as the projection onto $X^\lambda$ followed by the action of $u$. This way, one can identify the category $\mathcal{C}$ with the one of finite-dimensional {\it unital} $\dot{\mathbf{U}}$-modules, where unital means that all elements $1_\lambda$ act as $0$ but a finite number of them, and $ \sum_{\lambda\in P} 1_\lambda$ acts as the identity. It is proved in \cite[Section~29.5.1]{Lusztig}, that
\[
\Oo_q = \left\lbrace f\colon \dot{\mathbf{U}} \ra \mc(q) \, \bigg\vert \, \begin{matrix} \text{$f$ is $\mc(q)$-linear and vanishes on some}\\ \text{two-sided ideal of finite codimension of $\dot{\mathbf{U}}$} \end{matrix}\right\rbrace.
\]

There is an analogous realization of $\Oo_A$, of the form (see \cite[Sections~23.2 and~29.5.2]{Lusztig}, and~\cite{Lusztig3})
\[
\Oo_A = \left\lbrace f\colon \dot{\mathbf{U}}_A \ra A \, \bigg\vert \, \begin{matrix} \text{$f$ is $A$-linear and vanishes on some}\\ \text{two-sided ideal of finite corank of $\dot{\mathbf{U}}_A$} \end{matrix}\right\rbrace ,
\]
where $\dot{\mathbf{U}}_A$ is the $A$-subalgebra of $\dot{\mathbf{U}}$ generated by the elements \smash{$E_i^{(k)}1_{\lambda}$} and \smash{$F_i^{(k)}1_{\lambda}$}, for all $i\in \{1,\dots,m\}$, $k\in \mn$ and $\lambda \in P$. It is a $U_A^{\rm res}$-subbimodule of $\dot{\mathbf{U}}$, and the coproduct restricts to a map \smash{$\Delta \colon \dot{\mathbf{U}}_A \ra \dot{\mathbf{U}}_A {}^{\hat{\otimes}2}$}. The above identification of the category $\mathcal{C}$ with the one of finite-dimensional unital $\dot{\mathbf{U}}$-modules yields an identification of the category $\mathcal{C}_A$ of $U_A^{\rm res}$-modules of type $1$ with the category of \smash{$\dot{\mathbf{U}}_A$}-modules of finite rank.

The key advantage of this realization of $\Oo_A$ is that $\dot{\mathbf{U}}_A$ can be equipped with a canonical $A$-basis $\dot{\mathbf{B}}$. The construction of $\dot{\mathbf{B}}$ is described in \cite[Chapter~25]{Lusztig}. It relies on the Kashiwara--Lusztig {\it canonical basis} of $U_A^{\rm res}(\mathfrak{n}_-)$. This last basis, denoted by $\mathbf{B}^-$, is defined in \cite[Chapter~14]{Lusztig}, and \cite{Kashiwara1} (a review can be found in \cite[Chapter~14]{CP}). It enjoys the following nice properties. Denote by ${}^-\colon \mc(q) \ra \mc(q)$ the field involution such that $\overline{q} = q^{-1}$, and by ${}^-\colon U_q^{\rm ad} \ra U_q^{\rm ad}$ the homomorphism of $\mc$-algebras such that
\[\bar{E}_i = E_i,\qquad \bar{F}_i = F_i,\qquad \bar{K}_\lambda = K_{-\lambda},\qquad \overline{fx} = \bar{f}\bar{x}\]
for all $f\in \mc(q)$, $x\in U_q^{\rm ad}$ ($\bar{E}_i$ and $\bar{F}_i$ above, which will not appear elsewhere, should not be confused with the normalized elements in \eqref{normgen0}). The map ${}^-$ yields a $\mc$-algebra homomorphism~${{}^-\colon \dot{\mathbf{U}} \ra \dot{\mathbf{U}}}$. Then, we have
\begin{enumerate}\itemsep=0pt
\item[(1)] the elements of $\mathbf{B}^-$ are weight vectors under the adjoint action of $U_q^{\rm ad}(\mathfrak{h})$;
\item[(2)] for every $b\in \mathbf{B}^-$, $\bar{b} =b$;
\item[(3)] for every $b,b'\in \mathbf{B}^-$, $bb' = \sum_{b''\in \mathbf{B}^-} N^{bb'}_{b''} b''$ where $N^{bb'}_{b''}\in \mz\big[q,q^{-1}\big]$;
\item[(4)] for every $b,b'\in \mathbf{B}^-$, $\Delta(b) = \sum_{b',b''\in \mathbf{B}^-} C_{b'b''}^b b' \otimes b''$ where $C_{b'b''}^{b}\in \mz\big[q,q^{-1}\big]$;
\item[(5)] for every $\mu\in P^+$, denoting by $v_\mu$ the highest weight vector of the $U_A^{\rm res}$-module ${}_AV_\mu$, the elements $bv_\mu$ which are non-zero, where $b\in \mathbf{B}^-$, form an $A$-basis of ${}_AV_\mu$.
\end{enumerate}
When $\mathfrak{g}$ is simply laced, the coefficients \smash{$N^{bb'}_{b''}$} and \smash{$C_{b'b''}^b$} belong to $\mn\big[q,q^{-1}\big]$ \cite[Theorem~14.3.13]{Lusztig}. In the case of $\mathfrak{g}={\mathfrak{sl}_2}$, the elements of $\mathbf{B}^-$ are just the divided powers \smash{$F^{(k)}$}, $k\in \mn$. Formulas in terms of PBW basis elements are known also for $\mathfrak{g}=\mathfrak{sl}_3$ and $\mathfrak{sl}_4$, and an algorithm exists in the general case (see \cite{deGraaf} and the references therein).

Correspondingly to $\mathbf{B}^-$, the set $\mathbf{B}^+ = \omega(\mathbf{B}^-)$ is a basis of $U_A^{\rm res}(\mathfrak{n}_+)$, where $\omega\colon U_q^{\rm ad} \ra U_q^{\rm ad}$ is the ($\mc(q)$-linear) {\it Cartan automorphism}, defined by
\[\omega(E_i)=F_i , \qquad \omega(F_i) = E_i , \qquad \omega(K_i)=K_i^{-1}\]
for $i=1,\dots,m$. The triangular decomposition of $\dot{\mathbf{U}}$ implies that the elements $b^+1_\lambda b'{}^-$, where~${b^+ \in \mathbf{B}^+}$, $b'{}^- \in \mathbf{B}^-$ and $\lambda\in P$, form a basis of $\dot{\mathbf{U}}$. They form in fact an $A$-basis of~$\dot{\mathbf{U}}_A$, and its elements are fixed by the involution ${}^-\colon \dot{\mathbf{U}} {}\ra{} \dot{\mathbf{U}}$.

Lusztig has constructed another $A$-basis of $\dot{\mathbf{U}}_A$, denoted $\dot{\mathbf{B}}$, and called the {\it canonical basis} of~$\dot{\mathbf{U}}_A$. It satisfies numerous properties that we now review. Its elements are denoted by~${b\,\lozenge_\lambda\, b'}$, where $b, b' \in \mathbf{B}^-$ and $\lambda\in P$, and are related to the elements $b^+b'{}^-1_\lambda$, where $b^+ := \omega(b)$ and~${b'{}^-:=b'}$, by a specific trigonal change of basis with coefficients in~$A$. Although we always have $b^+1_\lambda$, $b'{}^-1_\lambda \in \dot{\mathbf{B}}$, to our knowledge explicit formulas of the elements of $\dot{\mathbf{B}}$ as linear combinations of elements $b^+1_\lambda b'{}^-$ or $b'{}^- 1_\lambda b^+$ are known only for $\mathfrak{g}={\mathfrak{sl}_2}$ or $\mathfrak{sl}_3$ (see \cite[Section~25.3]{Lusztig} and~\cite{Cui}). In the former case, identifying $P$ with $\mz$ and $Q$ with $2\mz$ the canonical basis $\dot{\mathbf{B}}$ is formed by the elements
\[
E^{(k)}1_{-n}F^{(l)}\quad \text{and}\quad F^{(l)}1_nE^{(k)}, \qquad k,l,n\in \mn,\qquad n\geq k+l,
\]
where $E^{(k)}1_{-n}F^{(l)} = F^{(l)}1_nE^{(k)}$ for $n=k+l$.

We are going to review Lusztig's construction of $\dot{\mathbf{B}}$, its canonical partition $\dot{\mathbf{B}} = \bigcup_{\lambda\in P_+} \dot{\mathbf{B}} [\lambda]$, the dual basis \smash{$\dot{\mathbf{B}} {}^*$}, and Kashiwara's approach to \smash{$\dot{\mathbf{B}} {}^*$} \cite{Kashiwara2,Kashiwara3}. The latter is stated in Theorem~\ref{canKashi} below. At first we need to recall the notions of based module and balanced triple; for details on these notions we refer to \cite[Chapter~27]{Lusztig} and \cite{Kashiwara2} (see also \cite{Kashiwara4}, \cite[Sections~3.15 and~3.16]{VY}, or \cite[Chapter~14]{CP} for overviews).

Denote by $\mathcal{A}_0\subset \mc(q)$ the ring of rational functions regular at $q=0$. By applying the involution ${}^-$, put $\mathcal{A}_\infty=\overline{\mathcal{A}}_0$.
Since $\mathcal{A}_0$ is the localization of $\mc[q]$ at $q=0$, we may regard $\mathcal{A}_\infty$ as the localization of $\mc\big[q^{-1}\big]$ at~${q=\infty}$.

Let us recall briefly the definition of crystal basis (see \cite{Kashiwara1}). Denote by $U_q^{\rm ad}(\mathfrak{g})_i$ the subalgebra of $U_q^{\rm ad}(\mathfrak{g})$ generated by $E_i$, $F_i$ and $K_i^{\pm 1}$; thus $U_q^{\rm ad}(\mathfrak{g})_i$ is isomorphic to $U_{q_i}({\mathfrak{sl}_2})$. Let $M$ be a~$U_q^{\rm ad}$-module of type $1$. Denote \smash{$M^\zeta$} the subspace of $M$ of weight $\zeta\in P$. For every $i=1,\dots, m$, we can regard $M$ as a $U_q^{\rm ad}(\mathfrak{g})_i$-module, so $M\cong \bigoplus_j V_{\lambda_j}$ for some simple $U_q^{\rm ad}(\mathfrak{g})_i$-modules $V_{\lambda_j}$. These being generated by primitive weight vectors, the PBW basis of $U_q^{\rm ad}(\mathfrak{g})_i$ yields
\[M = \bigoplus_{\zeta\in P} \bigoplus_{0\leq n\leq (\check{\alpha}_i, \zeta)} F_i^{(n)}\big({\rm Ker}(E_i)\cap M^\zeta\big).\]
The {\it Kashiwara operators} $\tilde{e}_i$, $\tilde{f}_i$ are the endomorphisms of $M$ defined by, for every $v\in {\rm Ker}(E_i)\cap M^\zeta$ and $0\leq n\leq (\check{\alpha}_i, \zeta)$,
\[
\tilde{f}_i\big(F_i^{(n)}v\big) = F_i^{(n+1)}v,\qquad \tilde{e}_i\big(F_i^{(n)}v\big) = F_i^{(n-1)}v.
\]
A {\it crystal basis of $M$ at $q=0$} consists of a pair $(\Ll, \mathcal{B})$, where
\begin{itemize}\itemsep=0pt
\item $\Ll$ is a free $\mathcal{A}_0$-sublattice of $M$ such that the canonical map $\Ll \bigotimes_{\mathcal{A}_0} \mc(q) \ra M$ is an isomorphism;
\item $\mathcal{B}$ is a basis of the $\mc$-vector space $\Ll/q\Ll$;
\item $ \Ll = \bigoplus_{\zeta \in P} \Ll^\zeta$ and \smash{$\mathcal{B} = \coprod_{\zeta \in P} \big(\mathcal{B}\cap \Ll^\zeta/q\Ll^\zeta\big)$}, where $\Ll^\zeta = \Ll\cap M^\zeta$;
\item for every $i=1,\dots, m$ the Kashiwara operators $\tilde{e}_i$, $\tilde{f}_i$ preserve $\Ll$, and the induced maps on~$\Ll/q\Ll$ send $\mathcal{B}$ into $\mathcal{B}\cup \{0\}$, and satisfy $b' = \tilde{f}_i(b)$ if and only if $b=\tilde{e}_i(b')$ for every~${b,b'\in B}$.
\end{itemize}
Crystal bases at $q=\infty$ are defined similarly, by replacing $\mathcal{A}_0$ with $\mathcal{A}_\infty$ and $q$ with $q^{-1}$.

 A {\it based module} consists of a pair $(M,B)$ where $M$ is a~$U_q^{\rm ad}$-module of type $1$ endowed with a~$\mc(q)$-basis $B$ such that the following conditions hold:
\begin{itemize}\itemsep=0pt
\item[(i)] For every weight $\zeta\in P$, the set $B\cap M^{\zeta}$ is a basis of the weight subspace $M^{\zeta}\subset M$.
\item[(ii)] The $A$-module ${}_A M$ generated by $B$ is stable under $U_A^{\rm res}$.

We will denote by $\Ll_M$ the $\mathcal{A}_0$-submodule of $M$ generated by $B$, and by $\bar{\Ll}_M$ the $\mathcal{A}_\infty$-submodule of $M$ generated by $B$.

\item[(iii)] The $\mc$-linear involution ${}^-\colon M\ra M$ defined by $\overline{fb} = \overline{f}b$ for all $f\in \mc(q)$ and $b\in B$ is compatible with the action of $U_q^{\rm ad}$ in the sense that $\overline{xm} = \bar{x}\bar{m}$ for all $x\in U_q^{\rm ad}$, $m\in M$.
\item[(iv)] The $\mathcal{A}_\infty$-submodule $\bar{\Ll}_M$ of $M$ together with the image of $B$ in $\bar{\Ll}_M/q^{-1}\bar{\Ll}_M$ forms a crystal basis of $M$ at $q=\infty$.
\end{itemize}
If $(M,B)$ is a based module, we will denote by $\overline{\mathcal{B}}$ the image of $B$ in $\bar{\Ll}_M/q^{-1}\bar{\Ll}_M$. From the notion of balanced triple that we recall now, denoting by $\mathcal{B}$ the image of $B$ in $\Ll_M/q\Ll_M$, we see that $(\Ll_M,\mathcal{B})$ is a crystal basis at $q=0$.

Indeed, consider more generally a $\mc(q)$-vector space $V$, finite-dimensional or not, a sub-$A$-module ${}_AV$, a sub-$\mathcal{A}_0$-module ${}_{\mathcal{A}_0}V$ and a sub-$\mathcal{A}_\infty$-module ${}_{\mathcal{A}_\infty}V$ satisfying the conditions (all isomorphisms being the canonical maps)
\[V\cong \mc(q) \bigotimes_A {}_AV ,\qquad V\cong \mc(q) \bigotimes_{\mathcal{A}_0} {}_{\mathcal{A}_0}V ,\qquad V\cong \mc(q) \bigotimes_{\mathcal{A}_\infty} {}_{\mathcal{A}_\infty}V.\]
Consider the $\mc$-vector space $E := {}_AV \cap {}_{\mathcal{A}_0}V \cap {}_{\mathcal{A}_\infty}V$. Then $({}_AV, {}_{\mathcal{A}_0}V,{}_{\mathcal{A}_\infty}V)$ is a {\it balanced triple} \cite{Kashiwara1,Kashiwara2} if the canonical maps
\begin{equation}\label{condbalanced}
A \bigotimes_\mc E \ra {}_AV,\qquad \mathcal{A}_0 \bigotimes_\mc E \ra {}_{\mathcal{A}_0}V ,\qquad \mathcal{A}_\infty \bigotimes_\mc E \ra {}_{\mathcal{A}_\infty}V
\end{equation}
are isomorphisms. Equivalently, $({}_AV, {}_{\mathcal{A}_0}V,{}_{\mathcal{A}_\infty}V)$ is balanced if and only if the canonical map $E \ra {}_{\mathcal{A}_0}V/q{}_{\mathcal{A}_0}V$ is an isomorphism, if and only if the canonical map $E \ra {}_{\mathcal{A}_\infty}V/q^{-1}{}_{\mathcal{A}_\infty}V$ is an isomorphism \cite[Lemma 2.1.1]{Kashiwara2}.

Given a based module $(M,B)$, the elements of $B$ are weight vectors and $\overline{b}=b$ for every~${b\in B}$. Also, if an element $m\in {}_A M$ satisfies $\overline{m}=m$ and $m\in B + q^{-1}\bar{\Ll}_M$, then $m\in B$ (see \cite[Section~27.1.5]{Lusztig} for details on this fact). It follows that the canonical quotient map
\begin{equation}\label{basedbalancediso} {}_A M \cap \Ll_M \cap \bar{\Ll}_M \ra \bar{\Ll}_M/q^{-1}\bar{\Ll}_M
\end{equation}
is an isomorphism of $\mc$-vector spaces. This provides another way of viewing based modules: by \eqref{basedbalancediso}, $\big({}_A M, \Ll_M,\bar{\Ll}_M\big)$ is a balanced triple, and by \eqref{condbalanced} the $A$-lattice ${}_A M$ is completely determined by the crystal base $\big(\bar{\Ll}_M,\overline{\mathcal{B}}\big)$. We will say that $\big(\bar{\Ll}_M,\overline{\mathcal{B}}\big)$ (or the corresponding crystal base $(\Ll_M,\mathcal{B})$ at $q=0$) is {\it melted into} the based module $(M,B)$.

We will indifferently apply the notion of based module to finite-dimensional unital $\dot{\mathbf{U}}$-modules, since they are equivalent to $U_q^{\rm ad}$-modules of type $1$.

Every module $V_\mu$, $\mu\in P^+$, supports a structure of based module (see \cite[Section~14.4.10]{Lusztig} and~\cite{Kashiwara1}); the corresponding basis, called {\it canonical basis} and that we will denote by $\underline{\mathbf{B}}_{\mu}$, is formed by the elements $bv_{\mu}\in {}_AV_\mu$ which are non-zero, where $b\in \mathbf{B}^-$ and $v_{\mu}$ is the canonical highest weight vector of $V_{\mu}$, corresponding to the coset of $1\in U_q^{\rm ad}(\mathfrak{n}_-)$ in the Verma module construction of $V_{\mu}$. Note that the involution $\ \bar{}\ \colon V_{\mu}\ra V_{\mu}$ defined by (iii) above is indeed an automorphism, for the space $V_\mu$ with action of $U_q^{\rm ad}$ defined by $x\cdot v := \bar{x}v$, for all $x\in U_q^{\rm ad}$, $v\in V_\mu$, has highest weight $\mu$, and is thus isomorphic to $V_\mu$ via the map $\ \bar{}\ $. The crystal base \smash{$\big(\Ll_\mu^{\rm low},\mathcal{B}_\mu^{\rm low}\big)$} at $q=0$ is formed by the $\mathcal{A}_0$-sublattice $\Ll_\mu^{\rm low}$ of $V_{\mu}$ generated by $\underline{\mathbf{B}}_{\mu}$ (which is eventually the same as the $\mathcal{A}_0$-sublattice generated by the vectors of the form \smash{$\tilde{f}_{i_1}\circ \dots \circ \tilde{f}_{i_k}(v_\mu)$}, where $i_1,\dots ,i_k \in \{1,\dots,m\}$), and $\mathcal{B}_\mu^{\rm low}$ is the set of non-zero images of these vectors in \smash{$\Ll_\mu^{\rm low}/q\Ll_\mu^{\rm low}$}.

There is the following uniqueness result \cite[Theorem 3]{Kashiwara1}.

\begin{teo}\label{Kashiuniqueteo} Let $M$ be a $U_q^{\rm ad}$-module of type $1$, and $(\Ll, \mathcal{B})$ a crystal base at $q=0$ of~$M$. Then there exists a $\mc(q)$-isomorphism $M\ra \bigoplus_j V_{\lambda_j}$ by which $(\Ll, \mathcal{B})$ is $\mathcal{A}_0$-isomorphic to \smash{$ \bigoplus_j \big(\Ll_{\lambda_j}^{\rm low}, \mathcal{B}_{\lambda_j}^{\rm low}\big)$}.
\end{teo}

The based modules form a category. Given based modules $(M,B)$ and $(M',B')$, a morphism of $U_q^{\rm ad}$-modules $f\colon M\ra M'$ is a morphism of based modules if{\samepage
\begin{itemize}\itemsep=0pt
\item[(a)] $f(b)\in B'\cup\{0\}$ for any $b\in B$;
\item[(b)] $B\cap {\rm Ker}(f)$ is a basis of ${\rm Ker}(f)$.
\end{itemize}}%

\pagebreak

\noindent
The direct sum of based modules $(M,B)$ and $(M',B')$ is a based module $(M \oplus M',B\cup B')$; and a submodule $M'$ of a~based module $(M,B)$ spanned over $\mc(q)$ by a subset $B'$ of $B$ forms a~based module $(M',B')$. The quotient module $M/M'$ together with the image of $B\setminus B'$ is then a based module.

The tensor product of based modules $(M,B)$, $(M',B')$ is also defined. Namely, consider the~$\mc$-linear map $\Psi\colon M\otimes M' \ra M\otimes M'$ defined by
\[
\Psi(m\otimes m') = \hat{R}^{-1}(\bar{m}\otimes \bar{m}'),
\] where $\hat{R} = \Theta^{-1}R$, see \eqref{Rmatfact} (note that, as we use the coproduct opposite to \cite{Lusztig} our quasi-$R$-matrix is $\hat{R}^{-1}$). It can be checked that $\Psi$ is an involution compatible with the action of $\dot{\mathbf{U}}$ in the sense of (iii) above in the definition of based module. Moreover, denote by $\Ll_{M,M'}$ the~$\mc[q^{-1}]$-submodule of $M\otimes M'$ spanned by the basis elements~$b\otimes b'$, where $b\in B$, $b'\in B'$. It is shown in \cite[Section~27.3]{Lusztig}, that for every pair $(b,b')\in B\times B'$ there is a unique element $b\,\lozenge\, b'\in \Ll_{M,M'}$ such that
\begin{itemize}\itemsep=0pt
\item[(a)] $\Psi(b\,\lozenge\, b')=b\,\lozenge\, b'$,
\item[(b)] $b\,\lozenge\, b' - b\otimes b'\in q^{-1} \mathcal{\Ll}_{M,M'}$.
\end{itemize}
Moreover, $B_{\lozenge} = \{b\,\lozenge\, b', b\in B, b'\in B'\}$ is a basis of $M\otimes M'$, a $\mc[q^{-1}]$-basis of $\Ll_{M,M'}$, a $\mc\big[q,q^{-1}\big]$-basis of the $\mc\big[q,q^{-1}\big]$-module ${}_A\Ll_{M,M'}$ of $M\otimes M'$ generated by the elements $b\otimes b'$, where $b\in B$, $b'\in B'$, and $(M\otimes M',B_{\lozenge})$ is a based module.

This construction of $B_{\lozenge}$ is associative. Since $(V_\mu, \underline{\mathbf{B}}_{\mu})$ is for every $\mu\in P_+$ a based module, it follows that any tensor product $M$ of a finite number of the simple modules $V_\mu$ is naturally a based module. The corresponding basis is called {\it the canonical basis} of $M$. These canonical basis have been computed explicitly in \cite{FK} in the case $\mathfrak{g}={\mathfrak{sl}_2}$.

Consider now the $U_q^{\rm ad}$-module ${}^{\omega}V_{\mu}$ with underlying space $V_{\mu}$, $\mu\in P_+$, and action defined by $x._{\omega} v := \omega(x)v$, for every \smash{$x\in U_q^{\rm ad}$} and $v\in V_{\mu}$ (as usual $\omega\colon U_q^{\rm ad} \ra U_q^{\rm ad}$ is the Cartan automorphism). Note that there are isomorphisms \smash{${}^{\omega}V_{\mu}\cong V_{-w_0(\mu)}\cong V_{\mu}^*$} (endowed with the standard left action of $U_q^{\rm ad}$). Let us denote by ${}^\omega v_\mu$ the vector $v_\mu$ regarded in ${}^{\omega}V_{\mu}$ (i.e., its canonical lowest weight vector), and by ${}^\omega\underline{\mathbf{B}}_{\mu}:=\{b._{\omega}{}^\omega v_{\mu},\, b\in \mathbf{B}^+\}\setminus \{0\}$ its canonical basis; note that ${}^\omega\underline{\mathbf{B}}_{\mu} = \{\omega(b){}v_{\mu},\, b\in \omega(\mathbf{B}^-)\}\setminus \{0\} = \{b v_{\mu},\, b\in \mathbf{B}^-\}\setminus \{0\} = \underline{\mathbf{B}}_{\mu}$. Then ${}^{\omega}V_{\mu'}\otimes V_{\mu''}$ has the canonical basis $\underline{\mathbf{B}}_{\mu',\mu''} := \{\underline{b}'\,\lozenge\, \underline{b}'', \, \underline{b}'\in {}^\omega\underline{\mathbf{B}}_{\mu'},\, \underline{b}''\in \underline{\mathbf{B}}_{\mu''}\}$. Since $\underline{b}'\,\lozenge\, \underline{b}''$ is canonically determined by the elements $b' ,b''\in \mathbf{B}^-$ such that $\underline{b}' = \omega(b')._{\omega} {}^\omega v_{\mu'}$, $\underline{b}'' = b''v_{\mu''}$, following Lusztig we denote it by~${(b' \,\lozenge\, b'')_{\mu',\mu''}}$.

 Denote by $v_{w_0(\mu)}$ the canonical lowest weight vector of $V_{\mu}$, and by ${}^\omega v_{w_0(\mu)}$ the vector $v_{w_0(\mu)}$ regarded in ${}^{\omega}V_{\mu}$. It is a crucial observation that ${}^\omega v_{w_0(\mu')}\otimes v_{w_0(\mu'')}$ is a cyclic vector of ${}^{\omega}V_{\mu'}\otimes V_{\mu''}$ (see, e.g., \cite[Proposition~23.3.6]{Lusztig}; note that \smash{${}^\omega v_{w_0(\mu')}\otimes v_{w_0(\mu'')}$} plays the role of $\xi_{-\mu'}\otimes \eta_{\mu''}:= {}^\omega v_{\mu'}\otimes v_{\mu''}$ in \cite{Lusztig}, because we use opposite coproducts on $U_q^{\rm ad}$ but the factors ${}^{\omega}V_{\mu'}$ and $V_{\mu''}$ are ordered in the same way).

 We can now state the definition of the canonical basis $\dot{\mathbf{B}}$ of $\dot{\mathbf{U}}$: each element $u$ of $\dot{\mathbf{B}}$ belongs to $\dot{\mathbf{U}}_A 1_\zeta$ for some (unique) $\zeta\in P$, and it is then uniquely determined by the property that, for every $\mu'$, $\mu''\in P^+$ such that $w_0(\mu''-\mu')=\zeta$, we have
\begin{equation}\label{udef}
u({}^\omega v_{w_0(\mu')}\otimes v_{w_0(\mu'')}) = (b' \,\lozenge\, b'')_{\mu',\mu''}\end{equation}
for some \smash{$(b' \,\lozenge\, b'')_{\mu',\mu''} \in \underline{\mathbf{B}}_{\mu',\mu''}$} \cite[Section~25.2]{Lusztig}. We write \smash{$u = b' \, \lozenge_\zeta\, b''$}, and as in \cite{Lusztig3} we denote by \smash{$\dot{\mathbf{B}}_{\mu',\mu''}$} the finite subset of \smash{$\dot{\mathbf{B}}$} which is in bijection with $\underline{\mathbf{B}}_{\mu',\mu''}$ under the map $u\mapsto u({}^\omega v_{w_0(\mu')}\otimes v_{w_0(\mu'')})$. So
\begin{equation}\label{Bpointunion}
\dot{\mathbf{B}} =\bigcup_{\mu',\mu''\in P_+} \dot{\mathbf{B}}_{\mu',\mu''}.
\end{equation}
Note in particular that $\dot{\mathbf{B}}$ is formed by weight vectors for the left and right action of $U_q^{\rm ad}(\mathfrak{h})$ (defined as usual by~\eqref{bimodstr}).

In a sense, one can view $\dot{\mathbf{U}}$ as the projective limit of an inverse system formed by the \smash{$\big(U_q^{\rm ad}\otimes U_q^{\rm ad}\big)$}-modules \smash{${}^{\omega}V_{\mu'}\otimes V_{\mu''}$}, where $\mu',\mu''\in P^+$; then $\dot{\mathbf{B}}$ is the basis resulting from the corresponding inverse system of basis $\{\dot{\mathbf{B}}_{\mu',\mu''}\}_{\mu',\mu''}$.

Lusztig has produced a partition of $\dot{\mathbf{B}}$ as follows. First, consider the situation of a based module $(M,B)$. For every $\lambda\in P_+$ denote by $M[\lambda]$ the sum of the simple submodules of $M$ isomorphic to $V_\lambda$ (i.e., its isotypical component). Set
\begin{equation}\label{decbaseisospace}
M[\geq\lambda] = \bigoplus_{\lambda'\geq \lambda} M[\lambda'].
\end{equation}
Then, for every base element $b\in B$ there is a unique $\lambda\in P_+$ such that $b\in M[\geq\lambda]$ and $\lambda$ is maximal with this property \cite[Section~27.2]{Lusztig}. Denote by $B[\lambda]$ the set of all $b\in B$ that give rise to $\lambda\in P_+$ in this way. Clearly, the sets $B[\lambda]$, $\lambda\in P_+$, form a partition of $B$.

Now, given $b\in \dot{\mathbf{B}}$, let $\zeta\in P$ be the unique weight such that $b \in \dot{\mathbf{U}}_A 1_\zeta$, and let ${\mu', \mu''\in P^+}$ be such that $w_0(\mu''-\mu')=\zeta$, and $(\check{\alpha}_i,\mu')$ is large enough for all $i=1,\dots,m$ so that $u({}^\omega v_{w_0(\mu')}\otimes v_{w_0(\mu'')})$ is non-zero. This element belongs to the canonical basis $\underline{\mathbf{B}}_{\mu',\mu''}$ of ${}^{\omega}V_{\mu'}\otimes V_{\mu''}$, and therefore to one of the subsets $\underline{\mathbf{B}}_{\mu',\mu''}[\lambda]$, for a unique $\lambda\in P_+$. It is a result that $\lambda$ does not depend on the choice of $\mu'$, $\mu''$ (see \cite[Section~29.1.1]{Lusztig}). Hence there is a well-defined map~${\dot{\mathbf{B}}\ra P_+}$,~$b\mapsto \lambda$. Denoting by $\dot{\mathbf{B}} [\lambda]$ the fiber of this map, we thus obtain a partition
\begin{equation}\label{partition} \dot{\mathbf{B}} = \coprod_{\lambda\in P_+}\dot{\mathbf{B}} [\lambda]. \end{equation}
The sets $\dot{\mathbf{B}} [\lambda]$ are called {\it $2$-sided cells}. They are finite sets and have the following remarkable properties. For every $\lambda\in P_+$ denote by~${\dot{\mathbf{U}} [\geq \lambda]}$ and $\dot{\mathbf{U}} [> \lambda]$ the subspaces of $\dot{\mathbf{U}}$ spanned by~\smash{$ \coprod_{\lambda'\geq \lambda} \dot{\mathbf{B}} [\lambda']$} and \smash{$ \coprod_{\lambda'> \lambda} \dot{\mathbf{B}}[\lambda']$} respectively. Then $\dot{\mathbf{U}} [\geq \lambda]$ (respectively~${\dot{\mathbf{U}} [> \lambda]}$) consists of the elements $u\in \dot{\mathbf{U}}$ such that if $u$ acts on \smash{$V_\mu$} by a non-zero linear map, then $\mu\geq \lambda$ (respectively~${\mu> \lambda}$) \cite[Lemmas~29.1.3 and~29.1.4]{Lusztig}. Both $\dot{\mathbf{U}} [\geq \lambda]$ and $\dot{\mathbf{U}} [> \lambda]$ are two-sided ideals of~$\dot{\mathbf{U}}$. Moreover, the algebra homomorphism $\pi_\lambda\colon \dot{\mathbf{U}} [\geq\lambda] \ra {\rm End}(V_\lambda)$ given by the $\dot{\mathbf{U}}$-module structure on $V_\lambda$ descends to an algebra and $U_q^{\rm ad}$-bimodule isomorphism (keeping the same notation) \cite[Proposition~29.2.2]{Lusztig}
\begin{equation}\label{pimuiso} \bar{\pi}_\lambda\colon\ \dot{\mathbf{U}} [\geq \lambda]/ \dot{\mathbf{U}} [> \lambda] \ra {\rm End}(V_\lambda) . \end{equation}

For instance, when $\mathfrak{g}={\mathfrak{sl}_2}$ the $2$-sided cell $\dot{\mathbf{B}} [n]$ associated to the simple $U_q^{\rm ad}({\mathfrak{sl}_2})$-module of type $1$ and dimension $n+1$ is the set of cardinality $(n+1)^2$ given by \cite[Section~29.4.3]{Lusztig}
\begin{equation}\label{2sidedcellsl2}
\dot{\mathbf{B}} [n] =\big\{E^{(k)}1_{-n}F^{(l)},\,n\geq k+l \big\} \cup \big\{F^{(l)}1_nE^{(k)},\,n\geq k+l\big\},
\end{equation}
with the identification $E^{(k)}1_{-n}F^{(l)} = F^{(l)}1_nE^{(k)}$ when $n=k+l$. As we are mainly interested in~$\Oo_A$, we are going to describe the dual partition of $\dot{\mathbf{B}} {}^*$, see Theorem~\ref{canKashi}. The duality with~\eqref{partition} is discussed after that theorem.

First, we follow the approach of Kashiwara \cite{Kashiwara2,Kashiwara3}. For every $\lambda\in P_+$, denote by $V_\lambda^r$ the dual space of $V_\lambda$ endowed with its natural structure of right $U_q^{\rm ad}$-module, defined by ${(f x)(v) = f(xv)}$
for every $f\in V_\lambda^r$, \smash{$x\in U_q^{\rm ad}$}, $v\in V_\lambda$.
Clearly, $V_\lambda^r$ is a simple module of highest weight $\lambda$.
 Let~\smash{$\varphi\colon U_q^{\rm ad} \ra U_q^{\rm ad}$} be the anti-automorphism of $\mc(q)$-algebra given by $\varphi(E_i) = F_i$, $\varphi(F_i) = E_i$, $\varphi(K_\lambda) = K_{\lambda}$. By using $\varphi$, any right $U_q^{\rm ad}$-module can be considered as a left $U_q^{\rm ad}$-module. In particular, by the Verma module construction of $V_\lambda$ it follows
\[V_\lambda^r \cong U_q^{\rm ad}\bigg/\bigg(\sum_{\mu\in P_+} \big(K_\mu -q^{(\lambda,\mu)}\big)U_q^{\rm ad}+ \sum_{i=1}^m E_i^{1+(\check{\alpha}_i,\lambda)}U_q^{\rm ad}\bigg),\]
and $\varphi$ affords an isomorphism of the right module $V_\lambda^r$ with the {\it left} module $V_\lambda$. We will denote by $f_\lambda$ the unique highest weight vector of $V_\lambda^r$ satisfying $\langle f_\lambda,v_\lambda\rangle = 1$.

The space $V_\lambda^r \otimes V_\lambda$ can be identified with ${\rm End}(V_\lambda)^*$, and thus acquires by duality a natural structure of $U_q^{\rm ad}$-bimodule \big(or equivalently left $U_q^{\rm ad} \otimes \big(U_q^{\rm ad}\big)^{\rm op}$-module\big); the left and right actions are given by
\begin{equation}\label{bimodaction}x(f \otimes v)y = fy \otimes xv
\end{equation}
for every $x,y\in U_q^{\rm ad}$, $f\in V_\lambda^r$, $v\in V_\lambda$. The space $V_\lambda^r \otimes V_\lambda$ also acquires by duality a natural ``upper" crystal structure over $U_q^{\rm ad} \otimes \big(U_q^{\rm ad}\big)^{\rm op}$, as we explain now. Denote by $\langle \ ,\ \rangle_\lambda \colon V_\lambda \times V_\lambda \ra \mc(q)$ the unique symmetric bilinear form such that
\begin{equation}\label{proppairinglambda}
\langle v_\lambda,v_\lambda\rangle_\lambda =1\qquad{\rm and}\qquad \langle \varphi(x)u,v\rangle_\lambda = \langle u, xv\rangle_\lambda
\end{equation}
for every $u,v\in V_\lambda$ and $x\in U_q^{\rm ad}$. Recall the crystal base $\big(\Ll_\mu^{\rm low},\mathcal{B}_\mu^{\rm low}\big)$ at $q=0$ introduced before Theorem \ref{Kashiuniqueteo}. In Kashiwara's terminology \cite{Kashiwara1, Kashiwara2}, the pair $\big(\Ll_{\lambda}^{\rm low},\mathcal{B}_\lambda^{\rm low}\big)$ is the {\it lower crystal base} of $V_\lambda$ at $q=0$. Applying the involution ${}^-\colon V_\lambda\ra V_\lambda$, one obtains the lower crystal base \smash{$\big(\overline{\Ll_{\lambda}^{\rm low}},\overline{\mathcal{B}_\lambda^{\rm low}}\big)$} at $q=\infty$. Because the canonical bases are determined by the crystal bases (see the discussion about \eqref{basedbalancediso}), we call $(V_\lambda,\underline{\mathbf{B}}_\lambda)$ the {\it lower} based module of $V_\lambda$, and $\underline{\mathbf{B}}_\lambda$ the {\it lower canonical basis} of $V_\lambda$.

Put
\begin{align}
& {}_AV_\lambda^{\rm up} := \{v\in V_\lambda, \langle v,{}_AV_\lambda\rangle_\lambda\subset A\},\qquad
\Ll_{\lambda}^{\rm up} := \big\{v\in V_\lambda, \langle v,\Ll_{\lambda}^{\rm low}\rangle_\lambda\subset \mathcal{A}_0\big\},\nonumber\\
&\overline{\Ll_{\lambda}^{\rm up}} := \big\{v\in V_\lambda, \langle v,\overline{\Ll_{\lambda}^{\rm low}}\rangle_\lambda\subset \mathcal{A}_\infty\big\}.\label{upmodule}
\end{align}
Then $\big({}_A\!V_\lambda^{\rm up},\Ll_{\lambda}^{\rm up},\overline{\Ll_{\lambda}^{\rm up}}\big)$ is a balanced triple \cite[Lemma 4.2.1]{Kashiwara2}. Denote by $\mathcal{B}_\lambda^{\rm up}\!$ the basis of~${\Ll_{\lambda}^{\rm up}\!/\!q\Ll_{\lambda}^{\rm up}}$ dual to $\mathcal{B}_\lambda^{\rm low}$ by the induced pairing $\langle\ ,\ \rangle_\lambda\colon \Ll_{\lambda}^{\rm up}/q\Ll_{\lambda}^{\rm up} \times \Ll_{\lambda}^{\rm low}/q\Ll_{\lambda}^{\rm low} \ra \mc$. The pair $\big(\Ll_{\lambda}^{\rm up},\mathcal{B}_\lambda^{\rm up}\big)$ is the {\it upper crystal base} of $V_\lambda$ at $q=0$. The weight spaces of the $\mathcal{A}_0$-modules $\Ll_{\lambda}^{\rm low}$ and $\Ll_{\lambda}^{\rm up}$ are related by
\begin{equation}\label{renormweightA} \big(\Ll_{\lambda}^{\rm up}\big)^\mu = q^{\frac{(\lambda,\lambda)}{2} - \frac{(\mu,\mu)}{2}} \big(\Ll_{\lambda}^{\rm low}\big)^\mu , \qquad \mu\in P. \end{equation}
Correspondingly, denoting $\big(\mathcal{B}_\lambda^{\rm up}\big)^\mu := \mathcal{B}_\lambda^{\rm up} \cap \big(\Ll_\lambda^{\rm up}\big)^\mu$ and $\big(\mathcal{B}_\lambda^{\rm low}\big)^\mu := \mathcal{B}_\lambda^{\rm low} \cap \big(\Ll_\lambda^{\rm low}\big)^\mu$, we have (see~\cite{Kashiwara1} and \cite[equation~(4.2.9)]{Kashiwara2})
\[\big(\mathcal{B}_\lambda^{\rm up}\big)^\mu = q^{\frac{(\lambda,\lambda)}{2} - \frac{(\mu,\mu)}{2}} \big(\mathcal{B}_\lambda^{\rm low}\big)^\mu.\]
The $A$-module ${}_AV_\lambda^{\rm up}$ is characterized by the following two properties \cite[equations~(4.2.10)--(4.2.12)]{Kashiwara2}:
\begin{align*}
\big({}_AV_\lambda^{\rm up}\big)^\lambda & = \mc\big[q,q^{-1}\big]v_\lambda,\qquad
\big({}_AV_\lambda^{\rm up}\big)^\mu = \bigl\{v\in V_\lambda \mid U_A^{\rm res}\big(\mathfrak{n}^+\big)^{\lambda-\mu} v \in \mc\big[q,q^{-1}\big]v_\lambda\bigr\},
\end{align*}
where $U_A^{\rm res}\big(\mathfrak{n}^+\big)^{\gamma} = \big\{u\in U_A^{\rm res}\big(\mathfrak{n}^+\big) \mid \forall \nu\in P, \, K_\nu uK_\nu^{-1} = q^{(\nu,\gamma)}u\big\}$. Denote by $\underline{\mathbf{B}}_\lambda^{\rm up}$ the inverse image of $\mathcal{B}_\lambda^{\rm up}$ by the isomorphism \smash{${}_AV_\lambda^{\rm up}\cap \Ll_{\lambda}^{\rm up}\cap \overline{\Ll_{\lambda}^{\rm up}} \ra \Ll_{\lambda}^{\rm up}/q\Ll_{\lambda}^{\rm up}$}. By \eqref{condbalanced}, the set $\underline{\mathbf{B}}_\lambda^{\rm up}$ is a~basis of ${}_AV_\lambda^{\rm up}$; we call it the {\it upper canonical basis} of $V_\lambda$. In the appendix, we describe in details the ${\mathfrak{sl}_2}$ case.

Similarly, the right module $V_\lambda^r$ with its canonical basis $\underline{\mathbf{B}}_\lambda^r{}=\{f_\lambda b,\, b\in \mathbf{B}^+\}\setminus \{0\}$ has the lower crystal base $\big(\Ll_{\lambda}^r{}^{\rm low},\mathcal{B}_\lambda^r{}^{\rm low}\big)$, and it supports a balanced triple $\big({}_AV_\lambda^r{}^{\rm up},\Ll_{\lambda} ^r{}^{\rm up},\overline{\Ll_{\lambda} ^r{}^{\rm up}}\big)$ defined again by duality. We denote by $(\Ll_{\lambda}^r{}^{\rm up},\mathcal{B}_\lambda^r{}^{\rm up})$ and $\underline{\mathbf{B}}_\lambda^r{}^{\rm up}$ the corresponding crystal base and upper canonical basis of $V_\lambda^r$, respectively.

It follows that $\big({}_AV_\lambda^r{}^{\rm up}\bigotimes_A {}_AV_\lambda^{\rm up},\Ll_{\lambda}^r{}^{\rm up} \bigotimes_{\mathcal{A}_0}\Ll_{\lambda}^{\rm up},\overline{\Ll_{\lambda}^r{}^{\rm up}}\bigotimes_{\mathcal{A}_\infty}\overline{\Ll_{\lambda}^{\rm up}}\big)$ is a balanced triple; equivalently~${V_\lambda^r \otimes V_\lambda}$ with the bimodule structure \eqref{bimodaction} and the basis $\underline{\mathbf{B}}_\lambda^r{}^{\rm up}\otimes \underline{\mathbf{B}}_\lambda^{\rm up}$ is a based $\big(U_q^{\rm ad} \otimes (U_q^{\rm ad})^{\rm op}\big)$-module.

Denote again by $\langle\cdot,\cdot \rangle \colon \Oq \times \dot{\mathbf{U}} \ra \mc(q)$ the pairing of $U_q^{\rm ad}$-bimodules induced by the canonical pairing $\langle\ ,\ \rangle \colon \Oq \times U_q^{\rm ad} \ra \mc(q)$, and let $\Phi_\lambda\colon V_\lambda^r \otimes V_\lambda \ra \Oo_q$, $\lambda\in P_+$, be the ``matrix coefficient" map, i.e.,
\begin{equation}\label{pairingOAKashi}
\langle \Phi_\lambda (f \otimes v),x\rangle = \langle f, xv\rangle_\lambda
\end{equation}
for every $f\in V_\lambda^r$, $x\in U_q^{\rm ad}$, $v\in V_\lambda$. The map $\Phi := \bigoplus_{\lambda\in P_+} \Phi_\lambda$ is an isomorphism of $U_q^{\rm ad}$-bimodules, so let us use it to identify $\Oo_q$ with $\bigoplus_{\lambda\in P_+} V_\lambda^r \otimes V_\lambda$ (which is the content of the Peter--Weyl decomposition \eqref{decompdirectOq}). Define
\begin{alignat*}{3}
& \Ll(\Oo_q) = \bigoplus_{\lambda\in P_+} \biggl(\Ll_{\lambda}^r {}^{\rm up} \bigotimes_{\mathcal{A}_0} \Ll_{\lambda}^{\rm up}\biggr),\qquad && \mathcal{B}(\Oo_q) := \coprod_{\lambda\in P_+} \mathcal{B}_\lambda^r {}^{\rm up} \otimes \mathcal{B}_\lambda^{\rm up},& \\
& \overline{\Ll}(\Oo_q) = \bigoplus_{\lambda\in P_+} \biggl(\overline{\Ll_{\lambda}^r {}^{\rm up}} \bigotimes_{\mathcal{A}_\infty} \overline{\Ll_{\lambda}^{\rm up}}\biggr),\qquad && \overline{\mathcal{B}}(\Oo_q) := \coprod_{\lambda\in P_+} \overline{\mathcal{B}_\lambda^r {}^{\rm up}} \otimes \overline{\mathcal{B}_\lambda^{\rm up}}.&
\end{alignat*}
\begin{teo}\label{canKashi} \quad
\begin{itemize}\itemsep=0pt
\item[$(i)$] The triple $\big(\Oo_A,\Ll(\Oo_q),\overline{\Ll}(\Oo_q)\big)$ is balanced. Therefore, denoting by $G$ the inverse of the canonical map $\Oo_A \cap \Ll(\Oo_q) \cap \overline{\Ll}(\Oo_q)\ra \Ll(\Oo_q)/q\Ll(\Oo_q)$, we have
\[\Oo_A = \bigoplus_{b\in \mathcal{B}(\Oo_q)} A G(b).\]
\item[$(ii)$] The basis $G(\mathcal{B}(\Oo_q)):=\{G(b), \, b\in \mathcal{B}(\Oo_q)\}$ coincides with the dual canonical basis $\dot{\mathbf{B}} {}^*$, i.e., the elements $a^* \in \Oo_A$, for every $a\in \dot{\mathbf{B}}$, defined by $a^* (a') = \delta_{a,a'}$ for every $a'\in \dot{\mathbf{B}}$. Therefore, \[\Oo_A = \bigoplus_{b\in \dot{\mathbf{B}}} Ab^*.\]
\end{itemize}
\end{teo}
The statement (i) is \cite[Theorem~1]{Kashiwara2}, and (ii) is \cite[Theorem 10.1 and Proposition 10.2.2]{Kashiwara3} and \cite[Section~29.5]{Lusztig}. The basis $G(\mathcal{B}(\Oo_q))=\dot{\mathbf{B}} {}^*$ is called the {\it global basis}, or {\it canonical basis}, of~$\Oo_q$. The proof of Theorem \ref{canKashi}\,(ii) in \cite{Kashiwara3} (see also \cite{Kashiwara4}) exhibits an isomorphism of crystals over $U_q^{\rm ad} \otimes \big(U_q^{\rm ad}\big){}^{\rm op}$,
\begin{equation}\label{isocrystal}
\psi\colon\ \mathcal{B}(\Oo_q) \ra \mathcal{B}\bigl(\dot{\mathbf{U}}\bigr),
\end{equation}
where \smash{$\big(\Ll\bigl(\dot{\mathbf{U}}\bigr),\mathcal{B}\bigl(\dot{\mathbf{U}}\bigr)\big)$} is the crystal base of $\dot{\mathbf{U}}$ associated to the canonical basis $\dot{\mathbf{B}}$. The isomorphism $\psi$ satisfies $\langle G(b),G(b')\rangle =\delta_{\psi(b),b'}$ for every $b\in \mathcal{B}(\Oo_q)$, $b'\in \mathcal{B}\bigl(\dot{\mathbf{U}}\bigr)$. The unit $1$ of $\Oo_A$ is~$(1_0)^*$; the constant structures of $\Oo_A$ are studied in \cite{Lusztig,Lusztig3}.

{\bf The canonical basis of $\boldsymbol{\Oo_A}$ when $\boldsymbol{\mathfrak{g}={\mathfrak{sl}_2}}$}. Denote by $a$, $b$, $c$, $d$ the matrix coefficients in the canonical basis $(v_+,v_-:=Fv_+)$ of $V_1$, the simple $U_q^{\rm ad}({\mathfrak{sl}_2})$-module of type $1$ and dimension two, read from the top left to the bottom right. In that case of $V_1$ the upper canonical basis~$\underline{\mathbf{B}}_{1}^r{}^{\rm up}$ and \smash{$\underline{\mathbf{B}}_{1}^{\rm up}$} coincide with the lower ones (this is not true in general, see Example~\ref{V2V20}). The basis~$\dot{\mathbf{B}} {}^*({\mathfrak{sl}_2})$ is formed by the monomials $c^sa^pb^r$ where $p,r,s\in \mn$, and $c^sd^pb^r$ where $p,r,s\in \mn$ and $p>0$; this is stated in \cite[Proposition~9.1.1]{Kashiwara2} (in \cite[Proposition~1.3]{DC-L}, similar monomials are shown to form an $A$-basis of $\Oo_A({\rm SL}_2)$, but without reference to the canonical basis; see the comments before \eqref{PWrel} below). More precisely, recall the $2$-sided cells~\eqref{2sidedcellsl2}. We verified by a~tedious though straightforward computation that we have the duality pairing
\begin{gather*}
\begin{split}
& \big\langle c^sd^pb^r, E^{(i)}1_{-k}F^{(j)} \big\rangle = \delta_{p+r+s,k}\delta_{r,i}\delta_{s,j} ,\qquad \big\langle c^sd^pb^r, F^{(j)}1_{k}E^{(i)} \big\rangle =0, \\
& \big\langle c^sa^pb^r, E^{(i)}1_{-k}F^{(j)} \big\rangle = 0 ,\qquad \big\langle c^sa^pb^r, F^{(j)}1_{k}E^{(i)} \big\rangle =\delta_{p+r+s,k}\delta_{r,i}\delta_{s,j}.
\end{split}
\end{gather*}
Therefore, \begin{align*}\dot{\mathbf{B}} [n]^* :={}& \{ c^sa^pb^r, \, p,r,s\in \mn,\, p+r+s=n\} \\&\cup \{ c^{s}d^{p}b^{r},\, p,r,s\in \mn,\, p>0,\, p+r+s=n\}.\end{align*}
A description of $\dot{\mathbf{B}}^*$ in the case of $\mathfrak{g} = {\mathfrak{sl}_n}$ can be found in \cite{Du}. Moreover, denote by $V_n$ the simple~\smash{$U_q^{\rm ad}({\mathfrak{sl}_2})$}-module of type $1$ and dimension $n+1$, by $(v_k)$ the canonical basis of $V_n$, by~$\big(v^k\big)$ the dual basis, and by \smash{$\pi_n\colon \dot{\mathbf{U}}({\mathfrak{sl}_2}) \rightarrow {\rm End}(V_n)$} the representation morphism. By using the above pairing, it is readily checked that for every $0\leq l,m\leq n$, we have
\begin{gather}
v^l(\pi_n( \cdot) \ v_m)\nonumber\\
\qquad= \sum_{\substack{0\leq i,j,k\\ i+j\leq k\leq n\\j-i=l-m}} \hspace*{-0.5cm}\delta_{-k,n-2(m+j)} \left[\begin{matrix} m+j \\ j \end{matrix}\right]_{q} \left[\begin{matrix}n-m+i-j \\ i \end{matrix}\right]_{q}\big(E^{(i)}1_{-k}F^{(j)}\big)^*\nonumber\\
\phantom{\qquad=}{} + \sum_{\substack{0\leq i,j,k\\ i+j< k\leq n
\\j-i=l-m}}\delta_{k,n-2(m-i)} \left[\begin{matrix} m-i+j \\ j \end{matrix}\right]_{q} \left[\begin{matrix} n-m+i \\ i \end{matrix}\right]_{q}\big(F^{(j)}1_{+k}E^{(i)}\big)^*.\label{formulesl2}
\end{gather}
In particular, we see in this case of $\mathfrak{g}={\mathfrak{sl}_2}$ that in general the matrix coefficients of simple~\smash{$U_A^{\rm res}$}-modules of type $1$ are not elements of the dual canonical basis $\dot{\mathbf{B}} {}^*$. Moreover, these matrix coefficients do not form a basis of $\Oo_A$. For instance, it follows from \eqref{formulesl2} that the matrix of matrix coefficients of $V_2$ has the following form:
\begin{equation}\label{matrixcoefV2}
\begin{pmatrix}a^2 & [2]_qab & b^2 \\ ca & [2]_qbc +1 & db \\ c^2 & [2]_q cd & d^2 \end{pmatrix}.
\end{equation}
The matrix coefficient $v_0^* \otimes v_0$ being equal to $[2]_qbc +1$, this shows $bc$ cannot be expressed as a~linear combination over $A$ of matrix coefficients of simple modules.

{\bf The refined Peter--Weyl theorem.} Let us discuss the $U_A^{\rm res}$-bimodule structure of $\Oo_A$, and its relation with the partition \eqref{partition}. For every $\lambda\in P_+$, put
\begin{equation}\label{AClambda}
{}_A \overset{\raisebox{-1.5pt}{\scriptsize$\smallbullet$}}{C} (\lambda) := \bigoplus_{b\in \dot{\mathbf{B}}[\lambda]} Ab^*
\end{equation}
and
\[\Oo_A(\leq \lambda) := \bigoplus_{\lambda'\leq \lambda} {}_A \overset{\raisebox{-1.5pt}{\scriptsize$\smallbullet$}}{C} (\lambda') ,\qquad \Oo_A(< \lambda) := \bigoplus_{\lambda'< \lambda} {}_A \overset{\raisebox{-1.5pt}{\scriptsize$\smallbullet$}}{C} (\lambda').\]
In particular, in the ${\mathfrak{sl}_2}$ case the $A$-module ${}_A\! \overset{\raisebox{-1.5pt}{\scriptsize$\smallbullet$}}{C} (n\varpi_1)$ has basis $\dot{\mathbf{B}} [n]^*$ given above, of cardinality~${(n+1)^2}$.

Recall that $\dot{\mathbf{U}} [\geq \lambda]$ and $\dot{\mathbf{U}} [> \lambda]$ are two-sided ideals of $\dot{\mathbf{U}}$, and the algebra (whence $U_q^{\rm ad}$-bimodule) isomorphism \smash{$\bar{\pi}_\lambda\colon \dot{\mathbf{U}} [\geq \lambda]/ \dot{\mathbf{U}} [> \lambda] \ra {\rm End}(V_\lambda)$} (see~\eqref{pimuiso}). In \cite[Section~29.3]{Lusztig}, Lusz\-tig groups this isomorphism and its properties under the general term of {\it refined Peter--Weyl theorem}. We are going to reinterpret it in terms of $\Oo_A$. First observe that

\begin{lem} The $A$-modules $\Oo_A(\leq \lambda)$ and $\Oo_A(< \lambda)$ are $U_A^{\rm res}$-bimodules, and the surjective map
\begin{gather}\label{dlambdamap}
d_\lambda\colon\ \Oo_A(\leq \lambda)\longrightarrow {\rm Hom}\bigl( \dot{\mathbf{U}}_A [\geq \lambda]/ \dot{\mathbf{U}}_A [> \lambda],A\bigr),\qquad \alpha\longmapsto\langle \alpha , \cdot \ \rangle
\end{gather}
descends to an isomorphism of $U_A^{\rm res}$-bimodules $\bar{d}_\lambda$ on $\Oo_A(\leq \lambda)/\Oo_A(< \lambda)$.
\end{lem}
\begin{proof} For every $\alpha\in \Oo_A(\leq \lambda)$, $x,y\in U_A^{\rm res}$, and $b\in \dot{\mathbf{B}} [\mu]$ with $\mu \nleqslant \lambda$, we have $xby \in \dot{\mathbf{U}}_A [\geq \mu]$. Since $\dot{\mathbf{U}}_A [\geq \mu]= \bigoplus_{\eta\geq \mu}\ A\dot{\mathbf{B}}[\eta]$ and $\eta\geq \mu$ implies $\eta \nleqslant \lambda$, it follows that $\langle xby, \alpha\rangle=0$, i.e., $(x\rhd \alpha \lhd y)(b)=0$. This shows $x\rhd \alpha \lhd y \in \Oo_A(\leq \lambda)$. The same proof applies as well to $\Oo_A(<\lambda)$, whence the first claim. Since $\dot{\mathbf{U}} [\geq \lambda]$ and $\dot{\mathbf{U}} [> \lambda]$ are two-sided ideals of $\dot{\mathbf{U}}$, $\dot{\mathbf{B}}$ is a basis of~$\dot{\mathbf{U}}_A$, and the $A$-modules $\dot{\mathbf{U}}_A [\geq \lambda]$ and $\dot{\mathbf{U}}_A [> \lambda]$ are spanned by \smash{$ \coprod_{\lambda'\geq \lambda} \dot{\mathbf{B}} [\lambda']$} and \smash{$ \coprod_{\lambda'> \lambda} \dot{\mathbf{B}} [\lambda']$}, both are two-sided ideals of $\dot{\mathbf{U}}_A$, and $\dot{\mathbf{U}}_A [\geq \lambda]/ \dot{\mathbf{U}}_A [> \lambda]$ inherits the quotient $U_A^{\rm res}$-bimodule structure. Clearly, the map $d_\lambda$ is well defined, it is a morphism of $U_A^{\rm res}$-bimodules, and its kernel contains $\Oo_A(< \lambda)$. Bijectivity of $\bar{d}_\lambda$ comes by comparing the cardinality of canonical bases: $\Oo_A(\leq \lambda)/\Oo_A(< \lambda)$ has the basis formed by the cosets of the elements of the basis $(\dot{\mathbf{B}}[\lambda])^*$ of~\smash{${}_A \overset{\raisebox{-1.5pt}{\scriptsize$\smallbullet$}}{C} (\lambda')$}, and $\dot{\mathbf{U}}_A [\geq \lambda]/ \dot{\mathbf{U}}_A [> \lambda]$ the basis formed by the cosets of the elements of $\dot{\mathbf{B}}[\lambda]$, all cosets being non-zero and pairwise distinct.
\end{proof}

Since $\dot{\mathbf{U}}_A$ preserves the canonical basis $\underline{\mathbf{B}}_\lambda$ of ${}_AV_\lambda$, $\bar{\pi}_\lambda$ descends to an isomorphism of $U_A^{\rm res}$-bimodules \smash{$\bar{\pi}_\lambda\colon \dot{\mathbf{U}}_A [\geq \lambda]/ \dot{\mathbf{U}}_A [> \lambda] \ra {\rm End}({}_A V_\lambda)$}. We thus get exact sequences of $U_A^{\rm res}$-bi\-mod\-ules%
\[\xymatrix{0 \ar[r] & \dot{\mathbf{U}}_A [> \lambda] \ar[r] & \dot{\mathbf{U}}_A [\geq \lambda] \ar[r]^{\bar{\pi}_\lambda\quad \quad } & {\rm End}({}_A V_\lambda) \longrightarrow 0}\]
and
\begin{equation}\label{exactsequence}
\xymatrix{0 \ar[r] & \Oo_A(< \lambda) \ar[r] & \Oo_A(\leq \lambda) \ar[rr]^{(\bar{\pi}_\lambda^{-1})^*\circ d_\lambda\quad \ } & & \left({\rm End}({}_A V_\lambda)\right)^* \ar[r] & 0.}
\end{equation}
They split as sequences of $A$-modules but not as sequences of bimodules. In fact,
\begin{align}
\left({\rm End}({}_A V_\lambda)\right)^*:={}& {\rm Hom}({\rm End}({}_A V_\lambda),A )\nonumber\\
\cong{} &{\rm Hom}({}_A^{\omega} V_\lambda \bigotimes_A {}_A V_\lambda,A ) = {}_AV_\lambda^{\rm up} \bigotimes_A ({}_A^{\omega}V_\lambda)^{\rm up},\label{identcoeffmatup}
\end{align}
 with the ``\ ${}^{\rm up}$\ '' structure defined in \eqref{upmodule}, and corresponding basis $\underline{\mathbf{B}}_\lambda^{\rm up} \otimes ({}^{\omega}\underline{\mathbf{B}}_\lambda)^{\rm up}$. Moreover, the exact sequence \eqref{exactsequence} shows that this $A$-module of matrix coefficients, regarded as an $A$-submodule of $\Oo_A$ by means of the coefficient map $\Phi := \bigoplus_{\lambda\in P_+} \Phi_\lambda$ (see \eqref{pairingOAKashi}), is contained in $\Oo_A(\leq \lambda)$. This for all $\lambda'\leq \lambda$ yields $ \bigoplus_{\lambda'\leq \lambda} \left({\rm End}({}_A V_{\lambda'})\right)^* \subset \Oo_A(\leq \lambda)$. Now, using the isomorphism~$\bar{\pi}_\lambda$, we get
\[{\rm rank}_A(\Oo_A(\leq \lambda)) = \sum_{\lambda'\leq \lambda} {\rm Card}\bigl( \dot{\mathbf{B}}[\lambda']\bigr) = \sum_{\lambda'\leq \lambda} {\rm rank}({}_AV_{\lambda'})^2\]
and therefore
\begin{equation}\label{rankOqleq}
{\rm dim}_{\mc(q)}\bigl(\Oo_A(\leq \lambda) \bigotimes_A \mc(q)\bigr) = \sum_{\lambda'\leq \lambda} {\rm dim}(V_{\lambda'})^2 = \sum_{\lambda'\leq \lambda} {\rm dim}((C(\lambda')),
\end{equation}
where as usual $C(\lambda')$ denotes the space of matrix coefficients of $V_{\lambda'}$ (see \eqref{defCmu}). It follows
\begin{equation}\label{tensorCOA}
\Oo_A(\leq \lambda) \bigotimes_A {\mathbb C}(q)=\bigoplus_{\lambda'\leq \lambda}C(\lambda'), \qquad\Oo_A(<\lambda) \bigotimes_A {\mathbb C}(q)=\bigoplus_{\lambda'< \lambda}C(\lambda').
\end{equation}
However, in general \smash{${}_A \overset{\raisebox{-1.5pt}{\scriptsize$\smallbullet$}}{C} (\lambda)\bigotimes_A {\mathbb C}(q)$} is not equal to $C(\lambda)$, \smash{${}_A \overset{\raisebox{-1.5pt}{\scriptsize$\smallbullet$}}{C} (\lambda)$} is not an $A$-sublattice of~$C(\lambda)$, and \smash{${}_A \overset{\raisebox{-1.5pt}{\scriptsize$\smallbullet$}}{C} (\lambda)$} is not a $U_A^{\rm res}$-bimodule (it is because of this discrepancy that we have introduced the dot notation ``\ \smash{${}^{\smallbullet}$}\ ''). For instance, we can see the first two facts in the case of $\mathfrak{g} = {\mathfrak{sl}_2}$, by inverting the system of identities \eqref{formulesl2} for all $0\leq l$, $m\leq n$ (or more simply by considering the identity $v_0^* \otimes v_0=[2]_qbc +1$ from \eqref{matrixcoefV2}). For the third fact, we have ${1_2 E\in \dot{\mathbf{B}} [2]}$ (see~\eqref{2sidedcellsl2}), so~${((1_2 E)^* \lhd E)(1_0) \!=\!\langle \Delta((1_2 E)^* ), E\otimes 1_0\rangle\!= \! \langle(1_2 E)^* , E1_0\rangle\!=\!\langle(1_2 E)^*, 1_2E\rangle\!=\!1}$ since $E1_0\! =\! 1_2E$. Therefore, \smash{$(1_2 E)^* \lhd E \notin {}_A \overset{\raisebox{-1.5pt}{\scriptsize$\smallbullet$}}{C} (2)$}.

From the formulas \eqref{formulesl2} and Appendix~\ref{lowupbases}, we can observe the isomorphism~\eqref{identcoeffmatup} in the case of $\mathfrak{g} = {\mathfrak{sl}_2}$. More simply, by projecting the matrix \eqref{matrixcoefV2} onto~$({\rm End}({}_A \!V_2))^*$ the entries are unchanged except the $(1,1)$ entry, which becomes $[2]_q bc$. All factors $[2]_q$ in the middle column disappear if one uses matrix coefficients in the upper canonical basis of~$V_2$, which is $v_0^{\rm up}:=v_0$, $v_1^{\rm up}:=[2]_q^{-1}v_1$, $v_2^{\rm up}:=v_2$ in the notations of \eqref{formulesl2}, since we have
\smash{$v^l(\pi_2(\cdot) \ v_m) = [\delta_{m,1}+1]_q \big\langle v_l^{\rm up}, \cdot v_m^{\rm up}\big\rangle$} for $l,m\in \{0,1,2\}$, where $\langle\ ,\ \rangle$ is the pairing \eqref{proppairinglambda}. Thus, in this particular example of \smash{$({\rm End}({}_A V_2))^*$} we see explicitly the identification of the basis $(\bar{\pi}_2^*)^{-1} \circ d_2(\dot{\mathbf{B}} [2]^*)$ and~${\underline{\mathbf{B}}_2^{\rm up} \otimes ({}^{\omega}\underline{\mathbf{B}}_2)^{\rm up}}$.

Summing up this discussion, the Lusztig refined Peter--Weyl theorem of \cite[Section~29.3]{Lusztig}, implies the following.

\begin{teo} As an $A$-module we have a direct sum decomposition
\begin{equation}\label{OAsurA} \Oo_A = \bigoplus_{\lambda\in P_+} {}_A \overset{\raisebox{-1.5pt}{\scriptsize$\smallbullet$}}{C} (\lambda),
\end{equation}
as $U_A^{\rm res}$-bimodules we have a $($directed by inclusion, and non direct$)$ sum
\begin{equation}\label{OAbimod}
\Oo_A = \sum_{\lambda\in P_+} \Oo_A(\leq \lambda),
\end{equation}
and the composition factors of $\Oo_A$ are the bimodules
\begin{equation}\label{OAfacteurs}
\left({\rm End}({}_A V_\lambda)\right)^* \cong \left({}_A^{\omega} V_\lambda \otimes {}_A V_\lambda \right)^*
\end{equation}
for every $\lambda\in P_+$, each of multiplicity $1$. \end{teo}
\begin{Remark}{\rm The above filtration and its composition factors appear in disguised manner as {\it good filtration} in \cite{APW} and \cite{Paradowski} (see also \cite{VV}).}
\end{Remark}
Because $\dot{\mathbf{B}}$ is formed by weight vectors for the left and right action of $U_q^{\rm ad}(\mathfrak{h})$ (see \eqref{Bpointunion}), the same is true of \smash{$\dot{\mathbf{B}}^*$} and \eqref{OAsurA} can thus be refined into a weight space decomposition
\begin{equation}\label{OAweightdecomp} \Oo_A = \bigoplus_{\mu,\nu\in P}\bigoplus_{\lambda\in P_+} \big({}_A \overset{\raisebox{-1.5pt}{\scriptsize$\smallbullet$}}{C} (\lambda)\big)_{\mu,\nu}.
\end{equation}

Now recall the property \eqref{Bpointunion}. Consider in particular the finite subsets $\dot{\mathbf{B}}_{0,\varpi_i}$ and \smash{$\dot{\mathbf{B}}_{\varpi_i,0}$} associated to the fundamental weights $\varpi_i$, $i=1,\dots, m$. The map \smash{$u\mapsto u({}^\omega v_0\otimes v_{w_0(\varpi_i)})$}, \smash{$u\in \dot{\mathbf{U}}$}, allows one to identify $\dot{\mathbf{B}}_{0,\varpi_i}$ with the canonical basis $\underline{\mathbf{B}}_{\varpi_i}$ of ${}^\omega V_{0}\otimes V_{\varpi_i} \cong V_{\varpi_i}$, and therefore with a uniquely determined finite subset $\mathbf{B}_{\varpi_i}$ of the canonical basis $\mathbf{B}^-$ of $U_q^{\rm ad}(\mathfrak{n}_-)$; similarly, one can identify \smash{$\dot{\mathbf{B}}_{\varpi_i,0}$} with a uniquely determined finite subset ${}^\omega \mathbf{B}_{\varpi_i}$ of the canonical basis $\mathbf{B}^+$ of $U_q^{\rm ad}(\mathfrak{n}_+)$. The elements of $\dot{\mathbf{B}}_{0,\varpi_i}$ and $\dot{\mathbf{B}}_{\varpi_i,0}$ are respectively of the form $b^-1_{\varpi_i}$ and $b^+1_{-\varpi_i}$, where $b^- \in \mathbf{B}_{\varpi_i}$ and $b^+\in {}^\omega \mathbf{B}_{\varpi_i}$, and we have (see \cite[Proposition~3.3 and Section~3.4]{Lusztig3}):

\begin{prop}\label{genOA} The algebra $\Oo_A$ is finitely generated. A system of generators is provided by the elements $a^*\in \dot {\mathbf{B}} {}^*$, where $a \in \bigcup_{i=1}^m (\dot{\mathbf{B}}_{0,\varpi_i}\cup \dot{\mathbf{B}}_{\varpi_i,0})$.\end{prop}
Note that the above system of generators of $\Oo_A$ has
$2 \sum_{i=1}^m \dim(V_{\varpi_i})$
 elements. In fact, recall that $\varphi\colon U_q^{\rm ad} \ra U_q^{\rm ad}$ is the anti-automorphism given by $\varphi(E_i) = F_i$, $\varphi(F_i) = E_i$, $\varphi(K_\lambda) = K_{\lambda}$. Denote by $v_{-\varpi_i}$ and $f_{-\varpi_i}$ the canonical lowest-weight vectors of the highest weight modules~${V_{-w_0(\varpi_i)}}$ and $V_{-w_0(\varpi_i)}^r$, respectively, and put the superscript `` ${}^{\rm up}$ '' for the upper canonical basis vectors.
\begin{lem} For every $b^- \in \mathbf{B}_{\varpi_i}$ and $b^+\in {}^\omega \mathbf{B}_{\varpi_i}$, we have
\begin{align}
 &(b^-1_{\varpi_i})^* = \Phi_{\varpi_i}((f_{\varpi_i} \varphi(b^-))^{\rm up} \otimes v_{\varpi_i}), \label{formgenOA1} \\
 &\big(b^+1_{-\varpi_i}\big)^* = \Phi_{-w_0(\varpi_i)}\big(\big(f_{-\varpi_i} \varphi\big(b^+\big)\big)^{\rm up} \otimes v_{-\varpi_i}\big). \label{formgenOA2} \end{align}
In other words, $(b^-1_{\varpi_i})^*$ and $\big(b^+1_{-\varpi_i}\big)^*$ are the matrix coefficients lying on the first and last columns of the matrix representations in the upper canonical bases of the spaces $V_{\varpi_i}$, ${i=1,\dots, m}$.
\end{lem}

\begin{proof} This can be checked by using the isomorphism \eqref{isocrystal}. The key observation is that
 \[\langle \Phi_\lambda(f_\lambda \otimes v_\lambda), 1_\mu\rangle = \langle f_\lambda , 1_\mu v_\lambda \rangle_\lambda = \delta_{\lambda,\mu}\]
 for every $\lambda\in P_+$, $\mu\in P$, and therefore $\Phi_\lambda(f_\lambda \otimes v_\lambda) = 1_\lambda^*$. Then the computation proceeds by using the equivariance of $\Phi$ under the action of \smash{$U_q^{\rm ad} \otimes (U_q^{\rm ad}){}^{\rm op}$}, the fact that $\langle\cdot,\cdot \rangle$ dualizes the bimodules structures on $\Oq$ and $\dot{\mathbf{U}}$, and the description of the associated Kashiwara operators on $\mathcal{B}(\Oo_q)$ and \smash{$\mathcal{B}\bigl(\dot{\mathbf{U}}\bigr)$}. Here is an alternative argument. By the very definition of the sets~\smash{$\dot{\mathbf{B}} [\lambda]$} we have \smash{$b^-1_{\varpi_i}\in \dot{\mathbf{B}} [\varpi_i]$}, \smash{$b^+1_{-\varpi_i}\in \dot{\mathbf{B}} [-\omega_0(\varpi_i)]$}. We wish to check if their duals~$(b^-1_{\varpi_i})^*$,~$(b^+1_{-\varpi_i})^*$ coincide with the elements of $\Oo_A$ on the right sides of \eqref{formgenOA1} and \eqref{formgenOA2}. As already noticed after~\eqref{exactsequence}, by the isomorphism $\Oo_A(\leq \lambda)/\Oo_A(<\lambda) \cong {\rm End}({}_AV_\lambda)^*$ every matrix coefficient of~${}_AV_\lambda$ belongs to $\Oo_A(\leq \lambda)$. Now, the $A$-modules $\Oo_A(\leq \varpi_i)$ and $\Oo_A(\leq -w_0(\varpi_i))$ are generated by~${\dot{\mathbf{B}}[\varpi_i]^*}$ and $\dot{\mathbf{B}} [-w_0(\varpi_i)]^*$, respectively. Because $((\bar{\pi}_{\lambda}^*)^{-1} \circ d_{\lambda})(\dot{\mathbf{B}} [\lambda]^*)$ coincides with~${\underline{\mathbf{B}}_{\lambda}^{\rm up} \otimes ({}^{\omega}\underline{\mathbf{B}}_{\lambda})^{\rm up}}$, the conclusion follows.
 \end{proof}

Note that the same argument implies that, for every $\lambda\in P_+$, any matrix coefficient of $V_\lambda$ in the upper canonical basis and vanishing on the elements of $\dot{\mathbf{B}} [\lambda']$ for $\lambda'<\lambda$ must belong to~$\dot{\mathbf{B}} [\lambda]^*$. For instance, in the ${\mathfrak{sl}_2}$ case, $\Oo_A(\leq 2)$ has canonical basis $\dot{\mathbf{B}} [0]^* \coprod \dot{\mathbf{B}} [2]^*$, so the matrix coefficients of $V_2$ vanishing on $1_0$ belong to $\dot{\mathbf{B}} [2]^*$. This can be observed in \eqref{matrixcoefV2}, using the comments in the paragraph before \eqref{OAsurA}.

Though the $A$-module ${}_AV_\mu \bigotimes_A {}_AV_\nu$ has no decomposition like \eqref{decomprep}, we can refine the map $C(\mu) \otimes C(\nu) \ra C(\mu+\nu)$ in \eqref{projOq} to an $A$-linear map defined on \smash{${}_A \overset{\raisebox{-1.5pt}{\scriptsize$\smallbullet$}}{C} (\mu) \bigotimes_A {}_A \overset{\raisebox{-1.5pt}{\scriptsize$\smallbullet$}}{C} (\nu)$}. Indeed, there is a unique injective morphism of $U_A^{\rm res}$-modules $\mathfrak{T}_{\mu,\nu}\colon{}_AV_{\mu+\nu}\ra {}_AV_\mu \bigotimes_A {}_AV_\nu$, which is given by $\mathfrak{T}_{\mu,\nu}(v_{\mu+\nu}) = v_{\mu} \otimes v_{\nu}$ \cite[Proposition 25.1.2\,(a)--(b)]{Lusztig}. It defines a morphism of based modules
\[(V_{\mu+\nu},\underline{\mathbf{B}}_{\mu+\nu}) \ra (V_\mu \otimes V_\nu ,\underline{\mathbf{B}}_{\mu} \,\lozenge\, \underline{\mathbf{B}}_\nu),
\]
where $\underline{\mathbf{B}}_{\mu} \,\lozenge\, \underline{\mathbf{B}}_\nu := \{b\,\lozenge\, b', b\in \underline{\mathbf{B}}_{\mu}, b'\in \underline{\mathbf{B}}_{\nu} \}$ \cite[Proposition~27.1.7]{Lusztig}. Hence, $\mathfrak{T}_{\mu,\nu}$ is a~split $A$-li\-near map, i.e., there exists a $A$-linear map $\mathfrak{S}_{\mu,\nu}\colon {}_AV_\mu \bigotimes_A {}_AV_\nu \ra {}_AV_{\mu+\nu} $ such that ${\mathfrak{S}_{\mu,\nu}\circ \mathfrak{T}_{\mu,\nu} ={\rm id}}$. Note that $\mathfrak{S}_{\mu,\nu}$ is not a $U_A^{\rm res}$-morphism. Similarly, the unique morphism of $U_A^{\rm res}$-modules \smash{${}^\omega\mathfrak{T}_{\mu,\nu}\colon{}_A^\omega V_{\mu+\nu}\ra {}_A^\omega V_\mu \bigotimes_A {}_A^\omega V_\nu$} is a split injection. Define \smash{$\rho_{\mu',\mu''}\colon \dot{\mathbf{U}}_A \ra {}_A^\omega V_{\mu'} \bigotimes_A {}_AV_{\mu''}$} by
\[
\rho_{\mu',\mu''}(u)= u\biggl({}^\omega v_{w_0(\mu')}\bigotimes_A v_{w_0(\mu'')}\biggr),\]
and $\rho_{\mu',\mu'',\nu',\nu''}\colon \dot{\mathbf{U}}_A{}^{ \hat{\otimes}2} \ra {}_A^\omega V_{\mu'} \bigotimes_A {}_AV_{\mu''} \bigotimes_A {}_A^\omega V_{\nu'} \bigotimes_A {}_AV_{\nu''}$ by
\[\rho_{\mu',\mu'',\nu',\nu''}(u)= u\biggl({}^\omega v_{w_0(\mu')}\bigotimes_A v_{w_0(\mu'')}\bigotimes_A {}^\omega v_{w_0(\nu')}\bigotimes_A v_{w_0(\nu'')}\biggr).\]
Define $\tau_{\mu',\mu'',\nu',\nu''}\colon {}_A^\omega V_{\mu'+\nu'} \bigotimes_A {}_AV_{\mu''+\nu''} \ra {}_A^\omega V_{\mu'} \bigotimes_A {}_AV_{\mu''} \bigotimes_A {}_A^\omega V_{\nu'} \bigotimes_A {}_AV_{\nu''}$ by
\[\tau_{\mu',\mu'',\nu',\nu''} = \big(1\otimes \hat{R}^{-1} \otimes 1\big)({}^\omega \mathfrak{T}_{\mu',\nu'} \otimes \mathfrak{T}_{\mu'',\nu''}).\]
It is an injective morphism of $U_A^{\rm res}$-modules. In \cite[Section 1.13]{Lusztig3}, Lusztig proved that $\tau_{\mu',\mu'',\nu',\nu''}$ is a split $A$-linear map (\cite{Lusztig3} uses $\hat{R}$ instead of $\hat{R}^{-1}$, since our coproducts on $U_q^{\rm ad}$ are opposite), and that it satisfies{\samepage
\begin{equation}\label{commutdiagLusztig}
\tau_{\mu',\mu'',\nu',\nu''} \rho_{\mu'+\mu'',\nu'+\nu''} = \rho_{\mu',\mu'',\nu',\nu''}\Delta,
\end{equation}
where $\Delta$ is the coproduct of $\dot{\mathbf{U}}_A$, see \eqref{defcoprod}.}

Now take $\mu := \mu'=\mu''$, $\nu:= \nu'=\nu'' \in P_+$, and put $\tau_{\mu,\nu} := \tau_{\mu,\mu,\nu,\nu}$. It follows from the classical decomposition~\eqref{decomprep} over $\mc(q)$, and \eqref{projOq} and \eqref{tensorCOA}, that the product of $\Oo_A$ yields a~map $m\colon \Oo_A(\leq \mu) \bigotimes_A \Oo_A(\leq \nu) \ra \Oo_A(\leq \mu+\nu)$.

Denote the projection map \smash{$p_{\mu+\nu}\colon \Oo_A(\leq \mu+\nu) \ra {}_A \overset{\raisebox{-1.5pt}{\scriptsize$\smallbullet$}}{C} (\mu+\nu)$}, define \smash{${}_A \overset{\smallbullet}{\tau}_{\mu,\nu}:= p_{\mu+\nu} \circ m$}, and put
\[\xymatrix{\pi_\lambda'\colon\ \Oo_A(\leq \lambda)\ar[r] & \Oo_A(\leq \lambda)/\Oo_A(< \lambda)\ar[rr]^{\quad (\bar{\pi}_\lambda^*)^{-1}\circ \bar{d}_\lambda} & & ({\rm End}({}_A V_\lambda) )^*},\]
where the first map is the quotient map. Consider the diagram
\[\xymatrix{ {}_A \overset{\raisebox{-1.5pt}{\scriptsize$\smallbullet$}}{C} (\mu) \otimes {}_A \overset{\raisebox{-1.5pt}{\scriptsize$\smallbullet$}}{C} (\nu) \ar[d]^{ \pi_\mu' \otimes \pi_\nu'} \ar[r]^{{}_A \overset{\smallbullet}{\tau}_{\mu,\nu}} & {}_A \overset{\raisebox{-1.5pt}{\scriptsize$\smallbullet$}}{C} (\mu+\nu) \ar[d]^{\pi_{\mu+\nu}'} \\ ({\rm End}({}_A V_\mu) )^* \otimes ({\rm End}({}_A V_\nu) )^* \ar[r]^{\hspace*{1cm} \tau_{\mu,\nu}^t} & ({\rm End}({}_A V_{\mu +\nu}) )^*,}\]
where $\tau_{\mu,\nu}^t$ is the transpose of Lusztig's map $\tau_{\mu,\nu}$.
\begin{prop} \label{teoLuzstigOAscinde} The map \smash{${}_A \overset{\smallbullet}{\tau}_{\mu,\nu} \colon {}_A \overset{\raisebox{-1.5pt}{\scriptsize$\smallbullet$}}{C} (\mu) \bigotimes_A {}_A \overset{\raisebox{-1.5pt}{\scriptsize$\smallbullet$}}{C} (\nu)\ra {}_A \overset{\raisebox{-1.5pt}{\scriptsize$\smallbullet$}}{C} (\mu+\nu)$} is split as an $A$-linear map and the above diagram is commutative.
\end{prop}
\begin{proof} The commutativity of the diagram comes from equation~\eqref{commutdiagLusztig}. The epimorphism $\pi_\lambda'$ is injective on \smash{${}_A \overset{\raisebox{-1.5pt}{\scriptsize$\smallbullet$}}{C} (\lambda)$}, and maps the canonical basis elements to the elements of the upper canonical basis $\underline{\mathbf{B}}_\lambda^{\rm up} \otimes ({}^{\omega}\underline{\mathbf{B}}_\lambda)^{\rm up}$. By Lusztig's results recalled above, the epimorphism $\tau_{\mu,\nu}^t$ splits as an $A$-linear map. Therefore, the same is true of \smash{${}_A \overset{\smallbullet}{\tau}_{\mu,\nu}$}.
\end{proof}

We stress that \smash{${}_A \overset{\smallbullet}{\tau}_{\mu,\nu}$} plays for $\Oo_A$ the same role as the map \eqref{projOq} for $\Oo_q$.

Finally, we consider for any $n\geq 1$ the invariant elements of $\Oo_A^{\otimes n}$ endowed with the action~${\rm coad}^r_n$ of $U_A^{\rm res}$, see~\eqref{actionprod} \big(recall that $\Ll_{0,n}=\Oo_q^{\otimes n}$ as $U_q^{\rm ad}$-module\big).

First note that, by definition, $\Oo_A(G^n)$ is the restricted dual of the Hopf algebra $U_A^{\rm res}\big(\mathfrak{g}^{\oplus n}\big)$, associated to its category of type $1$ modules. By ordering the summands of $\mathfrak{g}^{\oplus n}$ we get an isomorphism $U_A^{\rm res}\big(\mathfrak{g}^{\oplus n}\big)\cong U_A^{\rm res}(\mathfrak{g})^{\otimes n}$, and any type $1$ simple $U_A^{\rm res}(\mathfrak{g})^{\otimes n}$-module is isomorphic to ${V_{[\lambda]} := \bigotimes_{i=1}^n V_{\lambda_i}}$ endowed with the componentwise action, for some ${[\lambda] :=(\lambda_1,\dots, \lambda_n)\in P_+^n}$ (this is a~classical fact; see, e.g., \cite[Theorem~3.10.2]{Et}). Therefore, we have an isomorphism ${\Oo_A(G^n)\cong \Oo_A^{\otimes n}}$.
With the same notation $[\lambda]:=(\lambda_1,\dots, \lambda_n)\in P_+^n$, let us put
\begin{align*}
&{}_A \overset{\raisebox{-1.5pt}{\scriptsize$\smallbullet$}}{C} ([\lambda]) := \bigotimes_{i=1}^n {}_A \overset{\raisebox{-1.5pt}{\scriptsize$\smallbullet$}}{C} (\lambda_i) = \bigoplus_{b\in \bigotimes_{i=1}^n \dot{\mathbf{B}}[\lambda_i]^*} Ab,\\
&\Oo_A(\leq [\lambda]) := \bigotimes_{i=1}^n \Oo_A(\leq \lambda_i)=\bigoplus_{[\lambda'] \in P_+^n, \lambda_i'\leq \lambda_i} {}_A \overset{\raisebox{-1.5pt}{\scriptsize$\smallbullet$}}{C} ([\lambda']).
\end{align*}
We thus obtain a decomposition into based $(U_A^{\rm res}\otimes (U_A^{\rm res})^{\rm op})^{\otimes n}$-modules
\[
\Oo_A^{\otimes n} = \sum_{[\lambda] \in P_+^n} \Oo_A(\leq [\lambda]).
\]
Now ${\rm coad}_n^r = ({\rm coad}^r)^{\otimes n}\circ \Delta^{(n-1)}$ gives structures of $U_A^{\rm res}$-modules to $\Oo_A^{\otimes n}$ and $\Oo_A(\leq [\lambda])$. In order to make it a based module, we give it the ``$\lozenge$'' product of the canonical bases of the factors~${\Oo_A(\leq \lambda_i)}$, i.e.,
\[\dot{\mathbf{B}} [[\lambda]]^* := \lozenge_{i=1}^n \bigg(\coprod_{\lambda_i'\leq \lambda_i}\dot{\mathbf{B}} [\lambda_i']^*\bigg).\]
We thus obtain a decomposition into based $U_A^{\rm res}$-modules
\begin{equation}\label{decompositionOAcoad}
\Oo_A^{\otimes n} = \sum_{[\lambda] \in P_+^n} (\Oo_A(\leq [\lambda]), \dot{\mathbf{B}} [[\lambda]]^*),
\end{equation}
with composition factors $ \bigotimes_{i=1}^n \left({\rm End}({}_A V_{\lambda_i})\right)^*$. By the properties of ``$\lozenge$'' products of bases of based modules, the underlying $A$-module is
\begin{equation}\label{decompositionOA2}
\Oo_A^{\otimes n} = \bigoplus_{[\lambda] \in P_+^n} {}_A \overset{\raisebox{-1.5pt}{\scriptsize$\smallbullet$}}{C} ([\lambda]).
\end{equation}
Finally, we state the last property of based modules we need. Let $(M,B)$ be a based module. Recall the notations introduced around \eqref{decbaseisospace}. It is proved in \cite[Proposition 27.1.8]{Lusztig} that for every $\lambda\in P_+$ the submodule $M[\geq\lambda]$ is a sub-based module of~$M$, and that it has the basis
\begin{equation}\label{partbasemodule}
B\cap M[\geq\lambda] = \bigcup_{\lambda'\geq \lambda} B[\lambda'].
\end{equation}
Consider $M[\ne 0] := \bigoplus_{\lambda\ne 0} M[\lambda]$, the largest proper submodule of $M$ that contains no non-zero invariant element. Recall that the space of {\it coinvariants of} $M$ is
\begin{align*} M_{U_q^{\rm ad}} & = M/ M[\ne 0] = M/\mc(q)\big\{um-\varepsilon(u)m, \,m\in M, \,u\in U_q^{\rm ad}\big\}
\end{align*}
that is, the largest quotient of $M$ with trivial action, where $\varepsilon\colon U_q^{\rm ad}\ra \mc(q)$ is the counit. It follows from \eqref{partbasemodule} that $M[\ne 0]$ is a sub-based module of $M$, with the basis $\bigcup_{\lambda\ne 0} B[\lambda]$, and we have (this is, \cite[Proposition~27.2.6]{Lusztig}):
\begin{prop}\label{Linvbased} The quotient map $\pi\colon M \ra M_{U_q^{\rm ad}}$ is a morphism of based modules, where $M_{U_q^{\rm ad}}$ is endowed with the basis $B_{U_q^{\rm ad}} := \pi(B[0])$.
\end{prop}
Keeping the same notations, let ${}_AM \subset M$ be the $A$-module generated by $B$, and let ${{}_AM^* \subset M^*}$ be the $A$-module generated by $B^*$. They are $U_A^{\rm res}$-modules. Denote by $({}_A M^*)^{U_A^{\rm res}}$ the submodule of $U_A^{\rm res}$-invariant elements of ${}_A M^*$, regarded as a right module in the natural way.
\begin{lem} We have a direct sum decomposition of $A$-modules
\begin{equation}\label{decompAcoinv}
{}_A M^* = ({}_A M^*)^{U_A^{\rm res}} \bigoplus_A {}_A N,
\end{equation}
where ${}_A N\subset {}_A M^*$ is the $A$-submodule generated by $\bigcup_{\lambda\ne 0} B[\lambda]^*$.
\end{lem}
\begin{proof} By Proposition~\ref{Linvbased}, the transpose map \smash{$\pi^t\colon (M_{U_q^{\rm ad}})^* \ra M^*$} is a monomorphism mapping the dual basis \smash{$B_{U_q^{\rm ad}}^*$} to the subset $B[0]^*$ of $B^*$. The image of $\pi^t$ is \smash{$(M^*)^{U_q^{\rm ad}}$}. If we set \smash{${}_A M_{U_A^{\rm res}} = \pi({}_A M)$}, then \smash{$\pi^t(({}_A M_{U_A^{\rm res}})^*) = ({}_A M^*)^{U_A^{\rm res}}$} is generated by $B[0]^*$, which concludes the proof.
\end{proof}

Note that, since $B[0]$ is in general not invariant under the action of $U_A^{\rm res}$, ${}_A N$ need not be stable under this action.

We are now ready to draw consequences of this discussion and the previous results. As usual denote by \smash{$\big(\Oo_A^{\otimes n}\big)^{U_A^{\rm res}}$} the subspace of invariant elements of $\Oo_A^{\otimes n}$ for the action ${\rm coad}_n^r$. In the case $n=1$, it is just the center $\mathcal{Z}(\Oo_A)$.

\begin{teo}\label{teoLuzstigOAinv} \smash{$\big(\Oo_A^{\otimes n}\big)^{U_A^{\rm res}}$} is a direct summand of the $A$-module $\Oo_A^{\otimes n}$ for any $n\geq 1$.
\end{teo}

\begin{proof} By equation \eqref{decompositionOAcoad}, it is enough to show that for every $[\lambda]\in P_+^n$ the invariant elements of $\Oo_A(\leq [\lambda])$ form a direct summand, and these summands are compatible with non-empty intersections $\Oo_A(\leq [\lambda]) \cap \Oo_A(\leq [\lambda'])$. Using that $\Oo_A(G^n) \cong \Oo_A^{\otimes n}$ and viewing $P_+^n$ as the weight lattice of $G^n$, it is enough to prove these claims for $n=1$. Given $\lambda\in P_+$ put
\[P_{\lambda} = \{\lambda'\in P_+, \, \lambda' \nleqslant \lambda\},\]
and denote by $\dot{\mathbf{U}}_A [P_\lambda]$ the $A$-submodule of $\dot{\mathbf{U}}_A$ generated by $\coprod_{\lambda' \in P_\lambda} \dot{\mathbf{B}} [\lambda']$. Also, let us put \smash{$\dot{\mathbf{U}} [P_\lambda]= \dot{\mathbf{U}}_A [P_\lambda]\bigotimes_A \mc(q)$}.
The complement $P_+\setminus P_{\lambda}$ is finite, and if $\lambda'\in P_{\lambda}$ and $\lambda''\geq \lambda'$, then $\lambda''\in P_{\lambda}$. By the results of \cite[Section~29.2]{Lusztig}, $ \dot{\mathbf{U}} [P_\lambda]$ is a two-sided ideal, and the quotient algebra $\dot{\mathbf{U}} / \dot{\mathbf{U}} [P_\lambda]$ is finite-dimensional with unit the coset of \smash{$ \sum_{\lambda'\leq \lambda}1_{\lambda'}$}, and it is semisimple, isomorphic to \smash{$\bigoplus_{\lambda'\leq \lambda} {\rm End}(V_{\lambda'})$} (whereas $\dot{\mathbf{U}}_A / \dot{\mathbf{U}}_A [P_\lambda]$ has indecomposable modules, see Example~\ref{V2V20}).
It inherits from $\dot{\mathbf{U}}$ a canonical basis, formed by the non-zero cosets of elements of $\dot{\mathbf{B}}$, and with this basis $\dot{\mathbf{U}} / \dot{\mathbf{U}} [P_\lambda]$ is a based module for the right adjoint action~${\rm ad}^r$. Similarly as for~\eqref{dlambdamap}, we have a morphism of $U_A^{\rm res}$-modules
\[
\tilde{d}_\lambda\colon\ \Oo_A(\leq \lambda)\longrightarrow {\rm Hom} \big( \dot{\mathbf{U}}_A / \dot{\mathbf{U}}_A [P_\lambda],A\big),\qquad \alpha\longmapsto\langle \alpha , \cdot \ \rangle,\]
which is an isomorphism by \eqref{rankOqleq} and the computation \smash{$ {\rm dim}\bigl(\dot{\mathbf{U}} / \dot{\mathbf{U}} [P_\lambda]\bigr) = \sum_{\lambda'\leq \lambda} {\rm dim}(V_{\lambda'})^2$} in~\cite[Section~29.2]{Lusztig}. Applying Proposition \ref{Linvbased} and \eqref{decompAcoinv} to the based module ${M=\dot{\mathbf{U}} / \dot{\mathbf{U}} [P_\lambda]}$, we obtain that the invariant elements of $\Oo_A(\leq \lambda)$ form a direct summand. Finally, for any ${\lambda,\lambda'\!\in\! P_+}$ we have \smash{$\Oo_A(\leq \! \lambda) \cap \Oo_A(\leq \! \lambda')\cong {\rm Hom} \bigl( \dot{\mathbf{U}}_A \!/\bigl(\dot{\mathbf{U}}_A [P_\lambda]+ \dot{\mathbf{U}}_A [P_{\lambda'}]\bigr),A\bigr)$}.
 Applying Proposition \ref{Linvbased} and \eqref{decompAcoinv} to the based module $M:=\dot{\mathbf{U}} /(\dot{\mathbf{U}} [P_\lambda]+ \dot{\mathbf{U}} [P_{\lambda'}])$, we obtain that the invariant elements $({}_A M^*)^{U_A^{\rm res}}$ of $\Oo_A(\leq \lambda) \cap \Oo_A(\leq \lambda')$ form a direct $A$-summand. Since the latter is a~based $U_A^{\rm res}$-submodule of $ \Oo_A(\leq \lambda)$ and $\Oo_A(\leq \lambda')$, this summand is also a direct $A$-summand of $ \Oo_A(\leq \lambda)^{U_A^{\rm res}}$ and $\Oo_A(\leq \lambda')^{U_A^{\rm res}}$. This shows the $A$-modules $ \Oo_A(\leq \lambda)^{U_A^{\rm res}}$ for all $\lambda\in P_+$ match to form the $A$-summand $(\Oo_A)^{U_A^{\rm res}}$ of $\Oo_A$, and thus concludes the proof.
 \end{proof}

\begin{Remark} Let $(M,B)$, $(M', B')$ be based modules, with tensor product $(M\otimes M',B_{\lozenge})$, and $B_{\lozenge}[0]\subset B_{\lozenge}$ the subset in bijection with the canonical basis of the space of coinvariants~\smash{$(M\otimes M')_{U_q^{\rm ad}}$} (see Proposition \ref{Linvbased}). This subset is described in \cite[Proposition 27.3.8]{Lusztig} in terms of $B$ and $B'$. Since $\dot{\mathbf{U}} / \dot{\mathbf{U}} [P_{\lambda}]$ is semisimple with known summands, and the construction of the ``$\lozenge$'' product of canonical bases is associative, one can recursively compute the subset of the canonical basis of $\bigotimes_{i=1}^n \dot{\mathbf{U}} / \dot{\mathbf{U}} [P_{\lambda_i}]$ (endowed with the action dual to ${\rm coad}_n^r$) which is in bijection with the canonical basis of the space of coinvariants. Therefore, a complete (though highly nontrivial) characterization of the basis of \smash{$\big(\Oo_A^{\otimes n}\big)^{U_A^{\rm res}}$} can be obtained. Examples can be found in \cite[Section~27.3.10]{Lusztig}. In the case $\mathfrak{g}={\mathfrak{sl}_2}$, the canonical basis of the dual space~\smash{${\rm End}\big(V_1^{\otimes n}\big)^*$} has been identified in~\cite{FK} with the canonical basis of the Temperley--Lieb algebra $TL_n(q)$.
\end{Remark}

\begin{exa}\label{V2V20}{\rm The simplest case is already instructive. Namely, consider $V_1$ and $V_2$, the simple $U_q^{\rm ad}({\mathfrak{sl}_2})$-modules of type $1$ and dimension two and three.

On $V_1$, we have the lower canonical basis vectors $v_+$ and $v_-$, such that $Kv_+ = qv_+$, $Ev_+=0$, $v_- = Fv_+$. The canonical lower and upper bases of $V_1$ are both $\{v_+,v_-\}$. Using the relation~\eqref{udef}, we see that the elements of $\dot{\mathbf{B}}_{0,1}$ and $\dot{\mathbf{B}}_{1,0}$ are $1_1$, $F1_1$ and $1_{-1}$, $E1_{-1}$, respectively; the dual linear forms generate $\Oo_A({\rm SL}_2)$, they are the matrix coefficients $a$, $c$, $d$ and $b$ respectively. By \eqref{2sidedcellsl2}, we have $\dot{\mathbf{B}} [1] = \dot{\mathbf{B}}_{0,1}\coprod \dot{\mathbf{B}}_{1,0}$.

Next consider $V_2$. On $V_2$, we have the canonical highest weight vector $v_0$ of weight $2$, and lower canonical basis $\underline{\mathbf{B}}_{2}=\{v_0,v_1,v_2\}$, where $v_1=Fv_0$ and $v_2=F^{(2)}v_0$. We have ${\underline{\mathbf{B}}_{2}^{\rm up}=\big\{v_0,[2]_q^{-1}v_1,v_2\big\}}$ (see Appendix~\ref{lowupbases}). We can identify the ambient space of the right module~$V_2^r$ with that of~$V_2$; its highest weight vector is then~$v_0$, and its canonical lower and upper bases are $\underline{\mathbf{B}}_{2}^r=\{v_0,v_1,v_2\}$ and $\underline{\mathbf{B}}_{2}^r{}^{\rm up}=\big\{v_0,[2]_q^{-1}v_1,v_2\big\}$.

Consider now the module ${}^\omega V_1\otimes V_1$. We have
\[
\hat{R} = \sum_{n=0}^\infty \frac{\bigl(q-q^{-1}\bigr)^n}{[n]_q !} q^{n(n-1)/2} E^n\otimes F^n,
\]
 so the matrix of the involution $\Psi=\hat{R}^{-1}\circ \bar{}\ $ in the basis $v_+\otimes v_+, v_+\otimes v_-,v_-\otimes v_+,v_-\otimes v_-$ is
\[\big(\hat{R}^{-1}\circ \bar{}\ \big)_{{}^\omega V_1 ,V_1}=\begin{pmatrix} 1 &0 & 0 & 0\\0& 1& 0 &0\\0&0&1&0\\ q^{-1}-q &0&0& 1 \end{pmatrix}.\]
Therefore, the canonical basis $\underline{\mathbf{B}}_{1,1}$ is formed by the vectors $v_+ \,\lozenge\, v_+ =v_+\otimes v_+ + q^{-1} v_-\otimes v_-$ and $v_+ \,\lozenge\, v_- =v_+\otimes v_-$, $v_- \,\lozenge\, v_+ =v_-\otimes v_+ $, $v_- \,\lozenge\, v_- =v_-\otimes v_-$. Consider the partition $\underline{\mathbf{B}}_{1,1} = \underline{\mathbf{B}}_{1,1}[2] \cup \underline{\mathbf{B}}_{1,1}[0]$. We have $\underline{\mathbf{B}}_{1,1}[2] = \{v_- \,\lozenge\, v_+,v_+ \,\lozenge\, v_+, ,v_+ \,\lozenge\, v_-\}$, which is a basis of the three-dimensional submodule $W_2$ of $V_1\otimes V_1$. Since $\underline{\mathbf{B}}_{1,1}$ is an $A$-basis of ${}_A^\omega V_1\bigotimes_A {}_AV_1$, it follows that the epimorphism \smash{$\tau_{1,1}^t \colon {}_A\overset{\raisebox{-1.5pt}{\scriptsize$\smallbullet$}}{C}(1) \bigotimes_A {}_A\overset{\raisebox{-1.5pt}{\scriptsize$\smallbullet$}}{C}(1)\ra {}_A\overset{\raisebox{-1.5pt}{\scriptsize$\smallbullet$}}{C}(2)$} splits (see Proposition \ref{teoLuzstigOAscinde}). The vector $v_- \,\lozenge\, v_-$ is cyclic, so $\underline{\mathbf{B}}_{1,1}[0]=\{v_- \,\lozenge\, v_-\}$. By the definitions, we have $v_+ \,\lozenge\, v_+ = (1 \,\lozenge_0 \, 1)_{1,1}$, $v_+ \,\lozenge\, v_- = (1 \,\lozenge_0\, F)_{1,1}$, $v_- \,\lozenge\, v_+ = (F \,\lozenge_0\, 1)_{1,1}$, $v_- \,\lozenge\, v_- = (F \,\lozenge_0\, F)_{1,1}$, so the corresponding elements of \smash{$\dot{\mathbf{B}}_{1,1} \subset \dot{\mathbf{B}}$} are respectively $1_{0}$, $1_{-2}F$, $1_{2}E$, and $F1_{2}E = E1_{-2}F$.

The invariant submodule $W_0$ of ${}^\omega V_1\otimes V_1$ is generated by $v' = v_-\otimes v_- - q^{-1}v_+\otimes v_+$. The~$U_A^{\rm res}$-modules ${}_A^\omega V_1\bigotimes_A {}_AV_1$ and $W_2 \oplus W_0$ are not equal, though they are by extending scalars to $\mc(q)$. Indeed, we have
\[
v_+\otimes v_+= [2]_q^{-1}(qv_+ \,\lozenge\, v_+-v') \notin W_2 \oplus W_0.
\]
 The module of coinvariants is \smash{$({}^\omega V_1\otimes V_1)_{U_q^{\rm ad}} = \mc(q)\{\pi(v_- \otimes v_-)\}$},
where as usual $\pi\colon {}^\omega V_1\otimes V_1 \ra \smash{({}^\omega V_1\otimes V_1)_{U_q^{\rm ad}}}$ is the quotient map. The transpose map \smash{$\pi^t\colon (({}^\omega V_1\otimes V_1)_{U_q^{\rm ad}})^* \ra ({}^\omega V_1\otimes V_1)^*$} sends $(v_- \,\lozenge\, v_-)^*$ to the unique $U_q^{\rm ad}$-invariant linear map
\[
{\rm ev}_1\colon\ {}^\omega V_1\otimes V_1 \ra \mc(q)
\]
 such that ${\rm ev}_1(v_- \otimes v_-)=1$.

Note that, since elements of $\dot{\mathbf{U}}_A [\lambda > 2]$ act trivially on modules with all isotypical components of highest weight $\leq 2$, ${}_A^\omega V_1\bigotimes_A {}_AV_1$ is an indecomposable module over $\dot{\mathbf{U}}_A / \dot{\mathbf{U}}_A [\lambda > 2]$ (that is, \smash{$\dot{\mathbf{U}}_A / \dot{\mathbf{U}}_A [P_2]$} in the notations of Theorem~\ref{teoLuzstigOAinv}).}\end{exa}

\subsubsection[Some consequences on L\_{0,n}\^A and M\_{0,n}\^A]{Some consequences on $\boldsymbol{\Ll_{0,n}^A}$ and $\boldsymbol{\Mm_{0,n}^A}$}

Recall from Section \ref{defintform} the definition of the integral forms $\Ll_{0,n}^A$ and $\Mm_{0,n}^A$.
\begin{prop}\label{L0NAfgfree} $\Ll_{0,n}^A$ and $\Mm_{0,n}^A$ are free $A$-modules, and \smash{$\Mm_{0,n}^A$} is a direct summand of the $A$-module \smash{$\Ll_{0,n}^A$}. Moreover, \smash{$\Ll_{0,n}^A$} is a finitely generated ring.
\end{prop}
\begin{proof} Since \smash{$\Ll_{0,n}^A = \Oo_A^{\otimes n}$} as $U_A^{\rm res}$-modules, by \eqref{decompositionOA2} it has the basis \smash{$\bigcup_{[\lambda]\in P_+^n}\dot{\mathbf{B}} [[\lambda]]^*$}. Therefore, $\Ll_{0,n}^A$ is a free $A$-module. Since $A$ is a principal ideal domain, it follows that \smash{$\Mm_{0,n}^A$} is a~free $A$-submodule \cite[Appendix~2.2]{Lang}. By Theorem \ref{teoLuzstigOAinv}, there is a direct sum decomposition as $A$-module
\begin{equation}\label{decomposition OAn}
\Ll_{0,n}^A = \Mm_{0,n}^A \oplus {}_AN,
\end{equation} and the proof identifies a basis of $\Mm_{0,n}^A$ as a subset of $\bigcup_{[\lambda]\in P_+^n}\dot{\mathbf{B}} [[\lambda]]^*$.

Next, consider the question of finite generation. By the formula \eqref{prodL0nform}, it is enough to verify this for ${\mathcal L}_{0,1}^A$, but ${\mathcal L}_{0,1}^A= \Oo_A$ as an $A$-module, and $\Oo_A$ is finitely generated by the matrix coefficients of the fundamental $U_A^{\rm res}$-modules ${}_AV_{\varpi_k}$, $k\in\{1,\dots,m\}$ (see~\eqref{formgenOA1} and \eqref{formgenOA2}). Any monomials in these generators can be written as a $A$-linear combination of monomials in the same generators but with the product of $\Ll_{0,1}^A$, instead of the product $\star$. This follows from the integrality properties of the $R$-matrix, and the formula inverse to \eqref{mmtilde} (see in \cite[Section 3.3 and the formulas (4.6)--(4.8)]{BR}).
\end{proof}

\begin{Remark}\label{noHaarRoverA}\quad
\begin{itemize}\itemsep=0pt
\item[$(a)$] As noted in \eqref{decompAcoinv}, the $A$-module ${}_AN$ in the decomposition \eqref{decomposition OAn} is in general not a $U_A^{\rm res}$-module. Therefore, the $A$-linear projection map $\Rr_A\colon \Ll_{0,n}^A\ra \Mm_{0,n}^A$ such that ${\rm Ker}(\Rr_A)={}_AN$ is not a Reynolds operator, for it does not satisfy the identity $\Rr_A(\alpha \beta) = \alpha \Rr_A(\beta)$ for all $\alpha\in \Mm_{0,n}^A$, $\beta \in \Ll_{0,n}^A$.

\item[$(b)$] Recall \eqref{formuleRh}. In coherence with (a) above, there is no normalized Haar measure on $\Oo_A$ taking values in $A$. Indeed, by extending scalars over $\mc(q)$ it should otherwise coincide with the Haar measure $h\colon \Oo_q \ra \mc(q)$, but in the notations of Example \ref{V2V20} (see also the comments after \eqref{formulesl2}), since $h(v_0^* \otimes v_0)=0$ we have $h(bc)=-1/\big(q+q^{-1}\big)$, whence $h$ cannot be defined on $\Oo_A$.

\item[$(c)$] The Haar measure yields a well-defined $\mathcal{A}_0$-linear map $h\colon \Ll(\Oo_q) \ra \mathcal{A}_0$ (and analogously $\mathcal{A}_0$-linear and $\mathcal{A}_\infty$-linear maps $h\colon \Ll_{\lozenge}\big(\Oo_q^{\otimes n}\big) \ra \mathcal{A}_0$ and $\bar{h}\colon \bar{\Ll}_{\lozenge}\big(\Oo_q^{\otimes n}\big) \ra \mathcal{A}_\infty$ for any~$n\geq 1$, where \smash{$\big(\Ll_{\lozenge}\big(\Oo_q^{\otimes n}\big),\mathcal{B}[[\lambda]]^*\big)$} is the crystal basis at $q=0$ underlying the based $U_q^{\rm ad}$-module~\eqref{decompositionOAcoad}). Indeed, by \eqref{renormweightA} the lattice $\Ll_{\lambda}^r {}^{\rm up} \bigotimes_{\mathcal{A}_0} \Ll_{\lambda}^{\rm up}$ is generated by the matrix coefficients in the canonical bases of $V_\lambda^r$ and $V_\lambda$. Since the normalisation by powers of $q$ is vacuous on the trivial module $V_0^* \otimes V_0$, and $h$ vanishes on $V_\lambda^* \otimes V_\lambda$ for $\lambda \in P_+\setminus \{0\}$, the claim follows.
 \end{itemize}
\end{Remark}

\subsection{Perfect pairings}\label{intdual} We will need to restrict the morphisms $\Phi^+$, $\Phi^-$ in \eqref{phipm} on the integral forms $\Oo_A(B_+)$, $\Oo_A(B_-)$. We collect their properties in Theorem \ref{Phimap} and the discussion thereafter. In order to state it, we recall first a few facts about $R$-matrices and related pairings.

Recall that $\mathcal{C}_A$ is the category of $U_A^{\rm res}$-modules of type $1$. In \cite{Lusztig2,Lusztig}, Lusztig proved that $\mathcal{C}_A\bigotimes_A \mc[q^{\pm 1/D}]$ is braided and ribbon, with braiding given by the collection of endomorphisms
\[R = (R_{V,W})_{V,W\in {\rm Ob}(\mathcal{C}_A)}.\]
Actually, $R_{V,W}$ is represented by a matrix with coefficients in $q^{\mz/D}\mc\big[q^{\pm 1}\big]$ on the tensor product of the lower canonical bases of $V$ and $W$ (see \cite[Corollary~24.1.5]{Lusztig}).

This can be rephrased as follows in Hopf algebra terms. Denote by $\mathbb{U}_{\Gamma}$ the categorical completion of $\Gamma$, i.e., the Hopf algebra of natural transformations $F_{\mathcal{C}_A} \ra F_{\mathcal{C}_A}$, where ${F_{\mathcal{C}_A}\colon \mathcal{C}_A\ra A}$-$Mod_f$ is the forgetful functor towards the category $A$-$Mod_f$ of finite rank $A$-modules. Then $\mathbb{U}_{\Gamma}\bigotimes_A \mc[q^{\pm 1/D}]$ is quasi-triangular and ribbon with $R$-matrix
\[R\in \mathbb{U}_{\Gamma}^{\hat{\otimes} 2}\bigotimes_A \mc\bigl[q^{\pm 1/D}\bigr].\]
As in \eqref{Rcat}, we can write
 \[R^\pm=\sum_{(R)}R^\pm_{(1)}\otimes R^\pm_{(2)}.\]
There are pairings of Hopf algebras naturally related to the $R$-matrix $R$, considered as an element of $\mUq^{\hat{\otimes} 2}$. What follows is standard (see, e.g., \cite{KhorTol,KirRes,LS}), for details we refer to \cite[Proposition~3.73, Lemma~3.75, Theorem~3.92, Propositions~3.106 and 3.107]{VY}:
\begin{itemize}\itemsep=0pt
\item There is a unique pairing of Hopf algebras $\rho \colon U_q(\mathfrak{b}_-)^{\rm cop}\otimes U_q(\mathfrak{b}_+)\ra \mc\big(q^{1/D}\big)$ such that, for every $\alpha,\lambda\in P$ and $l,k\in U_q(\mathfrak{h})$,
\begin{gather}
\rho(K_\lambda,K_\alpha) = q^{(\lambda, \alpha)} ,\qquad \rho(F_i,E_j) = \delta_{i,j}\big(q_i-q_i^{-1}\big)^{-1},\nonumber \\
\rho(l,E_j) = \rho(F_i,k)=0.\label{defrho}
\end{gather}
\item The {\it Drinfeld pairing} $\tau \colon U_q(\mathfrak{b}_+)^{\rm cop}\otimes U_q(\mathfrak{b}_-)\ra \mc\big(q^{1/D}\big)$ is the bilinear map defined by $\tau(X,Y) = \rho(S(Y),X)$; it satisfies
\begin{gather}
\tau(K_\lambda,K_\alpha) = q^{-(\lambda, \alpha)} ,\qquad\tau(E_j,F_i) = - \delta_{i,j}\big(q_i-q_i^{-1}\big)^{-1} ,\nonumber\\
 \tau(l,F_i) = \tau(E_j,k)=0.\label{deftau}
\end{gather}
\item $\rho$ and $\tau$ are perfect pairings; this means that they yield {\it isomorphisms} of Hopf algebras $i_\pm \colon U_q(\mathfrak{b}_\pm) \ra \Oo_q(B_\mp)_{\rm op}$ (with coefficients {\it a priori} extended to $\mc\big(q^{1/D}\big)$, but see below) defined by, for every $X\in U_q(\mathfrak{b}_+)$, $Y\in U_q(\mathfrak{b}_-)$,
 \[\langle i_+(X), Y\rangle = \tau(S(X),Y),\qquad\langle i_-(Y) , X \rangle = \tau(X,Y).\]
Since $\Oo_q(B_\mp)_{\rm op}$ is equipped with the {\it inverse} of the antipode of $\Oo_q(B_\mp)$, which is induced by the antipode $S_{\Oq}$ of $\Oo_q$, it follows that $i_\pm \circ S = S_{\Oq}^{-1} \circ i_\pm$.
\item Denote by $p_\pm \colon \Oo_q(G) \ra \Oo_q(B_\pm)$ the canonical projection map, i.e., the Hopf algebra homomorphism dual to the inclusion map $U_q(\mathfrak{b}_\pm) \hookrightarrow U_q(\mathfrak{g})$. For every $\alpha,\beta \in \Oo_q(G)$, we have
\begin{equation}\label{Rpair}
\langle\alpha \otimes \beta , R \rangle = \tau \big(i_+^{-1}(p_-(\beta)), i_-^{-1}(p_+(\alpha))\big).
\end{equation}
\end{itemize}
Note that it is the use of weights $\alpha,\lambda\in P$ that forces the pairings $\rho$, $\tau$ to be defined over~$\mc\big(q^{1/D}\big)$, instead of $\mc(q)$. Then, let us consider the restrictions $\pi_q^+$ of $\rho$, and $\pi_q^-$ of $\tau$ defined by the formulas \eqref{defrho} and \eqref{deftau}, where now $\alpha\in Q$ and \smash{$k\in U_q^{\rm ad}(\mathfrak{h})$}. They take values in $\mc(q)$, and define pairings
\[\pi^{+}_q\colon\ U_q(\mathfrak{b}_-)^{\rm cop}\otimes U_q^{\rm ad}(\mathfrak{b}_+)\ra \mc(q),\qquad \pi^-_q\colon\ U_q(\mathfrak{b}_+)^{\rm cop}\otimes U_q^{\rm ad}(\mathfrak{b}_-)\ra \mc(q).\]
By the same arguments as for $\rho$ and $\tau$ (e.g., in \cite[Proposition 3.92]{VY}), it follows that $\pi^{\pm}_q$ are perfect pairings. Note also that $\pi^-_q =\kappa \circ \pi^+_q\circ (\kappa\otimes \kappa)$, where $\kappa\colon U_q \ra U_q$ is the $\mc$-linear automorphism extending \smash{${}^-\colon U_q^{\rm ad} \ra U_q^{\rm ad}$} in Section \ref{canbasemodqg}, so defined by
\begin{equation}\label{kappadef}
\kappa(E_i) = F_i,\qquad\kappa(F_i) = E_i,\qquad \kappa(K_\lambda) = K_{-\lambda},\qquad \kappa(q) = q^{-1}.
\end{equation}
In \cite{DC-L}, De Concini--Lyubashenko described integral forms of $\pi_q^\pm$ as follows. Denote by $m^*\colon \Oo_A\ra \Oo_A(B_+)\otimes \Oo_A(B_-)$ the map dual to the multiplication map $\Gamma(\mathfrak{b}_+)\otimes \Gamma(\mathfrak{b}_-)\ra \Gamma$, so $m^* = (p_+ \otimes p_-)\circ \Delta_{\Oo_A}$. Let $U_A(G^*)$ be the smallest $A$-subalgebra of $U_A(\mathfrak{b}_-)^{\rm cop}\otimes U_A(\mathfrak{b}_+)^{\rm cop}$ which contains the elements
\[1\otimes K_i^{-1}\bar{E}_i,\qquad \bar{F}_iK_i\otimes 1,\qquad L_i^{\pm 1}\otimes L_i^{\mp 1}, \qquad i=1,\dots,m,\]
and is stable under the diagonal action of $\mathcal{B}(\mathfrak{g})$. The reason for the notation $U_A(G^*)$ will be explained at the beginning of Section \ref{UOUlf}. Note that $U_A(G^*)$ is free over $A$, a Hopf subalgebra, and that a basis is given by the elements
\begin{equation}\label{basiselmt}
\bar{F}_{\beta_1}^{n_1}\cdots \bar{F}_{\beta_N}^{n_N}K_{n_1\beta_1+\dots+ n_N \beta_N}K_{\lambda}
\otimes K_{-\lambda}K_{-p_1\beta_1-\dots- p_N \beta_N}\bar{E}_{\beta_1}^{p_1}\cdots \bar{E}_{\beta_N}^{p_N},
\end{equation}
where $\lambda\in P$ and $n_1,\dots ,n_N, p_1,\dots ,p_N \in {\mathbb N}$.

Now, let $v$ be a lowest weight vector of the lowest weight $\Gamma$-module ${}_AV_{-\lambda}$, $\lambda \in P_+$. As after Theorem \ref{JLteo1}, denote by $v^*\in {}_AV_{-\lambda}^*$ the dual vector, and by $\psi_{-\lambda}\in \Oo_A$ the matrix coefficient defined by $\langle \psi_{-\lambda},x\rangle = v^*(xv)$ for every $x\in \Gamma$. Consider the maps $j_q^\pm\colon \Oo_q(B_\pm) \ra U_q(\mathfrak{b}_\mp)^{\rm cop}$ defined by
\[
\langle \alpha_+,X\rangle = \pi_q^+\big(j_q^+(\alpha_+),X\big),\qquad \langle \alpha_-,Y\rangle = \pi_q^-(j_q^-(\alpha_-),Y),
\]
where $\alpha_\pm\in \Oo_q(B_\pm)$, $X\in U_q^{\rm ad}(\mathfrak{b}_+)$, and $Y\in U_q^{\rm ad}(\mathfrak{b}_-)$.

The following theorem summarizes results proved in \cite[Sections~3 and~4]{DC-L}. Denote by~$\Oo_{\!A}\!\big[\psi_{-\rho}^{-1}\big]$ the localization of $\Oo_A$ by the element $\psi_{-\rho}$; this localization is well defined, for the set \smash{$\{\psi_{-\rho}^n\}_{n\in \mn}$} is a left and right multiplicative Ore subset of $\Oo_A$ (see Corollary \ref{Phi1int2} below for an analogous statement for $\Ll_{0,1}^A$). For the sake of clarity, let us spell out the correspondence of notations between statements: $\pi^+_q$, $\pi^-_q$, $U_q(\mathfrak{b}_\mp)^{\rm cop}$, $U_A(\mathfrak{b}_\mp)^{\rm cop}$, $\Oo_A(B_\pm)$, $U_A(G^*)$ and $\Phi$ are denoted in \cite{DC-L} respectively by $\pi''$, $\bar{\pi}''$, $U_q(\mathfrak{b}_\mp)_{\rm op}$, $R_q[B_\pm]''$, $R_q[B_\pm]$, $A''$ and $\mu''$ (the definition of~$j_A^\pm$ is implicit in \cite[Section 4.2]{DC-L}).

\begin{teo}\label{Phimap}\quad
\begin{itemize}\itemsep=0pt
\item[$(1)$] $\pi^\pm_q$ restricts to a perfect Hopf pairing between the unrestricted and restricted integral forms, \smash{$\pi^{\pm}_A\colon U_A(\mathfrak{b}_\mp)^{\rm cop}\otimes \Gamma(\mathfrak{b}_\pm)\ra A$}.
\item[$(2)$] \smash{$j_q^\pm$} yields an isomorphism of Hopf algebras \smash{$j_A^\pm\colon \Oo_A(B_\pm) \!\ra\! U_A(\mathfrak{b}_\mp)^{\rm cop}$}, satisfying $\langle \alpha_\pm,x_\pm\rangle\! = \pi^{\pm}_A\big(j_A^\pm(\alpha_\pm),x_\pm\big)$ for every $\alpha_\pm \in \Oo_A(B_\pm)$, $x_\pm\in \Gamma(\mathfrak{b}_\pm)$.
\item[$(3)$] The map $\Phi := \big(j_A^+\otimes j_A^-\big)\circ m^*\colon \Oo_A \ra U_A(G^*) \subset U_A(\mathfrak{b}_-)^{\rm cop}\otimes U_A(\mathfrak{b}_+)^{\rm cop}$ is an embedding of Hopf algebras, and it extends to an isomorphism $\Phi \colon \Oo_A\big[\psi_{-\rho}^{-1}\big] \ra U_A(G^*)$.
 \end{itemize}
\end{teo}

For our purposes, it is necessary to reformulate this result. Consider the morphisms of Hopf algebras $\Phi^\pm\colon \Oo_A(B_\pm) \ra U_A(\mathfrak{b}_\mp)^{\rm cop}$, $\alpha\mapsto (\alpha \otimes {\rm id})(R^\pm)$.
\begin{lem}\label{Phij} We have $\Phi^\pm = j_A^\pm$.
\end{lem}
\begin{proof} By definitions, for every \smash{$X\in U_q(\mathfrak{b}_+)^{\rm cop}$}, \smash{$Y\in U_q^{\rm ad}(\mathfrak{b}_-)$}, we have $\langle i_+(S^{-1}(X)), Y\rangle = \pi_q^-(X,Y)$, and similarly for every $X\in U_q^{\rm ad}(\mathfrak{b}_+)$, $Y\in U_q(\mathfrak{b}_-)^{\rm cop}$, we have $\langle i_-(S^{-1}(Y)), X\rangle = \pi_q^+(Y,X)$. By keeping these notations for $X$ and $Y$, we deduce \smash{$j_q^-(i_+(S^{-1}(X))) = X$} and \smash{$j_q^+(i_-(S^{-1}(Y))) = Y$}, i.e., \smash{$j_q^\pm = S\circ i_\mp^{-1}$}.
Because \smash{$S_{\Oq}^{-1} \circ i_\pm = i_\pm \circ S$}, it follows that
\begin{equation}\label{ijpmbis}
j_q^\pm\circ S_{\Oq} = S^{-1}\circ j_q^\pm .
\end{equation}
Also, for every $\alpha_-\in \Oo_q(B_-)$, we have
\begin{gather*}
\begin{split}
\big\langle \alpha_-, \Phi^+(i_-(Y))\big\rangle &= \langle i_-(Y) \otimes \alpha_-,R\rangle = \tau\big( i_+^{-1}(\alpha_-), Y\big)\\
& =
 \pi^-_q(j_q^-(S_{\Oq}(\alpha_-)),Y) = \langle \alpha_-,S(Y)\rangle,
\end{split}
\end{gather*}
where the first equality is by definition of $\Phi^+$ (see \eqref{phipm}), the second is \eqref{Rpair}, the third follows from \eqref{ijpmbis}, and the last from the definition of $j_q^-$. Similarly, for every $\alpha_+\in \Oo_q(B_+)$, we have
\begin{align*}
\langle \alpha_+, \Phi^-(i_+(X))\rangle & = \langle i_+(X) \otimes \alpha_+,R^-\rangle = \big\langle \alpha_+ \otimes S_{\Oq}^{-1}\circ i_+(X) , R\big\rangle = \langle \alpha_+ \otimes i_+(S(X)), R\rangle \\
& = \tau\big( S(X), i_-^{-1}(\alpha_+)\big) = \pi^+_q\big(S\big(i_-^{-1}(\alpha_+)\big), S(X)\big)\\
&= \pi^+_q\big(j_q^+(\alpha_+), S(X)\big) = \langle \alpha_+,S(X)\rangle.
\end{align*}
These computations imply $\Phi^\pm = S\circ i_\mp^{-1} = j_q^\pm$, and the result follows by taking integral forms.
\end{proof}

\subsection[Integral form and specialization of Phi\_n]{Integral form and specialization of $\boldsymbol{\Phi_n}$}\label{INTPHIn}
Recall the isomorphism of $U_q$-module algebras $\Phi_1\colon \Ll_{0,1} \ra U_q^{\rm lf}$, and that $U_A^{\rm lf} = U_A\cap U_q^{\rm lf}$. We have:

\begin{cor}\label{Phi1int} The map $\Phi_1$ affords an embedding of $U_A^{\rm res}$-module algebras $\Phi_1\colon \Ll_{0,1}^A \ra U_A^{\rm lf}$.
\end{cor}
\begin{proof} The only thing to be proved is that $\Phi_1(\Oo_A) \subset U_A^{\rm lf}$, since $\Ll_{0,1}^A = \Oo_A$ as $A$-module. But Lemma \ref{Phij} and \eqref{RSDmapdef} imply $\Phi_1 = m\circ \big({\rm id}\otimes S^{-1}\big)\circ \Phi$, and $\Phi$ maps $\Oo_A$ into $U_A(\mathfrak{b}_-)^{\rm cop}\otimes U_A(\mathfrak{b}_+)^{\rm cop}$ by Theorem \ref{Phimap}. The conclusion follows.
 \end{proof}

Let us denote
\[d=\psi_{-\rho}\in \Ll_{0,1}^A.\]
(The linear forms $\psi_{-\lambda}$ have been introduced before Theorem \ref{Phimap}.) When $\mathfrak{g}={\mathfrak{sl}_2}$ the element $d$ is one of the ``standard'' generators of $\Ll_{0,1}({\mathfrak{sl}_2})$ (see \eqref{rel01} below). In this case we have shown in \cite[Lemma 5.7]{BR} that $\Ll_{0,1}^A$ has a well-defined localization $\Ll_{0,1}^A\big[d^{-1}\big]$, and that \smash{$\Phi_1\colon \Ll_{0,1}^A\big[d^{-1}\big] \ra U_A^{\rm ad}=T_{2-}^{-1}U_A^{\rm lf}$} is an isomorphism of algebras. A generalization of these facts to any $\mathfrak{g}$ is provided by the following statement. As usual $\ell=K_{2\rho}$, the pivotal element.

\begin{cor}\label{Phi1int2}\quad
\begin{itemize}\itemsep=0pt
\item[$(1)$] The set $\{d^n\}_{n\in {\mathbb N}}$ is a left and right multiplicative Ore set in $\Ll_{0,1}^A$.
We can therefore define the localization $\Ll_{0,1}^A\big[d^{-1}\big]$.
\item[$(2)$] $\Phi_1$ extends to an embedding of $U_A^{\rm res}$-module algebras $\Phi_1\colon \Ll_{0,1}^A\big[d^{-1}\big] \ra U_A^{\rm lf}[\ell]$, and $U_A^{\rm lf}[\ell] = T_{2-}^{-1}U_A^{\rm lf}$.
\end{itemize}
\end{cor}
\begin{proof} (1) Because \smash{$\Ll_{0,1}^A$} has no nontrivial zero divisors, $d$ is a regular element. We have to show that for all \smash{$x\in \Ll_{0,1}^A$} there exists elements \smash{$y,y'\in\! \Ll_{0,1}^A$} and $d',d'' \in \!\{d^n\}_{n\in {\mathbb N}}$ such that~${xd'\!=\!dy}$ and $d''x=y'd$. In fact, $d'=d''=d$ in the present situation. Indeed by \eqref{Phi1psi-}, we have $\Phi_1(x)\Phi_1(d)=\Phi_1(x)K_{-2\rho}=K_{-2\rho}{\rm ad}^r(K_{2\rho})(\Phi_1(x))$, and ${\rm ad}^r(K_{2\rho})(\Phi_1(x))=\Phi_1({\rm coad}^r(K_{2\rho})(x))$. Therefore, the left Ore condition is satisfied with $y={\rm coad}^r(K_{2\rho})(x)$. Similarly, one finds $y'$.

(2) The first claim follows immediately from Corollary \ref{Phi1int} and $\Phi_1(d)=\ell^{-1}$, which is a~regular element of $U_A$. For the second claim, since \smash{$ K_{-2\rho}=\prod_{j=1}^m L_j^{-2}$}, localizing in $d$ we obtain
\[
 L_j^2=\prod_{k\not=j} L_k^{-2 }\Phi_1\big(d^{-1}\big)=\Phi_1\biggl(\prod_{k\not=j} \psi_{-\varpi_k}d^{-1}\biggr) \in \Phi_1\big( \Ll_{0,1}^A\big[d^{-1}\big] \big).
\]
Therefore, $T_{2-}^{-1}\!\subset \Phi_1\bigl( \Ll_{0,1}^A\big[d^{-1}\big] \bigr)$, which implies the assertion (2).
\end{proof}

We expect that the inclusion $\Phi_1(\Oo_A) \subset U_A^{\rm lf}$ is an equality, but have no proof yet. However, recall Joseph--Letzter's Theorem \ref{JLteo1}\,(1) and (2).
\begin{prop} \label{JLteo1overA} We have
\[
U_A = T_{2-}^{-1}U_A^{\rm lf}[T/T_{2}]=\Phi_1\big(\Ll_{0,1}^A\big[d^{-1}\big]\big)[T/T_{2}],\]
 and therefore \smash{$\Phi_1\colon \Ll_{0,1}^A\big[d^{-1}\big] \ra T_{2-}^{-1}U_A^{\rm lf}$} is an isomorphism. Moreover, \[\Phi_1(\Oo_A) = \bigoplus_{\lambda \in 2P_+} {\rm ad}^r(U_A^{\rm res})(K_{-\lambda}).\]
\end{prop}
\begin{proof} The inclusions $T\subset U_A$, $U_A^{\rm lf}\subset U_A$ and $\Phi_1\bigl(\Ll_{0,1}^A\big[d^{-1}\big]\bigr)\subset T_{2-}^{-1}U_A^{\rm lf}$ imply
\[\Phi_1\big(\Ll_{0,1}^A\big[d^{-1}\big]\big)[T/T_{2}]\subset T_{2-}^{-1}U_A^{\rm lf}[T/T_{2}]\subset U_A.\]
For the inverse inclusion, it is enough to show that any PBW basis vector of $U_A$ lies in $\Phi_1\big(\Ll_{0,1}^A\big[d^{-1}\big]\big)[T/T_{2}]$. This will follow at once if this is true of all root vectors $\bar{E}_{\beta_k}$, $\bar{F}_{\beta_k}$. Let us show this explicitly for the simple root vectors $\bar{E}_i$ and $\bar{F}_i$. For every positive root $\alpha$, define elements $\psi_{-\lambda}^{\alpha},\psi_{-\lambda}^{-\alpha}\in \Oo_A$ by the formulas
\[\langle \psi^{\alpha}_{-\lambda},x\rangle = v^*(xE_{\alpha} v),\qquad \langle \psi_{-\lambda}^{-\alpha},x\rangle = v^*(F_{\alpha} xv),\]
where $x\in \Gamma$. It is shown in \cite[Lemma~4.5]{DC-L} that
\begin{gather*}
\Phi(\psi_{-\lambda}) = K_{-\lambda} \otimes K_{\lambda},\qquad \Phi\big(\psi^{\alpha_i}_{-\varpi_j}\big) = - \delta_{i,j}q_i L_i^{-1} \otimes L_i K_i^{-1} \bar{E}_i, \\
 \Phi\big(\psi^{-\alpha_i}_{-\varpi_j}\big) =\delta_{i,j}q_i^{-1}\bar{F}_i K_i L_i^{-1}\otimes L_i .
\end{gather*}
\big(Note that the generators denoted by $E_i$ and $F_i$ in \cite{DC-L} are respectively $K_i^{-1}E_i$ and $F_iK_i$ in our notations, which explains the factors $q_i$, $q_i^{-1}$ in the formulas below; also $\kappa$ in \eqref{kappadef} maps $\bar{E}_i$, $\bar{F}_i$ to $-\bar{F}_i$, $-\bar{E}_i$, whence the sign for the expression of $\Phi\big(\psi^{\alpha_i}_{-\varpi_j}\big)$.\big) Since $\Phi_1 = m\circ \big({\rm id}\otimes S^{-1}\big)\circ \Phi$, we have
\begin{equation}
\label{formphi1K}
 \Phi_1(\psi_{-\lambda}) = K_{-2\lambda},\qquad \Phi_1\big(\psi^{\alpha_i}_{-\varpi_j}\big) =\delta_{i,j} L_i^{-2} \bar{E}_i ,\qquad \Phi_1\big(\psi^{-\alpha_i}_{-\varpi_j}\big) = \delta_{i,j}q_i^{-1}\bar{F}_i K_i L_i^{-2}.
\end{equation}
Therefore, \[\bar{E}_i, \bar{F}_i, L_i^{\pm 1} \in T_{2-}^{-1}\Phi_1\big(\Ll_{0,1}^A\big)[T/T_{2}] = \Phi_1\big(\Ll_{0,1}^A\big[d^{-1}\big]\big)[T/T_{2}].
\]
 These elements do not generate $U_A$; it is necessary to consider general root vectors. By the stability of $U_A(G^*)$ under $\mathcal{B}(\mathfrak{g})$ and the isomorphism \smash{$\Oo_A\big[\psi_{-\rho}^{-1}\big]\ra U_A(G^*)$} of Theorem \ref{Phimap}\,(3), for every positive root $\beta_k$, we have \smash{$1\otimes K_{\beta_k}^{-1}\bar{E}_{\beta_k}$} , \smash{$\bar{F}_{\beta_k}K_{\beta_k}\otimes 1 \in \Phi\big(\Oo_A\big[\psi_{-\rho}^{-1}\big]\big)=\Phi\big(\Ll_{0,1}^A\big[d^{-1}\big]\big)$}. Therefore, $\bar{F}_{\beta_k}K_{\beta_k}, S^{-1}\big(\bar{E}_{\beta_k}\big)K_{\beta_k} \in \Phi_1\big(\Ll_{0,1}^A\big[d^{-1}\big]\big)$, and $\bar{F}_{\beta_k}, S^{-1}(\bar{E}_{\beta_k})\in \Phi_1\big(\Ll_{0,1}^A\big[d^{-1}\big]\big)[T/T_{2}]$. The sets \smash{$S^{-1}\big(\bar{E}_{\beta_k}\big)U_A(\mathfrak{h})$} generate the subalgebra $U_A(\mathfrak{b}_+)$ of $U_A$ (in fact, let us quote that a formula of \smash{$S^{-1}(\bar{E}_{\beta_k})$} is given in~\cite{Xi}). From the triangular decomposition $U_A=U_A(\mathfrak{n}_-)U_A(\mathfrak{h})U_A(\mathfrak{n}_+)$, the inclusion $U_A\subset \Phi_1\big(\Ll_{0,1}^A\big[d^{-1}\big]\big)[T/T_{2}]$ follows, whence the equality too. In particular, $U_A$ is a free \smash{$\Phi_1\big(\Ll_{0,1}^A\big[d^{-1}\big]\big)$}-module with a basis formed by representatives of the cosets in $T/T_2$. By the uniqueness of this free decomposition, we find \smash{$\Phi_1\big(\Ll_{0,1}^A\big[d^{-1}\big]\big) = T_{2-}^{-1}U_A^{\rm lf}$}. Therefore, $\Phi_1$ in Corollary \ref{Phi1int2}\,(2) is surjective.

For the third claim, recall the isomorphism $\Phi_{1}\colon C(-w_0(\mu)) \ra {\rm ad}^r(U_q)(K_{-2\mu})$ (see \eqref{restphi1}), and that $\psi_{-\mu}$ is the matrix coefficient dual to the vector ${}^\omega v_{-\mu}\otimes v_{-\mu} \in {\rm End}_A(V_{-w_0(\mu)})$. This vector is cyclic by \eqref{udef}, so by equivariance $\Phi_1\colon {}_AC(-w_0(\mu)) \ra {\rm ad}^r(U_A^{\rm res})(K_{-2\mu})$ is an isomorphism of $U_A^{\rm res}$-modules. The second claim follows from this and \eqref{decompositionOA2} for $n=1$. \end{proof}

Recall from \eqref{Alekseevmap} the isomorphisms of $U_q$-module algebras \smash{$\Phi_n\colon \Ll_{0,n} \ra \big(U_q^{\otimes n}\big)^{\rm lf}$} and of algebras \smash{$\Phi_n\colon\Mm_{0,n}\rightarrow (U_q^{{\otimes} n})^{U_q}$}, and from \eqref{notambig0} the notations for specializations. Corollary \ref{Phi1int} can be extended to $\Phi_n$ as follows:

\begin{cor}\label{Phinint} The map $\Phi_n$ affords embeddings of module algebras \smash{$\Phi_n\colon \Ll_{0,n}^A \ra \big(U_A^{\otimes n}\big)^{\rm lf}$} and \smash{$\Phi_n\colon \Ll_{0,n}^\ep \ra \big(U^{\otimes n}\big)^{\rm lf}_\ep$}, $q=\epsilon'\in \mc^\times$.\end{cor}

\begin{proof} For the first claim, the only thing to prove is the inclusion $\Phi_n\big(\Ll_{0,n}^A\big) \subset U_A^{\otimes n}$. It follows from Corollary \ref{Phi1int} and the expression of $\Phi_n$ in terms of $\Phi_1$ and $R$-matrices (in particular, the fact that they preserve integrality, see \cite[Lemma~6.10]{BR}). For the specialization at $q=\epsilon'\in \mc^\times$, we have to justify that $\Phi_n$ is injective. One uses the fact, to be developed in Theorem~\ref{DCLteo1} below, that $\Phi \colon \Oo_\e \ra U_\e(G^*)$ is an embedding. The algebra $U_\e(G^*)$ has the basis elements \eqref{basiselmt}, and the map $m \circ \big({\rm id}\otimes S^{-1}\big)$ sends this basis to a free family of $U_\e$. Therefore, $\Phi_1\colon \Ll_{0,1}^\e\ra U_\e$ is injective. Since $\Phi_n$ differs from \smash{$\Phi_1^{\otimes n}$} by a linear isomorphism (induced by the conjugation action of $R$-matrices on the components \smash{${}_A\overset{\raisebox{-1.5pt}{\scriptsize$\smallbullet$}}{C}([\lambda])$} of $\Ll_{0,n}^A$ in~\eqref{decompositionOA2}, see \cite[equation~(6.10)]{BR}), $\Phi_n\colon\Ll_{0,n}^\e\to U_\epsilon^{\otimes n}$ is an embedding as well. \end{proof}

\begin{Remark} \quad
\begin{itemize}\itemsep=0pt
\item[$(1)$] It is a natural problem to determine the image of $\Phi_n$. One may expect that it would be~\smash{$\big(T_{2-}^{-1}U_A^{\rm lf}\big)^{\otimes n}$}, because this is true for $n=1$, as well as for any $n$ in the ${\mathfrak{sl}_2}$ case, as shown in \cite{BR}. Unfortunately, this is not so. This comes from the fact, e.g., for $n=2$, that the matrix elements of $R_{02} R_{01}R_{01}' R_{02}^{-1}$ do not belong to \smash{$\big(T_{2-}^{-1}U_A^{\rm lf}\big)^{\otimes 2}$} as can be shown by an explicit computation in the $\mathfrak{sl}(3)$ case.
\item[$(2)$] In the case of $\mathfrak{g} = {\mathfrak{sl}_2}$, we defined in \cite{BR} an algebra ${}_{\rm loc}\Ll_{0,n}^A$ generalizing \smash{$\Ll_{0,1}^A\big[d^{-1}\big]$} above, containing \smash{$\Ll_{0,n}^A$} as a subalgebra, and such that $\Phi_n$ extends to \smash{${}_{\rm loc}\Ll_{0,n}^A$} and yields an isomorphism \smash{$\Phi_n\colon {}_{\rm loc}\Ll_{0,n}^A\ra U_A^{\rm ad}({\mathfrak{sl}_2})^{\otimes n}$}. The definition of ${}_{\rm loc}\Ll_{0,n}^A$ involves elements \smash{$\xi^{(i)}\in \Ll_{0,n}^A$} ($i=1,\dots ,n$) such that
\smash{$\Phi_n\big(\xi^{(i)} \big)=\big(K^{-1}\big)^{(i)}\cdots\big(K^{-1}\big)^{(n)}$}. It may be of interest to study a~similar extension of $\Phi_n$ for general $\mathfrak{g}$.
\end{itemize}
\end{Remark}

\subsection[Structure theorems for U\_e and O\_e]{Structure theorems for $\boldsymbol{U_\e}$ and $\boldsymbol{\Oo_\e}$} \label{UOUlf}
As usual, we denote by $\e$ a primitive $l$-th root of unity, where $l$ is odd, and coprime to $3$ if $\mathfrak{g}$ has $G_2$-components.

Recall the subgroups $T_G$, $U_\pm$ and $B_\pm$ of $G$. Let $G^0 = B_+B_-$ (the {\it big cell} of $G$), and define the subgroup
\[G^* = \big\{\big(u_+t,u_-t^{-1}\big), \,t\in T_G,\, u_\pm \in U_\pm\big\}\subset B_+^{\rm op}\times B_-^{\rm op},\]
where $B_\pm^{\rm op}$ is the group $B_\pm$ with opposite multiplication. The group $G^*$ can be naturally identified with the Poisson--Lie dual of $G$ with its standard structure.

Recall also that there is an injective homomorphism $\gamma_q^{-1}\circ h_q\colon \mathcal{Z}(U_q) \ra U_q(\mathfrak{h})$, defined by means of the quantum Harish-Chandra homomorphism (see, e.g., \cite[Section 9.1.C]{CP}, or \cite[Section 3.13]{VY}). The image of $\gamma_q^{-1}\circ h_q$ is the set \smash{$U_q(\mathfrak{h})^{\tilde{W}}$} of invariant elements under $\tilde{W}$, the subgroup of $W\ltimes P_2^*$ generated by the conjugates $\sigma W \sigma$ of $W$ by elements $\sigma\in P_2^*$. Here, $P_2^*$ is the group of homomorphisms $P\ra \mz/2\mz$, and the semidirect product $W\ltimes P_2^*$ acts on $U_q(\mathfrak{h})$ by the standard action of the Weyl group $W$, and by the action of $P_2^*$ given by $\sigma\cdot K_\lambda := \sigma(\lambda)K_\lambda$.

Consider the inverse map \smash{$h_q^{-1}\circ \gamma_q \colon U_q(\mathfrak{h})^{\tilde{W}} \ra \mathcal{Z}(U_q)$}. The elements of the domain and target, when expanded in the PBW basis, have coefficients in $\mc(q)$. It was shown in \cite[Section 21.1]{DC-K-P2} that if an element of \smash{$U_q(\mathfrak{h})^{\tilde{W}}$} has no coefficient with a pole at $q=\e$, then its image by $h_q^{-1}\circ \gamma_q$ has no coefficient with a pole at $q=\e$. We therefore have a well-defined injection
\[U_\e(\mathfrak{h})^{\tilde{W}}\ra \mathcal{Z}(U_\e).\]
We denote its image by $\mathcal{Z}_1(U_\e)$. For instance, when $U_\e = U_\e({\mathfrak{sl}_2})$, $\mathcal{Z}_1(U_\e)$ is the polynomial algebra generated by the Casimir element $\Omega=\big(\e-\e^{-1}\big)^2FE+\e K+\e^{-1}K^{-1}$.

Denote by $\mathcal{Z}_0(U_\e)\subset U_\e$ the smallest subalgebra containing the elements $E_i^l$, $F_i^l$, $K_\alpha^{l}$, for $i\in \{1,\dots,m\}$, $\alpha\in P$, and stable under $\mathcal{B}(\mathfrak{g})$; it is also the subalgebra generated by \smash{$E_{\beta_k}^l$},~\smash{$F_{\beta_k}^l$},~$L_i^{\pm l}$, for $k\in \{1,\dots, N\}$ and $i\in \{1,\dots,m\}$ \cite[Section~18]{DC-K-P2}. We will denote by $\mathcal{Z}_0(U_\e(\mathfrak{n}_-))$, $\mathcal{Z}_0(U_\e(\mathfrak{h}))$ and $\mathcal{Z}_0(U_\e(\mathfrak{n}_+))$ the subalgebras of $\mathcal{Z}_0(U_\e)$ generated by the elements \smash{$F_{\beta_k}^l$}, $K_{\lambda}^l$ ($\lambda\in P$), and~$E_{\beta_k}^l$, respectively. In \cite[Sections~1.8, 3.3 and 3.8]{DC-K} and \cite[Theorem~14.1 and Sections~20--21]{DC-K-P2}, the following results are proved:

\begin{teo} \label{DCKteo1}\quad
\begin{itemize}\itemsep=0pt
\item[$(1)$] $U_\e$ has no nontrivial zero divisors, $\mathcal{Z}_0(U_\e)$ is a central Hopf subalgebra of $U_\e$, and $U_\e$ is a free $\mathcal{Z}_0(U_\e)$-module of rank $l^{\dim \mathfrak{g}}$. Moreover, the classical fraction algebra $Q(U_\e) = Q(\mathcal{Z}(U_\e))\bigotimes_{\mathcal{Z}(U_\e)} U_\e$ is a central simple algebra of PI degree $l^{N}$, and $U_\e$ is a maximal order of $Q(U_\e)$.
\item[$(2)$] ${\rm Maxspec}(\mathcal{Z}_0(U_\e))$ is a group isomorphic to $G^*$ above, and the multiplication map yields an isomorphism $\mathcal{Z}_0(U_\e)\bigotimes_{\mathcal{Z}_0(U_\e)\cap \mathcal{Z}_1(U_\e)} \mathcal{Z}_1(U_\e) \ra \mathcal{Z}(U_\e)$.
\end{itemize}
\end{teo}

By this theorem, the dimension of $Q(U_\e)$ over its center $Q(\mathcal{Z}(U_\e))$ is $l^{2N}$, and its dimension over $Q(\mathcal{Z}_0(U_\e))$ is $l^{\dim \mathfrak{g}} = l^{m+2N}$. Therefore, the field $Q(\mathcal{Z}(U_\e))$ is an extension of $Q(\mathcal{Z}_0(U_\e))$ of degree~$l^m$.

Note that, because $\mathcal{Z}_0(U_\e)$ is an affine and commutative algebra, the maximal spectrum ${\rm Maxspec}(\mathcal{Z}_0(U_\e)$), viewed as the set of characters of $\mathcal{Z}_0(U_\e)$, acquires by duality a structure of affine algebraic group. Thus, the first claim of (2) in the theorem means precisely that this group can be identified with $G^*$. See, for instance, \cite[Section~7.2.1]{BR} for an explicit description in the~${\mathfrak{sl}_2}$ case.

In addition, ${\rm Maxspec}(\mathcal{Z}_0(U_\e)$) and $G^*$ have natural Poisson structures which correspond one to the other under the isomorphism of (2), and we have the following identifications (see~\cite[Section 21.2]{DC-K-P2}). The dual isomorphism $\mathcal{O}(G^*)\!\ra\! \mathcal{Z}_0(U_\e)$ identifies $\Oo(T_G)$ with ${\mathcal{Z}_0(U_\e)\! \cap \! U_\e(\mathfrak{h})\! =\! \mc[lP]}$, where as usual $U_\e(\mathfrak{h}) = U_A(\mathfrak{h}) \bigotimes_A \mc_\e$. Therefore, we can identify $\mc[P]$ with $\Oo\big(\tilde T_G\big)$, the coordinate ring of the $l^m$-fold covering space \smash{$\tilde T_G \ra T_G$}. The quantum Harish-Chandra isomorphism identifies $\mathcal{Z}_1(U_\e)$ with \smash{$\mc[2P]^W \cong \Oo\big(\tilde T_G/(2)\big)^W$}, where we denote by $(2)$ the subgroup of $2$-torsion elements in $\tilde T_G$. Consider the map
\[\sigma\colon\ B_+\times B_-\longrightarrow G^0,\qquad (b_+,b_-)\longmapsto b_+b_-^{-1}.\]
The restriction of $\sigma$ to $G^*$ is an unramified covering map of degree $2^m$.
Composing $\sigma\colon G^*\ra G^0$ with the quotient map under conjugation, $G^0 \hookrightarrow G \ra G/ /G$, we get dually an embedding of~${\Oo(G/ /G)=\mathcal{O}(G)^G}$ in $\mathcal{O}(G^*)$. Collecting these observations, we see that the isomorphism of Theorem \ref{DCKteo1}\,(2) affords identifications
\[\mathcal{Z}_0(U_\e)\cap \mathcal{Z}_1(U_\e) \cong \mathcal{O}(G)^G\]
as a subalgebra of $\mathcal{Z}_0(U_\e)\cong \Oo(G^*)$, and
\[\mathcal{Z}_0(U_\e)\cap \mathcal{Z}_1(U_\e) = \mc[2lP]^W \cong \Oo\bigl(\tilde T_G/(2l)\bigr)^W\cong \Oo(T_G/(2))^W\]
as a subalgebra of $\mathcal{Z}_1(U_\e)\cong \Oo\big(\tilde T_G/(2)\big)^W$.

We will use the following obvious though crucial fact. Note that $U_A^{\rm ad}$ is naturally a subalgebra of $U_A^{\rm res}$, and therefore acts on $U_\e^{\rm res}$-modules. Denote by \smash{$\mathcal{Z}_0\big(U_A^{\rm ad}\big)\subset U_A^{\rm ad}$} the subalgebra generated by the elements $\bar{E}_{\beta_k}^l$, $\bar {F}_{\beta_k}^l$, $K_i^{\pm l}$, for $k\in \{1,\dots, N\}$ and $i\in \{1,\dots,m\}$.
\begin{lem}\label{acttrivial} For every $U_A^{\rm res}$-module $V$ of type $1$, the action of $\mathcal{Z}_0\big(U_A^{\rm ad}\big)$ on the specialization $V_\e := V\bigotimes_A \mc_\e$ is trivial.
\end{lem}

\begin{proof} This comes from \smash{$E_i^l = [l]_{q_i}! E_i^{(l)}$}, \smash{$F_i^l = [l]_{q_i}! F_i^{(l)}$} and the fact that $K_i$ acts on $V$ by powers of $q_i$. Specializing to $q=\e$ ends the proof. \end{proof}

A result similar to Theorem \ref{DCKteo1} holds true for $\Oo_\e$. Namely, take the specializations at $q=\e$ in Theorem~\ref{Phimap}. Denote by $\mathcal{Z}_0(U_\e(G^*))$ the subalgebra of $U_\e(G^*)$ generated by the elements ($k\in \{1,\dots, N\}, i\in \{1,\dots, m\}$)
\[1\otimes K_{-l \beta_k}E_{\beta_k}^l ,\qquad F_{\beta_k}^l K_{l\beta_k}\otimes 1,\qquad L_i^{\pm l}\otimes L_i^{\mp l}.\]
It is a central Hopf subalgebra. Recall that the coordinate ring $\Oo(G)$ can be identified as a~Hopf algebra with $U(\mathfrak{g})^\circ$, where as usual $U(\mathfrak{g})^\circ$ denotes the restricted dual of the enveloping algebra~$U(\mathfrak{g})$ over $\mc$. In \cite[Section 6]{DC-L}, De Concini--Lyubashenko introduced an epimorphism of Hopf algebras $\eta\colon\Gamma_\e\ra U(\mathfrak{g})$ (essentially a version of Lusztig's ``Frobenius'' epimorphism in~\cite{Lusztig2}), defined by
\begin{gather}
 \eta\big(E_i^{(p)}\big) = \begin{cases} \dfrac{e_i^{p/l}}{(p/l)!} & {\rm if}\ l \ {\rm divides}\ p,\\ 0 & {\rm otherwise},\end{cases} \qquad \eta\big(F_i^{(p)}\big) = \begin{cases} \dfrac{f_i^{p/l}}{(p/l)!} & {\rm if}\ l \ {\rm divides}\ p,\\ 0 & {\rm otherwise},\end{cases} \nonumber\\
 \eta(K_i)=1 , \qquad \eta((K_i;p)_{q_i}) = \begin{cases} \dfrac{h_i(h_i-1)\cdots (h_i-(p/l)+1)}{(p/l)!} & {\rm if}\ l \ {\rm divides}\ p,\\ 0 & {\rm otherwise},\end{cases}\label{etadef}
\end{gather}
where $p\in \mn$, and $e_i$, $f_i$ and $h_i$, $i\in \{1,\dots,m\}$, denote the standard generators of $U(\mathfrak{g})$. The kernel of $\eta$ is generated by the elements $E_i$, $F_i$, $K_i-1$, and $(K_i;p)_{q_i}$ where $l$ does not divide $p$. Put
\begin{equation}\label{Z0Odef}
\mathcal{Z}_0(\Oo_\e) := \eta^*(\Oo(G)),
\end{equation}
where $\eta^*\colon U(\mathfrak{g})^\circ \ra \Gamma_\e^\circ$ is the monomorphism dual to $\eta$. Let us define special matrix coefficients, analogous to those introduced in Theorem \ref{Phimap}. Denote by $v_{\varpi_i}$ and $v_{w_0(\varpi_i)}$ a highest weight vector and a lowest weight vector of the $\Gamma$-module ${}_AV_{\varpi_i}$. Denote also by $v_{w_0(\varpi_i)}^*$ and $v_{\varpi_i}^*$ a~highest and lowest weight vector of the dual module $\Gamma$-module ${}_AV_{\varpi_i}^* \cong {}_AV_{-w_0(\varpi_i)}$. Define the matrix coefficients $b_{\varpi_i}, c_{\varpi_i}\in \Oo_A$ by
\[b_{\varpi_i}(x) = v_{\varpi_i}^*(xv_{w_0(\varpi_i)}) ,\qquad c_{\varpi_i}(x) = v_{w_0(\varpi_i)}^*(xv_{\varpi_i})\]
for all $x\in \Gamma$. We consider them as elements of $\Oo_\e$. Denote by $\mathcal{Z}_1(\Oo_\e)$ the subalgebra of $\Oo_\e$ generated by the elements $b_{\varpi_i}^kc_{\varpi_i}^{l-k}$ for
$1 \leq i\leq m$ and $0\leq k\leq l$.

\begin{teo} \label{DCLteo1} \quad
\begin{itemize}\itemsep=0pt
\item[$(1)$] $\mathcal{Z}_0(\Oo_\e)$ is a central Hopf subalgebra of $\Oo_\e\!\subset\! \Gamma_\e^\circ$, and $Q(\mathcal{Z}(\Oo_\e))$ is an extension of $Q(\mathcal{Z}_0(\Oo_\e))$ of degree $l^m$.
\item[$(2)$] $\psi_{-l\rho}\in \mathcal{Z}_0(\Oo_\e)$, and $\mathcal{Z}_0(\Oo_\e)$ is generated by matrix coefficients of irreducible $\Gamma$-modules of highest weight $l\lambda$, $\lambda\in P_+$. Moreover, the multiplication map yields an isomorphism
 \[
 \mathcal{Z}_0(\Oo_\e)\bigotimes_{\mathcal{Z}_0(\Oo_\e)\cap \mathcal{Z}_1(\Oo_\e)} \mathcal{Z}_1(\Oo_\e) \ra \mathcal{Z}(\Oo_\e),
 \]
 and the map $\Phi$ in Theorem {\rm\ref{Phimap}} affords an algebra embedding $\mathcal{Z}_0(\Oo_\e) \ra \mathcal{Z}_0(U_\e(G^*))$ and algebra isomorphisms $\mathcal{Z}_0(\Oo_\e)\big[\psi_{-l\rho}^{-1}\big] \ra \mathcal{Z}_0(U_\e(G^*))$ and $\Oo_\e\big[\psi_{-l\rho}^{-1}\big] \ra U_\e(G^*)$.
\item[$(3)$] $\Oo_\e$ has no nontrivial zero divisors, and it is a free $\mathcal{Z}_0(\Oo_\e)$-module of rank $l^{\dim \mathfrak{g}}$. Moreover, the classical fraction algebra \smash{$Q(\Oo_\e) = Q(\mathcal{Z}(\Oo_\e))\bigotimes_{\mathcal{Z}(\Oo_\e)} \Oo_\e$} is a central simple algebra of~PI degree $l^{N}$, and $\Oo_\e$ is a maximal order of $Q(\Oo_\e)$.
 \end{itemize}
\end{teo}
For the proof, see \cite{DC-L}: Proposition~6.4 for the first claim of (1) (where $\mathcal{Z}_0(\Oo_\e)$ and~$\mathcal{Z}_0(U_\e(G^*))$ are denoted $F_0$ and $A_0$ respectively), the appendix of Enriquez and \cite{Enr} for the second claim of~(1) and~(2), Propositions 6.4 and~6.5 for the other claims of~(2), Theorem~7.2 (where $\Oo_\e$ is shown to be projective over $\mathcal{Z}_0(\Oo_\e)$) and~\cite{BGS} (which provides the additional K-theoretic arguments to deduce that $\Oo_\e$ is free), or \cite[Remark~2.18\,(b)]{AnGa}, for the second claim of (3), and Corollary~7.3 and Theorem~7.4 for the third claim. The fact that $\Oo_\e$ has no nontrivial zero divisors follows from the embedding $\Oo_\e \ra U_\e(G^*)$ via $\Phi$.

As above for $U_\e$, it follows directly from~(3) that $Q(\mathcal{Z}(\Oo_\e))$ has degree $l^m$ over $Q(\mathcal{Z}_0(\Oo_\e))$. For a more complete description of $\mathcal{Z}(\Oo_\e)$ we refer to~\cite{Enr} and Enriquez' appendix in~\cite{DC-L}, as well as~\cite{BGOe}.

We do not know a basis of $\Oo_\e$ over $\mathcal{Z}_0(\Oo_\e)$ for general $G$, but see~\cite{DRZ} for the case of~${\rm SL}_2$. We will recall the known results in this case of ${\rm SL}_2$ before Lemma~\ref{L01sl2}.

Finally, there is a natural action of the braid group $\mathcal{B}(\mathfrak{g})$ on $\Oo_\e$, that we will use. Namely, let~${n_i\in N(T_G)}$ be a representative of the reflection $s_i\in W=N(T_G)/T_G$ associated to the simple root $\alpha_i$. In \cite{Soi,SV}, Soibelman--Vaksman introduced functionals $t_i\colon \mathcal{O}_q\ra \mc(q)$ which quantize the elements $n_i$. They correspond dually to generators of the quantum Weyl group of~$\mathfrak{g}$; in the appendix, we recall their main properties, in particular, they map $\Oo_A$ to $A$ (see also~\cite[Section~8.2]{CP}, and \cite{DC-L,KhorTol,KirRes,LS,SV}). Denote by $\lhd$ the natural right action of functionals on~$\Oo_A$, namely (using Sweedler's notation)
\[\alpha\lhd h = \sum_{(\alpha)}h(\alpha_{(1)}) \alpha_{(2)}\]
for every $\alpha\in \Oo_A$ and $h\in \Oo_A\ra A$. Let us identify $\mathcal{Z}_0(\Oo_\e)$ with $\Oo(G)$ by means of~\eqref{Z0Odef}. We have \cite[Proposition~7.1]{DC-L}:

\begin{prop} \label{funcDCL} The maps $\lhd t_i$ on $\Oo_\e$ preserve $\mathcal{Z}_0(\Oo_\e)$, and satisfy $(f\lhd t_i)(a) = f(n_ia)$ and $(f\star\alpha)\lhd t_i = (f\lhd t_i)(\alpha\lhd t_i)$ for every $f\in \mathcal{Z}_0(\Oo_\e)$, $a\in G$, $\alpha\in \mathcal{O}_\e$.
\end{prop}

We provide an alternative, non computational, proof of this result in Appendix~\ref{actO}.

\section{Noetherianity and finiteness}\label{HNsec}
In this section, we prove Theorem~\ref{HN}. Recall that by Noetherian we mean right and left Noetherian. We begin with

\begin{teo} \label{LonNoeth} The algebras ${\mathcal L}_{0,n}$, $\Ll_{0,n}^A$ and $\Ll_{0,n}^\ep$, $\ep\in \mc^\times$, are Noetherian.
\end{teo}

By Proposition \ref{L0NAfgfree}, each of the algebras in this theorem is finitely generated.

Theorem \ref{LonNoeth} for ${\mathcal L}_{0,1}$ and any $\mathfrak{g}$ follows immediately from Joseph--Letzter's Theorem \ref{JLteo1}, claim (3), by identifying ${\mathcal L}_{0,1}$ with \smash{$U_q^{\rm lf}$} via $\Phi_1$. The method of proof uses filtration arguments. An alternative proof in the case of $\mathfrak{sl}(n)$, which works also for ${\mathcal L}_{0,1}^A$, was obtained by Domokos--Lenagan in \cite{DomLen}, by exhibiting special sequences of generators of \smash{${\mathcal L}_{0,1}^A$} satisfying {\it polynormal} relations, as we define now.

\begin{defi}[{see \cite[Proposition 3.133]{VY}}] Let $R$ be a Noetherian Abelian ring, and~$B$ a~finitely generated $R$-algebra with product $\circ$. We call polynormal a set of relations between generators $u_1,\dots,u_M$ of~$B$, of the form{\samepage
\begin{equation}\label{echangegrad++}
u_i\circ u_j - q_{ij} u_j\circ u_i = \sum_{s=1}^{j-1}\sum_{t=1}^M \bigl(\alpha_{ij}^{st} u_s\circ u_t + \beta_{ij}^{st} u_t\circ u_s\bigr)
\end{equation}
for all $1\leq j<i\leq M$, where $\alpha_{ij}^{st}, \beta_{ij}^{st} \in R$, and the elements $q_{ij}\in R$ are invertible.}
\end{defi}
Note that this definition is more restrictive than the more standard one, e.g., in~\cite[Definition~II.4.1]{BG}. If such a set of relations exists in $B$, then $B$ can be endowed with an algebra filtration such that the associated graded algebra is a quotient of a skew-polynomial algebra~\cite[Proposition I.8.17]{BG}. By classical results, we have (see, e.g., \cite[Theorems~1.2.9, 1.6.9 and Examples~1.6.11]{MC-R}, or~\cite[Lemmas~3.130--3.131]{VY}):

\begin{teo} If the algebra filtration is well founded, then $B$ is a Noetherian ring.
\end{teo}

In \cite{DomLen}, Theorem \ref{LonNoeth} is also proved for any $n\geq 1$ in the case of $\mathfrak{g}={\mathfrak{sl}_2}$ by considering~$\Ll_{0,n}^A({\mathfrak{sl}_2})$ as an iterated overring of ${\mathcal L}_{0,1}({\mathfrak{sl}_2})$.

The proof of Theorem \ref{LonNoeth} that we develop for any $\mathfrak{g}$ and $n\geq 1$ is also based on polynormal relations. In our proof, the generating set of $\Ll_{0,n}$ that we will consider is evident, as they are matrix coefficients in the modules $V_{\varpi_k}$, $k\in\{1,\dots,m\}$; the task is then to exhibit a set of polynormal relations between them, that hold in a certain graded algebra associated to $\Ll_{0,n}$. Indeed, as explained above this will imply that the graded algebra is Noetherian, and that $\Ll_{0,n}$ is Noetherian as well. In the case of $\Ll_{0,n}^A$, the proof is formally similar, but it needs the use of canonical bases discussed in Section \ref{canbasemodqg}.

\begin{proof}[Proof of Theorem \ref{LonNoeth}] First, we develop the proof for $\Ll_{0,n}$, and then for $\Ll_{0,n}^A$; the result for%
\[\Ll_{0,n}^\ep\!= \Ll_{0,n}^A/(q-\ep)\Ll_{0,n}^A
\]
 follows immediately by lifting ideals by the quotient map \smash{$\Ll_{0,n}^A \!\ra \Ll_{0,n}^\ep$}.

We adapt the proof of Theorem \ref{JLteo1}\,(3) given in \cite[Theorem 3.137]{VY}. Let us begin by recalling these arguments. In doing this, let us stress that \cite{VY} takes on ${\mathcal O}_q$ and $\Ll_{0,1}$ the product opposite to ours, and below in \eqref{echangegrad} and \eqref{echangegrad+} we respect their convention.

As usual, let $C(\mu)$ be the vector space generated by the matrix coefficients of $V_\mu$, the simple~$U_q^{\rm ad}$-module of highest weight $\mu \in P_+$. Denote by $C(\mu)_\lambda \subset C(\mu)$ the subspace of weight $\lambda$ for the left coregular action of $U_q({\mathfrak h})$; so $\alpha\in C(\mu)_\lambda$ if \smash{$K_\nu\rhd \alpha = q^{(\nu,\lambda)}\alpha$}, $\nu\in P$. Consider the semigroup
\[\Lambda=\{(\mu, \lambda)\in P_+\times P, \, \lambda \;\text{is a weight of}\;V_\mu\}.\]
Recall that the partial order $\preceq$ on $P$ is defined by $\mu\preceq \mu'$ if and only if $\mu'-\mu\in D^{-1}Q_+$. Define~$\preceq$ on $\Lambda$ by: $(\mu,\lambda)\preceq (\mu', \lambda')$ if and only if $\mu'-\mu\in D^{-1}Q_+$ and $\lambda'-\lambda\in D^{-1}Q_+$. If~${(\mu,\lambda)\preceq (\mu', \lambda')}$ and $(\mu,\lambda)\ne (\mu', \lambda')$, we write $(\mu,\lambda)\prec (\mu', \lambda')$. Since ${\mathcal L}_{0,1}$ and ${\mathcal O}_q$ are isomorphic vector spaces, we have $ {\mathcal L}_{0,1}=\bigoplus_{\mu\in P_+} C(\mu)=\bigoplus_{(\mu,\lambda)\in \Lambda} C(\mu)_\lambda$. Consider the family of subspaces
\begin{align*}
& {\mathcal F}_2^{\mu,\lambda} :=\bigoplus_{(\mu',\lambda')\preceq (\mu,\lambda)} C(\mu')_{\lambda'} , \qquad{\mathcal F}_2^{\prec \mu,\lambda}:=\bigoplus_{(\mu',\lambda')\prec (\mu,\lambda)} C(\mu')_{\lambda'} , \qquad (\mu,\lambda)\in \Lambda.
\end{align*}
We have
\begin{equation}\label{filtrationMatthieu}
{\mathcal L}_{0,1} = \bigcup_{(\mu,\lambda)\in \Lambda} {\mathcal F}_2^{\mu,\lambda}.
\end{equation}
Indeed, clearly
\[ {\mathcal L}_{0,1} = \sum_{(\mu,\lambda)\in \Lambda} {\mathcal F}_2^{\mu,\lambda},
\]
 so \eqref{filtrationMatthieu} follows from the following fact: for every $(\mu,\lambda), (\mu',\lambda')\in \Lambda$, the element $(\mu'',\lambda''):= (\mu+\mu',\lambda+\lambda')$ is such that
\[{\mathcal F}_2^{\mu,\lambda}+ {\mathcal F}_2^{\mu',\lambda'} \subset {\mathcal F}_2^{\mu'',\lambda''}.\]
Note that in general, since $Q_+\nsubseteq P_+$ (but $P_+\subset D^{-1}Q_+$), it is not true that there exists an element $(\mu'',\lambda'')$ satisfying such an inclusion if one replaces $\preceq$ with the standard ``product'' partial order $\leq$ on $\Lambda$, defined by $(\mu,\lambda)\leq (\mu', \lambda')$ if and only if $\mu'-\mu\in Q_+$ and $\lambda'-\lambda\in Q_+$. Note also that $\preceq$ is finer than $\leq$, in the sense that if $\mu \leq \mu'$, then $\mu\preceq \mu'$. Again, this would not be true if we had replaced $D^{-1}Q_+$ by $P_+$ in the definition of $\preceq$.

The family \smash{${\mathcal F}_2:=\big\{{\mathcal F}_2^{\mu,\lambda}\big\}_{(\mu,\lambda)\in \Lambda}$} is a filtration of the vector space $ {\mathcal L}_{0,1}$, which is clearly well founded (i.e., every subset of $\Lambda$ contains a minimal element, or equivalently any decreasing infinite sequence of elements in $\Lambda$ is eventually constant).

Consider the associated graded vector space \smash{${\rm Gr}_{\mathcal F_2}({\mathcal L}_{0,1}) := \bigoplus_{(\mu,\lambda)} {\mathcal F}_2^{\mu,\lambda}/{\mathcal F}_2^{\prec \mu,\lambda}$}. By identifying an element $x\in C(\mu)_\lambda$ with its coset \smash{$\bar{x}\in {\mathcal F}_2^{\mu,\lambda}/{\mathcal F}_2^{\prec \mu,\lambda}$}, we get an equality of vector spaces $ {\rm Gr}_{{\mathcal F}_2}({\mathcal L}_{0,1})= \bigoplus_{(\mu,\lambda)\in \Lambda} C(\mu)_\lambda$. Now, one has the following facts:

(i) Taking the product in $\Ll_{0,1}$, we have
\begin{equation}\label{gradalg01}
\alpha \beta\in {\mathcal F}_2^{\mu_1+\mu_2,\lambda_1+\lambda_2}\qquad \mathrm{for}\quad \alpha \in C(\mu_1)_{\lambda_1},\quad \beta \in C(\mu_2)_{\lambda_2}.
\end{equation}
This follows from \eqref{decomprep} and the fact that, for every $\nu\in P_+$ and every summand of the formula~\eqref{mmtilde}, denoting by $-r\in -Q_+$ the weight of the $R$-matrix component $R_{(2)}$ we have
\begin{gather*}
K_\nu\rhd \bigl((R_{(2')}{S}(R_{(2)}) \rhd \alpha) \star (R_{(1')}\rhd \beta \lhd R_{(1)})\bigr) \\
\qquad= q^{(\nu,\lambda_1+\lambda_2-r)}(R_{(2')}{S}(R_{(2)}) \rhd \alpha) \star (R_{(1')}\rhd \beta \lhd R_{(1)}).
\end{gather*}
(Details of a similar computation are given below \eqref{prodL0nform-duplicate}.) It follows from~\eqref{gradalg01} that~$\mathcal{F}_2$ is an algebra filtration of $\Ll_{0,1}$, and then ${\rm Gr}_{\mathcal F_2}({\mathcal L}_{0,1})$ is a graded algebra.

(ii) Denote by $\alpha \circ \beta$ the product in $ {\rm Gr}_{{\mathcal F}_2}({\mathcal L}_{0,1})$ of $\alpha,\beta\in \Ll_{0,1}$. The space $C(\mu_1+\mu_2)$ has multiplicity one in $C(\mu_1)\otimes C(\mu_2)$ (again by~\eqref{decomprep}), therefore if $\alpha \in C(\mu_1)_{\lambda_1}$ and $\beta \in C(\mu_2)_{\lambda_2}$, then $\alpha \circ \beta$ is the projection of $\alpha \beta$ onto $C(\mu_1+\mu_2)_{\lambda_1+\lambda_2}$.
 Denote by $\bar{\star}$ the product $\star$ of $\Oo_q$ followed by the projection onto the component $C(\mu + \nu)$. Then, we have
\begin{equation}\label{filtration=}
C(\mu) \circ C(\nu) = C(\mu)\ \bar{\star} \ C(\nu) = C(\mu + \nu).
\end{equation}
This follows from the formula \eqref{mmtilde}, and the fact that it is given by an invertible twist of the product $\star$.

(iii) For every $\mu\in P_+$, fix a basis of weight vectors \smash{$e_1^\mu,\dots,e_{d(\mu)}^\mu$} of $V_\mu$. Denote by $e^1_\mu,\dots,e^{d(\mu)}_\mu\!\in V_\mu^*$ the dual basis, and by $w\big(e_i^\mu\big)$ the weight of $e_i^\mu$. Consider the matrix coefficients ${}_{\mu}\phi^i_j(x) :=\smash{ e^i_\mu\big(\pi_V(x)\big(e_j^\mu\big)\big)}$, $x\in U_q$. By using the formula \eqref{mmtilde} and the explicit form of the $R$-matrix, one can check that
\begin{align}
{}_{{\mu}}\phi^i_j \circ {}_{{\nu}}\phi^k_l & = \sideset{}{'}\sum_{j',l'} c_{j',l'}^{ikjl}\ {}_{{\mu}}\phi^i_{j'} \ \bar{\star} \ {}_{{\nu}}\phi^k_{l'}\nonumber \\ & = q^{(w(e_j^\mu),w(e_l^\nu)-w(e_k^\nu))} {}_{{\mu}}\phi^i_{j} \ \bar{\star} \ {}_{{\nu}}\phi^k_{l} + \sideset{}{'}\sum_{\substack{j',l'\\ j'\ne j,\, l'\ne l}} d_{j',l'}^{ikjl}\ {}_{{\mu}}\phi^i_{j'} \circ {}_{{\nu}}\phi^k_{l'},\label{prodcirc1}
\end{align}
where \smash{$\sum^{'}_{j',l'}$} is the sum over indices with weights satisfying
\[
w\big(e_j^\mu\big) + w\big(e_l^\nu\big) = w\big(e_{j'}^\mu\big) + w\big(e_{l'}^\nu\big), \qquad w\big(e_{j'}^\mu\big)\leq w\big(e_{j}^\mu\big)
\qquad \text{and} \qquad w\big(e_{l'}^\nu\big)\geq w\big(e_{l}^\nu\big),
\]
 and the coefficient \smash{$c_{j,l}^{ikjl}$}, equal to \smash{$q^{(w(e_j^\mu),w(e_l^\nu)-w(e_k^\nu))}$}, is computed from the term $\Theta$ in the $R$-matrix factorization \eqref{Rmatfact}. In general, all the coefficients~\smash{$c_{j',l'}^{ikjl}$} and \smash{$d_{j',l'}^{ikjl}$} belong to $\mc(q)$ (see \cite[Proposition 4.1]{BR}); in particular \smash{$q^{(w(e_j^\mu),w(e_l^\nu)-w(e_k^\nu))}\in q^\mz$} since $w(e_l^\nu)-w(e_k^\nu)\in Q$. The second equality follows by repeated use of the first and \eqref{filtration=}. Similarly, by using \eqref{mmtildebis} one gets
\begin{align*}
{}_{{\nu}}\phi^k_l \circ {}_{{\mu}}\phi^i_j & = \sideset{}{'}\sum_{i',k'} e_{i',k'}^{kilj}\ {}_{{\mu}}\phi^{i'}_{j} \ \bar{\star} \ {}_{{\nu}}\phi^{k'}_{l}\\
 & = q^{(w(e_i^\mu),w(e_k^\nu)-w(e_l^\nu))} {}_{{\mu}}\phi^i_{j} \ \bar{\star} \ {}_{{\nu}}\phi^k_{l} + \sideset{}{'}\sum_{\substack{i',k'\\ i'\ne i,\, k'\ne k}} e_{i',k'}^{kilj}\ {}_{{\mu}}\phi^{i'}_{j} \ \bar{\star}\ {}_{{\nu}}\phi^{k'}_{l}\\
 & = q^{(w(e_i^\mu),w(e_k^\nu)-w(e_l^\nu))} {}_{{\mu}}\phi^i_{j} \ \bar{\star} \ {}_{{\nu}}\phi^k_{l} + \sideset{}{'}\sum_{\substack{i',k',j',l'\\ i'\ne i,\, k'\ne k}} f_{i',k'}^{kilj}\ {}_{{\mu}}\phi^{i'}_{j'} \ \circ\ {}_{{\nu}}\phi^{k'}_{l'},
\end{align*}
where \smash{$e_{i',k'}^{kilj}$}, \smash{$f_{i',k'}^{kilj}\in \mc(q)$}, and \smash{$\sum^{'}_{i',k'}$} is the sum over indices with weights satisfying \begin{gather*}
w\big(e_i^\mu\big) + w(e_k^\nu) = w\big(e_{i'}^\mu\big) + w(e_{k'}^\nu), \qquad w\big(e_{i'}^\mu\big)\leq w\big(e_{i}^\mu\big),\\
w(e_{k'}^\nu)\geq w(e_{k}^\nu) , \qquad e_{i,k}^{kilj} = q^{(w(e_i^\mu),w(e_k^\nu)-w(e_l^\nu))}.
\end{gather*}
 The third equality comes from the second and \eqref{prodcirc1}; the sum is over indices with weights satisfying
 \begin{gather*}
 w\big(e_i^\mu\big) + w(e_k^\nu) = w\big(e_{i'}^\mu\big) + w(e_{k'}^\nu), \\
 w\big(e_{i'}^\mu\big)< w\big(e_{i}^\mu\big), \qquad w(e_{k'}^\nu)> w(e_{k}^\nu), \qquad w\big(e_{j'}^\mu\big)\leq w\big(e_{j}^\mu\big), \qquad w(e_{l'}^\nu)\geq w(e_{l}^\nu).
 \end{gather*} By eliminating the leading term \smash{${}_{{\mu}}\phi^i_{j} \ \bar{\star} \ {}_{{\nu}}\phi^k_{l}$}, one deduces
\begin{align}\label{commutphiklij}
&{}_{{\nu}}\phi^k_l \circ {}_{{\mu}}\phi^i_j - q_{ijkl} \ {}_{{\mu}}\phi^i_j \circ {}_{{\nu}}\phi^k_l = \sideset{}{'}\sum_{\substack{i',k',j',l'\\ i'\ne i,\, k'\ne k}} f_{i',k'}^{kilj}\ {}_{{\mu}}\phi^{i'}_{j'} \ \circ\ {}_{{\nu}}\phi^{k'}_{l'}-\sideset{}{'}\sum_{\substack{j',l'\\ j'\ne j,\, l'\ne l}} q_{ijkl} d_{j',l'}^{ikjl}\ {}_{{\mu}}\phi^i_{j'} \circ {}_{{\nu}}\phi^k_{l'},\!\!\!
\end{align}
where $q_{ijkl} = q^{(w(e_j^\mu)+ w(e_i^\mu),w(e_k^\nu)-w(e_l^\nu))}$.

(iv) We can always reorder the weight vectors $e_1^\mu,\dots,e_{d(\mu)}^\mu$ so that $w\big(e_i^\mu\big) > w\big(e_j^\mu\big)$ implies $i < j$; then \eqref{commutphiklij} reads
\begin{align}
{}_{{\nu}}\phi^k_l \circ {}_{{\mu}}\phi^i_j - q_{ijkl} \ {}_{{\mu}}\phi^i_j\circ {}_{{\nu}}\phi^k_l ={}& \sum_{r=i}^{d(\mu)}\sum_{s=1}^k \sum_{u=1}^{l-1} \sum_{v=j+1}^{d(\mu)} \delta^{ijkl}_{rsuv} \ {}_{\mu}\phi^r_v \circ {}_{{\nu}}\phi^s_u \nonumber\\
 & - \sum_{r=i+1}^{d(\mu)}\sum_{s=1}^{k-1} q_{ijkl}\gamma_{rs}^{ijkl} \ {}_{\mu}\phi^r_j \circ {}_{{\nu}}\phi^s_l, \label{echangegrad}
\end{align}
where \smash{$\gamma_{rs}^{ijkl}, \delta^{ijkl}_{rsuv}\in \mc(q)$} are such that $\gamma_{rs}^{ijkl}=0$ unless $w(e_r^\mu) < w(e_i^\mu)$ and $w(e_s^\nu) > w(e_k^\nu)$, and \smash{$\delta^{ijkl}_{rsuv}=0$} unless $w(e_u^\nu) > w(e_l^\nu)$, $w(e_v^\mu) < w(e_j^\mu)$, $w(e_r^\mu) \leq w(e_i^\mu)$ and $w(e_s^\nu) \geq w(e_k^\nu)$. Now, from \eqref{echangegrad} one can extract a defining set of polynormal relations for $ {\rm Gr}_{{\mathcal F}_2}({\mathcal L}_{0,1})$, as in \eqref{echangegrad++}. Indeed, like ${\mathcal L}_{0,1}$ the algebra $ {\rm Gr}_{{\mathcal F}_2}({\mathcal L}_{0,1})$ is generated by the matrix coefficients ${}_{{\varpi_k}} \phi_i^j$ of the fundamental representations $V_{\varpi_k}$. One can list these matrix coefficients, say $M$ in number, in an ordered sequence $u_1,\dots,u_M$ such that the following condition holds: if $w(e_k^{\varpi_s}) < w(e_i^{\varpi_r})$, or~${w(e_k^{\varpi_s}) = w(e_i^{\varpi_r})}$ and $w(e_l^{\varpi_s}) < w(e_j^{\varpi_r})$, then $u_a:={}_{{\varpi_r}} \phi^i_j$ and $u_b:={}_{{\varpi_s}} \phi^k_l$ satisfy $b<a$. Then denoting \smash{${}_{{\mu}}\phi^i_j$}, \smash{${}_{{\nu}}\phi^k_l$} in \eqref{echangegrad} by $u_j$, $u_i$, respectively, and assuming $u_j<u_i$, one finds that all terms~${u_s:= {}_{\mu}\phi^r_v}$, ${}_{\mu}\phi^r_j$ in the sums are $< u_j$. Therefore, for all $1\leq j<i\leq M$ it takes the form
\begin{equation}\label{echangegrad+}
u_i\circ u_j - q_{ij} u_j\circ u_i = \sum_{s=1}^{j-1}\sum_{t=1}^M \alpha_{ij}^{st} u_s\circ u_t
\end{equation}
for some $q_{ij}\in q^{\mz}$ and $\alpha_{ij}^{st}\in \mc(q)$. As explained after \eqref{echangegrad++}, it follows that $ {\rm Gr}_{{\mathcal F}_2}({\mathcal L}_{0,1})$ is a~Noetherian ring, and since the filtration $\mathcal{F}_2$ is well founded, it implies that $\Ll_{0,1}$ is Noetherian too.

We are going to extend all these facts to $\Ll_{0,n}$, $n>1$. First, we need to refine the filtration~${\mathcal F}_2$ on $\Ll_{0,1}$. Consider the action of $U_q({\mathfrak h})$ on $C(\mu)_{\lambda}$ given by
\begin{equation}\label{coadweightnew}
K_\nu.\alpha := {\rm coad}\big(K_\nu^{-1}\big)(\alpha),\qquad \nu\in P, \quad \alpha \in C(\mu)_{\lambda}.
\end{equation}
Denote by $C(\mu)_{\lambda,\gamma} \subset C(\mu)_{\lambda}$ the subspace of weight $\gamma$ for
this action; so $\alpha \in C(\mu)_{\lambda,\gamma}$ if $K_\nu.\alpha= q^{(\nu,\gamma)}\alpha$.
Consider the semigroup
\[\Lambda_P=\{(\mu, \lambda,\gamma)\in P_+\times P^2, \, \lambda \;\text{is a weight of}\;V_\mu \;\text{for} \; \rhd ,\, \gamma \;\text{is a weight of}\;V_\mu \;\text{for}\ . \}\]
with the partial order $(\mu,\lambda,\gamma)\preceq (\mu', \lambda',\gamma')$ if and only if $\mu'-\mu,\lambda'-\lambda,\gamma'-\gamma\in D^{-1}Q_+$. Define
\begin{gather*}
[\Lambda_P] =\bigl\lbrace([\mu],[\lambda],[\gamma])\in P_+^n\times P^n\times P^n \\
\hphantom{[\Lambda_P] =\bigl\lbrace([\mu],[\lambda],[\gamma]}{} \mid (\mu_i,\lambda_i,\gamma_i)\in \Lambda_P,\ [\mu] = (\mu_i)_{i=1}^n, [\lambda] = (\lambda_i)_{i=1}^n, [\gamma]=(\gamma_i)_{i=1}^n\bigr\rbrace .
\end{gather*}
Let us put the following lexicographic order on $ [\Lambda_P]$, starting from the tail: $([\mu'],[\lambda'],[\gamma'])\preceq ([\mu], [\lambda],[\gamma])$ if $(\mu_n',\lambda_n',\gamma_n')\prec (\mu_n, \lambda_n,\gamma_n)$, or $(\mu_n,\lambda_n,\gamma_n)= (\mu'_n, \lambda'_n,\gamma'_n)$ and $(\mu_{n-1}',\lambda_{n-1}',\gamma_{n-1}')\prec (\mu_{n-1}, \lambda_{n-1},\gamma_{n-1}),\dots$, or $(\mu_k,\lambda_k,\gamma_k)= (\mu'_k, \lambda'_k,\gamma'_k)$ for all $1< k\leq n$ and $(\mu_1',\lambda_1',\gamma_1')\preceq (\mu_1, \lambda_1,\gamma_1)$. (As usual, we write $([\mu'],[\lambda'],[\gamma'])\prec ([\mu], [\lambda],[\gamma])$ for $([\mu'],[\lambda'],[\gamma'])\preceq ([\mu], [\lambda],[\gamma])$ and $([\mu'],[\lambda'],[\gamma'])\ne ([\mu], [\lambda],[\gamma])$.)

Now recall that $\Ll_{0,n} = {\mathcal L}_{0,1}^{\otimes n} = {\mathcal O}_q^{\otimes n}$ as vector spaces. For every $([\mu],[\lambda],[\gamma])\in[\Lambda_P]$, consider the subspace $C([\mu])_{[\lambda],[\gamma]}\subset \Ll_{0,n}$ defined by
\begin{align*}
&C([\mu]) =C(\mu_1)\otimes \dots \otimes C(\mu_n),\qquad
C([\mu])_{[\lambda],[\gamma]}=C(\mu_1)_{\lambda_1,\gamma_1}\otimes \dots \otimes C(\mu_n)_{\lambda_n,\gamma_n}.
\end{align*}
Then $ {\mathcal L}_{0,n}=\bigoplus_{[\mu]\in P_+^n} C({[\mu]})$ and $ C({[\mu]})=\bigoplus_{([\lambda],[\gamma])} C([\mu])_{[\lambda],[\gamma]}$. For every
$([\mu],[\lambda],[\gamma])\in [\Lambda_P]$ define
\begin{align}
&{\mathcal F}_3^{[\mu],[\lambda],[\gamma]} =\bigoplus_{([\mu'],[\lambda'],[\gamma'])\preceq ([\mu],[\lambda],[\gamma])} C([\mu'])_{[\lambda'],[\gamma']},\label{F2grad}\\
&{\mathcal F}_3^{\prec [\mu],[\lambda],[\gamma]} =\bigoplus_{([\mu'],[\lambda'],[\gamma'])\prec ([\mu],[\lambda],[\gamma])} C([\mu'])_{[\lambda'],[\gamma']}.\notag
\end{align}
Clearly, $\Ll_{0,n}$ is the union of the subspaces ${\mathcal F}_3^{[\mu],[\lambda],[\gamma]}$ over all $([\mu],[\lambda],[\gamma])\in [\Lambda_P]$, so these form a vector space filtration of $\Ll_{0,n}$. Let us denote it ${\mathcal F}_3$, and define
\[{\rm Gr}_{{\mathcal F}_3}({\mathcal L}_{0,n})_{[\mu],[\lambda],[\gamma]} = {\mathcal F}_3^{[\mu],[\lambda],[\gamma]}/{\mathcal F}_3^{\prec [\mu],[\lambda],[\gamma]}.\]
This space is canonically identified with $C({[\mu]})_{[\lambda],[\gamma]}$, so the graded vector space associated to~$\mathcal{F}_3$~is%
\begin{equation}\label{decompgrad}
{\rm Gr}_{{\mathcal F}_3}({\mathcal L}_{0,n}) = \bigoplus_{([\mu],[\lambda],[\gamma])\in [\Lambda_P]} {\rm Gr}_{{\mathcal F}_3}({\mathcal L}_{0,n})_{[\mu],[\lambda],[\gamma]} = \bigoplus_{([\mu],[\lambda],[\gamma])\in [\Lambda_P]} C({[\mu]})_{[\lambda],[\gamma]}.
\end{equation}
 We claim that ${\mathcal F}_3$ is an algebra filtration with respect to the product of $\Ll_{0,n}$, and therefore~${{\rm Gr}_{{\mathcal F}_3}({\mathcal L}_{0,n})}$ is a graded algebra.

For notational simplicity, let us prove it for $n=2$, the general case being strictly similar. Recall the $R$-matrix factorization \eqref{Rmatfact}. Take tuples $([\mu],[\lambda],[\gamma]) =((\mu_1,\mu_2),(\lambda_1,\lambda_2),(\gamma_1,\gamma_2))$ and $([\mu'],[\lambda'],[\gamma'])=((\mu_1',\mu_2'),(\lambda_1',\lambda_2'),(\gamma_1',\gamma_2'))$ in $[\Lambda_P]$, and elements $\alpha\otimes \beta \in C([\mu])_{[\lambda],[\gamma]}$ and~${\alpha'\otimes \beta'\in C([\mu'])_{[\lambda'],[\gamma']}}$. Recall from \eqref{prodL0nform} that the product of $\Ll_{0,2}$ is given by the formula
\begin{gather}
(\alpha\otimes \beta) (\alpha' \otimes \beta') \nonumber\\
\qquad= \sum_{(R^1),\dots, (R^4)} \alpha \bigl(S(R^3_{(1)}R^4_{(1)})\rhd \alpha' \lhd R^1_{(1)}R^2_{(1)}\bigr)\otimes \bigl( S(R^1_{(2)}R^3_{(2)}) \rhd \beta \lhd R^2_{(2)}R^4_{(2)}\bigr)\beta'.\label{prodL0nform-duplicate}
\end{gather}
For every $\nu\in P$ and any of the components \smash{$R^1_{(2)},\dots,R^4_{(2)}$}, denoting by $-r_j\in -Q_+$ the weight of $R^j_{(2)}$, we have
\begin{gather*}
K_\nu \rhd \big( S\big(R^1_{(2)}R^3_{(2)}\big) \rhd \beta \lhd R^2_{(2)}R^4_{(2)}\big) \\
 \qquad= \sum_{(\beta)} \beta_{(1)}\big(R^2_{(2)}R^4_{(2)}\big)\big(K_\nu S\big(R^1_{(2)}R^3_{(2)}\big)\rhd \beta_{(2)}\big)\\
 \qquad= q^{-(\nu,r_1+r_3)}\sum_{(\beta)} \beta_{(1)}\big(R^2_{(2)}R^4_{(2)}\big)\big(S\big(R^1_{(2)}R^3_{(2)}\big)K_\nu \rhd \beta_{(2)}\big)\\
 \qquad= q^{(\nu,\lambda_2-r_1-r_3)}\sum_{(\beta)} \beta_{(1)}\bigl(R^2_{(2)}R^4_{(2)}\bigr)\big(S\big(R^1_{(2)}R^3_{(2)}\big)\rhd \beta_{(2)}\big)\\
 \qquad= q^{(\nu,\lambda_2-r_1-r_3)} \big( S\big(R^1_{(2)}R^3_{(2)}\big) \rhd \beta \lhd R^2_{(2)}R^4_{(2)}\big).
\end{gather*}
By similar computations for the action ${\rm coad}(K_\nu^{-1})$, and for all terms in the right-hand side of \eqref{prodL0nform-duplicate}, and using \eqref{gradalg01} componentwisely, we find that
\[\alpha \big(S\big(R^3_{(1)}R^4_{(1)}\big)\rhd \alpha' \lhd R^1_{(1)}R^2_{(1)}\big)\otimes \big( S\big(R^1_{(2)}R^3_{(2)}\big) \rhd \beta \lhd R^2_{(2)}R^4_{(2)}\big)\beta' \in \Ff_3^{[\mu]+[\mu'], [\lambda''],[\gamma'']},\]
where
\begin{align*}
&\lambda'' = (\lambda_1+\lambda_1'+r_3+r_4,\lambda_2+\lambda_2'-r_1-r_3),\\
& \gamma'' = (\gamma_1+\gamma_1'+r_1+r_2+r_3+r_4,\gamma_2+\gamma_2'-r_1-r_2-r_3-r_4).
\end{align*}
Since $r_1+r_2+r_3+r_4=0$ implies $r_1=r_2=r_3=r_4=0$, by the order we have put on $[\Lambda_P]$, we deduce
 \[(\alpha\otimes \beta)(\alpha'\otimes \beta') \in \Ff_3^{[\mu]+[\mu'],[\lambda]+[\lambda'],[\gamma]+[\gamma']}.\]
Note that the filtration $\mathcal{F}_3$, taking the action \eqref{coadweightnew} into account, is crucial for this argument to work. Similar arguments work for any $n\geq 2$. This proves that ${\rm Gr}_{{\mathcal F}_3}({\mathcal L}_{0,n})$ is a graded algebra. We denote its product by $\circ_n$.

Next, we show that \eqref{filtration=} implies the analogous property for the product $\circ_n$. For simplicity of notations let us again assume that $n=2$. Recall that the product $\circ_2$ is defined on homogeneous elements \smash{$\overline{\alpha\otimes \beta} \in {\rm Gr}_{{\mathcal F}_3}({\mathcal L}_{0,n})_{[\mu],[\lambda]}$} and \smash{$\overline{\alpha'\otimes \beta'} \in {\rm Gr}_{{\mathcal F}_3}({\mathcal L}_{0,n})_{[\mu'],[\lambda']}$} by \[\overline{\alpha\otimes \beta}\circ_n\overline{\alpha'\otimes \beta'} = (\alpha\otimes \beta)(\alpha'\otimes \beta') + {\mathcal F}_3^{\prec [\mu+\mu'],[\lambda+\lambda']}.\]
Clearly, \eqref{filtration=} gives $(C(\mu_1)\circ C(\mu_1'))\otimes (C(\mu_2)\circ C(\mu_2')) = C([\mu + \mu'])$, and \eqref{prodL0nform-duplicate} gives
 \[C([\mu]) \circ_n C([\mu']) \subset (C(\mu_1)\circ C(\mu_1'))\otimes (C(\mu_2)\circ C(\mu_2')).\]
The converse inclusion holds true as well, as one can see by expressing, reciprocally, the (componentwise) product of $\Ll_{0,1}^{\otimes n}$ in terms of the product of $\Ll_{0,n}$ via the formula \eqref{prodL0nform3}. In conclusion,%
\[
C([\mu]) \circ_n C([\mu']) = C([\mu + \mu']).
\]
We are left to show that \eqref{echangegrad} generalizes to $\Ll_{0,n}$. First, note that for every $1\leq a\leq n$ the embedding $\mathfrak{i}_a\colon \Ll_{0,1}\ra \Ll_{0,n}$ in \eqref{plgtL01a} is a morphism of the filtered algebras $(\Ll_{0,1},\mathcal{F}_2)$ and $(\Ll_{0,n},\mathcal{F}_3)$, in the sense that
\[\mathfrak{i}_a\big({\mathcal F}_2^{\mu,\lambda}\big) \subset \sum_{\gamma \in P} {\mathcal F}_3^{[\mu_a],[\lambda_a],[\gamma_a]},
\]
where by definition $[\mu_a] = (0,\dots,0,\mu,0,\dots,0)$ with $\mu$ on the $a$-th entry, and similarly $[\lambda_a] = (0,\dots,0,\lambda,0,\dots,0)$ and $[\gamma_a] = (0,\dots,0,\gamma,0,\dots,0)$. Therefore, the relation \eqref{echangegrad} yields in ${\rm Gr}_{{\mathcal F}_3}({\mathcal L}_{0,n})$ similar relations between elements of the form (matrix coefficient)$\otimes 1$, or $1 \otimes$(matrix coefficient).

We now consider the case of ``mixed'' products. We give the details when $n=2$, the general case being similar. Let us write the twist $F$ in \eqref{prodL0nform2} as
\[F = \sum_{(F)} F_{(1)}\otimes F_{(2)} = \sum_{(F)} F_{(1)1}\otimes F_{(1)2}\otimes F_{(2)1}\otimes F_{(2)2},\]
that is, we set $F_{(1)1} := R^2_{(2)}R^4_{(2)}$, $F_{(1)2} := R^1_{(2)}R^3_{(2)}$, $F_{(2)1} := R^1_{(1)}R^2_{(1)}$, $F_{(2)2} := R^3_{(1)}R^4_{(1)}$. Put~${d(\mu):= \dim(V_{\mu})}$, $\mu\in P_+$, and
\[\Delta^{(2)}\big({}_{{\mu_2}}\phi^{k_2}_{l_2}\big) = \sum_{p,s=1}^{{ d(\mu_2)}}{}_{{\mu_2}}\phi^{k_2}_{p}\otimes {}_{{\mu_2}}\phi^{p}_{s}\otimes {}_{{\mu_2}}\phi^{s}_{l_2},\qquad \Delta^{(2)}\big({}_{{\mu_1'}}\phi^{k_1'}_{l_1'}\big) = \sum_{p',s'=1}^{{ d(\mu_1')}}{}_{{\mu_1'}}\phi^{k_1'}_{p'}\otimes {}_{{\mu_1'}}\phi^{p'}_{s'}\otimes {}_{{\mu_1'}}\phi^{s'}_{l_1'}.\]
From \eqref{prodL0nform-duplicate}, one obtains
\begin{align}
\big(1 \otimes {}_{{\mu_2}}\phi^{k_2}_{l_2}\big) \big({}_{{\mu_1'}}\phi^{k_1'}_{l_1'} \otimes 1\big)
 ={}& \sum_{(F)} \sum_{p,s=1}^{\scriptscriptstyle{ d(\mu_2)}} \sum_{p',s'=1}^{\scriptscriptstyle{ d(\mu_1')}} \big({}_{{\mu_1'}}\phi^{p'}_{s'} \big({}_{{\mu_1'}}\phi^{k_1'}_{p'}(F_{(2)1}){}_{{\mu_1'}}\phi^{s'}_{l_1'}(S(F_{(2)2}))\big)\big) \nonumber\\
&\otimes \big({}_{{\mu_2}}\phi^{p}_{s} \big({}_{{\mu_2}}\phi^{k_2}_{p}(F_{(1)1}){}_{{\mu_2}}\phi^{s}_{l_2}(S(F_{(1)2}))\big)\big).\label{identcoefnew}
\end{align}
It is immediate that
\[{}_{{\mu_1'}}\phi^{p'}_{s'} \otimes {}_{{\mu_2}}\phi^{p}_{s}\in C(\mu_1')_{w(e_{s'}^{\mu_1'}),w(e_{s'}^{\mu_1'})-w(e_{p'}^{\mu_1'})}\otimes C(\mu_2)_{w(e_{s}^{\mu_2}),w(e_{s}^{\mu_2})-w(e_{p}^{\mu_2})}.\]
As in (iv) above, for every $\mu\in P_+$ we order the weight vectors $e_1^\mu,\dots,e_m^\mu$ so that $w\big(e_i^\mu\big) > w\big(e_j^\mu\big)$ implies $i < j$. With such an ordering the factorization \smash{$R=\Theta\hat{R}$} (see~\eqref{Rmatfact}) implies
\[
{}_{{\mu_2}}\phi^{k_2}_{p}(F_{(1)1}){}_{{\mu_2}}\phi^{s}_{l_2}(S(F_{(1)2}))=0 \qquad \text{unless $k_2\geq p$ and $s\geq l_2$},
\]
and
 \[
 {}_{{\mu_1'}}\phi^{k_1'}_{p'}(F_{(2)1}){}_{{\mu_1'}}\phi^{s'}_{l_1'}(S(F_{(2)2}))=0\qquad \text{unless $k_1'\leq p'$ and $s'\leq l_1'$}.
 \]
Since $s>l_2$, we have $w\big(e_{s}^{\mu_2}\big)\leq w\big(e_{l_2}^{\mu_2}\big)$, and if $w\big(e_{s}^{\mu_2}\big) < w\big(e_{l_2}^{\mu_2}\big)$, then \smash{${}_{{\mu_2}}\phi^{p}_{s} \in \mathcal{F}_2^{< \mu_2,w(e_{l_2}^{\mu_2})}$}. In this last situation, the summands \smash{${}_{{\mu_1'}}\phi^{p'}_{s'} \otimes {}_{{\mu_2}}\phi^{p}_{s}$} in the sum above vanish in ${\rm Gr}_{\mathcal{F}_3}(\Ll_{0,2})$. In order to find all the non-zero summands, we have to consider also the weights with respect to the action~\eqref{coadweightnew}. Since $k_2\geq p$ implies $w\big(e_{k_2}^{\mu_2}\big)\leq w\big(e_{p}^{\mu_2}\big)$, we have~${w\big(e_{s}^{\mu_2}\big)-w\big(e_{p}^{\mu_2}\big) \leq w\big(e_{l_2}^{\mu_2}\big)-w\big(e_{k_2}^{\mu_2}\big)}$. Therefore, the summands which are non-zero in ${{\rm Gr}_{\mathcal{F}_3}(\Ll_{0,2})}$ have both weights $w\big(e_{s}^{\mu_2}\big) = w\big(e_{l_2}^{\mu_2}\big)$ and $w\big(e_{p}^{\mu_2}\big) = w\big(e_{k_2}^{\mu_2}\big)$. Doing similarly with the weights of~\smash{${}_{{\mu_1'}}\phi^{p'}_{s'}$}, we find that also \smash{$w\big(e_{s'}^{\mu_1'}\big) = w\big(e_{l_1'}^{\mu_1'}\big)$} and \smash{$w\big(e_{p'}^{\mu_1'}\big) = w\big(e_{k_1'}^{\mu_1'}\big)$}. When all these conditions on weights are satisfied, the corresponding components $F_{(1)1}$, $F_{(1)2}$, $F_{(2)1}$, $F_{(2)2}$ have zero weight. Therefore, the sum reduces to
\begin{align*}
&\sum_{(F)} {}_{{\mu_2}}\phi_{k_2}^{k_2}(F_{(1)1}){}_{{\mu_2}}\phi_{l_2}^{l_2}(S(F_{(1)2})){}_{{\mu_1'}}\phi_{k_1'}^{k_1'}(F_{(2)1}){}_{{\mu_1'}}\phi_{l_1'}^{l_1'}(S(F_{(2)2}) \\ & \qquad =\big\langle {}_{{\mu_2}}\phi_{k_2}^{k_2}\otimes {}_{{\mu_2}}\phi_{l_2}^{l_2}\otimes {}_{{\mu_1'}}\phi_{k_1'}^{k_1'}\otimes {}_{{\mu_1'}}\phi_{l_1'}^{l_1'}, \Theta_{13}\Theta_{14}^{-1}\Theta_{24}\Theta_{23}^{-1}\big\rangle = q^{(w(e_{k_2}^{\mu_2})-w(e_{l_2}^{\mu_2}), w(e_{k_1'}^{\mu_1'}) - w(e_{l_1'}^{\mu_1'}) )}.
\end{align*}
Denoting by $q'_{k_2l_2k_1'l_1'}$ this scalar, it follows
\begin{align*}
\big(1 \otimes {}_{{\mu_2}}\phi^{k_2}_{l_2}\big) \circ_2 \big({}_{{\mu_1'}}\phi^{k_1'}_{l_1'} \otimes 1\big) & = q'_{k_2l_2k_1'l_1'} \ {}_{{\mu_1'}}\phi^{k_1'}_{l_1'} \otimes {}_{{\mu_2}}\phi^{k_2}_{l_2} = q'_{k_2l_2k_1'l_1'} \big({}_{{\mu_1'}}\phi^{k_1'}_{l_1'} \otimes 1\big)\circ_2 \big(1 \otimes {}_{{\mu_2}}\phi^{k_2}_{l_2}\big).
\end{align*}
This is the relation analogous to \eqref{echangegrad} for mixed products in ${\rm Gr}_{\mathcal{F}_3}(\Ll_{0,2})$.

Recall that in \eqref{echangegrad+} we denoted by $u_1,\dots,u_M$ the ordered list of matrix coefficients ${}_{\varpi_k}\phi_{i}^{j}$. Let us order in a lexicographic way the elements $u_i\otimes u_j$, i.e., as a sequence \smash{$u_1^{(2)},\dots,u_{M^2}^{(2)}$} such that the following condition holds: if \smash{${}_{\varpi_{l'}}\phi_{s'}^{t'} < {}_{\varpi_{k'}}\phi_{i'}^{j'}$}, or \smash{${}_{\varpi_{l'}}\phi_{s'}^{t'} = {}_{\varpi_{k'}}\phi_{i'}^{j'}$} and \smash{${}_{\varpi_l}\phi_{s}^{t} < {}_{\varpi_k}\phi_{i}^{j}$}, then~\smash{$u_a^{(2)}:={}_{\varpi_k}\phi_{i}^{j} \otimes {}_{\varpi_{k'}}\phi_{i'}^{j'}$} and \smash{$u_b^{(2)}:= {}_{\varpi_l}\phi_{s}^{t} \otimes {}_{\varpi_{l'}}\phi_{s'}^{t'}$} satisfy \smash{$u_b^{(2)} < u_a^{(2)}$}. Then, for this ordering the polynormal relations \eqref{echangegrad+} hold true for all \smash{$1\leq u_j^{(2)} < u_i^{(2)}\leq M^2$}. As described after \eqref{echangegrad++}, it follows that ${\rm Gr}_{\mathcal{F}_3}(\Ll_{0,n})$ is Noetherian. The filtration $\mathcal{F}_3$ being well founded, it implies that~$\Ll_{0,n}$ is Noetherian too.

Finally, we consider the $A$-algebra $\Ll_{0,n}^A$, and prove it is Noetherian. We proceed in exactly the same way as for $\Ll_{0,n}$, changing the generators and replacing key arguments of the steps (i)--(iv) by the corresponding results over $A$. Let us describe these modifications step by step.

First, consider the case $n=1$. Recall the $A$-lattices \smash{${}_A \overset{\raisebox{-1.5pt}{\scriptsize$\smallbullet$}}{C} (\lambda)$} (see \eqref{AClambda}), and the decomposition \eqref{OAweightdecomp} of $\Oo_A$ into weight subspaces. In particular, have a decomposition into weight subspaces for the left coregular action,
\[{}_A \overset{\raisebox{-1.5pt}{\scriptsize$\smallbullet$}}{C} (\lambda) = \bigoplus_{\lambda'\in P} {}_A \overset{\raisebox{-1.5pt}{\scriptsize$\smallbullet$}}{C} (\lambda)_{\lambda'}.\]
Define
\[{}_A\mathcal{F}_2^{\mu,\lambda} := \bigoplus_{(\mu',\lambda')\preceq (\mu,\lambda)} {}_A \overset{\raisebox{-1.5pt}{\scriptsize$\smallbullet$}}{C} (\mu')_{\lambda'}.\]
Recall that every $A$-module of matrix coefficients $({\rm End}({}_A V_\mu))^*$, $\mu\in P_+$, is contained in $\Oo_A(\leq \mu)$, and by inverting over $\mc(q)$ the corresponding linear triangular system between basis elements, and using that the order relation $\preceq$ is finer than $\leq$, we obtain an inclusion
\[
 \bigoplus_{\mu'\preceq \mu}\ {}_A \overset{\raisebox{-1.5pt}{\scriptsize$\smallbullet$}}{C} (\mu') \subset \bigoplus_{\mu'\preceq \mu} \ C(\mu')
\]
 (see \eqref{exactsequence}--\eqref{tensorCOA}). It follows that \smash{${}_A\mathcal{F}_2^{\mu,\lambda} = \mathcal{F}_2^{\mu,\lambda} \cap \Oo_A$}, and therefore, like $\mathcal{F}_2$ the family${}_A\mathcal{F}_2 :=\smash{\big\{{}_A\mathcal{F}_2^{\mu,\lambda}\big\}_{(\mu,\lambda)\in \Lambda}}$ is a well-founded filtration of $\Oo_A$. Put \smash{${}_A\mathcal{F}_2^{\prec \mu,\lambda} = \mathcal{F}_2^{\prec \mu,\lambda} \cap \Oo_A$}, and consider the graded $A$-module \smash{${\rm Gr}_{ {}_A\mathcal{F}_2}\big({\mathcal L}_{0,1}^A\big)$} associated to ${}_A\mathcal{F}_2$. By \eqref{OAsurA}--\eqref{OAfacteurs} and the fact that~${\Oo_A=\Ll_{0,1}^A}$ as an $A$-module, we have the $A$-module decomposition
\[{\rm Gr}_{ {}_A\mathcal{F}_2}\big({\mathcal L}_{0,1}^A\big)=\bigoplus_{(\mu,\lambda)\in \Lambda} {}_AC(\mu)_\lambda,\]
where ${}_AC(\mu)_\lambda$ is the submodule of weight $\lambda$ (for the left coregular action) of
\[{}_AC(\mu) := \left({\rm End}({}_A V_\mu\right))^*.\]
Then, we can proceed as before. By step (i), we deduce that ${}_A\mathcal{F}_2$ is an algebra filtration of $\Ll_{0,1}^A$. By Proposition~\ref{teoLuzstigOAscinde}, the $A$-module \smash{${}_A \overset{\raisebox{-1.5pt}{\scriptsize$\smallbullet$}}{C} (\mu_1+\mu_2)$} has multiplicity one in \smash{${}_A \overset{\raisebox{-1.5pt}{\scriptsize$\smallbullet$}}{C} (\mu_1)\otimes {}_A \overset{\raisebox{-1.5pt}{\scriptsize$\smallbullet$}}{C} (\mu_2)$}. In fact, by step (ii), ${}_AC(\mu_1+\mu_2)$ has multiplicity one in ${}_AC(\mu_1) \bigotimes_A {}_AC(\mu_2)$, so exactly in the same way as \eqref{filtration=}, we obtain in ${\rm Gr}_{ {}_A\mathcal{F}_2}\big({\mathcal L}_{0,1}^A\big)$ the equality
\[
{}_AC(\mu) \circ{}_AC(\nu) = {}_AC(\mu + \nu).
\]
In step (iii), we fixed a basis of each space $C(\mu)$, consisting of a set of matrix coefficients~$\big\{{}_\mu\phi^i_j\big\}$ with respect to dual basis of weight vectors of the modules $V_\mu$ and $V_\mu^*$. In step (iv), the basis elements of $V_\mu$ and $V_\mu^*$ were ordered by means of the weights, and we used the fact that the matrix coefficients in the spaces $C(\varpi_1),\dots,C(\varpi_m)$ form a generating set of the algebra $ {\rm Gr}_{\mathcal{F}_2}({\mathcal L}_{0,1})$. The only property of the matrix coefficients used in the computations was that they are weight vectors for the left coregular action (and later, in the case $n>1$, for the action~\eqref{coadweightnew}).

We can proceed exactly in the same manner by working with the $A$-modules of matrix coefficients ${}_AC(\mu)$. If one wishes to work at the lever of $\Oo_A$, recall that any set of generators of~$\Oo_A$ generates $\Ll_{0,1}^A$ as well (see the proof of Proposition~\ref{L0NAfgfree}). Then, one can replace the basis~\smash{$\big\{{}_\mu\phi^i_j\big\}$} of each space $C(\mu)$ with the canonical basis~\smash{$\dot{\mathbf{B}} [\mu]^*$} of \smash{${}_A \overset{\raisebox{-1.5pt}{\scriptsize$\smallbullet$}}{C} (\mu)$}, and take the generating set of~$\Oo_A$ formed by the elements in \smash{$\dot{\mathbf{B}} [\varpi_i]^*$}, $i=1,\dots,m$ (see Proposition~\ref{genOA} and the comments thereafter). By the integrality properties satisfied by the $R$-matrix and the twists, all the computations in the proof of steps (iii) and (iv) can be done using such basis elements, and eventually take place over $A$ (see \cite[Propositions 4.10 and 6.9]{BR}). Therefore, we obtain a relation like~\eqref{echangegrad+} with coefficients \smash{$\alpha_{ij}^{st}\in A$}. Since $A$ is a Noetherian ring, again this proves~${ {\rm Gr}_{ {}_A\mathcal{F}_2}\big({\mathcal L}_{0,1}^A\big)}$, whence ${\mathcal L}_{0,1}^A$, are Noetherian rings.

This being done, the adaptation of the proof when $n>1$ is immediate. The filtration $\mathcal{F}_3$ has to be replaced with \smash{${}_A\mathcal{F}_3:=\big\{{}_A{\mathcal F}_3^{[\mu],[\lambda],[\gamma]}\big\}_{([\mu],[\lambda],[\gamma])}$}, where \smash{${}_A{\mathcal F}_3^{[\mu],[\lambda],[\gamma]}$} is the $A$-module defined by
\[
{}_A{\mathcal F}_3^{[\mu],[\lambda],[\gamma]}=\bigoplus_{([\mu'],[\lambda'],[\gamma'])\preceq ([\mu],[\lambda],[\gamma])} {}_A \overset{\raisebox{-1.5pt}{\scriptsize$\smallbullet$}}{C} ([\mu'])_{[\lambda'],[\gamma']},
\]
where
\[
{}_A \overset{\raisebox{-1.5pt}{\scriptsize$\smallbullet$}}{C} ([\mu])_{[\lambda],[\gamma]} ={}_A \overset{\raisebox{-1.5pt}{\scriptsize$\smallbullet$}}{C} (\mu_1)_{\lambda_1,\gamma_1}\bigotimes_A \dots \bigotimes_A {}_A \overset{\raisebox{-1.5pt}{\scriptsize$\smallbullet$}}{C} (\mu_n)_{\lambda_n,\gamma_n},
\]
 and \smash{${}_A \overset{\raisebox{-1.5pt}{\scriptsize$\smallbullet$}}{C} (\mu)_{\lambda,\gamma}$} is the subspace of \smash{${}_A \overset{\raisebox{-1.5pt}{\scriptsize$\smallbullet$}}{C} (\mu)_{\lambda}$} of weight $\gamma$ for the action~\eqref{coadweightnew}. Then the proof proceeds in exactly the same way, replacing in~\eqref{identcoefnew} and all subsequent computations the matrix coefficients by the generators of~$\Oo_A$ provided by Proposition~\ref{genOA}. This concludes the proof.
\end{proof}

\begin{teo}\label{MonNoeth} The algebra $\Mm_{0,n} = {\mathcal L}_{0,n}^{U_q}$ is Noetherian and generated over ${\mathbb C}(q)$ by a finite number of elements.
\end{teo}
Our method of proof follows closely that of the Hilbert--Nagata theorem (see \cite{DC}). Let us recall one version of this theorem. Let $K$ be an arbitrary field, $\mathfrak{A}$ a commutative algebra over~$K$ finitely generated by elements $a_1,\dots,a_n$, and $G$ a group of algebra automorphisms of $\mathfrak{A}$.
\begin{teo} If the action of $G$ on $\mathfrak{A}$ is completely reducible on finite-dimensional representations, then the ring $\mathfrak{A}^G$ of invariants of $\mathfrak{A}$ with respect to $G$ is Noetherian and a finitely generated algebra over $K$.
\end{teo}
We recall here the main steps of the proof that we will adapt in order to prove Theorem \ref{MonNoeth}:
\begin{enumerate}\itemsep=0pt
\item[(a)] From the complete reducibility of the action of $G$ on $\mathfrak{A}$, one can define a linear map
\[R\colon \ \mathfrak{A}\rightarrow \mathfrak{A}^G\]
namely the projection onto the space of invariant elements along the sum of nontrivial isotypical components of $\mathfrak{A}$. This linear map is the Reynolds operator; we already discussed it in \eqref{Reynoldsdef} in the case of $U_q$ acting on $\Ll_{0,n}$. By the same arguments we developed there, it satisfies $R(hf)=hR(f)$ for every $f\in \mathfrak{A}$, $h\in \mathfrak{A}^G$.

\item[(b)] Let $I$ be an ideal of $\mathfrak{A}^G$. Then $I=R(\mathfrak{A} I)=\mathfrak{A} I\cap \mathfrak{A}^G$. Because $\mathfrak{A}I$ is an ideal of~$\mathfrak{A}$, and $\mathfrak{A}$ is Noetherian, there exist elements $b_1,\dots , b_s$, that can be chosen in $I\subset \mathfrak{A}^G$, such that~${\mathfrak{A}I= \mathfrak{A} b_1+\dots + \mathfrak{A} b_s}$. Since $I=R(\mathfrak{A} I)=R(\mathfrak{A}b_1+\dots + \mathfrak{A}b_s)=\mathfrak{A}^Gb_1+\dots + \mathfrak{A}^G b_s$, $I$~is finitely generated over $\mathfrak{A}^G$. Therefore, $\mathfrak{A}^G$ is Noetherian.

\item[(c)] Let $\mathfrak{B}$ be an algebra graded over $\mn$ (for simplicity of notations): $ \mathfrak{B}=\bigoplus_{n=0}^{+\infty }\mathfrak{B}_n$, with ${\mathfrak{B}_m.\mathfrak{B}_n \subset \mathfrak{B}_{m+n}}$. The augmentation ideal of $\mathfrak{B}$ is $ \mathfrak{B}^+:=\bigoplus_{n=1}^{+\infty }\mathfrak{B}_n$. If $\mathfrak{B}^+$ is a Noetherian ideal of $\mathfrak{B}$, then~$\mathfrak{B}$ is a finitely generated algebra over $\mathfrak{B}_0$. This is \cite[Lemma~2.4.5]{Sp} (in that statement~$\mathfrak{B}$ is commutative, but this hypothesis is not necessary for the proof).

\item[(d)] Assume that $\mathfrak{A}^G$ is graded over $\mn$ (for simplicity of notations): $ \mathfrak{A}^G=\bigoplus_{n=0}^{+\infty }\mathfrak{A}^G_n$ with ${\mathfrak{A}^G_0=K}$. Then $ \mathfrak{A}^G{}^+=\bigoplus_{n=1}^{+\infty }\mathfrak{A}^G_n$ is an ideal of $\mathfrak{A}^G$, which is Noetherian by (b) above. Applying~$(c)$, we deduce that $\mathfrak{A}^G$ is a finitely generated algebra over~$K$.
\end{enumerate}

\begin{proof}[Proof of Theorem \ref{MonNoeth}] Consider the filtration $\Ff$ of ${\mathcal L}_{0,n}$ by the subspaces
\[{\mathcal F}^{[\mu]} =\bigoplus_{[\mu']\preceq [\mu]} C([\mu']), \qquad \mu\in P_+^n,\]
where $P_+^n$ is given the lexicographic partial order induced from $[\Lambda]$. It is easily seen that $\mathcal{F}$ is an algebra filtration: the coregular actions $\rhd$, $\lhd$ fix globally each component $C(\mu)$ of $\Ll_{0,1}$, so the claim follows from~\eqref{mmtilde},~\eqref{prodL0nform} and the fact that $C(\mu)\star C(\nu) \subset C(\mu+\nu)$ for all $\mu,\nu\in P_+$. Denote by ${\rm Gr}_{\mathcal F}({\mathcal L}_{0,n})$ the corresponding graded algebra. As a vector space, we have
\begin{equation}\label{Grvect}
{\rm Gr}_{{\mathcal F}}({\mathcal L}_{0,n})= \bigoplus_{[\mu]\in P_+^n} C([\mu]).
\end{equation}
Because each space $C([\mu])$ is stabilized by the coadjoint action of $U_q$, \eqref{Grvect} has a key advantage on the refined decomposition \eqref{decompgrad}. Indeed, since ${\mathcal L}_{0,n}$ is a $U_q$-module algebra, the action of $U_q$ is well defined on ${\rm Gr}_{\mathcal F}({\mathcal L}_{0,n})$ and gives it a structure of $U_q$-module algebra. As vector spaces, we have
\[
{\rm Gr}_{{\mathcal F}}({\mathcal L}_{0,n})^{U_q}=\bigoplus_{[\mu]\in P_+^n} C([\mu])^{U_q}.
\]
Now we can adapt the different steps (a)--(d) recalled above:
\begin{enumerate}\itemsep=0pt 
\item[(a$'$)] The action of $U_q$ on ${\rm Gr}_{\mathcal F}({\mathcal L}_{0,n})$ is completely reducible. This follows from \eqref{Grvect} and the fact that the spaces $C(\mu)$ are finite-dimensional and thus completely reducible $U_q$-modules. We can therefore define the Reynolds operator as in (a),
\[R\colon\ {\rm Gr}_{{\mathcal F}}({\mathcal L}_{0,n})\ra {\rm Gr}_{{\mathcal F}}({\mathcal L}_{0,n})^{U_q}.\]
\item[(b$'$)] ${\rm Gr}_{{\mathcal F}}({\mathcal L}_{0,n})$ is Noetherian, because \eqref{Grvect} shows it is filtered by ${\mathcal F}_3$, and the associated graded algebra ${\rm Gr}_{{\mathcal F}_3}({\rm Gr}_{{\mathcal F}}({\mathcal L}_{0,n}))={\rm Gr}_{{\mathcal F}_3}({\mathcal L}_{0,n})$ is Noetherian by Theorem \ref{LonNoeth}. As in~(b), we deduce that \smash{${\rm Gr}_{{\mathcal F}}({\mathcal L}_{0,n})^{U_q}$} is Noetherian. But \smash{${\rm Gr}_{{\mathcal F}}({\mathcal L}_{0,n})^{U_q}={\rm Gr}_{{\mathcal F}}\big({\mathcal L}_{0,n}^{U_q}\big)$}, which implies that~\smash{${\mathcal L}_{0,n}^{U_q}$} is Noetherian.
\item[(c$'$)] (and (d$'$)) Then we can apply the steps (c)--(d). As a result ${\rm Gr}_{{\mathcal F}}({\mathcal L}_{0,n})^{U_q}$ is finitely generated, say by~$k$ non-zero elements $\bar{x}_1,\dots,\bar{x}_k$, which we may assume homogeneous.
\item[(e$'$)] We can now deduce that $\Ll_{0,n}^{U_q}$ is generated by elements $x_i$ with leading terms the
$\bar{x}_i$'s. Indeed, let \smash{$x\in \Ll_{0,n}^{U_q}$}, and $[\mu]\in P_+^n$ such that \smash{$x\in {\mathcal F}^{[\mu]}\!\setminus \! {\mathcal F}^{\prec[\mu]}$}, where \smash{${\mathcal F}^{\prec [\mu]} := \bigoplus_{[\mu']\prec [\mu]}\! C([\mu'])$}. In~\smash{${\rm Gr}_{{\mathcal F}}({\mathcal L}_{0,n})^{U_q}_{[\mu]}\! =\!{\mathcal F}^{[\mu]}/ {\mathcal F}^{\prec [\mu]}$}, we have
\[
 \bar{x}\! =\! \sum_{(i_1,\dots ,i_k)\in I} \lambda_{(i_1,\dots ,i_k)} \bar{x}_1^{i_1}\cdots\bar{x}_k^{i_k}
 \]
 for some finite set ${I\!\subset\! \mn^k}$, scalars \smash{$\lambda_{(i_1,\dots ,i_k)}\in \mc(q)$}, and monomials \smash{$\bar{x}_1^{i_1}\cdots\bar{x}_k^{i_k}$} of degree $[\mu]$. By definition of the product in~\smash{${\rm Gr}_{{\mathcal F}}({\mathcal L}_{0,n})^{U_q}$}, \[
 \bar{x}_1^{i_1}\cdots\bar{x}_k^{i_k} = x_1^{i_1}\cdots x_k^{i_k} + {\mathcal F}^{\prec[\mu]},
 \]
 so \smash{$x_1^{i_1}\cdots x_k^{i_k} \in {\mathcal F}^{[\mu]}\setminus {\mathcal F}^{\prec [\mu]}$}, whence $\smash{\bar{x}_1^{i_1}\cdots\bar{x}_k^{i_k}} = \smash{\overline{x_1^{i_1}\cdots x_k^{i_k}}}$ and
 \[
 x- \sum_{(i_1,\dots ,i_k)\in I} \lambda_{(i_1,\dots ,i_k)} x_1^{i_1}\cdots x_k^{i_k} \in {\mathcal F}^{\prec [\mu]}.
 \]
 The conclusion follows by decreasing induction on $[\mu]$, since at last we terminate at ${{\mathcal F}^{[0]}\cong \mc(q)}$.
\end{enumerate}

By combining the steps (a$'$) to (e$'$), we get that $\Mm_{0,n}$ is a Noetherian and finitely generated ring.
\end{proof}

\begin{Remark} \quad
\begin{itemize}\itemsep=0pt
\item[(1)] Because ${\mathcal L}_{0,1}^{U_q}$ is the center of ${\mathcal L}_{0,1}$, (e$'$) proves it is finitely generated. Of course this follows also from the isomorphism \smash{$\Ll_{0,1}\cong U_q^{\rm lf}$} and the fact that the center of \smash{$U_q^{\rm lf}$} is the center of~$U_q$ (by Theorem~\ref{JLteo1}), plus the well-known description of the latter.
\item[(2)] In the ${\mathfrak{sl}_2}$ case the filtration $\mathcal{F}$ on \smash{${\mathcal L}_{0,n}^{U_q}$} should be related via the Wilson loop isomorphism (defined in \cite[Section 8.2]{BR}) to the filtration of skein algebras of spheres with $n+1$ punctures used in \cite{PS}.
 \end{itemize}
\end{Remark}

\section[Proof of Theorem 1.2]{Proof of Theorem \ref{Llibre}}\label{ORDER} As usual we let $\e$ be a primitive $l$-th root of unity with $l$ odd and $l>d_i$ for all $i\in \{1,\dots,m\}$. We now consider the specialization $\Ll_{0,n}^\e$ of $\Ll_{0,n}$ at $q=\e$, defined in Section \ref{defintform}. Recall the isomorphism of algebras $\eta^*\colon \Oo(G)\ra \mathcal{Z}_0(\Oo_\e)$ (see \eqref{etadef}), and that $\Ll_{0,n}^\e = \Oo_\e^{\otimes n}$ as a vector space. Consider the linear subspace of $\Ll_{0,n}^\e$ defined by $\mathcal{Z}_0\bigl(\Ll_{0,n}^\e\bigr) := \mathcal{Z}_0(\Oo_\e)^{\otimes n}$.
This space is naturally a subalgebra of $\Oo_\e^{\otimes n}$ (endowed with the componentwise product $\star$). In fact, we also have the following.

\begin{prop} \label{Z0L0n} \quad
\begin{itemize}\itemsep=0pt
\item[$(1)$] $\mathcal{Z}_0\bigl(\Ll_{0,n}^\e\bigr)$ is a central subalgebra of the algebra $\Ll_{0,n}^\e$, and the $\mathcal{Z}_0\bigl(\Ll_{0,n}^\e\bigr)$-modules $\Ll_{0,n}^\e$ and \smash{$\Oo_\e^{\otimes n}$}, with actions defined by the respective products of these algebras, do coincide.
\item[$(2)$] $\Ll_{0,n}^\e$ is a free $\mathcal{Z}_0\bigl(\Ll_{0,n}^\e\bigr)$-module of rank $l^{n.\dim \mathfrak{g}}$.
\item[$(3)$] \smash{$\big(\eta^*{}^{-1}\big)^{\otimes n}\colon \mathcal{Z}_0\bigl(\Ll_{0,n}^\e\bigr)\ra\Oo(G)^{\otimes n}$} is an isomorphism of algebras, and $\mathcal{Z}_0\bigl(\Ll_{0,n}^\e\bigr)$ is a Noetherian ring.
\item[$(4)$] The $\mathcal{Z}_0\bigl(\Ll_{0,n}^\e\bigr)$-module $\Ll_{0,n}^\e$ is finite and Noetherian. Therefore, $\Ll_{0,n}^\e$ is a Noetherian ring.
\end{itemize}
\end{prop}
Note that the proof we give in (4) of the fact that $\Ll_{0,n}^\e$ is Noetherian is independent from the proof of Theorem \ref{LonNoeth}.

\begin{proof} (1) Let us show that $\mathcal{Z}_0\bigl(\Ll_{0,n}^\e\bigr)$ is a central subalgebra of $\Ll_{0,n}^\e$. In the case $n=1$, the formula \eqref{mmtilde} implies that $\alpha \beta = \alpha\star \beta$ for all $\alpha\in \mathcal{Z}_0(\Oo_\e)$ and $\beta\in \Ll_{0,1}^\e$. Indeed, by \eqref{mmtilde} we have 
\begin{align*}
\alpha\beta & = \sum_{(R),(R)}(R_{(2')}{S}(R_{(2)}) \rhd \alpha) \star (R_{(1')}\rhd \beta \lhd R_{(1)})\\
& = \sum_{(R),(R),(\alpha), (\beta)} \alpha_{(1)} \star ( \beta_{(1)}(R_{(1)}\alpha_{(3)}(S(R_{(2)}))\beta_{(3)}(R_{(1')}\alpha_{(2)}(R_{(2')}))\beta_{(2)}),
\end{align*}
where all components $\alpha_{(1)}, \alpha_{(2)}, \alpha_{(3)}\in \mathcal{Z}_0(\Oo_\e)$, since $\mathcal{Z}_0(\Oo_\e)$ is a Hopf subalgebra of $\Oo_\e$. But
\[
 \sum_{(R)} R_{(1)}\alpha_{(3)}(S(R_{(2)})) = S^{-1}(\Phi^{-}(S_{\Oo_\e}(\alpha_{(3)})))\in \mathcal{Z}_0(U_\epsilon),
\]
 since $\Phi^-(S_{\Oo_\e}(\alpha_{(3)})) \in \mathcal{Z}_0(U_\epsilon)$ by Theorem \ref{DCLteo1}\,(2). Similarly, \smash{$ \sum_{(R)} R_{(1')}\alpha_{(2)}(R_{(2')})\in \mathcal{Z}_0(U_\epsilon)$}. In general, these elements belong to~$\mathcal{Z}_0(U_\epsilon)$ and not ${\mathcal Z}_0\big(U_\epsilon^{\rm ad}\big)$ because of the ``diagonal'' factor $\Theta$ of the $R$-matrix in \eqref{Rmatfact}. By Lemma~\ref{acttrivial}, \smash{${\mathcal Z}_0\big(U_A^{\rm ad}\big)$} acts by the trivial character $\varepsilon$ (the counit) on specializations of~$\Gamma$-modules. The action of ${\mathcal Z}_0(U_A)$ is the counit $\varepsilon$ multiplied with some powers of $\epsilon^{1/D}$.
However, \cite[Propositions 4.1 and 4.10]{BR} show that such powers of \smash{$\epsilon^{1/D}$} eventually disappear in the sum above; this is because the sum can be rewritten in terms of copies of the quasi $R$-matrix~$\hat{R}$ in~\eqref{Rmatfact} and the pivotal element $\ell$, instead of copies of $R$. Therefore, \begin{equation}\label{produitZ0eps}
\alpha\beta = \sum_{(\alpha), (\beta)} \alpha_{(1)} \star ( \varepsilon(\beta_{(1)})\varepsilon(\alpha_{(3)})\varepsilon(\beta_{(3)})\varepsilon(\alpha_{(2)})\beta_{(2)}) = \alpha \star \beta.
\end{equation}
This shows $\Ll_{0,1}^\e$ and $\Oo_\e$ coincide as modules over
$\mathcal{Z}_0\bigl(\Ll_{0,1}^\e\bigr)= \mathcal{Z}_0(\Oo_\e)$.
 Next, we show that the subalgebras \smash{$\mathcal{Z}_0(\Oo_\e)^{(a)}$} are central in $\Ll_{0,n}^\e$ for all $a=1,\dots,n$. This fact will conclude the proof that $\Ll_{0,n}^\e$ and $\Oo_\e^{\otimes n}$ coincide as $\mathcal{Z}_0\bigl(\Ll_{0,n}^\e\bigr)$-modules, because the subalgebras~\smash{$\mathcal{Z}_0(\Oo_\e)^{(a)}$} generate the space $\mathcal{Z}_0\bigl(\Ll_{0,n}^\e\bigr)$ in $\bigl(\Ll_{0,1}^\e\bigr)^{\otimes n}$, and hence in $\Ll_{0,n}^\e$ (this follows from the comment before~\eqref{prodL0nform2}).

In order to show that \smash{$\mathcal{Z}_0(\Oo_\e)^{(a)}$} is central in $\Ll_{0,n}^\e$ for all $a=1,\dots,n$, it is enough to show \smash{$\mathcal{Z}_0(\Oo_\e)^{(a)}$} commutes with the elements of $\Ll_{0,n}^\e$ supported by the tensor factors \smash{$\bigl(\Ll_{0,1}^\e\bigr)^{(b)}$} with $b\ne a$. Since \smash{$(\alpha)^{(a)}\otimes (\beta)^{(b)} = ((\alpha)^{(a)} \otimes 1)(1\otimes (\beta)^{(b)})$} by the definition, we have to show that \smash{$\big(1\!\otimes\! (\beta)^{(b)}\big)\big((\alpha)^{(a)} \!\otimes \!1\big) = (\alpha)^{(a)}\!\otimes\! (\beta)^{(b)}$ whenever $\alpha\!\in \!\mathcal{Z}_0(\Oo_\e)$}. We have (denoting \smash{$ \sum_{(\alpha),(\alpha),(\alpha),(\alpha)}$} by \smash{$\sum_{(\alpha)^4}$}, \smash{$ \Delta(\alpha_{(1)}) = \sum_{(\alpha)} \alpha_{(1)(1)} \otimes \alpha_{(1)(2)}$} etc.):
\begin{align*}
\big(1\otimes (\beta)^{(b)}\big) \big((\alpha)^{(a)} \otimes 1\big) ={}&\sum_{\scriptscriptstyle(R^i)} \big(S\big(R^3_{(1)}R^4_{(1)}\big)\rhd \alpha \lhd R^1_{(1)}R^2_{(1)}\big)^{(a)} \\
&\otimes \big(S\big(R^1_{(2)}R^3_{(2)}\big) \rhd \beta \lhd R^2_{(2)}R^4_{(2)}\big)^{(b)}
\\
={}&\sum_{\scriptscriptstyle(R^i),(\alpha)^4,(\beta)^2} (\alpha_{(2)})^{(a)}\otimes (\beta_{(2)})^{(b)}\\
& \times \beta_{(1)}\big(\alpha_{(1)(2)}\big(R^2_{(1)}\big) R^2_{(2)} \alpha_{(3)(1)}\big(S\big(R^4_{(1)}\big)\big) R^4_{(2)}\big)\\
& \times \beta_{(3)}\big(\alpha_{(3)(2)}\big(R^3_{(1)}\big)R^3_{(2)} \alpha_{(1)(1)}\big(R^1_{(1)}\big) S\big(R^1_{(2)}\big)\big).
\end{align*}
By Theorem \ref{DCLteo1}\,(2), it follows that
\[
\alpha_{(1)(2)}\big(R^2_{(1)}\big) R^2_{(2)} = \Phi^+(\alpha_{(1)(2)}) \in \mathcal{Z}_0(U_\epsilon),
\]
 and similarly
\[
\alpha_{(3)(1)}\big(S\big(R^4_{(1)}\big)\big) R^4_{(2)}, \alpha_{(3)(2)}\big(R^3_{(1)}\big)R^3_{(2)}, \alpha_{(1)(1)}\big(R^1_{(1)}\big) S\big(R^1_{(2)}\big) \in \mathcal{Z}_0(U_\epsilon).
\]
 Denote by $z$ any such element; ${\mathcal Z}_0\big(U_\epsilon^{\rm ad}\big)$ acts by the trivial character (the counit $\varepsilon$) on specializations of $\Gamma$-modules. By using \cite[Proposition 6.2]{BR}, arguing as above \eqref{produitZ0eps}, we obtain that the expression of $z$ in terms of the corresponding $\alpha_{(i)(j)}$ involves $\varepsilon(z)= \varepsilon(\alpha_{(i)(j)})$ only \big(no root~$\epsilon^{1/D}$\big). It follows
\begin{gather*}
 \beta_{(1)}\big(\alpha_{(1)(2)}\big(R^2_{(1)}\big) R^2_{(2)} \alpha_{(3)(1)}\big(S\big(R^4_{(1)}\big)\big) R^4_{(2)}\big) \\
 \qquad = \varepsilon(\alpha_{(1)(2)}\alpha_{(3)(1)})\beta_{(1)}(1) = \varepsilon(\alpha_{(1)(2)})\varepsilon(\alpha_{(3)(1)})\varepsilon(\beta_{(1)}),\\
\beta_{(3)}\big(\alpha_{(3)(2)}\big(R^3_{(1)}\big)R^3_{(2)} \alpha_{(1)(1)}\big(R^1_{(1)}\big) S\big(R^1_{(2)}\big)\big) = \varepsilon(\alpha_{(3)(2)})\varepsilon(\alpha_{(1)(1)})\varepsilon(\beta_{(3)}).
\end{gather*}
Therefore, $\big(1\otimes (\beta)^{(b)}\big)\big((\alpha)^{(a)} \otimes 1\big) = (\alpha)^{(a)}\otimes (\beta)^{(b)}$. It follows that $\Ll_{0,n}^\e = \Oo_\e^{\otimes n}$ as modules over $\mathcal{Z}_0\bigl(\Ll_{0,n}^\e\bigr)$; for instance when $n=2$, given $\alpha',\beta'\in \mathcal{Z}_0\bigl(\Ll_{0,1}^\e\bigr)$ we have $(\alpha'\otimes \beta')(\alpha\otimes \beta) = (\alpha'\otimes 1)(1\otimes \beta')(\alpha\otimes 1)(1\otimes \beta)$ immediately by \eqref{prodL0nform}, and $(1\otimes \beta')(\alpha\otimes 1) = \alpha\otimes \beta' = (\alpha\otimes 1)(1\otimes \beta')$ as above. Then $(\alpha'\otimes \beta')(\alpha\otimes \beta) = \alpha'\alpha\otimes \beta'\beta$. In particular, $\mathcal{Z}_0\bigl(\Ll_{0,n}^\e\bigr)$ is a central subalgebra of~$\Ll_{0,n}^\e$.

(2) Since $\Ll_{0,n}^\e$ and $\Oo_\e^{\otimes n}$ coincide as modules over $\mathcal{Z}_0\big(\Ll_{0,n}^\e\big)=\mathcal{Z}_0(\Oo_\e^{\otimes n})$, the claim follows from Theorem~\ref{DCLteo1}, that is, from \cite[Theorem~7.2]{DC-L}, which shows that $\Oo_\e$ is a finitely generated projective module of rank $l^{\dim \mathfrak{g}}$ over $\mathcal{Z}_0(\Oo_\e)$, and from the arguments of \cite{BGS} (using that $K_0(\mathcal{O}(G))=\mz$ by~\cite{Marlin}), which imply that this module is free. Alternatively, it follows from the fact that $\Oo_\e$ is a cleft extension of $\Oo(G)$ (see \cite[Remark~2.18\,(b)]{AnGa}, and \cite[Section~3.2]{BC}).

(3) The linear isomorphism $\big(\eta^*{}^{-1}\big)^{\otimes n}\colon \mathcal{Z}_0\bigl(\Ll_{0,n}^\e\bigr)\ra\Oo(G)^{\otimes n}$ is an isomorphism of algebras because $\mathcal{Z}_0\bigl(\Ll_{0,n}^\e\bigr)$ is central in $\Ll_{0,n}^\e$. It implies that $\mathcal{Z}_0\bigl(\Ll_{0,n}^\e\bigr)$ is a Noetherian ring, since $\Oo(G)^{\otimes n} = \Oo(G^n)$ and $G^n$ is an affine algebraic variety.

(4) The fact that $\Ll_{0,n}^\e$ is a finitely generated $\mathcal{Z}_0\bigl(\Ll_{0,n}^\e\bigr)$-module follows from (2); an alternative proof of this fact will be provided at the end of the proof of Theorem \ref{modZ0}. Since $\Ll_{0,n}^\e$ is finite over $\mathcal{Z}_0\bigl(\Ll_{0,n}^\e\bigr)$ and $\mathcal{Z}_0\bigl(\Ll_{0,n}^\e\bigr)$ is Noetherian, $\Ll_{0,n}^\e$ is a Noetherian $\mathcal{Z}_0\bigl(\Ll_{0,n}^\e\bigr)$-module (e.g., by \cite[Proposition~6.5]{AM}). It follows that $\Ll_{0,n}^\e$ is a Noetherian ring (e.g., by \cite[Chapter~1, Section~1.3]{MC-R}).\looseness=1
\end{proof}

Recall that we denote \smash{$U_\e^{\rm lf} = U_A^{\rm lf}\bigotimes_A \mc_\e$} (see \eqref{notambig0}), and $\mathcal{Z}_0(U_\e)\subset U_\e$ is the central polynomial subalgebra generated by \smash{$E_{\beta_k}^l$}, \smash{$F_{\beta_k}^l$}, \smash{$L_i^{\pm l}$}, for $k\in \{1,\dots, N\}$ and $i\in \{1,\dots,m\}$. Since \smash{$\Phi_1\colon \Ll_{0,1}^\e \ra U^{\rm lf}_\e$} is an embedding of algebras (see Corollary \ref{Phinint}), it identifies $\mathcal{Z}_0\bigl(\Ll_{0,1}^\e\bigr)$ with a~central subalgebra of \smash{$U^{\rm lf}_\e$}. Put \smash{$\mathcal{Z}_0\big(U^{\rm lf}_\e\big):= \Phi_1(\mathcal{Z}_0\bigl(\Ll_{0,1}^\e\bigr))$}.
Recall Theorem \ref{JLteo1}, Proposition~\ref{JLteo1overA}, and let \smash{$T^{(l)}$}, \smash{$T_{2-}^{(l)}$} and $T_{2}^{(l)}$ be the subsets of $T$, $T_{2-}$ and~$T_{2}$ formed by the elements $K_{\lambda l}$ with $\lambda\in P$, $\lambda\in -2P_+$ and $\lambda\in 2P$, respectively.

\begin{prop}\label{JLad2} We have \smash{$U_\e = T_{2-}^{-1}U_\e^{\rm lf}[T/T_{2}] = \Phi_1\big(\Ll_{0,1}^\e\big[d^{-1}\big]\big)[T/T_{2}]$}, and therefore the map \smash{$\Phi_1\colon \Ll_{0,1}^\e\big[d^{-1}\big] \ra T_{2-}^{-1}U_\e^{\rm lf}$} is an isomorphism.

 Moreover, \smash{$\mathcal{Z}\big(U_\e^{\rm lf}\big) = U_\e^{\rm lf} \cap \mathcal{Z}(U_\e)$}, and
\[\mathcal{Z}_0(U_\e) = T_{2-}^{(l)-1}\mathcal{Z}_0\big(U_\e^{\rm lf}\big)\big[T^{(l)}/T_{2}^{(l)}\big],\qquad
\mathcal{Z}(U_\e) = T_{2-}^{(l)-1}\mathcal{Z}\big(U_\e^{\rm lf}\big)\big[T^{(l)}/T_{2}^{(l)}\big].\]
\end{prop}

\begin{proof} The first claim follows immediately from Proposition \ref{JLteo1overA} by specialization at ${q=\e}$. For the second claim, the inclusion $U_\e^{\rm lf} \cap \mathcal{Z}(U_\e) \subset \mathcal{Z}\big(U_\e^{\rm lf}\big)$ is clear, and for the converse inclusion we only have to show that the elements of \smash{$\mathcal{Z}\big(U_\e^{\rm lf}\big)$} commute with $T$. They commute with~\smash{$T_2\subset U_\e^{\rm lf}$}, so the conjugation action by elements of $T$ on \smash{$\mathcal{Z}\big(U_\e^{\rm lf}\big)$} has order at most $2$. Let~\smash{$x\in \mathcal{Z}\big(U_\e^{\rm lf}\big)$} with decomposition $x= \sum_i c_i x_i$ with all $c_i\in \mc$ and $x_i$ PBW basis vectors, and let $\lambda\in P$. We have~${K_\lambda x K_{-\lambda} = \sum_i c_i q(x_i) x_i}$, where $q(x_i) \in \e^\mz$ satisfies $q(x_i)^2 = 1$ for all $i$. Because $\e$ has odd order the only possibility is $q(x_i) = 1$, whence $K_\lambda x K_{-\lambda} = x$. The conclusion follows.

The inclusion $\mathcal{Z}_0\big(U_\e^{\rm lf}\big) \subset \mathcal{Z}_0(U_\e)$ follows from the definition $\mathcal{Z}_0\bigl(\Ll_{0,1}^\e\bigr) = \mathcal{Z}_0(\Oo_\e)$, the formula~${\Phi_1 = m\circ \big({\rm id} \otimes S^{-1}\big)\circ \Phi}$, and the fact that $\Phi$ affords an embedding $\mathcal{Z}_0(\Oo_\e)\ra \mathcal{Z}_0(U_\e(G^*))$ (see Theorem \ref{DCLteo1}\,(2)). Since \smash{$T^{(l)}\subset \mathcal{Z}_0(U_\e)$}, we obtain
\[
T_{2-}^{(l)-1}\mathcal{Z}_0\big(U_\e^{\rm lf}\big)\big[T^{(l)}/T_{2}^{(l)}\big]\subset \mathcal{Z}_0(U_\e).
\]
 The proof of the converse inclusion is similar to that in Proposition~\ref{JLteo1overA}. The isomorphism $\smash{\mathcal{Z}_0(\Oo_\e)\big[\psi_{-l\rho}^{-1}\big]}\ra \mathcal{Z}_0(U_\e(G^*))$ of Theorem \ref{DCLteo1}\,(2) implies
\[
F_{\beta_k}^lK_{\beta_k}^l\otimes 1, 1\otimes K_{\beta_k}^{-l}E_{\beta_k}^l\in \Phi\big(\mathcal{Z}_0(\Oo_\e)\big[\psi_{-l\rho}^{-1}\big]\big)
\]
 for every positive root $\beta_k$. Since $\psi_{-l\rho} = \Phi_1^{-1}(K_{-2l\rho})=\psi_{-\rho}^l$ (the $l$-th power of $\psi_{-\rho}$ in $\Ll_{0,1}^\e$), and \[
 \Phi_1\big(\mathcal{Z}_0\bigl(\Ll_{0,1}^\e\bigr)\big[\psi_{-\rho}^{-l}\big]\big) = T_{2-}^{(l)-1}\mathcal{Z}_0\big(U_\e^{\rm lf}\big),
 \] it follows that
\[
F_{\beta_k}^lK_{\beta_k}^l, S^{-1}\bigl(E_{\beta_k}^l\bigr) K_{\beta_k}^l\in T_{2-}^{(l)-1}\mathcal{Z}_0\big(U_\e^{\rm lf}\big).
\]
 Hence $F_{\beta_k}^l, S^{-1}\big(E_{\beta_k}^l\big)\in T_{2-}^{(l)-1}\mathcal{Z}_0\big(U_\e^{\rm lf}\big)\big[T^{(l)}/T_{2}^{(l)}\big]$. The sets $S^{-1}\big(E_{\beta_k}^l\big)\mathcal{Z}_0(U_\e(\mathfrak{h}))$, $k=1,\dots,N$, generate the subalgebra $\mathcal{Z}_0(U_\e(\mathfrak{b}_+))$ of $\mathcal{Z}_0(U_\e)$, so from the triangular decomposition $\mathcal{Z}_0(U_\e)=\mathcal{Z}_0(U_\e(\mathfrak{n}_-))\mathcal{Z}_0(U_\e(\mathfrak{h}))\mathcal{Z}_0(U_\e(\mathfrak{n}_+))$ this proves the inclusion \smash{$\mathcal{Z}_0(U_\e)\subset T_{2-}^{(l) -1}\mathcal{Z}_0\big(U_\e^{\rm lf}\big)\big[T^{(l)}/T_{2}^{(l)}\big]$}.
 From the isomorphism
 \[
 \mathcal{Z}_0(U_\e)\bigotimes_{\mathcal{Z}_0(U_\e)\cap \mathcal{Z}_1(U_\e)} \mathcal{Z}_1(U_\e) \ra \mathcal{Z}(U_\e)
 \]
 (see Theorem \ref{DCKteo1}), and the fact that $\mathcal{Z}(U_q) \subset U_q^{\rm lf}$ \big(whence $\mathcal{Z}_1(U_\e)\subset \mathcal{Z}\big(U_\e^{\rm lf}\big)$\big), the equality \smash{$\mathcal{Z}(U_\e) = T_{2-}^{(l)-1}\mathcal{Z}\big(U_\e^{\rm lf}\big)\big[T^{(l)}/T_{2}^{(l)}\big]$} follows at once.
 \end{proof}

\begin{Remark}\label{Z0viaUe} Let us explain how this can be used to give an interpretation of the isomorphism $\mathcal{Z}_0\bigl(\Ll_{0,1}^\e\bigr)\cong \Oo(G)$. Recall the notations introduced around Theorem \ref{DCKteo1}. Since $G^* = U_+T_GU_-$, we have $\Oo(G^*) = \Oo(U_+)\Oo(T_G)\Oo(U_-)$, and the map $\sigma$ yields an identification
\begin{equation}\label{OG0}
\Oo\big(G^0\big) = \Oo(U_+)\Oo(T_G/(2))\Oo(U_-).
\end{equation}
We can identify $\Oo\big(G^0\big)$ with the subalgebra $(\sigma_{\vert G^*})^*\big(\Oo\big(G^0\big)\big)\subset\Oo(G^*)$. Consider the exterior power $V=\wedge^N \mathfrak{g}$ endowed with the action $\wedge^N {\rm Ad}$ of $G$. Put on $\mathfrak{g}$ a basis consisting of one element~$e_\alpha$ per root space $\mathfrak{g}_\alpha$, along with a basis of $\mathfrak{h}$. Let $v\in V$ be the exterior power of the~$e_\alpha$'s for $\alpha$ negative, and $v^*$ a dual vector such that $v^*(v)=1$ and $v^*$ vanishes on a $T_G$-invariant complement of $v$. It is classical that $G\setminus G^0$ has defining equation~${\delta(g)=0}$, where $\delta$ is the matrix coefficient $\delta(g) = v^*(\pi_V(g)v)$ (see, e.g., \cite[p.~174]{Hum}). Hence $\Oo\big(G^0\big) =\Oo(G)\big[\delta^{-1}\big]$. On $G^0$ we have $\delta(u_+tu_-) = \chi_{-2\rho}(t)$, where $\chi_{-2\rho}$ is the character of $T_G$ associated to the root $-2\rho$. Now we can make the connection with $U_\e$. The isomorphism $\mathcal{Z}_0(U_\e)\cong\Oo(G^*)$ of Theorem~\ref{DCKteo1}\,(2) identifies $\mathcal{Z}_0(U_\e(\mathfrak{h}))=\mc\big[T^{(l)}\big]$ with $\Oo(T_{G})$ by mapping $K_{\lambda l}$ to the character of $T_G$ associated to~$\lambda$. Therefore, it maps \smash{$\mc\big[T_2^{(l)}\big]$ to $\Oo(T_G/(2))$}, and \smash{$T_{2-}^{(l)-1}\mathcal{Z}_0\big(U_\e^{\rm lf}\big)$} to $\Oo\big(G^0\big)$ by \eqref{OG0} and Proposition~\ref{JLad2}. Since $\Oo\big(G^0\big) = \Oo(G)\big[\delta^{-1}\big]$ and~\smash{$T_{2-}^{(l)-1}\mathcal{Z}_0\big(U_\e^{\rm lf}\big) = \mathcal{Z}_0\big(U_\e^{\rm lf}\big)\big[\ell^l\big]$}, it follows that~\smash{$\mathcal{Z}_0\big(U_\e^{\rm lf}\big)$} and $\Oo(G)$ coincide after localization by~$\ell^l$ and $\delta$ respectively. By using the Bruhat decomposition of $G$ as in~\eqref{BruhatactionX} in the proof of Theorem~\ref{modZ0} below, one can deduce \smash{$\mathcal{Z}_0\big(U_\e^{\rm lf}\big)\cong \Oo(G)$}, whence~${\mathcal{Z}_0\bigl(\Ll_{0,1}^\e\bigr)\cong \Oo(G)}$ by injectivity of~$\Phi_1$.
\end{Remark}

Let us make the following observation. Since \smash{$\Ll_{0,n}^\e = \Ll_{0,n}^A \bigotimes_A \mc_\e$}, with \smash{$\Ll_{0,n}^A = \Oo_A^{\otimes n}$} as an $A$-module, and a generating system of \smash{$\Oo_A^{\otimes n}$} is also a generating system of $\Ll_{0,n}^A$, it follows from Proposition~\ref{genOA} and the identities \eqref{formgenOA1}--\eqref{formgenOA2} that $\Ll_{0,n}^\e$ is generated by elements of the form~${\alpha_1 \otimes \dots \otimes \alpha_n}$, where $\alpha_1,\dots , \alpha_n$ belong to the set $C_{\rm gen}$ of matrix coefficients lying on the first and last columns of the matrix representations of $U_A^{\rm res}$ in the canonical bases of the modules~${{}_AV_{\varpi_i}}$,~$i=1,\dots, m$. Denote by $\alpha^{\star k}$, $k\in \mathbb{N}$, the $k$-th power of an element~${\alpha\in \Oo_A}$.
\begin{lem} For all $\alpha\in C_{\rm gen}$, $\alpha^{\star l}\in \mathcal{Z}_0\bigl(\Ll_{0,1}^\e\bigr)$.
\end{lem}

\begin{proof} Recall that the Frobenius epimorphism $\eta\colon U_A^{\rm res} \bigotimes_A \mc_\e \ra U(\mathfrak{g})$ in \eqref{etadef} has kernel the ideal $I$ generated by the elements $E_i$, $F_i$, $K_i-1$, and $(K_i;p)_{q_i}$ where $l$ does not divide~$p$, ${i=1,\dots,m}$. It follows that an element of $\Oo_\e$ belongs to $\mathcal{Z}_0(\Oo_\e)=\eta^*(\Oo(G))$ if and only if it vanishes on~$I$. But this is immediate to check for the elements of the form~$\alpha^{\star l}$ with $\alpha\in C_{\rm gen}$, using that~$K_i$ is grouplike and the pure summands of $\Delta(E_i)$ and $\Delta(F_i)$ have one component equal to $1$ or $K_i^{\pm 1}$ and the other component equal to $E_i$ or $F_i$. For instance,
\[
\psi_{\varpi_i}^{\star l}(K_i-1) = \psi_{\varpi_i}(K_i)^l-1 = \e^{l(\alpha_i,\varpi_i)} - 1=0.
\]
 Similarly, for every $\alpha \in C_{\rm gen}$, we find
 \[
\alpha^{\star l}(E_i) =\alpha^{\otimes l}\big(\Delta^{(l)}(E_i)\big)=0 \qquad \text{and} \qquad \alpha^{\star l}(F_i)=\alpha^{\star l}(K_i-1)=0.\tag*{\qed}
\]\renewcommand{\qed}{}
\end{proof}

We need below explicit descriptions of the centers of $\Oo_\e({\rm SL}_2)$ and $\Ll_{0,1}^\e({\mathfrak{sl}_2})$ and their $\mathcal{Z}_0$-subalgebras. Denote by $a$, $b$, $c$, $d$ the standard generators of $\mathcal{O}_q({\rm SL}_2)$, i.e., the matrix coefficients in the basis of weight vectors $v_0$, $v_1=F.v_0$ of the $2$-dimensional irreducible representation $V_1$ of $U_q({\mathfrak{sl}_2})$. As above, denote by $x^{\star k}$, $k\in \mathbb{N}$, the $k$-th power of an element $x\in \Oo_A({\rm SL}_2)$. The algebra $\mathcal{O}_A({\rm SL}_2)$ is generated by $a$, $b$, $c$, $d$; the monomials $a^{\star i } \star b^{\star j } \star d^{\star r }$ and $a^{\star i } \star c^{\star k } \star d^{\star r }$, $i,j,k,r\in \mn, k>0$, form an $A$-basis of $\mathcal{O}_A({\rm SL}_2)$. The algebra $\mathcal{Z}_0(\Oo_\e({\rm SL}_2))$ is generated by $a^{\star l}$,~$b^{\star l}$,~$c^{\star l}$,~$d^{\star l}$; the monomials $a^{\star i l} \star b^{\star j l} \star d^{\star r l}$ and $a^{\star i l} \star c^{\star k l} \star d^{\star r l}$ form a basis of $\mathcal{Z}_0(\Oo_\e({\rm SL}_2))$, and $\mathcal{Z}(\Oo_\e({\rm SL}_2))$ is generated by $\mathcal{Z}_0(\Oo_\e({\rm SL}_2))$ and the elements $b^{\star (l-k)}\star c^{\star k}$, $k=0,\dots,l$ (see~\cite[Proposition~1.4 and the appendix]{DC-L}). We have the relation
\begin{equation}\label{PWrel} a^{\star l}\star d^{\star l}-b^{\star l}\star c^{\star l} =1
\end{equation}
and the Frobenius isomorphism of Parshall--Wang (see \cite[Chapter~7]{PW}) coincides with the map
\[{\rm Fr}_{\rm PW}\colon\ \Oo({\rm SL}_2) \ra \mathcal{Z}_0(\Oo_\e({\rm SL}_2))\]
induced by $\eta^*$; it sends the standard generators $\underline{a}$, $\underline{b}$, $\underline{c}$, $\underline{d}$ of $\Oo({\rm SL}_2) = \mathcal{O}_1({\rm SL}_2)$ respectively to~$a^{\star l}$,~$b^{\star l}$,~$c^{\star l}$,~$d^{\star l}$. Finally, we have seen that $\Oo_\e({\rm SL}_2)$ is a free $\mathcal{Z}_0(\Oo({\rm SL}_2))$-module of rank~$l^3$ (see Theorem \ref{DCLteo1}\,(3)). In \cite{DRZ}, it is shown that a basis of this module is formed by the monomials~$a^mb^nc^{s'}$ and $b^nc^{s''}d^r$, with the integers $m$, $n$, $r$, $s'$, $s''$ in the range
\begin{equation}\label{basisSL2}
1\leq m \leq l-1,\qquad 0 \leq n,r \leq l-1,\qquad m \leq s' \leq l-1,\qquad 0 \leq s'' \leq l-r-1.
\end{equation}
Now consider $\Ll_{0,1}^A({\mathfrak{sl}_2})$. Recall that $\Ll_{0,1}^A=\Oo_A$ as $U_A$-modules. The algebra $\Ll_{0,1}^A({\mathfrak{sl}_2})$ is also generated by $a$, $b$, $c$, $d$; a set of defining relations is (see \cite[Section 5]{BR}):
\begin{gather}
ad=da ,\qquad db = q^2 bd, \qquad cd = q^2dc, \qquad
 ab-ba = -\big(1-q^{-2}\big)bd,\nonumber \\ cb-bc = \big(1-q^{-2}\big)\big(da-d^2\big),\qquad ac-ca=\big(1-q^{-2}\big)dc,\qquad ad-q^2bc=1. \label{rel01}
\end{gather}
The element $\omega:=qa+q^{-1}d$ is central. Let $T_k$, $k\in \mn$, be such that $T_k(x)/2$ is the $k$-th Chebyshev polynomial of the first type in the variable $x/2$. We have (see \cite[Proposition 7.2]{BR}, for the generalization to $\Ll_{0,n}^\e({\mathfrak{sl}_2})$):
\begin{lem}\label{L01sl2} Let $\mathcal{I}$ be the ideal of $\mc\big[\omega,b^l,c^l,d^l\big]$ generated by $\big(T_l(\omega)-d^l\big)d^l-b^lc^l-1$, we have
\[
\mathcal{Z}(\Ll_{0,1}^\e({\mathfrak{sl}_2})) = \mc\big[\omega,b^l,c^l,d^l\big]/\mathcal{I}\qquad \text{and} \qquad \mathcal{Z}_0(\Ll_{0,1}^\e({\mathfrak{sl}_2})) = \mc\big[T_l(\omega),b^l,c^l,d^l\big]/\mathcal{I}.
\]
\end{lem}
Here $b^l$, $c^l$, $d^l$ are the $l$-th powers of $b$, $c$, $d$ computed using the product of $\Ll_{0,1}^A({\mathfrak{sl}_2})$, not the product $\star$ of $\mathcal{Z}_0(\Oo_\e({\rm SL}_2))$. The above generator of $\mathcal{I}$ can be interpreted as a determinant, and~$\omega$ as a quantum trace on $V_1$. The following has also been observed in \cite{KQ}.

\begin{lem}\label{identZF} Viewed as elements of $\Oo_A({\rm SL}_2)$, $T_l(\omega)-d^l = a^{\star l}$ and $x^l = x^{\star l}$, $x\in \{b,c,d\}$.
\end{lem}
\begin{proof} Let $\alpha$ and $\varpi$ be the simple root and fundamental weight of ${\mathfrak{sl}_2}$. In the notations of~\eqref{formphi1K}, we have $b = \psi^{-\alpha}_{-\varpi}$, $c=\psi^\alpha_{-\varpi}$, $d=\psi_{-\varpi}$; the formulas of $\Phi$ give
\[
\Phi_1\big(b^{\star l}\big) = \big(q-q^{-1}\big)^lF^l , \qquad \Phi_1\big(c^{\star l}\big) = \big(q-q^{-1}\big)^lE^lK^{-l},
\qquad \Phi_1\big(d^{\star l}\big) = K^{-l}.
\]
 These coincide respectively with $\Phi_1\big(b^{l}\big)$, $\Phi_1(c^{l})$, $\Phi_1\big(d^{l}\big)$ (see \cite[equation~(5.3)]{BR}). By passing to the localization $\Oo_A({\rm SL}_2)\big[d^{-1}\big]$, and using Par\-shall--Wang's relation \eqref{PWrel}, one deduces easily
\[
\Phi_1\big(a^{\star l}\big) = K^l + \big(q-q^{-1}\big)^{2l}F^lE^l = T_l(\Omega)-K^{-l},
\]
 where \smash{$\Omega=\big(\e-\e^{-1}\big)^2FE+\e K+\e^{-1}K^{-1}$} is the Casimir element, and $T_l(x)/2$ is the $l$-th Chebyshev polynomial of the first type in the variable $x/2$. We have $\Phi_1(\omega) = \Omega$, so $\Phi_1\big(a^{\star l}\big) = T_l(\omega)-d^l$. The conclusion follows from the injectivity of $\Phi_1$.
\end{proof}

This lemma proves that we have a commutative diagram
\[\xymatrix{
 \Oo({\rm SL}_2) \ar[r]^{{\rm Fr}_{\rm PW}\ \ \ } \ar[rd]_{\rm Fr} & \mathcal{Z}_0(\Oo_\e({\rm SL}_2)) \ar[d]\ar@{^{(}->}[r] & \Oo_\e({\rm SL}_2) \ar[d] \\
 & \quad \quad \mathcal{Z}_0(\Ll_{0,1}^\e({\mathfrak{sl}_2})) \ar@{^{(}->}[r] & \Ll_{0,1}^\e({\mathfrak{sl}_2}),}\]
where ${\rm Fr}_{\rm PW}$ is Parshall--Wang's Frobenius isomorphism recalled above, ${\rm Fr}$ is the Frobenius isomorphism introduced in \cite[Definition~7.1]{BR}, and the vertical arrows are the identifications as vector spaces (the middle one proved by Proposition~\ref{Z0L0n}).
\begin{Remark}\label{changebasis}
By Lemma \ref{L01sl2}, the monomials $T_l(\omega)^{i} b^{jl} d^{rl}$ and $T_l(\omega)^{i} c^{kl} d^{rl}$, for $i,j,k,r\in \mn$ and~${k>0}$, form an $A$-basis of $\mathcal{Z}_0(\Ll_{0,1}^\e({\mathfrak{sl}_2}))$. It is straightforward (though cumbersome) to express these basis elements in terms of the basis elements \smash{$a^{\star i l} \star b^{\star j l} \star d^{\star r l}$} and \smash{$a^{\star i l} \star c^{\star k l} \star d^{\star r l}$} of $\mathcal{Z}_0(\Oo_\e({\rm SL}_2))$, and conversely; this can be done by using Lemma \ref{identZF}, the formula \eqref{mmtilde} or the inverse one (expressing $\star$ in terms of the product of $\Ll_{0,1}$, see \cite[equation~(4.6)]{BR}), and the formula of the coproduct $\Delta\colon \Ll_{0,1}^\e({\mathfrak{sl}_2})) \ra \Ll_{0,2}^\e({\mathfrak{sl}_2}))$ restricted to $\mathcal{Z}_0(\Ll_{0,1}^\e({\mathfrak{sl}_2}))$ (given in \cite[Lemma~7.5]{BR}).
\end{Remark}

Since \smash{$\Ll_{0,1}^A=\Oo_A$} as an $A$-module, the functionals $t_i$ in Proposition \ref{funcDCL} can be seen as maps $t_i\colon \Ll_{0,1}^A \ra A$. Though the algebra structures of $\Oo_\e$ and $\Ll_{0,1}^\e$ are very different, $\Ll_{0,1}^\e$ satisfies a~result analogous to Proposition \ref{funcDCL}:
\begin{prop}\label{functL01} The maps $\lhd t_i$ preserve $\mathcal{Z}_0\bigl(\Ll_{0,1}^\e\bigr)$, and they satisfy $(f\lhd t_i)(a) = f(n_ia)$ and~${(f\alpha)\lhd t_i = (f\lhd t_i)(\alpha\lhd t_i)}$ for every $f\in \mathcal{Z}_0\bigl(\Ll_{0,1}^\e\bigr)$, $a\in G$, $\alpha\in \Ll_{0,1}^\e$.
\end{prop}

\begin{proof} The first two claims follow from Proposition \ref{funcDCL} and the definition ${\mathcal{Z}_0\bigl(\Ll_{0,1}^\e\bigr) = \mathcal{Z}_0(\Oo_\e)}$.

The last claim follows from the case $\mathfrak{g}={\mathfrak{sl}_2}$, as in the proof of \cite[Proposition~7.1]{DC-L}. In fact, it is enough to show that $t(fg) = t(f)t(g)$ for every $f\in \mathcal{Z}_0(\Ll_{0,1}^\e({\mathfrak{sl}_2})$, $g\in \Ll_{0,1}^\e({\mathfrak{sl}_2})$; for completeness we explain this in Appendix~\ref{actO}, see~\eqref{fundrelt}. A word of caution is needed: the proof of \eqref{fundrelt} uses that $\Delta\colon \Oo_\e \ra \Oo_\e \otimes \Oo_\e$ is a morphism of algebras. The analogous property for $\Ll_{0,1}^\e$ is that $\Delta$ yields a morphism of algebras $\Delta\colon \Ll_{0,1}^\e \ra \Ll_{0,2}^\e$. Since the algebra structure of $\Ll_{0,2}^\e$ is not the product one on $\Ll_{0,1}^\e\otimes \Ll_{0,1}^\e$, it is not true in general that
\[ \sum_{(f),(g)} (f_{(1)}\otimes f_{(2)})(g_{(1)} \otimes g_{(2)}) = \sum_{(f),(g)} f_{(1)}g_{(1)}\otimes f_{(2)}g_{(2)}
\]
 for every $f,g\in \Ll_{0,1}^\e$. However, it holds whenever $f\in \mathcal{Z}_0\bigl(\Ll_{0,1}^\e\bigr)$, since $\Delta(\mathcal{Z}_0\bigl(\Ll_{0,1}^\e\bigr)) \subset \mathcal{Z}_0\bigl(\Ll_{0,1}^\e\bigr) \otimes \mathcal{Z}_0\bigl(\Ll_{0,1}^\e\bigr)$ and therefore $f_{(2)}\in \mathcal{Z}_0\bigl(\Ll_{0,1}^\e\bigr) = \mathcal{Z}_0(\Oo_\e)$ commutes in $\Ll_{0,2}^\e$ with any ${g_{(1)}\in \Ll_{0,1}^\e = \Oo_\e}$.

It is enough to prove the identity $t(fg) = t(f)t(g)$ when $f$ ranges in a set of generators of the algebra $\mathcal{Z}_0(\Ll_{0,1}^\e({\mathfrak{sl}_2}))$. So one can take $f$ among, say, \smash{$T_l(\omega)-d^l = a^{\star l}$} and \smash{$x^l = x^{\star l}$}, $x\in \{b,c,d\}$ (using Lemma \ref{L01sl2}). By \eqref{mmtilde} and Proposition \ref{tidstar}, we have
\[t(fg) = \sum_{(R),(R)} t\left(R_{(2')}{S}(R_{(2)}) \rhd f\right) t\left(R_{(1')}\rhd g \lhd R_{(1)}\right).\]
Expanding coproducts and using that $R^{-1} = (S\otimes {\rm id})(R)$, we deduce
\begin{align*}\label{tprod1}
 t(fg) & = \sum_{(f),(R),(R)} t\big(f_{(1)}\big) \left\langle f_{(2)} , R_{(2')}{S}(R_{(2)}) \right\rangle t\left(R_{(1')}\rhd g \lhd R_{(1)}\right)\\
 & = \sum_{(f),(R),(R)} t\big(f_{(1)}\big) t(\left\langle f_{(2)} , R_{(2')}\right\rangle R_{(1')}\rhd g \lhd \left\langle f_{(3)} ,S(R_{(2)}) \right\rangle R_{(1)})\\
 & = \sum_{(f)} t(f_{(1)}) t\big( S^{-1}(\Phi^-(f_{(2)}))\rhd g \lhd S^{-2}(\Phi^-(f_{(3)}))\big)\\
 & = \sum_{(f)} t(f_{(1)}) \big\langle g, S^{-2}(\Phi^-(f_{(3)}))\underline{w} S^{-1}(\Phi^-(f_{(2)})) \big\rangle\\
 & = \sum_{(f)} t(f_{(1)}) \varepsilon\big( S^{-2}(\Phi^-(f_{(3)})) \big) \varepsilon\big( S^{-1}(\Phi^-(f_{(2)})) \big) t(g),
\end{align*}
where $\underline{w}\in\mathbb{U}_\Gamma$ is the quantum Weyl group element dual to $t$ (see Appendix~\ref{QW}), and in the last equality we used that $\Phi^-$ maps $\mathcal{Z}_0(\Oo_\e)$ into $\mathcal{Z}_0(U_\e)$ (see Theorem \ref{DCLteo1}\,(2)), which acts on specializations of $\Gamma$-modules by the trivial character (the counit) $\varepsilon \colon U_\e \ra \mc$. \mbox{By~\eqref{tform1}--\eqref{tform2}}, we have $t\big(a^{\star l}\big)=t\big(d^{\star l}\big)=0$ and $t\big(b^{\star l}\big)=1$, $t\big(c^{\star l}\big)=-1$. Now the computation of $t(fg)$ follows easily. For instance, taking $f=b^l=b^{\star l}$, by using $\Delta\big(b^{\star l}\big) = a^{\star l}\otimes b^{\star l}+b^{\star l}\otimes d^{\star l}$ and~\smash{$\Delta\big(d^{\star l}\big) = c^{\star l}\otimes b^{\star l}+d^{\star l}\otimes d^{\star l}$}, we get
\begin{align*}
 t(b^lg) ={}& \varepsilon\big( S^{-2}\big(\Phi^-\big(b^{\star l}\big)\big) \big) \varepsilon\big( S^{-1}\big(\Phi^-\big(c^{\star l}\big)\big) \big) t(g) + \varepsilon\big( S^{-2}\big(\Phi^-\big(d^{\star l}\big)\big) \big) \varepsilon\big( S^{-1}\big(\Phi^-\big(d^{\star l}\big)\big) \big) t(g).
 \end{align*}
Since $b^{\star l}\in \Oo_\e(U_+)$, $\Phi^-\big(b^{\star l}\big)=0$. Also, it is immediate from the definition of $\Phi^-$ that $\Phi^-\big(d^{\star l}\big) = \Phi^-(d)^l = L^{l}$; alternatively, one can bypass this computation by observing that $\Phi^-$ sets an isomorphism from $\Oo_\e(T_G) = \Oo_\e(B_+)\cap \Oo_\e(B_-)$ to $\mc\big[L^{\pm 1}\big] = U_\e(\mathfrak{b}_+)\cap U_\e(\mathfrak{b}_-)$, mapping a generator $d$ to $L$ or $L^{-1}$. We have $\varepsilon\big(L^l\big)=1$, and therefore
\begin{align*}
t\big(b^lg\big) = t(g) = t\big(b^l\big)t(g).
 \end{align*}
The other cases $f=T_l(\omega)-d^l, c^l,d^l$ are similar.
\end{proof}

\begin{teo} \label{modZ0} $\Ll_{0,n}^\e$ is a free $\mathcal{Z}_0\bigl(\Ll_{0,n}^\e\bigr)$-module of rank $l^{n.\dim \mathfrak{g}}$, and \smash{$\bigl(\Ll_{0,n}^\e\bigr)^{U_\e}$} is a Noetherian ring and a finite, whence Noetherian, $\mathcal{Z}_0\bigl(\Ll_{0,n}^\e\bigr)$-module.
\end{teo}

\begin{proof} We already proved the first claim in Proposition \ref{Z0L0n}, and that $\Ll_{0,n}^\e$ is a Noetherian $\mathcal{Z}_0\bigl(\Ll_{0,n}^\e\bigr)$-module. For the second claim, it follows that the $\mathcal{Z}_0\bigl(\Ll_{0,n}^\e\bigr)$-submodule \smash{$\bigl(\Ll_{0,n}^\e\bigr)^{U_\e}$} is necessarily finitely generated. But $\mathcal{Z}_0\bigl(\Ll_{0,n}^\e\bigr)$ being Noetherian, \smash{$\bigl(\Ll_{0,n}^\e\bigr)^{U_\e}$} is then a Noetherian $\mathcal{Z}_0\bigl(\Ll_{0,n}^\e\bigr)$-module and a Noetherian ring.

For the sake of clarity, let us provide a self-contained proof of the first claim, not invoking directly \cite{BGS,DC-L} or \cite{AnGa,BC}, but applying the same arguments directly to $\Ll_{0,n}^\e$. Since $\Ll_{0,n}^\e$ and $\Ll_{0,1}^{\otimes n}$ coincide as modules over $\mathcal{Z}_0\bigl(\Ll_{0,n}^\e\bigr)=\mathcal{Z}_0\bigl(\Ll_{0,1}^\e\bigr)^{\otimes n}$ by Proposition \ref{Z0L0n}, the result will follow from the case $n=1$. Then we argue in four steps. First, using Theorem \ref{JLteo1} we show that a certain localization of $\Ll_{0,1}^\e$ is a free module of rank \smash{$l^{\dim \mathfrak{g}}$}. Then, assuming that $\Ll_{0,1}^\e$ is finitely generated and projective, we explain why it has constant rank \smash{$l^{\dim \mathfrak{g}}$} (this is very classical). Thirdly, we prove that $\Ll_{0,1}^\e$ is finitely generated and projective as in \cite[Theorem 7.2]{DC-L}. Finally, we obtain that it is a free module as in~\cite{BGS}.

Recall Proposition \ref{JLad2}: $U_\e$ is a free $\Phi_1\big(\Ll_{0,1}^\e\big[d^{-l}\big]\big)$-module of rank $2^m$ \big(note that $\Ll_{0,1}^\e\big[d^{-l}\big]=\Ll_{0,1}^\e\big[d^{-1}\big]$\big), $\mathcal{Z}_0(U_\e)$ is free over
\[
T_{2-}^{(l)-1}\mathcal{Z}_0\big(U_\e^{\rm lf}\big) = \Phi_1\big(\mathcal{Z}_0\bigl(\Ll_{0,1}^\e\bigr)\big[d^{-l}\big]\big)
\]
 of rank $2^m$. Since $U_\e$ is also free of rank \smash{$l^{\dim \mathfrak{g}}$} over $\mathcal{Z}_0(U_\e)$ (see Theorem \ref{DCKteo1}\,(1)), it is free over \smash{$\Phi_1\big(\mathcal{Z}_0\bigl(\Ll_{0,1}^\e\bigr)\big[d^{-l}\big]\big)$} of rank \smash{$2^ml^{\dim \mathfrak{g}}$}. The decomposition being unique, \smash{$\Phi_1\big(\Ll_{0,1}^\e\big[d^{-l}\big]\big)$} is free of rank \smash{$l^{\dim \mathfrak{g}}$} over \smash{$\Phi_1\big(\mathcal{Z}_0\bigl(\Ll_{0,1}^\e\bigr)\big[d^{-l}\big]\big)$}, and injectivity of $\Phi_1$ implies that $\Ll_{0,1}^\e\big[d^{-l}\big]$ is free of rank~\smash{$l^{\dim \mathfrak{g}}$} over \smash{$\mathcal{Z}_0\bigl(\Ll_{0,1}^\e\bigr)\big[d^{-l}\big]$}.

Assume now that $\Ll_{0,1}^\e$ is finitely generated and projective. Let us show that its rank is~\smash{$l^{\dim \mathfrak{g}}$}. The localization $\bigl(\Ll_{0,1}^\e\bigr)_P$ of $\Ll_{0,1}^\e$ at any prime ideal $P$ of $\mathcal{Z}_0\bigl(\Ll_{0,1}^\e\bigr)$ is a free module over $\mathcal{Z}_0\bigl(\Ll_{0,1}^\e\bigr)_P$ \cite[Proposition~2.12.15]{Rowen}; the ranks of such modules are finite in number \cite[Proposition~2.12.20]{Rowen}. If these ranks are all equal, then, by definition, it is the rank of $\Ll_{0,1}^\e$ over $\mathcal{Z}_0\bigl(\Ll_{0,1}^\e\bigr)$. This happens if $\mathcal{Z}_0\bigl(\Ll_{0,1}^\e\bigr)$ has no nontrivial (i.e., $\ne 1$) idempotent~\cite[Corollary~2.12.23]{Rowen}, which is the case since it has no nontrivial zero divisors. To compute the rank, suppose $P$ does not contain~\smash{$d^l=\psi_{-\rho}^l$}. Such ideals $P$ are in $1$-$1$ correspondence with the prime ideals of $\mathcal{Z}_0\bigl(\Ll_{0,1}^\e\bigr)\big[d^{-l}\big]$ by the natural ring monomorphism \smash{$\mathcal{Z}_0\bigl(\Ll_{0,1}^\e\bigr) \ra \mathcal{Z}_0\bigl(\Ll_{0,1}^\e\bigr)\big[d^{-l}\big]$}. The set~${S = \mathcal{Z}_0\bigl(\Ll_{0,1}^\e\bigr)\setminus P}$ is multiplicatively closed, and we have also a ring morphism \smash{$\mathcal{Z}_0\bigl(\Ll_{0,1}^\e\bigr)\big[d^{-l}\big] \ra S^{-1}\mathcal{Z}_0\bigl(\Ll_{0,1}^\e\bigr)$}, which is also an injection (there are no zero divisors in $\mathcal{Z}_0\bigl(\Ll_{0,1}^\e\bigr)$, whence in $S$). Then
\[
\bigl(\Ll_{0,1}^\e\bigr)_P = S^{-1}\Ll_{0,1}^\e = \Ll_{0,1}^\e\big[d^{-l}\big] \bigotimes_{\mathcal{Z}_0\bigl(\Ll_{0,1}^\e\bigr)[d^{-l}]} S^{-1}\mathcal{Z}_0\bigl(\Ll_{0,1}^\e\bigr)
\]
shows that $\bigl(\Ll_{0,1}^\e\bigr)_P$ has over $\mathcal{Z}_0\bigl(\Ll_{0,1}^\e\bigr)_P = S^{-1}\mathcal{Z}_0\bigl(\Ll_{0,1}^\e\bigr)$ the same rank $l^{\dim \mathfrak{g}}$ as $\Ll_{0,1}^\e\big[d^{-l}\big]$ over $\mathcal{Z}_0\bigl(\Ll_{0,1}^\e\bigr)\big[d^{-l}\big]$. This proves our claim.

In order to show that $\Ll_{0,1}^\e$ is finitely generated and projective over $\mathcal{Z}_0\bigl(\Ll_{0,1}^\e\bigr)$, it is enough to show it is finite locally free, i.e., there are elements \smash{$d_i\in \mathcal{Z}_0\bigl(\Ll_{0,1}^\e\bigr)$} such that the localizations~${\Ll_{0,1}^\e\big[d_i^{-1}\big]}$ are finite free $\mathcal{Z}_0\bigl(\Ll_{0,1}^\e\bigr)\big[d_i^{-1}\big]$-modules, and ${\rm Maxspec}\bigl(\mathcal{Z}_0\bigl(\Ll_{0,1}^\e\bigr)\bigl)$ is covered by the open sets $U(d_i)$ made of the ideals not containing $d_i$ (see \cite[Lemma~77.2]{StackP}).

We have seen above that $\Ll_{0,1}^\e\big[d^{-l}\big]$ is free of rank $l^{\dim \mathfrak{g}}$ over $\mathcal{Z}_0\bigl(\Ll_{0,1}^\e\bigr)\big[d^{-l}\big]$. By Remark~\ref{Z0viaUe}, \smash{$\mathcal{Z}_0\bigl(\Ll_{0,1}^\e\bigr)\big[d^{-l}\big]\cong \mathcal{Z}_0\big(U_\e^{\rm lf}\big)\big[\ell^l\big]$} is isomorphic to $\Oo\big(G^0\big)$, and $\Oo\big(G^0\big) = \Oo(G)\big[\delta^{-1}\big]$. Now, given $w\in W$ with a reduced expression $s_{i_1}\cdots s_{i_k}$, put $t_w = t_{i_1}\cdots t_{i_k}$. Let $w$ be represented by $n_w = n_{i_1}\cdots n_{i_k}$ in $N(T_G)$. By Proposition~\ref{functL01}, we have $(f\lhd t_w)(x) = f(n_wx)$ for every $f\in \mathcal{Z}_0\bigl(\Ll_{0,1}^\e\bigr)$, $x\in G$. Then
\begin{equation}\label{BruhatactionX}
\mathcal{Z}_0\bigl(\Ll_{0,1}^\e\bigr)\big[d^{-l}\big]\lhd t_w \cong \Oo\big(n_{w}^{-1}G^0\big) \cong \Oo(G)\big[(\delta\lhd t_w)^{-1}\big].
\end{equation}
If $b_1,\dots,b_r$ \big($r:= l^{\dim \mathfrak{g}}$\big) is a basis of $\Ll_{0,1}^\e\big[d^{-l}\big]$ over $\mathcal{Z}_0\bigl(\Ll_{0,1}^\e\bigr)\big[d^{-l}\big]$, then $\Ll_{0,1}^\e\big[d^{-l}\big]\lhd t_w$ is free over $\mathcal{Z}_0\bigl(\Ll_{0,1}^\e\bigr)\big[(d\lhd t_w)^{-l}\big]\cong \Oo\big(n_{w}^{-1}G^0\big)$ with basis $b_1\lhd t_w,\dots,b_r\lhd t_w$. Consider the Bruhat decomposition of $G$: any $g\in G$ can be written in the form $g=b_1n b_2$, where $b_1,b_2\in B_-$, $n\in W$. Hence $g = nn^{-1}b_1nb_2\in nB_+B_- = nG^0$, and therefore
\[G = \bigcup_{w\in W} (B_- n_w B_- )= \bigcup_{w\in W} \big(n_w G^0\big).\]
For every $w\in W$, put $d^l_w := d^{l}\lhd t_w$.
Under the isomorphism of $\mathcal{Z}_0\bigl(\Ll_{0,1}^\e\bigr)$ with $\Oo(G)$, we thus get that ${\rm Maxspec}\bigl(\mathcal{Z}_0\bigl(\Ll_{0,1}^\e\bigr)\bigr)$ is covered by the open sets \smash{$U\big(d^l_w\big) \cong n_w G^0$}, and $\Ll_{0,1}^\e\big[d_w^{-l}\big]$ is finite free over $\mathcal{Z}_0\bigl(\Ll_{0,1}^\e\bigr)\big[d_w^{-l}\big]$. Therefore, $\Ll_{0,1}^\e$ is finitely generated and projective over $\mathcal{Z}_0\bigl(\Ll_{0,1}^\e\bigr)$.

Finally, let us explain why $\Ll_{0,1}^\e$ is free over $\mathcal{Z}_0\bigl(\Ll_{0,1}^\e\bigr)$, following the arguments of~\cite{BGS}. Let~$R$ be a commutative Noetherian ring, put $X={\rm Maxspec}(R)$, and let $P$ be an $R$-module. Denote by~$R_I$,~$P_I$ the localizations of $R$, $P$ at a maximal ideal $I\in X$. Define the {\it f-rank} of~$P$ as f-rank$ (P) = \inf_{I\in X} \{$ f-rank${}_{R_I}(P_I)\}$, where f-rank${{}_{R_I}(P_I) = \sup\bigl\{ r\in \mathbb{N}, R_I^{\otimes r}\subset P_I\bigr\}\in \mn\cup\{+\infty\}}$ (i.e., the maximal dimension of a free summand of $P_I$). Bass' Cancellation theorem asserts that if~$P$ is projective and f-rank$(P)> {\rm dim}(X)$, and $P\oplus Q \cong M\oplus Q$ for some $R$-modules~$Q$ and~$M$ such that~$Q$ is finitely generated and projective, then $P\cong M$ (see \cite[Section~IV.3.5, p.~167 and~p.~170]{Bass}, taking~${A=R}$, or \cite[Section~11.7.13]{MC-R}). Let us apply this to ${R=\Oo(G)}$ and $P=\Ll_{0,1}^\e$. We proved above that f-rank\smash{${}_{R_I}(P_I) = l^{\dim \mathfrak{g}}$}, a constant, and we have \smash{${l^{\dim \mathfrak{g}}> \dim \mathfrak{g} = \dim (G)}$}. By a~result of Marlin~\cite{Marlin}, $G$ being semisimple and simply connected the Grothendieck ring $K_0(\Oo(G)\!)$ is isomorphic to $\mz$. Therefore, $\Ll_{0,1}^\e\oplus Q \cong \Oo(G)^r$ for some free $\Oo(G)$-module $Q$ and~${r\in \mathbb{N}}$. Then Bass' cancellation implies $\Ll_{0,1}^\e$ is free over $\mathcal{Z}_0(\Ll_{0,1})\cong \Oo(G)$.
\end{proof}

\section{Proof of Theorem \ref{degree}}\label{DEGREEsec}
We begin with the following lemma, interesting by itself.
\begin{lem}\label{ZfiniteZ0} $\mathcal{Z}\bigl(\Ll_{0,n}^\e\bigr)$ is a finite $\mathcal{Z}_0\bigl(\Ll_{0,n}^\e\bigr)$-module and a Noetherian ring. Therefore, the ring $\mathcal{Z}\bigl(\Ll_{0,n}^\e\bigr)$ is integral over $\mathcal{Z}_0\bigl(\Ll_{0,n}^\e\bigr)$.
\end{lem}
\begin{proof} We know by Proposition \ref{Z0L0n} that $\mathcal{Z}_0\bigl(\Ll_{0,n}^\e\bigr)$ is a Noetherian ring, and $\Ll_{0,n}^\e$ is a finite Noetherian $\mathcal{Z}_0\bigl(\Ll_{0,n}^\e\bigr)$-module. Therefore, the submodule $\mathcal{Z}\bigl(\Ll_{0,n}^\e\bigr)$ is finitely generated. Being finite over $\mathcal{Z}_0\bigl(\Ll_{0,n}^\e\bigr)$, it is necessarily a Noetherian ring (e.g., by \cite[Proposition~7.2]{AM}).

Let $x\in \mathcal{Z}\bigl(\Ll_{0,n}^\e\bigr)$. The $\mathcal{Z}_0\bigl(\Ll_{0,n}^\e\bigr)$-submodule $\mathcal{Z}_0\bigl(\Ll_{0,n}^\e\bigr)[x]$ of $\Ll_{0,n}^\e$ is finitely generated by the same argument. Using the fact that an element $x$ is integral over $\mathcal{Z}_0\bigl(\Ll_{0,n}^\e\bigr)$ if and only if~${\mathcal{Z}_0\bigl(\Ll_{0,n}^\e\bigr)[x]}$ is a finitely generated $\mathcal{Z}_0\bigl(\Ll_{0,n}^\e\bigr)$-module (e.g., by~\cite[Proposition 5.1]{AM}), this proves the last claim.
\end{proof}

We will use the following notations. Let $A$ be a ring with no nontrivial zero divisors. The center $Z=Z(A)$ is a commutative integral domain. We denote by $Q(Z)$ its field of fractions, and put
\[Q(A):=Q(Z)\bigotimes_Z A.\]
It is an algebra over its center $Q(Z)$. Since $\Ll_{0,n}^\e$ has no nontrivial zero divisors \cite[Pro\-position~6.30]{BR}, we can take $A=\Ll_{0,n}^\e$, or \smash{$A=\bigl(\Ll_{0,n}^\e\bigr)^{U_\e}$}.

By the lemma, $\mathcal{Z}\bigl(\Ll_{0,n}^\e\bigr)$ is finite over $\mathcal{Z}_0\bigl(\Ll_{0,n}^\e\bigr)$, so the ring \smash{$\mathcal{Z}\bigl(\Ll_{0,n}^\e\bigr) \bigotimes_{\mathcal{Z}_0(\Ll_{0,n}^\e)} Q\bigl(\mathcal{Z}_0\bigl(\Ll_{0,n}^\e\bigr)\bigr)$} is a~field. Necessarily it coincides with $Q(\mathcal{Z}\bigl(\Ll_{0,n}^\e\bigr))$, and therefore
\begin{equation}\label{QZ0}
Q\bigl(\Ll_{0,n}^\e\bigr) = Q\bigl(\mathcal{Z}\bigl(\Ll_{0,n}^\e\bigr)\bigr) \bigotimes_{\mathcal{Z}(\Ll_{0,n}^\e)} \Ll_{0,n}^\e = Q\bigl(\mathcal{Z}_0\bigl(\Ll_{0,n}^\e\bigr)\bigr) \bigotimes_{\mathcal{Z}_0(\Ll_{0,n}^\e)} \Ll_{0,n}^\e.
\end{equation}
Recall that we denote by $N$ the number of positive roots of $\mathfrak{g}$.
\begin{teo} \label{LisCSA} $Q\bigl(\Ll_{0,n}^\e\bigr)$ is a division algebra and a central simple algebra of PI degree~$l^{Nn}$.
\end{teo}

\begin{proof} It follows from \eqref{QZ0} and Theorem \ref{modZ0} that $Q\bigl(\Ll_{0,n}^\e\bigr)$ is a vector space of dimension~\smash{$l^{n.\dim\mathfrak{g}}$} over $Q\bigl(\mathcal{Z}_0\bigl(\Ll_{0,n}^\e\bigr)\bigr)$, and therefore has finite dimension over its center $Q(\mathcal{Z}\bigl(\Ll_{0,n}^\e\bigr))$. Because $\Ll_{0,n}^\e$ has no nontrivial divisors \cite[Proposition 6.30]{BR} and $Q\bigl(\Ll_{0,n}^\e\bigr)$ is finite-dimensional over $Q\bigl(\mathcal{Z}\bigl(\Ll_{0,n}^\e\bigr)\bigr)$, $Q\bigl(\Ll_{0,n}^\e\bigr)$ is a division algebra, whence a central simple algebra. By classical theory (see, e.g., \cite[Section~13.3.5]{MC-R}, or \cite[Corollary~2.3.25]{Rowen}), there is a finite extension $\mathbb{F}$ of~$Q\bigl(\mathcal{Z}\bigl(\Ll_{0,n}^\e\bigr)\bigr)$, a {\it splitting field}, such that
\[\mathbb{F} \bigotimes_{Q(\mathcal{Z}(\Ll_{0,n}^\e))} Q\bigl(\Ll_{0,n}^\e\bigr) = M_{d}(\mathbb{F}),\]
where $d\in \mathbb{N}$, the PI degree of $Q\bigl(\Ll_{0,n}^\e\bigr)$, satisfies
\begin{equation}\label{factdegree}
d^2 = \bigl[Q\bigl(\Ll_{0,n}^\e\bigr):Q\bigl(\mathcal{Z}\bigl(\Ll_{0,n}^\e\bigr)\bigr)\bigr] = \frac{[Q\bigl(\Ll_{0,n}^\e\bigr):Q\bigl(\mathcal{Z}_0\bigl(\Ll_{0,n}^\e\bigr)\bigr)]}{\bigl[Q(\mathcal{Z}\bigl(\Ll_{0,n}^\e\bigr)): Q\bigl(\mathcal{Z}_0\bigl(\Ll_{0,n}^\e\bigr)\bigr)\bigr]}.
\end{equation}
We have to show $d^2=l^{2nN}$. We will obtain this equality by proving firstly that $d^2\geq l^{2nN}$, and then $d^2\leq l^{2nN}$.

In order to show that $d^2\geq l^{2nN}$, it is enough to exhibit an irreducible representation $V$ of~$\Ll_{0,n}^\e$ of dimension \smash{$k:=l^{nN}$}. Indeed, the representation map $\rho_V\colon \Ll_{0,n}^\e\ra {\rm End}_{\mc}(V)$ being surjective, given basis elements $v_1,\dots, v_{k^2}\in {\rm End}(V)$, and elements $\alpha_1,\dots,\alpha_{k^2}\in \Ll_{0,n}^\e$ such that~${\rho(\alpha_i) = v_i}$ for every $i\in \bigl\{1,\dots,k^2\bigr\}$, necessarily $\alpha_1,\dots,\alpha_{k^2}$ form a free family of $Q\bigl(\Ll_{0,n}^\e\bigr)$. For, if there was a nontrivial relation $ \sum_i z_i\alpha_i=0$, with $z_i\in Q\bigl(\mathcal{Z}_0\bigl(\Ll_{0,n}^\e\bigr)\bigr)$, by clearing denominators and then applying the representation map $\rho_V$, we would get a nontrivial relation in~${\rm End}_{\mc}(V)$ between~${v_1,\dots, v_{k^2}}$.

Now, by Theorem \ref{DCKteo1}\,(1) (see \cite[Section 20]{DC-K-P2}), the dimension of a generic irreducible representation space of $U_\e$ is $l^N$. Because $U_\e = T_{2-}^{-1}U_\e^{\rm lf}[T/T_{2}]$ by Proposition \ref{JLad2}, an irreducible representation of $U_\e$ yields an irreducible representation of \smash{$U_\e^{\rm lf}$}. Moreover, the tensor product of $n$ irreducible representation spaces of \smash{$U_\e^{\rm lf}$} of dimension \smash{$l^N$} is an irreducible representation space of \smash{$\big(U_\e^{\rm lf}\big)^{\otimes n}$} of dimension \smash{$l^{nN}$} (see, e.g., \cite[Theorem 3.10.2]{Et}). Applying the linear isomorphism \smash{$\psi_n = \Phi_n \circ \big(\Phi_1^{-1}\big)^{\otimes n}$} in \eqref{Phipsi} thus provides an irreducible representation of $\Ll_{0,n}^\e$ of dimension~\smash{$l^{nN}$}.

It remains to show $d^2\leq l^{2nN}$, which by \smash{$\bigl[Q\bigl(\Ll_{0,n}^\e\bigr):Q\bigl(\mathcal{Z}_0\bigl(\Ll_{0,n}^\e\bigr)\bigr)\bigr]= l^{n(2N+m)}$} is equivalent to~${\bigl[Q(\mathcal{Z}\bigl(\Ll_{0,n}^\e\bigr)):Q\bigl(\mathcal{Z}_0\bigl(\Ll_{0,n}^\e\bigr)\bigr)\bigr] \geq l^{mn}}$. For this, it is enough to exhibit an extension of $Q\bigl(\mathcal{Z}_0\bigl(\Ll_{0,n}^\e\bigr)\bigr)$ contained in $Q\bigl(\mathcal{Z}\bigl(\Ll_{0,n}^\e\bigr)\bigr)$ and of degree $l^{mn}$. There is a very natural one, which we denote by~\smash{$Q\big(\hat{\mathcal{Z}}_{0}\bigl(\Ll_{0,n}^\e\bigr)\big)$} and is constructed as follows. Consider for every $\lambda\in P_+$ the matrices
\[
M_\lambda:=\bigl({}_{{}_AV_{\lambda}}\phi_{e_k}^{e_l}\bigr)_{k,l}\in {\rm End}({}_AV_{\lambda})\otimes\Ll_{0,n}^A,\qquad M_\lambda^{(i)}:=\big(\bigl({}_{{}_AV_{\lambda}}\phi_{e_k}^{e_l}\bigr)^{(i)}\big)_{k,l}\in {\rm End}({}_AV_{\lambda})\otimes\Ll_{0,n}^A,
\]
where $i=1,\dots,n$, and as usual ${}_{{}_AV_{\lambda}}\phi_{e_k}^{e_l}$ is a matrix coefficient of ${}_AV_{\lambda}$, $\{e_k\}$ the canonical basis of ${}_AV_{\lambda}$, and \smash{$\bigl({}_{V_{\lambda}}\phi_{e_k}^{e_l}\bigr)^{(i)} := 1^{\otimes (i-1)}\otimes {}_{V_{\lambda}}\phi_{e_k}^{e_l}\otimes 1^{\otimes (n-i)}$}. Set
\[
{}_{\lambda}\omega := Tr(\pi_{V_{\lambda}}(\ell)M_\lambda),\qquad {}_{\lambda}\omega^{(i)} := Tr\big(\pi_{V_{\lambda}}(\ell)M_\lambda^{(i)}\big),
\]
where $Tr$ is the standard trace on ${\rm End}(V_{\lambda})$. Clearly, ${}_{\lambda}\omega\in \Ll_{0,1}^A$, ${}_{\lambda}\omega^{(i)}\in \Ll_{0,n}^A$. By \cite[Propositions~4.8 and~6.24]{BR}, the family of elements \smash{$ \prod_{i=1}^n {}_{\lambda_i}\omega^{(i)}$}, where $\lambda_1,\dots, \lambda_n\in P_+$, is a basis of~$\mathcal{Z}(\Ll_{0,n})$; moreover the Alekseev map $\Phi_n$ affords an isomorphism from $\mathcal{Z}(\Ll_{0,n})$ to $\mathcal{Z}(U_q)^{\otimes n}$, and $\Phi_n\big({}_{\lambda}\omega^{(i)}\big) = (\Phi_1({}_{\lambda}\omega))^{(i)}$. For $n=1$, specializing $q$ to $\e$ it follows
\begin{equation}\label{Z1omegalambda}
\mathcal{Z}_1(U_\e) = {\rm Vect} \lbrace \Phi_1({}_{\lambda}\omega), \, \lambda\in P_+ \rbrace ,
\end{equation}
where $\mathcal{Z}_1(U_\e)$ is defined before Theorem \ref{DCKteo1}. Then, for every $i=1,\dots,n$ define
\[\mathcal{Z}_{0,(i)}\bigl(\Ll_{0,n}^\e\bigr) := \mathcal{Z}_0\bigl(\Ll_{0,n}^\e\bigr)\big[\big\{{}_{\lambda}\omega^{(i)},\lambda\in P_+\big\}\big]
\]
and let $\hat{\mathcal{Z}}_{0}\bigl(\Ll_{0,n}^\e\bigr)\subset \mathcal{Z}\bigl(\Ll_{0,n}^\e\bigr)$ be the algebra generated by $\mathcal{Z}_{0,(1)}\bigl(\Ll_{0,n}^\e\bigr),\dots,\mathcal{Z}_{0,(n)}\bigl(\Ll_{0,n}^\e\bigr)$. The fields~${Q\bigl(\mathcal{Z}_{0,(i)}\bigl(\Ll_{0,n}^\e\bigr)\bigr)}$ are $n$ linearly disjoint extensions of $Q\bigl(\mathcal{Z}_0\bigl(\Ll_{0,n}^\e\bigr)\bigr)$, so
\[\big[Q\big(\hat{\mathcal{Z}}_{0}\bigl(\Ll_{0,n}^\e\bigr)\big):Q\bigl(\mathcal{Z}_0\bigl(\Ll_{0,n}^\e\bigr)\bigr)\big] = \prod_{i=1}^n \bigl[Q\bigl(\mathcal{Z}_{0,(i)}\bigl(\Ll_{0,n}^\e\bigr)\bigr):Q\bigl(\mathcal{Z}_0\bigl(\Ll_{0,n}^\e\bigr)\bigr)\bigr].\]
Now, by Proposition \ref{JLad2}, we know that $\Phi_1$ affords isomorphisms $Q\bigl(\mathcal{Z}_0\bigl(\Ll_{0,1}^\e\bigr)\bigr)\cong Q\big(\mathcal{Z}_0\big(U_\e^{\rm lf}\big)\big)$ and~\smash{$Q(\mathcal{Z}\bigl(\Ll_{0,1}^\e\bigr))\cong Q\big(\mathcal{Z}\big(U_\e^{\rm lf}\big)\big)$}, and moreover
\begin{equation}\label{extQUverslf}
Q(\mathcal{Z}_0(U_\e)) = Q\big(\mathcal{Z}_0\big(U_\e^{\rm lf}\big)\big)\big(T^{(l)}/T_{2}^{(l)}\big), \qquad Q(\mathcal{Z}(U_\e)) = Q\big(\mathcal{Z}\big(U_\e^{\rm lf}\big)\big)\big(T^{(l)}/T_{2}^{(l)}\big).
\end{equation}
Computing via the field embedding $\Phi_1^{\otimes n}\colon Q\big(\hat{\mathcal{Z}}_{0}\bigl(\Ll_{0,n}^\e\bigr)\big) \ra Q\bigl(\mathcal{Z}\bigl(U_\e^{\otimes n}\bigr)\bigr)$, we deduce
\begin{gather*}
\bigl[Q\bigl(\mathcal{Z}_{0,(i)}\bigl(\Ll_{0,n}^\e\bigr)\bigr) :Q\bigl(\mathcal{Z}_0\bigl(\Ll_{0,n}^\e\bigr)\bigr)\bigr]\\
\qquad
 = \bigl[\Phi_1^{\otimes n}\bigl(Q\bigl(\mathcal{Z}_{0,(i)}\bigl(\Ll_{0,n}^\e\bigr)\bigr)\bigr):\Phi_1^{\otimes n}\bigl(Q\bigl(\mathcal{Z}_0\bigl(\Ll_{0,n}^\e\bigr)\bigr)\bigr)\bigr]\\
\qquad = \big[Q\big(\mathcal{Z}_0\big(U_\e^{\rm lf}\big)^{\otimes n}\big)\big[\big\{(\Phi_1({}_{\lambda}\omega))^{(i)},\,\lambda\in P_+, \, i=1,\dots,n\big\}\big]:Q\big(\mathcal{Z}_0\big(U_\e^{\rm lf}\big)^{\otimes n}\big)\big]\\
\qquad = \big[Q\big(\mathcal{Z}_0(U_\e)^{\otimes n}\big)\big[\big\{(\Phi_1({}_{\lambda}\omega))^{(i)},\,\lambda\in P_+, \, i=1,\dots,n\big\}\big]:Q\big(\mathcal{Z}_0(U_\e)^{\otimes n}\big)\big]=l^m.
\end{gather*}
The second and third equalities follow from \eqref{extQUverslf} and the properties of $\Phi_1$ recalled before it, and the last equality follows from Theorem \ref{DCLteo1}\,(2) and \eqref{Z1omegalambda}. As a result, we have
\[
\big[Q\big(\hat{\mathcal{Z}}_{0}\bigl(\Ll_{0,n}^\e\bigr)\big):Q\bigl(\mathcal{Z}_0\bigl(\Ll_{0,n}^\e\bigr)\bigr)\big]=l^{mn},
\]
 whence \[
 \bigl[Q\bigl(\mathcal{Z}\bigl(\Ll_{0,n}^\e\bigr)\bigr):Q\bigl(\mathcal{Z}_0\bigl(\Ll_{0,n}^\e\bigr)\bigr)\bigr]\geq l^{mn}.
 \]
Since $\bigl[Q\bigl(\Ll_{0,n}^\e\bigr):Q\bigl(\mathcal{Z}_0\bigl(\Ll_{0,n}^\e\bigr)\bigr)\bigr] = l^{n(m+2N)}$, by~\eqref{factdegree} we obtain \smash{$d^2 \leq l^{2nN}$}, which concludes the proof.
\end{proof}
\begin{Remark} It follows $\bigl[Q\bigl(\mathcal{Z}\bigl(\Ll_{0,n}^\e\bigr)\bigr):Q\bigl(\mathcal{Z}_0\bigl(\Ll_{0,n}^\e\bigr)\bigr)\bigr]=l^{mn}$ by the degree computation above, whence \smash{$Q\bigl(\mathcal{Z}\bigl(\Ll_{0,n}^\e\bigr)\bigr) = Q\big(\hat{\mathcal{Z}}_{0}\bigl(\Ll_{0,n}^\e\bigr)\big)$}. In~\cite{BFR2}, we prove that \smash{$\mathcal{Z}\bigl(\Ll_{0,n}^\e\bigr) = \hat{\mathcal{Z}}_{0}\bigl(\Ll_{0,n}^\e\bigr)$}.
\end{Remark}
\begin{teo}\label{central simpleQ} $Q\big(\bigl(\Ll_{0,n}^\e\bigr)^{U_\e}\big)$, $n\geq 2$, is a division algebra and a central simple algebra of PI degree \smash{$l^{N(n-1)-m}$}.
\end{teo}

\begin{proof} The center of \smash{$\bigl(\Ll_{0,n}^\e\bigr)^{U_\e}$} contains $\mathcal{Z}\bigl(\Ll_{0,n}^\e\bigr)$, so the finite-dimensionality of $Q\bigl(\Ll_{0,n}^\e\bigr)$ over $Q\bigl(\mathcal{Z}\bigl(\Ll_{0,n}^\e\bigr)\bigr)$ implies the finite-dimensionality of \smash{$Q\big(\bigl(\Ll_{0,n}^\e\bigr)^{U_\e}\big)$} over its center. Since it has no non-zero divisors, this proves \smash{$Q\big(\bigl(\Ll_{0,n}^\e\bigr)^{U_\e}\big)$} is a division algebra.

Now denote by $\Delta^{(n)}\colon \Oo_\e\ra \Oo_\e^{\otimes n}$, $n\geq 2$, the $n$-fold coproduct, i.e., $\Delta^{(2)} := \Delta$, the standard coproduct of $\Oo_\e$, and \smash{$\Delta^{(n)} := ({\rm id}\otimes \Delta)\circ \Delta^{(n-1)}$} for $n\geq 3$. Identifying $\Ll_{0,n}^\e$ with $\Oo_\e^{\otimes n}$ as a~vector space, we consider \smash{$\Delta^{(n)}$} as a map \smash{$\Delta^{(n)}\colon \Ll_{0,1}^\e\ra \Ll_{0,n}^\e$}. It is an algebra morphism \cite[Proposition~6.18]{BR}, injective because \smash{$\big(\varepsilon^{\otimes (n-1)} \otimes {\rm id}\big)\Delta^{(n)} = {\rm id}$}. Then it extends uniquely to the fraction algebra $Q\bigl(\Ll_{0,1}^\e\bigr)$. As noted above, \smash{$Q\bigl(\Ll_{0,1}^\e\bigr) = Q\bigl(\mathcal{Z}_0\bigl(\Ll_{0,1}^\e\bigr)\bigr) \!\bigotimes_{\mathcal{Z}_0(\Ll_{0,1}^\e)}\! \Ll_{0,1}^\e$}. Since~${\mathcal{Z}_0\bigl(\Ll_{0,1}^\e\bigr)=\mathcal{Z}_0(\Oo_\e)}$ is a~Hopf subalgebra of $\Oo_\e$ \cite[Proposition~6.4]{DC-L}, \smash{$\Delta^{(n)}$} maps $\mathcal{Z}_0\bigl(\Ll_{0,1}^\e\bigr)$ to $\mathcal{Z}_0\bigl(\Ll_{0,1}^\e\bigr)^{\otimes n}$. Then, extending the scalars of \smash{$\Delta^{(n)}\bigl(Q\bigl(\Ll_{0,1}^\e\bigr)\bigr)$} by the field $Q\bigl(\mathcal{Z}\bigl(\Ll_{0,n}^\e\bigr)\bigr)$, consider the algebra
\begin{align*}
 Q_{\mathcal{Z}}\big(\Delta^{(n)}\bigl(\Ll_{0,1}^\e\bigr)\big) :={}& Q\bigl(\mathcal{Z}\bigl(\Ll_{0,n}^\e\bigr)\bigr) \bigotimes_{\Delta^{(n)}(\mathcal{Z}_0 (\Ll_{0,1}^\e ))} \Delta^{(n)}\bigl(\Ll_{0,1}^\e\bigr)\\
 ={}& Q\bigl(\mathcal{Z}\bigl(\Ll_{0,n}^\e\bigr)\bigr) \bigotimes_{\Delta^{(n)}(Q(\mathcal{Z}_0 (\Ll_{0,1}^\e )))} \Delta^{(n)}\bigl(Q\bigl(\Ll_{0,1}^\e\bigr)\bigr) \\
 ={}& Q\bigl(\mathcal{Z}\bigl(\Ll_{0,n}^\e\bigr)\bigr) \bigotimes_{\Delta^{(n)}(Q(\mathcal{Z}_0 (\Ll_{0,1}^\e )))} \Delta^{(n)}\bigl(Q\bigl(\mathcal{Z}\bigl(\Ll_{0,1}^\e\bigr)\bigr)\bigr) \\ &\bigotimes_{\Delta^{(n)}(Q(\mathcal{Z} (\Ll_{0,1}^\e )))} \Delta^{(n)}\bigl(Q\bigl(\Ll_{0,1}^\e\bigr)\bigr).
 \end{align*}
By Proposition \ref{LisCSA}, $\Delta^{(n)}\bigl(Q\bigl(\Ll_{0,1}^\e\bigr)\bigr)$ is a $\Delta^{(n)}\bigl(Q\bigl(\mathcal{Z}\bigl(\Ll_{0,1}^\e\bigr)\bigr)\bigr)$-central simple algebra. The left factor is a field, so \smash{$Q_{\mathcal{Z}}\big(\Delta^{(n)}\bigl(\Ll_{0,1}^\e\bigr)\big)$} is a central simple algebra over it (see, e.g., \cite[Theorem~1.7.27]{Rowen}, or~\cite[Lemma~4.9]{StackP2}). Note that the left factor can also be written as
\[\tilde{Q}\bigl(\mathcal{Z}\bigl(\Ll_{0,n}^\e\bigr)\bigr):= Q\bigl(\mathcal{Z}\bigl(\Ll_{0,n}^\e\bigr)\bigr) \bigotimes_{\Delta^{(n)}(\mathcal{Z}_0(\Ll_{0,1}^\e))} \Delta^{(n)}\bigl(\mathcal{Z}\bigl(\Ll_{0,1}^\e\bigr)\bigr)\]
for it contains $\tilde{Q}\bigl(\mathcal{Z}\bigl(\Ll_{0,n}^\e\bigr)\bigr)$, it is contained in its fraction field, and $\tilde{Q}\bigl(\mathcal{Z}\bigl(\Ll_{0,n}^\e\bigr)\bigr)$ is a field because $\mathcal{Z}\bigl(\Ll_{0,1}^\e\bigr)$ is finite over $\mathcal{Z}_0\bigl(\Ll_{0,1}^\e\bigr)$ and has no nontrivial zero divisors. Note that
\[\big[\tilde{Q}(\mathcal{Z}\bigl(\Ll_{0,n}^\e\bigr)):Q\bigl(\mathcal{Z}\bigl(\Ll_{0,n}^\e\bigr)\bigr)\big] =l^m. \]
We proved in \cite[Proposition 6.19]{BR} that the ring \smash{$\big(\Ll_{0,n}^A\big)^{U_A}$} is the centralizer of $\Delta^{(n)}\big(\Ll_{0,1}^A\big)$ in~\smash{$\Ll_{0,n}^A$}; the same arguments show that \smash{$\bigl(\Ll_{0,n}^\e\bigr)^{U_\e}$} is the centralizer of \smash{$\Delta^{(n)}\bigl(\Ll_{0,1}^\e\bigr)$} in $\Ll_{0,n}^\e$. So the algebra
\[Q\big(\bigl(\Ll_{0,n}^\e\bigr)^{U_\e}\big) := Q\bigl(\mathcal{Z}\bigl(\Ll_{0,n}^\e\bigr)\bigr) \bigotimes_{\mathcal{Z}(\Ll_{0,n}^\e)} \bigl(\Ll_{0,n}^\e\bigr)^{U_\e}
\]
is the centralizer of $Q_{\mathcal{Z}}\big(\Delta^{(n)}\bigl(\Ll_{0,1}^\e\bigr)\big)$ in $Q\bigl(\Ll_{0,n}^\e\bigr)$. Since the latter is simple, we can apply the double centralizer theorem (see, e.g., \cite[Theorem 7.1.9]{Rowen}, or \cite[Theorem 7.1]{StackP2}): \smash{$Q\big(\bigl(\Ll_{0,n}^\e\bigr)^{U_\e}\big)$} is a simple algebra, we have
\begin{align*}
\big[Q\big(\bigl(\Ll_{0,n}^\e\bigr)^{U_\e}\big):Q\bigl(\mathcal{Z}\bigl(\Ll_{0,n}^\e\bigr)\bigr)\big] & = \frac{\bigl[Q\bigl(\Ll_{0,n}^\e\bigr):Q\bigl(\mathcal{Z}\bigl(\Ll_{0,n}^\e\bigr)\bigr)\bigr]} {\bigl[Q_{\mathcal{Z}}\big(\Delta^{(n)}\bigl(\Ll_{0,1}^\e\bigr)\big):Q\bigl(\mathcal{Z}\bigl(\Ll_{0,n}^\e\bigr)\bigr)\bigr]} = l^{2nN-(2N+m)},
\end{align*}
and the centralizer of $Q\big(\bigl(\Ll_{0,n}^\e\bigr)^{U_\e}\big)$ is $Q_{\mathcal{Z}}\big(\Delta^{(n)}\bigl(\Ll_{0,1}^\e\bigr)\big)$. In particular, $Q\big(\bigl(\Ll_{0,n}^\e\bigr)^{U_\e}\big)$ has center \smash{$Q\big(\bigl(\Ll_{0,n}^\e\bigr)^{U_\e}\big)\cap Q_{\mathcal{Z}}\big(\Delta^{(n)}\bigl(\Ll_{0,1}^\e\bigr)\big)$}, which is easily shown to be \smash{$\tilde{Q}\bigl(\mathcal{Z}\bigl(\Ll_{0,n}^\e\bigr)\bigr)$}. It then follows
\begin{align*}
\big[Q\big(\bigl(\Ll_{0,n}^\e\bigr)^{U_\e}\big):\tilde{Q}\bigl(\mathcal{Z}\bigl(\Ll_{0,n}^\e\bigr)\bigr)\big] & =\frac{\big[Q\big(\bigl(\Ll_{0,n}^\e\bigr)^{U_\e}\big):Q\bigl(\mathcal{Z}\bigl(\Ll_{0,n}^\e\bigr)\bigr)\big]}{\big[\tilde{Q}\bigl(\mathcal{Z}\bigl(\Ll_{0,n}^\e\bigr)\bigr) :Q\bigl(\mathcal{Z}\bigl(\Ll_{0,n}^\e\bigr)\bigr)\big]} \\ & = l^{2nN-(2N+m)} . l^{-m} = l^{2(N(n-1)-m)}.\end{align*}
Therefore, $Q\big(\bigl(\Ll_{0,n}^\e\bigr)^{U_\e}\big)$ is a central simple algebra of PI degree $l^{N(n-1)-m}$.
\end{proof}

\appendix
\section[Low and up crystal structures in the sl\_2 case]{Low and up crystal structures in the $\boldsymbol{{\mathfrak{sl}_2}}$ case}\label{lowupbases}
Let $k\in {\mathbb N}$, and denote by $V_k$ the simple $U_q^{\rm ad}({\mathfrak{sl}_2})$ module of dimension $k+1$. It has a basis $v_0,\dots, v_k$ such that
\begin{gather*}
K.v_j=q^{k-2j}v_j,\qquad
F.v_j=[j+1]_qv_{j+1}\quad \mathrm{if}\quad j<k,\qquad F.v_k=0,\\
 E.v_j=[k-j+1]_q v_{j-1}\quad \mathrm{if}\quad j>0,\qquad E.v_0=0.
\end{gather*}
This basis defines the full $A$-sublattice
${}_AV_k$, which is left invariant by $U_A^{\rm res}$, and we have
\[ F^{(a)}.v_j=\left[\begin{matrix} j+a \\ a \end{matrix}\right]_{q} v_{j+a},\qquad E^{(a)}.v_j=\left[\begin{matrix} k-j+a \\ a \end{matrix}\right]_{q} v_{j-a}, \qquad a\geq 0.
\]

The action of the Kashiwara operator $\tilde{e}$, $\tilde{f}$ on $V_k$ are given by $\tilde{f}(v_j)=v_{j+1}$, $\tilde{e}(v_j)=v_{j-1}$.

The crystal basis $\big({\mathcal L}^{\rm low}, {\mathcal B}^{\rm low}\big)$ at $q=0$ is formed by the $\Aa_0$-sublattice ${\mathcal L}^{\rm low}$ generated by~${v_0,\dots ,v_k}$, and ${\mathcal B}^{\rm low}$ by the images $\bar{v}_0,\dots ,\bar{v}_k $ of these vectors in ${\mathcal L}^{\rm low}/q {\mathcal L}^{\rm low}$.

The bilinear form $\langle\,\ \rangle_k$ defined by (\ref{proppairinglambda}) is easily computed
\[
\langle v_i, v_j\rangle_k =\big\langle F^{(i)}.v_0, F^{(j)}.v_0\big\rangle_k= \big\langle v_0, E^{(i)}F^{(j)}.v_0\big\rangle_k=
\left[\begin{matrix} k\\ i \end{matrix}\right]_{q}\delta_{i,j}.
\]

By definition,
\[
{}_AV_k^{\rm up}=\{v \in V_k,\, \langle v, {}_AV_k\rangle_k\subset A\}=\bigoplus_{j=0}^k A v_j^{\rm up},
\]
where
\[v_j^{\rm up}=\left[\begin{matrix} k\\ j \end{matrix}\right]_{q}^{-1} v_j.\]
The upper crystal basis $({\mathcal L}^{\rm up}, {\mathcal B}^{\rm up})$ at $q=0$ is formed by the $\Aa_0$-sublattice ${\mathcal L}^{\rm up}$ generated by~${v_0^{\rm up},\dots ,v_k^{\rm up}}$, and ${\mathcal B}^{\rm up}$ by the images $\bar{v}_0^{\rm up},\dots ,\bar{v}_k^{\rm up} $ of these vectors in ${\mathcal L}^{\rm up}/q {\mathcal L}^{\rm up}$.

Using that $[n]_q\in q^{1-n}(1+q\Aa_0)$, we obtain
\[\left[\begin{matrix} k\\ j \end{matrix}\right]_{q}\in q^{j^2-kj}(1+ q\Aa_0).\]
As a result, we get $\bar{v}_j^{\rm up}=q^{kj-j^2}\bar{v}_j$, which is exactly the relation (\ref{renormweightA}) relating the low and up crystal bases, with $\lambda=k\varpi_1$, $\mu=(k-2j)\varpi_1$.

\section{Quantum Weyl group}\label{QW}
We recall some of the formulas of \cite{BRT}. Let $e_q(z)$ be the formal power series in $z$ with coefficients in $\mathbb{C}(q)$ defined by
\[
e_q(z)=\sum_{n=0}^{+\infty}\frac{z^n}{(n)_{q^2}!}.
\]
We first consider the case of $\mathfrak{g}= {\mathfrak{sl}_2}$. As explained in \cite[Section 3]{BR}, the Cartan element~${H\in \mathfrak{g}}$ defines an element of ${\mathbb U}_q({\mathfrak{sl}_2})$. Viewed as elements of ${\mathbb U}_q(\mathfrak{sl}_2)$ we have $L\!=\!q^{H/2}$. The series~${\Theta=q^{H\otimes H/2}}$ defines an element of \smash{${\mathbb U}_q({\mathfrak{sl}_2})^{\hat{\otimes} 2}$}, its image under multiplication being \smash{$q^{H^2/2}$}. The $R$-matrix can be expressed as
$R=\Theta \hat{R}$ where $\hat{R}=e_{q^{-1}}\big(\big(q-q^{-1}\big) E\otimes F\big)$ is a well defined element of \smash{${\mathbb U}_q (\mathfrak{sl}_2)^{\hat{\otimes}2}$}.
Consider the Lusztig \cite{Lusztig2} braid group automorphism of $U_q({\mathfrak{sl}_2})$, defined~by%
\begin{equation}\label{Tdef}
T(L)=L^{-1},\qquad T(E)=- FK^{-1},\qquad T(F)=-KE.
\end{equation}
For every $x \in U_q({\mathfrak{sl}_2})$ it satisfies: \smash{$\Delta(T(x))=\hat{R}^{-1} (T\otimes T)(\Delta(x))\hat{R}$}.
Define the quantum Weyl group element $\hat{w}\in {\mathbb U}_q({\mathfrak{sl}_2})$ by Saito's formula \cite{Saito}:
\begin{equation}\label{Saitoform}
\hat{w}=e_{q^{-1}}(F)q^{-H^2/4}e_{q^{-1}}(-E)q^{-H^2/4}e_{q^{-1}}(F)q^{-H/2}.
\end{equation}
For every $x \in U_q({\mathfrak{sl}_2})$, it satisfies
\begin{align}
&T(x) =\hat{w}x\hat{w}^{-1},\qquad
\Delta(\hat{w}) =\hat{R}^{-1}(\hat{w}\otimes\hat{w}),\label{relTw}\\
 &\hat{w}^2 =q^{H^2/2}\xi \theta,\label{wnorm}
\end{align}
where $\theta\in {\mathbb U}_q({\mathfrak{sl}_2})$ is the ribbon element, and $\xi \in {\mathbb U}_q({\mathfrak{sl}_2})$ is the central group element whose value on the type $1$ simple module $V_k$ of $U_q^{\rm ad}({\mathfrak{sl}_2})$ of dimension $k+1$ is the scalar endomorphism~${(-1)^{k} id_{V_k}}$.

In order to compare our setting to the one of \cite{DC-L}, we need an explicit formula of $\hat{w}$.
Using the basis $v_j$ of $V_k$ of Appendix~\ref{lowupbases}, \eqref{Tdef},~\eqref{relTw} and \eqref{wnorm}, we obtain
\begin{equation}\label{wact}
\hat{w}{v}_j=(-1)^j q^{-j(k-j-1)-k}{v}_{k-j}.
\end{equation}
In \cite{DC-L}, another quantum Weyl group element $\underline{w}$ is defined. It is dual to the Vaksman--Soibelman functional $t\colon \Oo_q({\rm SL}_2)\ra \mc(q)$ of \cite{Soi,SV}, that is, $t(\alpha)=\langle \alpha , \underline{w} \rangle$ for all $\alpha \in \Oo_q({\rm SL}_2)$. By comparing \eqref{wact} with the formulas defining the action of~$t$ in \cite[Section~1.7]{DC-L}, we find $\underline{w}=\xi \hat{w}K$ and the basis vectors $w_r^p$ of~\cite{DC-L}, where $p\in (1/2)\mn$ and $r\in\{-p,-p+1,\dots,p-1,p\}$, are related to the vectors $v_j$ above as follows:
${v}_j=\lambda_j w_r^p$,
where $k=2p$, $j=p-r$, $\lambda_0=1$, $\lambda_1=[k]q^{-k}$, and
\[\lambda_j=\frac{[k]!}{[j]![k-(j-2)]!}q^{j(j+1)-j(k+2)}, \qquad j\geq 2.\]
Explicit formulas of the evaluation of $t$ on basis vectors of $\Oo_q({\rm SL}_2)$ can be computed. We get%
\begin{gather}
 t\big(\tilde{a}^{\star m}\star \tilde{b}^{\star n}\star \tilde{d}^{\star p}\big)=\delta_{m,p}q^{-np}\prod_{i=1}^{p}\big(1-q^{-2i}\big),\label{tform1}\\
 t\big(\tilde{a}^{\star m}\star \tilde{c}^{\star n}\star \tilde{d}^{\star p}\big)=(-1)^n\delta_{m,p}q^{-n(p+1)}\prod_{i=1}^{p}\big(1-q^{-2i}\big)\label{tform2},
\end{gather}
where $
\tilde{a} = a$, $\tilde{b} = q b$, $\tilde{c} = q^{-1}c$, $\tilde{d}=d$ and as usual $a$, $b$, $c$, $d$ are the standard generators of~$\mathcal{O}_q({\rm SL}_2)$, i.e., the matrix coefficients in the basis of weight vectors $v_0$, $v_1$ of the $2$-dimensional irreducible representation $V_1$ of $U_q({\mathfrak{sl}_2})$ such that $K.v_0=qv_0$ and $v_1=F.v_0$. Here we have introduced the generators $\tilde{a},\dots , \tilde{d}$ to facilitate the comparison with the formulas in \cite{DC-L}; these generators come naturally in their setup because they use different generators $E_i$ and $F_i$ of $U_q(\mathfrak{g})$, which in our notations can be written respectively as $K_i^{-1}E_i$ and $F_iK_i$.

The formulas \eqref{tform1}--\eqref{tform2} can be shown by two independent methods. The first uses a~definition of $t$ as a $GNS$ state associated to an infinite-dimensional representation of $\Oo_q({\rm SL}_2)$, as recalled in \cite[Section 1.6]{DC-L}. The second is to write, e.g.,
\begin{equation}\label{twformcop}
t\big(\tilde{a}^{\star m}\star \tilde{b}^{\star n}\star \tilde{d}^{\star p}\big) = \big\langle \tilde{a}^{\otimes m}\otimes \tilde{b}^{\otimes n}\otimes \tilde{d}^{\otimes p}, \Delta^{(m+n+p)}(\underline{w})\big\rangle
\end{equation}
and to use explicit expressions of $\Delta^{(m+n+p)}(\underline{w})$ when represented on $V_1^{\otimes (m+n+p)}$. In general, one can check that
\begin{align*}
\Delta^{(n)}(\hat{\omega}) ={}&\big( \Delta^{(n-1)} \otimes {\rm id}\big)\big(\widehat{R}^{-1}\big)\big(\big( \Delta^{(n-2)} \otimes {\rm id}\big)\big(\widehat{R}^{-1}\big) \otimes {\rm id}\big)\cdots \big(\big( \Delta \otimes {\rm id}\big)\big(\widehat{R}^{-1}\big) \otimes {\rm id}^{\otimes (n-3)}\big)\\ &\times\big(\widehat{R}^{-1} \otimes {\rm id}^{\otimes (n-2)}\big)\hat{\omega}^{\otimes n}.
\end{align*}
By \eqref{wact} or \eqref{tform1}--\eqref{tform2}, we see that $\hat{w}$ (or $\underline{w}$) and $t$ are well defined on the integral forms,
\[\hat{w} \in \mathbb{U}_{\Gamma},\qquad t \colon\ \Oo_A({\rm SL}_2)\ra A.\]
We now consider the case where $\mathfrak{g}$ is of rank $m\geq 2$. To each simple root $\alpha_i$, $1\leq i\leq m$, is associated the subalgebra of $U_q$ generated by $E_i$, $F_i$, $L_i$, $L_i^{-1}$. It is a copy of $U_{q_i}({\mathfrak{sl}_2})$, where~${q_i = q^{d_i}}$. Let $\hat{w}_i$ be the corresponding quantum Weyl group element in $\mUq = \mUq(\mathfrak{g})$, defined by Saito's formula \eqref{Saitoform}, replacing $H$, $E$, $F$ by $H_i$, $E_i$ and $F_i$. Also, denote by $\nu_i \colon \mathcal{O}_q \ra \mathcal{O}_{q_i}({\rm SL}_2)$ the projection map dual to the inclusion $U_{q_i}({\mathfrak{sl}_2}) \bigotimes_{\mc(q_i)} \mc(q)\hookrightarrow U_q$, and put $t_i = t \circ \nu_i$. Let $\underline{w}_i$ be the corresponding quantum Weyl group element in $\mUq$, i.e., $t_i(\alpha) = \langle \alpha,\underline{w}_i\rangle$ for all $\alpha\in \Oo_q$. On integral forms they yield well-defined elements $\hat{w}_i,\underline{w}_i \in \mathbb{U}_{\Gamma}$ and $t_i\colon \Oo_A \ra A$ (see \cite[Proposition~5.1]{DC-L}, and \cite{Lusztig3} for a different construction). They satisfy the defining relations of the braid group~$\mathcal{B}(\mathfrak{g})$ of $\mathfrak{g}$ \cite{KirRes}:
\begin{align*}
&\hat{w}_i\hat{w}_j\hat{w}_i = \hat{w}_j\hat{w}_i\hat{w}_j\qquad \text{if}\quad a_{ij}a_{ji}=1,\\
& (\hat{w}_i\hat{w}_j)^k = (\hat{w}_j\hat{w}_i)^k \qquad \text{for}\quad k=1,2,3 \quad \text{if}\quad a_{ij}a_{ji}=0,2,3,
\end{align*}
and similarly by replacing $\hat{w}_i$ with $\underline{w}_i$, or with $t_i$ (see \cite{Soi} for the latter). The Weyl group $W=W(\mathfrak{g}) = N(T_G)/T_G$ is generated by the reflections $s_i$ associated to the simple roots $\alpha_i$. Denote by $n_i\in N(T_G)$ a representative of $s_i$. Let $w\in W$ and denote by $w=s_{i_1}\dots s_{i_k}$ a reduced expression. Because of the braid group relations the elements $\hat{w}=\hat{w}_{i_1}\cdots \hat{w}_{i_k}$, $\underline{w}=\underline{w}_{i_1}\cdots \underline{w}_{i_k}$ and the functional $t_w=t_{i_1}\cdots t_{i_k}$
do not depend on the choice of reduced expression. The Lusztig~\cite{Lusztig2} braid group automorphism $T_w \colon \Gamma \ra \Gamma$ associated to $w$ satisfies (see~\cite{DC-L})
\[T_w(x)=\hat{w}x\hat{w}^{-1}, \qquad x \in \Gamma.\]
Let $w_0$ be the longest element in $W$. We have
\begin{equation}\label{Deltaw}
\Delta(\hat{w}_0)=\hat{R}^{-1}(\hat{w}_0\otimes\hat{w}_0),
\end{equation}
where as usual $R={\Theta}{\hat R}$.

\section[Regular action on O\_e]{Regular action on $\boldsymbol{\Oo_\e}$}\label{actO}

The following result is proved in \cite[Section 1.10]{DC-L}. For completeness, let us give a (different) proof. Recall from \eqref{Z0Odef} that we may identify $\mathcal{Z}_0(\Oo_\e)$ with $\Oo(G)$.

\begin{prop}\label{tidstar} For every $f\in \mathcal{Z}_0(\Oo_\e)$, $g\in {\mathcal O}_\epsilon$, we have
\begin{gather}
t_i(f)=f(n_i), \label{tonF0} \\
t_i(f\star g)=t_i(f)t_i(g).\label{tfstarg}
\end{gather}
\end{prop}
\begin{proof} It is sufficient to prove the results for ${\rm SL}_2$ because $\nu_i\colon {\mathcal O}_\epsilon \rightarrow
{\mathcal O}_{\epsilon_i} ({\rm SL}_2)$ is a morphism of Hopf algebras and $\nu_i(\mathcal{Z}_0(\Oo_\e))\subset \mathcal{Z}_0(\Oo_{\e_i}({\rm SL}_2))$. In this case, \eqref{tonF0} can be proved by using~\mbox{\eqref{tform1}--\eqref{tform2}}, evaluating $t$ on basis elements of $\mathcal{Z}_0(\Oo_\e({\rm SL}_2))$ as is done in \cite[Lemma~1.5\,(a)]{DC-L}. Such a basis is formed by monomials like in \eqref{tform1}--\eqref{tform2}, with all exponents divisible by $l$; then for instance
\[t\big(\tilde{a}^{\star m l}\star \tilde{b}^{\star n l}\star \tilde{d}^{\star p l}\big)=\delta_{p,0}\delta_{m,0} = \underline{a}^m\underline{b}^n\underline{d}^p (n),
\]
where $\underline{a},\dots,\underline{d}$ are the generators of $\Oo(G) = \Oo_1(G)$ corresponding to $a,\dots,d$, and we take
\[n=\left(\begin{matrix} \hphantom{-}0 & 1 \\ -1 & 0\end{matrix}\right)\]
as representative of the reflection $s$ generating the Weyl group $W({\mathfrak{sl}_2})$. Here is an alternative proof of \eqref{tonF0}: \eqref{tfstarg} shows that $t$ is a homomorphism on $\mathcal{Z}_0(\Oo_\e({\rm SL}_2))$, so by proving \eqref{tfstarg} at first one is reduced to check \eqref{tonF0} on the generators $a^{\star l},\dots,d^{\star l}$, which is easy by means of \eqref{twformcop} and \eqref{Deltaw}.

We provide a proof of (\ref{tfstarg}) that we find more conceptual than the one in \cite[Lemma~1.5\,(b)]{DC-L} (which uses again \eqref{tform1}--\eqref{tform2}). As above, let us denote $\underline{w}=\xi \hat{w}K$. For any $f,g\in {\mathcal O}_\epsilon$, we have
\begin{align*}
 t(f\star g)& =(f\otimes g)(\Delta(\underline{w}))=(f\otimes g)\big(\widehat{R}^{-1} (\underline{w} \otimes \underline{w})\big)=\sum_{(\widehat{R}^{-1}) }f\big( \big( \widehat{R}^{-1}\big)_{(1) } \underline{w} \big) g\big( \big( \widehat{R}^{-1}\big)_{(2) } \underline{w}\big)\\
&=\sum_{(\widehat{R}^{-1}),(f) }f_{(1)}\big( \big( \widehat{R}^{-1}\big)_{(1) } \big) f_{(2)}(\underline{w} ) g\big( \big( \widehat{R}^{-1}\big)_{(2) } \underline{w}\big)=\sum_{(f) }f_{(2)}(\underline{w} ) g\big( (f_{(1)}\otimes {\rm id} )\big(\widehat{R}^{-1}\big) \underline{w}\big).
\end{align*}
Assume now $f\in \mathcal{Z}_0(\Oo_\e({\rm SL}_2))$. Since $\mathcal{Z}_0(\Oo_\e({\rm SL}_2))$ is a Hopf subalgebra of ${\mathcal O}_{\epsilon}({\rm SL}_2)$, we have \linebreak ${f_{(1)}\in \mathcal{Z}_0(\Oo_\e({\rm SL}_2))}$. From Theorem \ref{DCLteo1}\,(2), we deduce
\[(f_{(1)}\otimes {\rm id} )\big(\widehat{R}^{-1}\big) \in U_\epsilon (\mathfrak{n}_-)\cap \mathcal{Z}_0\big(U_\epsilon^{\rm ad}\big).\]
Denote by $z$ this element. Note that from its expression we have $\epsilon(z)= \epsilon(f_{(1)})$. Now
 $g( z\underline{w})=\sum_{(g)}g_{(1)}(z)g_{(2)}(\underline{w})$, but $g_{(1)}$ is a linear combination of matrix elements of $\Gamma$-modules, on which ${\mathcal Z}_0\big(U_\epsilon^{\rm ad}\big)$ acts by the trivial character. Therefore,
\[ g( z\underline{w})=\sum_{(g)}\epsilon(z)g_{(1)}(1)g_{(2)}(\underline{w})=\epsilon(z)g(\underline{w})= \epsilon(f_{(1)})g(\underline{w}),\]
and eventually
 \[t(f\star g)=\sum_{(f) }f_{(2)}(\underline{w}) \epsilon(f_{(1)}) g(\underline{w})=t(f)t(g).\]
This concludes the proof.
\end{proof}

For the sake of completeness, let us show how this result implies:

\begin{proof}[Proof of Proposition \ref{funcDCL} (i.e., {\cite[Proposition 7.1]{DC-L}})] We have $ f\lhd t_i = \sum_{(f)} t_i(f_{(1)})f_{(2)}$, $f\in \mathcal{Z}_0(\Oo_\e)$. Since $\mathcal{Z}_0(\Oo_\e)$ is a Hopf subalgebra of ${\mathcal O}_{\epsilon}$, $f_{(2)} \in \mathcal{Z}_0(\Oo_\e)$ and therefore the maps $\lhd t_i\colon \Oo_\e \ra \Oo_\e$ preserve $\mathcal{Z}_0(\Oo_\e)$. Moreover, $ (f\lhd t_i)(a) = \sum_{(f)} f_{(1)}(n_i)f_{(2)}(a) = f(n_ia)$, $a\in G$, by \eqref{tonF0}.

It remains to show that $(f\star \alpha)\lhd t_i = (f\lhd t_i)(\alpha\lhd t_i)$ for every $f\in \mathcal{Z}_0(\Oo_\e)$, $\alpha\in \mathcal{O}_\e$. We have
\begin{align}
 (f\star g)\lhd t_i &= \sum_{(f\star g)} t_i\left((f\star g)_{(1)}\right)(f\star g)_{(2)}= \sum_{(f),(g)} t_i\left(f_{(1)}\star g_{(1)}\right)f_{(2)}\star g_{(2)}\nonumber\\ & = \sum_{(f),(g)} t\left(\nu_i(f_{(1)})\nu_i(g_{(1)})\right)f_{(2)}\star g_{(2)}\nonumber\\
& =\sum_{(f),(g)} t\left(\nu_i(f_{(1)})\right)t\left(\nu_i(g_{(1)})\right)f_{(2)}\star g_{(2)},\label{fundrelt}
\end{align}
using that $\nu_i$ is a homomorphism in the third equality, and \eqref{tfstarg} in the last one. The result is just $(f\lhd t_i)(g\lhd t_i)$.
\end{proof}

\subsection*{Acknowledgements}We are grateful to M.~Faitg for many valuable discussions on the subject, especially concerning the filtration arguments in the proof of Theorem~\ref{HN}, and the use of the partial order $\preceq$ in the proof of Theorem~\ref{HN}. We also thank K.A.~Brown for pointing out the references \cite{AnGa} and~\cite{BC} (see the comments before Theorem~\ref{Llibre}). We thank also the referees for their comments and suggestions, which greatly improved the clarity of our work.


\pdfbookmark[1]{References}{ref}
\LastPageEnding


\begin{thebibliography}{99}
\footnotesize\itemsep=0pt

\bibitem{A}
Alekseev A.Yu., Integrability in the {H}amiltonian {C}hern--{S}imons theory,
 \textit{St.~Petersburg Math.~J.} \textbf{6} (1995), 241--253,
 \href{https://arxiv.org/abs/hep-th/9311074}{arXiv:hep-th/9311074}.

\bibitem{AGS1}
Alekseev A.Yu., Grosse H., Schomerus V., Combinatorial quantization of the
 {H}amiltonian {C}hern--{S}imons theory.~{I}, \href{https://doi.org/10.1007/BF02099431}{\textit{Comm. Math. Phys.}}
 \textbf{172} (1995), 317--358, \href{https://arxiv.org/abs/hep-th/9403066}{arXiv:hep-th/9403066}.

\bibitem{AGS2}
Alekseev A.Yu., Grosse H., Schomerus V., Combinatorial quantization of the
 {H}amiltonian {C}hern--{S}imons theory.~{II}, \href{https://doi.org/10.1007/BF02101528}{\textit{Comm. Math. Phys.}}
 \textbf{174} (1996), 561--604, \href{https://arxiv.org/abs/hep-th/9408097}{arXiv:hep-th/9408097}.

\bibitem{AG3}
Alekseev A.Yu., Schomerus V., Representation theory of {C}hern--{S}imons
 observables, \href{https://doi.org/10.1215/S0012-7094-96-08519-1}{\textit{Duke Math.~J.}} \textbf{85} (1996), 447--510,
 \href{https://arxiv.org/abs/q-alg/9503016}{arXiv:q-alg/9503016}.

\bibitem{APW}
Andersen H.H., Polo P., Wen K.X., Representations of quantum algebras,
 \href{https://doi.org/10.1007/BF01245066}{\textit{Invent. Math.}} \textbf{104} (1991), 1--59.

\bibitem{AnGa}
Andruskiewitsch N., Garc\'{\i}a G.A., Quantum subgroups of a simple quantum
 group at roots of $1$, \href{https://doi.org/10.1112/S0010437X09003923}{\textit{Compos. Math.}} \textbf{145} (2009), 476--500,
 \href{https://arxiv.org/abs/0707.0070}{arXiv:0707.0070}.

\bibitem{AM}
Atiyah M.F., Macdonald I.G., Introduction to commutative algebra, Vol.~64,
 \href{https://doi.org/10.1201/9780429493638}{Addison-{W}esley Publishing Co.}, Reading, Mass., 1969.

\bibitem{B}
Baseilhac S., Quantum coadjoint action and the {$6j$}-symbols of {$U_q{\rm
 sl}_2$}, in Interactions between {H}yperbolic {G}eometry, {Q}uantum
 {T}opology and {N}umber {T}heory, \textit{Contemp. Math.}, Vol. 541, \href{https://doi.org/10.1090/conm/541/10681}{American
 Mathematical Society}, Providence, RI, 2011, 103--143, \href{https://arxiv.org/abs/1101.3440}{arXiv:1101.3440}.

\bibitem{BB-3}
Baseilhac S., Benedetti R., Quantum hyperbolic invariants of 3-manifolds with
 {${\rm PSL}(2,\mathbb C)$}-characters, \href{https://doi.org/10.1016/j.top.2004.02.001}{\textit{Topology}} \textbf{43} (2004),
 1373--1423, \href{https://arxiv.org/abs/math.GT/0306280}{arXiv:math.GT/0306280}.

\bibitem{BB-2}
Baseilhac S., Benedetti R., Classical and quantum dilogarithmic invariants of
 flat {${\rm PSL}(2,\mathbb C)$}-bundles over 3-manifolds, \href{https://doi.org/10.2140/gt.2005.9.493}{\textit{Geom.
 Topol.}} \textbf{9} (2005), 493--569, \href{https://arxiv.org/abs/math.GT/0306283}{arXiv:math.GT/0306283}.

\bibitem{BB-1}
Baseilhac S., Benedetti R., Quantum hyperbolic geometry, \href{https://doi.org/10.2140/agt.2007.7.845}{\textit{Algebr. Geom.
 Topol.}} \textbf{7} (2007), 845--917, \href{https://arxiv.org/abs/math.GT/0611504}{arXiv:math.GT/0611504}.

\bibitem{BB0}
Baseilhac S., Benedetti R., The {K}ashaev and quantum hyperbolic link
 invariants, \textit{J.~G\"okova Geom. Topol. GGT} \textbf{5} (2011), 31--85,
 \href{https://arxiv.org/abs/1101.1851}{arXiv:1101.1851}.

\bibitem{BB1}
Baseilhac S., Benedetti R., Analytic families of quantum hyperbolic invariants,
 \href{https://doi.org/10.2140/agt.2015.15.1983}{\textit{Algebr. Geom. Topol.}} \textbf{15} (2015), 1983--2063,
 \href{https://arxiv.org/abs/1212.4261}{arXiv:1212.4261}.

\bibitem{BB2}
Baseilhac S., Benedetti R., Non ambiguous structures on 3-manifolds and quantum
 symmetry defects, \href{https://doi.org/10.4171/QT/101}{\textit{Quantum Topol.}} \textbf{8} (2017), 749--846,
 \href{https://arxiv.org/abs/1506.01174}{arXiv:1506.01174}.

\bibitem{BB3}
Baseilhac S., Benedetti R., On the quantum {T}eichm\"uller invariants of fibred
 cusped 3-manifolds, \href{https://doi.org/10.1007/s10711-017-0315-0}{\textit{Geom. Dedicata}} \textbf{197} (2018), 1--32,
 \href{https://arxiv.org/abs/1704.05667}{arXiv:1704.05667}.

\bibitem{BFR}
Baseilhac S., Faitg M., Roche P., Unrestricted quantum moduli
 algebras,~{III}:~surfaces of arbitrary genus and skein algebras,
 \href{https://arxiv.org/abs/2302.00396}{arXiv:2302.00396}.

\bibitem{BFR2}
Baseilhac S., Faitg M., Roche P., Structure and representations of quantum
 moduli and ${\mathfrak g}$-skein algebras at roots of unity, {i}n
 preparation.

\bibitem{BR}
Baseilhac S., Roche P., Unrestricted quantum moduli algebras.~{I}.~{T}he case
 of punctured spheres, \href{https://doi.org/10.3842/SIGMA.2022.025}{\textit{SIGMA}} \textbf{18} (2022), 025, 78~pages,
 \href{https://arxiv.org/abs/1912.02440}{arXiv:1912.02440}.

\bibitem{Bass}
Bass H., Algebraic {$K$}-theory, \textit{Math. Lect. Note Ser.}, W.A.~Benjamin, Inc.,
 New York, 1968.

\bibitem{Bau1}
Baumann P., Another proof of {J}oseph and {L}etzter's separation of variables
 theorem for quantum groups, \href{https://doi.org/10.1007/BF01237175}{\textit{Transform. Groups}} \textbf{5} (2000),
 3--20.

\bibitem{BBG}
Beliakova A., Blanchet C., Geer N., Logarithmic {H}ennings invariants for
 restricted quantum {$\mathfrak{sl}(2)$}, \href{https://doi.org/10.2140/agt.2018.18.4329}{\textit{Algebr. Geom. Topol.}}
 \textbf{18} (2018), 4329--4358, \href{https://arxiv.org/abs/1705.03083}{arXiv:1705.03083}.

\bibitem{B-BZ-J}
Ben-Zvi D., Brochier A., Jordan D., Quantum character varieties and braided
 module categories, \href{https://doi.org/10.1007/s00029-018-0426-y}{\textit{Selecta Math. (N.S.)}} \textbf{24} (2018),
 4711--4748, \href{https://arxiv.org/abs/1606.04769}{arXiv:1606.04769}.

\bibitem{BW}
Bonahon F., Wong H., Quantum traces for representations of surface groups in
 {${\rm SL}_2(\mathbb C)$}, \href{https://doi.org/10.2140/gt.2011.15.1569}{\textit{Geom. Topol.}} \textbf{15} (2011),
 1569--1615, \href{https://arxiv.org/abs/1003.5250}{arXiv:1003.5250}.

\bibitem{BW1}
Bonahon F., Wong H., Representations of the {K}auffman bracket skein
 algebra~{I}: invariants and miraculous cancellations, \href{https://doi.org/10.1007/s00222-015-0611-y}{\textit{Invent. Math.}}
 \textbf{204} (2016), 195--243, \href{https://arxiv.org/abs/1206.1638}{arXiv:1206.1638}.

\bibitem{BC}
Brown K.A., Couto M., Affine commutative-by-finite {H}opf algebras,
 \href{https://doi.org/10.1016/j.jalgebra.2020.12.039}{\textit{J.~Algebra}} \textbf{573} (2021), 56--94, \href{https://arxiv.org/abs/1907.10527}{arXiv:1907.10527}.

\bibitem{BG}
Brown K.A., Goodearl K.R., Lectures on algebraic quantum groups, \textit{Adv. Courses Math. CRM Barcelona}, \href{https://doi.org/10.1007/978-3-0348-8205-7}{Birkh\"auser}, Basel, 2002.

\bibitem{BGOe}
Brown K.A., Gordon I., The ramifications of the centres: quantised function
 algebras at roots of unity, \href{https://doi.org/10.1112/plms/84.1.147}{\textit{Proc. London Math. Soc.}} \textbf{84}
 (2002), 147--178, \href{https://arxiv.org/abs/math.RT/9912042}{arXiv:math.RT/9912042}.

\bibitem{BGS}
Brown K.A., Gordon I., Stafford J.T., {$\mathcal{O}_\varepsilon[G]$} is a free
 module over {$\mathcal{O}[G]$}, \href{https://arxiv.org/abs/math.QA/0007179}{arXiv:math.QA/0007179}.

\bibitem{BR1}
Buffenoir E., Roche P., Two-dimensional lattice gauge theory based on a quantum
 group, \href{https://doi.org/10.1007/BF02099153}{\textit{Comm. Math. Phys.}} \textbf{170} (1995), 669--698,
 \href{https://arxiv.org/abs/hep-th/940512}{arXiv:hep-th/9405126}.

\bibitem{BR2}
Buffenoir E., Roche P., Link invariants and combinatorial quantization of
 {H}amiltonian {C}hern--{S}imons theory, \href{https://doi.org/10.1007/BF02101008}{\textit{Comm. Math. Phys.}}
 \textbf{181} (1996), 331--365, \href{https://arxiv.org/abs/q-alg/9507001}{arXiv:q-alg/9507001}.

\bibitem{BRT}
Buffenoir E., Roche P., Terras V., Quantum dynamical coboundary equation for
 finite dimensional simple {L}ie algebras, \href{https://doi.org/10.1016/j.aim.2007.02.001}{\textit{Adv. Math.}} \textbf{214}
 (2007), 181--229, \href{https://arxiv.org/abs/math.QA/0512500}{arXiv:math.QA/0512500}.

\bibitem{Bul}
Bullock D., A finite set of generators for the {K}auffman bracket skein
 algebra, \href{https://doi.org/10.1007/PL00004727}{\textit{Math.~Z.}} \textbf{231} (1999), 91--101.

\bibitem{BFK}
Bullock D., Frohman C., Kania-Bartoszy\'nska J., Topological interpretations of
 lattice gauge field theory, \href{https://doi.org/10.1007/s002200050471}{\textit{Comm. Math. Phys.}} \textbf{198} (1998),
 47--81, \href{https://arxiv.org/abs/q-alg/9710003}{arXiv:q-alg/9710003}.

\bibitem{Ca}
Caldero P., \'El\'ements ad-finis de certains groupes quantiques,
 \textit{C.~R.~Acad. Sci. Paris S\'er.~I Math.} \textbf{316} (1993), 327--329.

\bibitem{CP}
Chari V., Pressley A., A guide to quantum groups, Cambridge University Press,
 Cambridge, 1994.

\bibitem{CL}
Costantino F., L\^e T.T.Q., Stated skein algebras of surfaces, \href{https://doi.org/10.4171/jems/1167}{\textit{J.~Eur.
 Math. Soc.}} \textbf{24} (2022), 4063--4142, \href{https://arxiv.org/abs/1907.11400}{arXiv:1907.11400}.

\bibitem{Cui}
Cui W., Canonical bases of modified quantum algebras for type {$A_2$},
 \href{https://doi.org/10.1142/S021949881850113X}{\textit{J.~Algebra Appl.}} \textbf{17} (2018), 1850113, 27~pages,
 \href{https://arxiv.org/abs/1208.5531}{arXiv:1208.5531}.

\bibitem{DRZ}
Dabrowski L., Reina C., Zampa A., {$A({\rm SL}_q(2))$} at roots of unity is a
 free module over {$A({\rm SL}(2))$}, \href{https://doi.org/10.1023/A:1007601131002}{\textit{Lett. Math. Phys.}} \textbf{52}
 (2000), 339--342, \href{https://arxiv.org/abs/math.QA/0004092}{arXiv:math.QA/0004092}.

\bibitem{DC-K}
De~Concini C., Kac V.G., Representations of quantum groups at roots of {$1$},
 in Operator {A}lgebras, {U}nitary {R}epresentations, {E}nveloping {A}lgebras,
 and {I}nvariant {T}heory ({P}aris, 1989), \textit{Progr. Math.}, Vol.~92,
 Birkh\"auser, Boston, MA, 1990, 471--506.

\bibitem{DC-K-P1}
De~Concini C., Kac V.G., Procesi C., Quantum coadjoint action, \href{https://doi.org/10.2307/2152754}{\textit{J.~Amer.
 Math. Soc.}} \textbf{5} (1992), 151--189

\bibitem{DC-L}
De~Concini C., Lyubashenko V., Quantum function algebra at roots of {$1$},
 \href{https://doi.org/10.1006/aima.1994.1071}{\textit{Adv. Math.}} \textbf{108} (1994), 205--262.

\bibitem{DC-K-P2}
De~Concini C., Procesi C., Quantum groups, in {$D$}-modules, {R}epresentation
 {T}heory, and {Q}uantum {Q}roups ({V}enice, 1992), \textit{Lecture Notes in
 Math.}, Vol. 1565, \href{https://doi.org/10.1007/BFb0073466}{Springer}, Berlin, 1993, 31--140.

\bibitem{deGraaf}
de~Graaf W.A., Constructing canonical bases of quantized enveloping algebras,
 \href{https://doi.org/10.1080/10586458.2002.10504683}{\textit{Experiment. Math.}} \textbf{11} (2002), 161--170.

\bibitem{dRGGP-MR}
De~Renzi M., Gainutdinov A.M., Geer N., Patureau-Mirand B., Runkel I.,
 $3$-dimensional {TQFT}s from non-semisimple modular categories,
 \href{https://doi.org/10.1007/s00029-021-00737-z}{\textit{Selecta Math.~(N.S.}} \textbf{28} (2022), 42, 60~pages,
 \href{https://arxiv.org/abs/1912.02063}{arXiv:1912.02063}.

\bibitem{dRGP-M}
De~Renzi M., Geer N., Patureau-Mirand B., Renormalized Hennings Invariants and
 $2+1$-{TQFT}s, \href{https://doi.org/10.1007/s00220-018-3187-8}{\textit{Comm. Math. Phys.}} \textbf{362} (2020), 855--907,
 \href{https://arxiv.org/abs/1707.08044}{arXiv:1707.08044}.

\bibitem{DC}
Dieudonn\'e J.A., Carrell J.B., Invariant theory, old and new, \href{https://doi.org/10.1016/0001-8708(70)90015-0}{\textit{Adv. Math.}} \textbf{4} (1970), 1--80.

\bibitem{DomLen}
Domokos M., Lenagan T.H., Quantized trace rings, \href{https://doi.org/10.1093/qmathj/hah051}{\textit{Q.~J.~Math.}}
 \textbf{56} (2005), 507--523, \href{https://arxiv.org/abs/math.QA/0407053}{arXiv:math.QA/0407053}.

\bibitem{Dr}
Drinfeld V.G., On almost cocommutative {H}opf algebras, \textit{Leningrad
 Math.~J.} \textbf{1} (1990), 321--342.

\bibitem{Du}
Du J., Global {IC} bases for quantum linear groups, \href{https://doi.org/10.1016/S0022-4049(96)00045-X}{\textit{J.~Pure Appl.
 Algebra}} \textbf{114} (1996), 25--37.

\bibitem{Enr}
Enriquez B., Le centre des alg\'ebres de coordonn\'ees des groupes quantiques
 aux racines {$p^\alpha$}-i\'emes de l'unit\'e, \href{https://doi.org/10.24033/bsmf.2242}{\textit{Bull. Soc. Math.
 France}} \textbf{122} (1994), 443--485.

\bibitem{Et}
Etingof P., Golberg O., Hensel S., Liu T., Schwendner A., Vaintrob D., Yudovina
 E., Introduction to representation theory, \textit{Stud. Math. Libr.},
 Vol.~59, \href{https://doi.org/10.1090/stml/059}{American Mathematical Society}, Providence, RI, 2011,
 \href{https://arxiv.org/abs/0901.0827}{arXiv:0901.0827}.

\bibitem{Faitg}
Faitg M., Mapping class groups, skein algebras and combinatorial quantization,
 Ph.D.~Thesis, {M}ontpellier de {U}niversit\'e, 2019, \href{https://arxiv.org/abs/1910.04110}{arXiv:1910.04110}.

\bibitem{Faitg2}
Faitg M., Projective representations of mapping class groups in combinatorial
 quantization, \href{https://doi.org/10.1007/s00220-019-03470-z}{\textit{Comm. Math. Phys.}} \textbf{377} (2020), 161--198,
 \href{https://arxiv.org/abs/1812.00446}{arXiv:1812.00446}.

\bibitem{Faitg4}
Faitg M., Holonomy and (stated) skein algebras in combinatorial quantization,
 \textit{Quantum Topol.}, {t}o appear, \href{https://arxiv.org/abs/2003.08992}{arXiv:2003.08992}.

\bibitem{FR}
Fock V.V., Rosly A.A., Poisson structure on moduli of flat connections on
 {R}iemann surfaces and the {$r$}-matrix, in Moscow {S}eminar in
 {M}athematical {P}hysics, \textit{Amer. Math. Soc. Transl. Ser.}, Vol. 191,
 \href{https://doi.org/10.1090/trans2/191/03}{American Mathematical Society}, Providence, RI, 1999, 67--86,
 \href{https://arxiv.org/abs/math.QA/9802054}{arXiv:math.QA/9802054}.

\bibitem{FK}
Frenkel I.B., Khovanov M.G., Canonical bases in tensor products and graphical
 calculus for {$U_q(\mathfrak{sl}_2)$}, \href{https://doi.org/10.1215/S0012-7094-97-08715-9}{\textit{Duke Math.~J.}} \textbf{87}
 (1997), 409--480.

\bibitem{FKL}
Frohman C., Kania-Bartoszynska J., L\^e T., Unicity for representations of the
 {K}auffman bracket skein algebra, \href{https://doi.org/10.1007/s00222-018-0833-x}{\textit{Invent. Math.}} \textbf{215} (2019),
 609--650, \href{https://arxiv.org/abs/1707.09234}{arXiv:1707.09234}.

\bibitem{GJS19}
Ganev I., Jordan D., Safronov P., The quantum {F}robenius for character
 varieties and multiplicative quiver varieties, \href{https://doi.org/10.4171/JEMS/1427}{\textit{J.~Eur. Math. Soc.}}, {t}o appear, \href{https://arxiv.org/abs/1901.11450}{arXiv:1901.11450}.

\bibitem{Hum}
Humphreys J.E., Linear algebraic groups, \textit{Grad. Texts Math.}, Vol.~21,
 \href{https://doi.org/10.1007/978-1-4684-9443-3}{Springer}, New York, 1975.

\bibitem{JW}
Jordan D., White N., The center of the reflection equation algebra via quantum
 minors, \href{https://doi.org/10.1016/j.jalgebra.2019.08.038}{\textit{J.~Algebra}} \textbf{542} (2020), 308--342,
 \href{https://arxiv.org/abs/1709.09149}{arXiv:1709.09149}.

\bibitem{Jos}
Joseph A., Quantum groups and their primitive ideals, \textit{Ergeb. Math.
 Grenzgeb.~(3)}, Vol.~29, \href{https://doi.org/10.1007/978-3-642-78400-2}{Springer}, Berlin, 1995.

\bibitem{JL}
Joseph A., Letzter G., Local finiteness of the adjoint action for quantized
 enveloping algebras, \href{https://doi.org/10.1016/0021-8693(92)90157-H}{\textit{J.~Algebra}} \textbf{153} (1992), 289--318.

\bibitem{JL2}
Joseph A., Letzter G., Separation of variables for quantized enveloping
 algebras, \href{https://doi.org/10.2307/2374984}{\textit{Amer.~J. Math.}} \textbf{116} (1994), 127--177.

\bibitem{KK}
Karuo H., Korinman J., Classification of semi-weight representations of reduced
 stated skein algebras, \href{https://arxiv.org/abs/2303.09433}{arXiv:2303.09433}.

\bibitem{Kashiwara1}
Kashiwara M., On crystal bases of the {$Q$}-analogue of universal enveloping
 algebras, \href{https://doi.org/10.1215/S0012-7094-91-06321-0}{\textit{Duke Math.~J.}} \textbf{63} (1991), 465--516.

\bibitem{Kashiwara2}
Kashiwara M., Global crystal bases of quantum groups, \href{https://doi.org/10.1215/S0012-7094-93-06920-7}{\textit{Duke Math.~J.}}
 \textbf{69} (1993), 455--485.

\bibitem{Kashiwara3}
Kashiwara M., Crystal bases of modified quantized enveloping algebra,
 \href{https://doi.org/10.1215/S0012-7094-94-07317-1}{\textit{Duke Math.~J.}} \textbf{73} (1994), 383--413.

\bibitem{Kashiwara4}
Kashiwara M., On crystal bases, in Representations of {G}roups ({B}anff, {AB},
 1994), \textit{CMS Conf. Proc.}, Vol.~16, American Mathematical Society,
 Providence, RI, 1995, 155--197.

\bibitem{KhorTol}
Khoroshkin S.M., Tolstoy V.N., Universal {$R$}-matrix for quantized
 (super)algebras, \href{https://doi.org/10.1007/BF02102819}{\textit{Comm. Math. Phys.}} \textbf{141} (1991), 599--617.

\bibitem{KirRes}
Kirillov A.N., Reshetikhin N., {$q$}-{W}eyl group and a multiplicative formula
 for universal {$R$}-matrices, \href{https://doi.org/10.1007/BF02097710}{\textit{Comm. Math. Phys.}} \textbf{134} (1990),
 421--431.

\bibitem{KS}
Klimyk A., Schm\"udgen K., Quantum groups and their representations, \textit{Texts Monogr. Phys.}, \href{https://doi.org/10.1007/978-3-642-60896-4}{Springer}, Berlin, 1997.

\bibitem{KLNY}
Kolb S., Lorenz M., Nguyen B., Yammine R., On the adjoint representation of a
 {H}opf algebra, \href{https://doi.org/10.1017/S0013091520000358}{\textit{Proc. Edinb. Math. Soc.}} \textbf{63} (2020),
 1092--1099, \href{https://arxiv.org/abs/1905.03020}{arXiv:1905.03020}.

\bibitem{K0}
Korinman J., Unicity for representations of reduced stated skein algebras,
 \href{https://doi.org/10.1016/j.topol.2020.107570}{\textit{Topology Appl.}} \textbf{293} (2021), 107570, 28~pages,
 \href{https://arxiv.org/abs/2001.00969}{arXiv:2001.00969}.

\bibitem{K1}
Korinman J., Finite presentations for stated skein algebras and lattice gauge
 field theory, \href{https://doi.org/10.2140/agt.2023.23.1249}{\textit{Algebr. Geom. Topol.}} \textbf{23} (2023), 1249--1302,
 \href{https://arxiv.org/abs/2012.03237}{arXiv:2012.03237}.

\bibitem{KQ}
Korinman J., Quesney A., Classical shadows of stated skein representations at
 odd roots of unity, \href{https://arxiv.org/abs/1905.03441}{arXiv:1905.03441}.

\bibitem{Ko}
Kostant B., Groups over {$Z$}, in Algebraic {G}roups and {D}iscontinuous
 {S}ubgroups ({P}roc. {S}ympos. {P}ure {M}ath., {B}oulder, {C}olo., 1965),
 American Mathematical Society, Providence, RI, 1966, 90--98.

\bibitem{Lang}
Lang S., Algebraic structures, Addison-Wesley Publishing Co., Reading, Mass.,
 1967.

\bibitem{Le}
L\^e T.T.Q., Triangular decomposition of skein algebras, \href{https://doi.org/10.4171/QT/115}{\textit{Quantum
 Topol.}} \textbf{9} (2018), 591--632, \href{https://arxiv.org/abs/1609.04987}{arXiv:1609.04987}.

\bibitem{LeSik}
L\^e T.T.Q., Sikora A.S., Stated ${\rm SL}(n)$-skein modules and algebras,
 \href{https://arxiv.org/abs/2201.00045}{arXiv:2201.00045}.

\bibitem{LY}
L\^e T.T.Q., Yu T., Quantum traces and embeddings of stated skein algebras into
 quantum tori, \href{https://doi.org/10.1007/s00029-022-00781-3}{\textit{Selecta Math.~(N.S.)}} \textbf{28} (2022), 66, 48~pages,
 \href{https://arxiv.org/abs/2012.15272}{arXiv:2012.15272}.

\bibitem{LS}
Levendorskii S.Z., Soibelman Y.S., Some applications of the quantum
 {W}eyl groups, \href{https://doi.org/10.1016/0393-0440(90)90013-S}{\textit{J.~Geom. Phys.}} \textbf{7} (1990), 241--254.

\bibitem{Lusztig2}
Lusztig G., Quantum groups at roots of~{$1$}, \href{https://doi.org/10.1007/BF00147341}{\textit{Geom. Dedicata}}
 \textbf{35} (1990), 89--113.

\bibitem{Lusztig}
Lusztig G., Introduction to quantum groups, \textit{Mod. Birkh\"auser Class.},
 Vol. 110, \href{https://doi.org/10.1007/978-0-8176-4717-9}{Birkh\"auser}, Boston, MA, 1993.

\bibitem{Lusztig3}
Lusztig G., Study of a $\mathbb{Z}$-form of the coordinate ring of a reductive
 group, \href{https://doi.org/10.1090/S0894-0347-08-00603-6}{\textit{J.~Amer. Math. Soc.}} \textbf{22} (2009), 739--769.

\bibitem{LM}
Lyubashenko V., Majid S., Braided groups and quantum {F}ourier transform,
 \href{https://doi.org/10.1006/jabr.1994.1165}{\textit{J.~Algebra}} \textbf{166} (1994), 506--528.

\bibitem{Maj}
Majid S., Braided matrix structure of the {S}klyanin algebra and of the quantum
 {L}orentz group, \href{https://doi.org/10.1007/BF02096865}{\textit{Comm. Math. Phys.}} \textbf{156} (1993), 607--638,
 \href{https://arxiv.org/abs/hep-th/9208008}{arXiv:hep-th/9208008}.

\bibitem{Marlin}
Marlin R., Anneaux de {G}rothendieck des vari\'et\'es de drapeaux,
 \href{https://doi.org/10.24033/bsmf.1831}{\textit{Bull. Soc. Math. France}} \textbf{104} (1976), 337--348.

\bibitem{MC-R}
McConnell J.C., Robson J.C., Noncommutative {N}oetherian rings, \textit{Grad.
 Stud. Math.}, Vol.~30, \href{https://doi.org/10.1090/gsm/030}{American Mathematical Society}, Providence, RI, 2001.

\bibitem{Meus}
Meusburger C., Kitaev lattice models as a {H}opf algebra gauge theory,
 \href{https://doi.org/10.1007/s00220-017-2860-7}{\textit{Comm. Math. Phys.}} \textbf{353} (2017), 413--468,
 \href{https://arxiv.org/abs/1607.01144}{arXiv:1607.01144}.

\bibitem{Mu}
Murakami J., Generalized {K}ashaev invariants for knots in three manifolds,
 \href{https://doi.org/10.4171/QT/86}{\textit{Quantum Topol.}} \textbf{8} (2017), 35--73, \href{https://arxiv.org/abs/1312.0330}{arXiv:1312.0330}.

\bibitem{Paradowski}
Paradowski J., Filtrations of modules over the quantum algebra, in Algebraic
 {G}roups and their {G}eneralizations: {Q}uantum and {I}nfinite-{D}imensional
 {M}ethods ({U}niversity {P}ark, {PA}, 1991), \textit{Proc. Sympos. Pure
 Math.}, Vol.~56, \href{https://doi.org/10.1090/pspum/056.2/1278729}{American Mathematical Society}, Providence, RI, 1994,
 93--108.

\bibitem{PW}
Parshall B., Wang J.P., Quantum linear groups, \href{https://doi.org/10.1090/memo/0439}{\textit{Mem. Amer. Math. Soc.}}
 \textbf{89} (1991), vi+157~pages.

\bibitem{PS}
Przytycki J.H., Sikora A.S., Skein algebras of surfaces, \href{https://doi.org/10.1090/tran/7298}{\textit{Trans. Amer.
 Math. Soc.}} \textbf{371} (2019), 1309--1332, \href{https://arxiv.org/abs/1602.07402}{arXiv:1602.07402}.

\bibitem{RSTS}
Reshetikhin N.Yu., Semenov-Tian-Shansky M.A., Quantum {$R$}-matrices and
 factorization problems, \href{https://doi.org/10.1016/0393-0440(88)90018-6}{\textit{J.~Geom. Phys.}} \textbf{5} (1988), 533--550
 (1989).

\bibitem{RT}
Reshetikhin N.Yu., Turaev V.G., Invariants of {$3$}-manifolds via link
 polynomials and quantum groups, \href{https://doi.org/10.1007/BF01239527}{\textit{Invent. Math.}} \textbf{103} (1991),
 547--597.

\bibitem{Rowen}
Rowen L.H., Ring theory, Academic Press, Inc., Boston, MA, 1991.

\bibitem{Saito}
Saito Y., P{BW} basis of quantized universal enveloping algebras, \href{https://doi.org/10.2977/prims/1195166130}{\textit{Publ.
 Res. Inst. Math. Sci.}} \textbf{30} (1994), 209--232.

\bibitem{Soi}
Soibelman Y.S., Algebra of functions on a compact quantum group and its
 representations, \textit{Leningrad Math.~J.} \textbf{2} (1990), 161--268.

\bibitem{Sp}
Springer T.A., Invariant theory, \textit{Lect. Notes in Math.}, Vol. 585,
 \href{https://doi.org/10.1007/BFb0095644}{Springer}, Berlin, 1977.

\bibitem{StackP}
{The Stacks Project Authors}, Commutative {A}lg, {T}he {S}tack {P}roject,
 {C}hapter~10, available at \url{https://stacks.math.columbia.edu}.

\bibitem{StackP2}
{The Stacks Project Authors}, Brauer groups, {T}he {S}tack {P}roject,
 {C}hapter~11, available at \url{https://stacks.math.columbia.edu}.

\bibitem{SV}
Vaksman L.L., Soibelman Y.S., Algebra of functions on the quantum group ${\rm
 SU}(2)$, \href{https://doi.org/10.1007/BF01077623}{\textit{Funct. Anal. Appl.}} \textbf{22} (1988), 170--181.

\bibitem{VV}
Varagnolo M., Vasserot E., Double affine {H}ecke algebras at roots of unity,
 \href{https://doi.org/10.1090/S1088-4165-2010-00384-2}{\textit{Represent. Theory}} \textbf{14} (2010), 510--600,
 \href{https://arxiv.org/abs/math.RT/0603744}{arXiv:math.RT/0603744}.

\bibitem{VY}
Voigt C., Yuncken R., Complex semisimple quantum groups and representation
 theory, \textit{Lect. Notes in Math.}, Vol. 2264, \href{https://doi.org/10.1007/978-3-030-52463-0}{Springer}, Cham, 2020.

\bibitem{Wang}
Wang Z., On stated ${\rm SL}(n)$-skein modules, \href{https://arxiv.org/abs/2307.10288}{arXiv:2307.10288}.

\bibitem{Wi}
Witten E., Topological quantum field theory, \href{https://doi.org/10.1007/BF01223371}{\textit{Comm. Math. Phys.}}
 \textbf{117} (1988), 353--386.

\bibitem{Xi}
Xi N.H., Root vectors in quantum groups, \href{https://doi.org/10.1007/BF02564506}{\textit{Comment. Math. Helv.}}
 \textbf{69} (1994), 612--639.

\end{thebibliography}
\end{document}